\documentstyle[12pt]{article}

\newcommand{\numero}[1]{
\addtocounter{section}{1}
\begin{center}{\bf \thesection .\
#1\vspace{-.1in}}\end{center}
\setcounter{subsection}{0}
\setcounter{lemma}{0}\indent}

\newcommand{\subnumero}[1]{
\pagebreak[1]\begin{center}{\em #1}\nopagebreak\end{center} }

\newcommand{\eop}{\hfill $\Box$\vspace{.1in}}

\newtheorem{lemma}{Lemma}[section]
\newtheorem{theorem}[lemma]{Theorem}

\newtheorem{corollary}[lemma]{Corollary}

\newtheorem{proposition}[lemma]{Proposition}

\newcommand{\zz}{{\bf Z}}
\newcommand{\nn}{{\bf N}}

\newcommand{\Cc}{{\cal C}}

\newcommand{\Ff}{{\cal F}}
\newcommand{\Gg}{{\cal G}}
\newcommand{\Qq}{{\cal Q}}

\newcommand{\Dd}{{\cal D}}
\newcommand{\Aa}{{\cal A}}
\newcommand{\Bb}{{\cal B}}
\newcommand{\Mm}{{\cal M}}

\newcommand{\Uu}{{\cal U}}
\newcommand{\Ll}{{\cal L}}

\newcommand{\Ii}{{\cal I}}
\newcommand{\Jj}{{\cal J}}
\newcommand{\Nn}{{\cal N}}
\newcommand{\Pp}{{\cal P}}
\newcommand{\Rr}{{\cal R}}

\newcommand{\Xx}{{\cal X}}
\newcommand{\Yy}{{\cal Y}}
\newcommand{\Zz}{{\cal Z}}

\newcommand{\mylabel}[1]{\label{#1}}

\begin{document}

\section*{Flexible sheaves}

Carlos SIMPSON\newline
Laboratoire Emile Picard,\newline UMR 5580, CNRS\newline
Universit\'e Paul Sabatier\newline
31062 Toulouse CEDEX

\begin{center}
{\bf Introduction}
\end{center}

In the early 1970's was introduced a notion of {\em homotopy-coherent diagram}
(Segal \cite{SegalTopology}, Leitch \cite{Leitch}, Vogt \cite{Vogt1} and Mather
\cite{Mather}). This is a way of using cubical homotopies to deal at once with
all of the higher homotopies of coherence involved in a diagram of spaces which
``commutes up to homotopy''. This notion was subsequently studied by Vogt
\cite{Vogt1} \cite{Vogt2} and others (for example Hardie and Kamps
\cite{HardieKamps}---see also the references therein), and especially by Cordier
and Porter  \cite{Cordier1} \cite{Cordier2} \cite{Cordier3}, \cite{CoPo1}
\cite{CoPo2} \cite{CoPoBook} \cite{CoPo3} \cite{CoPo4} \cite{CordierPorter},
\cite{Porter1} \cite{Porter2}.

In {\tt v1} of the present paper, unaware of all of the references cited
above, I looked at homotopy-coherent diagrams of spaces over a Grothendieck
site $\Xx$, calling them ``flexible functors'' or ``flexible presheaves''. The
present revision is designed to correct the problem of inadequate references;
however, we retain the terminology ``flexible functor''.

One of the main results which one shows is that a homotopy-coherent diagram can
be replaced by an equivalent strictly commutative diagram. This
``rectification'' was done by Segal in \cite{SegalTopology} at the same time as
he introduced the notion of homotopy-coherent diagram. A similar theorem is that
of Dwyer and Kan \cite{DwyerKanEquivs}. This is also discussed by Cordier and
Porter \cite{CoPo2}. See \S 3 below for our version. (In fact, the whole
story dates back to SGA 1 where a fibered category is replaced by an equivalent
split one.)

As noticed early on by Vogt \cite{Vogt1}, the advantage of the
homotopy-coherent  point of view over that of strictly-commutative diagrams, is
that it allows one to get a hold on the space of morphisms between two diagrams
$R$ and $T$, in an explicit way: a morphism from $R$ to $T$ is just a
homotopy-coherent diagram (flexible presheaf) on $\Xx \times I$ which restricts
to $R$ and $T$ on $\Xx \times \{ 0\}$ and $\Xx \times \{ 1\}$. These morphism
spaces don't have a natural composition but are naturally organized into a
``Segal category'', i.e. a simplicial space satisfying the condition that the
{\em Segal maps}  are equivalences (see \cite{SegalTopology}, \cite{Adams},
\cite{DwyerKanSmith}, \cite{effective}). A version of this observation occurs in
Vogt \cite{Vogt1}, where the morphism sets are not topologized but where the
resulting simplicial set satisfies the {\em restricted Kan condition}
(\cite{BoardmanVogt}). Vogt proved that the resulting category of morphisms
up to
homotopy, is equivalent to the homotopy category of diagrams obtained by
inverting the level-wise homotopy equivalences \cite{Vogt1}, see Cordier-Porter
\cite{CoPo2}.

One of the main steps in Vogt's argument, also taken up by Cordier and Porter
\cite{CoPo2}, was to show invertibility of level-wise weak equivalences as
homotopy-coherent maps. Below, we show a refined version, that the space of
homotopy-coherent diagrams indexed by the category $\overline{I}$ having two
isomorphic objects $0,1$, is contractible. (It is not clear whether this
contractibility was known to Vogt; in any case it is clearly not present
in \cite{CoPo2}.) Crucial to our treatment of the Segal condition for the
morphism spaces, and this  invertibility, is an obstruction theory first
developed by Cooke \cite{Cooke} and Dwyer-Kan-Smith \cite{DwyerKanSmith}. Our
treatment of this obstruction theory is at the end of \S 1. The applications to
morphisms are in \S 2.

What seems (according my present state of knowledge of the literature, at least)
to be the main new contribution of the present paper is a generalization of
Vogt's theorem to the case of a site $\Xx$ with nontrivial Grothendieck
topology. We show how to obtain the morphisms in Illusie's derived category
\cite{Illusie1}
(obtained from the category of simplicial presheaves on $\Xx$ by inverting the
Illusie weak equivalences) as homotopy classes of homotopy-coherent morphisms.
For this, it is natural to introduce the notion of {\em flexible sheaf} (\S 4).
This is a diagram which satisfies an additional descent condition, analogous to
the usual descent condition for sheaves of sets, or similarly analogous to the
descent condition for $1$-stacks. In \S 4 we treat this condition and start to
define an operation $T\mapsto HT$ related to the process of taking the
associated sheaf. In particular, $T$ is a flexible sheaf if and only if
$T\rightarrow HT$ is an (object-by-object) equivalence. After some further
materiel in \S 5 and the start of \S 6 generalizing standard constructions for
spaces (e.g. homotopy group sheaves, fiber products), we get to the
construction of the flexible sheaf associated to an $n$-truncated flexible
presheaf $T$, given as $H^{n+2}T$. This is the analogue of the classical
construction of the sheaf associated to a presheaf of sets (the case $n=0$), as
obtained by applying a natural operation twice.
In \S 7 we get to the statement of the analogue of Vogt's theorem: up to
homotopy, morphisms from $R$ to $T$ in Illusie's derived category can be
represented as homotopy-coherent morphisms from $R$ to $T$ under the
condition that $T$ be a flexible sheaf.

Starting from \S 8 on, we treat some standard things from topology in this
``relative'' case. In \S 10 we draw the relationship between the objects we
study here and classical $1$-stacks: this is the case where the values of the
diagram $T$ are of the form $K(\pi , 1)$.

In the original version (1993) there were several unfinished chapters at the
end of the paper heading towards some applications of these ideas e.g.
to nonabelian de Rham  cohomology. These chapters were left out in the
electronic
preprint {\tt v1}, because they were unfinished (and partly wrong). These things
may now be found in \cite{kobe}, \cite{SantaCruz}, \cite{RelativeLie},
\cite{geometricN}, and the reader in quest of motivation is refered there.

Another recent preprint concerning these questions is \cite{descente} (joint
with A. Hirschowitz). There we present a more modern approach based on closed
model categories. The condition which we call below ``flexible sheaf'' is the
same as the condition of being a {\em stack} in \cite{descente}. The
generalization of Vogt's theorem in \S 7 below is essentially the same as the
result of \cite{descente} Lemme 9.2 stating that a simplicial presheaf is
fibrant for the site $\Xx$ if and only if it is a stack (i.e. flexible sheaf in
the terminology of the present paper), and fibrant for the category $\Xx$ with
the coarse topology.

In a similar vein, the version of invertibility of homotopy
equivalences which we give below (cf \cite{Vogt1}, \cite{CoPo2}) occurs in a
somewhat different guise as Theorem 2.5.1 of \cite{limits}.

\bigskip

I apologize to all of the authors whose work was not correctly referred to in
the first preprint version {\tt v1}. I hope that the references in the present
version are substantially correct and complete; but if anybody knows of other
references or comments which need to be added, I would appreciate it if they
could let me know.

\newpage

\setcounter{section}{0}

\numero{Flexible functors}

We define what we call a ``flexible functor''. This is what is known in the
literature as a ``homotopy-coherent functor''. This fundamental notion is due
(c. 1973) to G. Segal \cite{SegalTopology}, R. Leitch \cite{Leitch}, R. Vogt
\cite{Vogt1} \cite{Vogt2}, M. Mather \cite{Mather} and has been extensively
investigated by these authors, by J.-M. Cordier and T. Porter \cite{Cordier1}
\cite{Cordier2} \cite{Cordier3} \cite{Porter1} \cite{Porter2}
\cite{CoPo1}  \cite{CoPo2} \cite{CoPo3} \cite{CordierPorter},
Edwards and Hastings \cite{EdwardsHastings},
W. Dwyer and D. Kan
\cite{DwyerKanEquivs} \cite{DwyerKanSmith} and others. The approaches of Segal
and Leitch \cite{SegalTopology} \cite{Leitch}, Vogt \cite{Vogt1}, and Dwyer-Kan
\cite{DwyerKan1} were compared by Cordier in \cite{Cordier1} \cite{Cordier2} and
shown to be the same. Our approach is the same as these.

A {\em continuous category} $\Cc$ is a category enriched over the
category $Top$ of topological spaces \cite{Kelly}, i.e. a category in
the category of topological spaces, where the space of objects is given the
discrete topology. In other words, it consists of a set of objects, and for any
two objects a space of morphisms $\Cc (X,Y)$, such that composition is
given by a
continuous map, and containing the identity. There is an obvious notion of
continuous functor between two continuous categories, including the condition
that the image of the identity is the identity.

A {\em continuous semi-category} $\Cc$ is a semi-category (that is, an
object like a category but without the existence of identity morphisms
presupposed) in
the category of topological spaces, again where the space of objects is a
discrete set. In other words, it consists of a set of objects, and for any two
objects a space of morphisms $\Cc (X,Y)$, such that composition is given by a
continuous map. There is an obvious notion of continuous functor between two
continuous semi-categories, or from a continuous semi-category to a continuous
category. This is different from the notion of continuous functor between two
continuous categories, in that the latter notion includes the condition of
respecting the identities.

We keep the natural structure of topological space inherited by the set of
functors $Hom (\Cc , \Dd )$ between two continuous semi-categories.

The category $Top$ of topological spaces has a structure of continuous
category: the mapping spaces are given the compact-open topology.

The following definition underlies the Segal-Leitch-Vogt approach to
homotopy coherent diagrams \cite{Vogt1}, \cite{Leitch}, see also
\cite{SegalTopology}.  Suppose $\Xx$ is a category. We define a continuous
semi-category $M\Xx$ to be the category whose objects are the objects of $\Xx$,
but where the set of morphisms from $X$ to $Y$ is replaced by a {\em topological
space} of morphisms, $M(X,Y)$ (denoted $M\Xx (X,Y)$ if there is risk of
confusion), defined as follows. The points of $M(X,Y)$ consist of pairs $(\phi
,t)$ where $\phi =(\phi _1,\ldots , \phi _n)$ is a composable sequence of
morphisms in $\Xx$ with $X_0=Y$ and $X_n=X$, and $t=(t_1,\ldots , t_{n-1})\in
[0,1]^{n-1}$. The points are identified by an equivalence relation: if $t_i=0$
then $(\phi ,t)\sim (\phi ',t')$ where $t_i$ is removed to obtain $t'$ and $\phi
_i$ and $\phi _{i+1}$ are composed to obtain $\phi '$. The space is given the
obvious topology which is the quotient of the usual topology on the disjoint
union of the hypercubes, modulo the equivalence relation. Composition is defined
by concatenating the sequences $\phi$, and by concatenating the sequences $t$
with a $t_j=1$ added in between. This gives an associative composition where the
law $M(X,Y)\times M(Y,Z)$ is a continuous map.

There is a distinguished object $(1_X)$ in $M(X,X)$, which does not act as
the identity, however. If $\Cc$ is a continuous category, a {\em flexible
functor} (or ``homotopy coherent diagram'' \cite{SegalTopology}
\cite{Leitch} \cite{Vogt1} \cite{Mather} \cite{Cordier1} \cite{Cordier2})
$$
T:\Xx \rightarrow \Cc
$$
is a continuous functor (of
semicategories) from $M\Xx$ to $\Cc$ such that $T(1_X)=1_{T_X}$. In particular,
for $\Cc =Top$ we obtain the notion of {\em flexible functor} $T\rightarrow
Top$. In principle we may consider covariant or contravariant functors (note
that $M(\Xx ^o)=(M\Xx )^o$). For our purposes, we make the convention that we
always speak of {\em contravariant flexible functors} and drop the word
contravariant. We can also use the term {\em flexible presheaf} (of topological
spaces) for a flexible functor $T\rightarrow Top$.

We give a more explicit description (which follows directly from the
definition of $M\Xx$). See \cite{Leitch} \cite{SegalTopology}.
A flexible functor $T:\Xx \rightarrow Top$ consists of the following data.
For each object $X\in \Xx$, a topological space $T_X$; and for each
$n$-tuple of morphisms $(\varphi _1,\ldots , \varphi _n)$ with $\varphi
_i:X_i\rightarrow X_{i-1}$ (for objects $X_0,\ldots ,X_n$), a morphism
$$
T(\varphi _1,\ldots , \varphi _n): [0,1]^{n-1}\times T_{X_0}\rightarrow
T_{X_n} $$
denoted $T(\varphi _1,\ldots , \varphi _n; t_1,\ldots ,
t_{n-1}):T_{X_0}\rightarrow T_{X_n}$,
such that for each $i=1,\ldots , n-1$,
$$
T(\varphi _1,\ldots , \varphi _n; t_1,\ldots , t_{i-1}, 0, t_{i+1}, \ldots
, t_{n-1}) =
$$
$$
T(\varphi _1,\ldots , \varphi _i\varphi _{i+1},\ldots , \varphi _n;
t_1,\ldots , t_{i-1}, t_{i+1}, \ldots , t_{n-1}) $$
and
$$
T(\varphi _1,\ldots , \varphi _n; t_1,\ldots , t_{i-1}, 1, t_{i+1}, \ldots
, t_{n-1}) =
$$
$$
T(\varphi _{i+1},\ldots , \varphi _n; t_{i+1}, \ldots , t_{n-1}) T(\varphi
_1,\ldots , \varphi _i; t_1,\ldots , t_{i-1}). $$
Finally, $T(1_X)=1_{T_X}$.
Note that for $n=1$, $T(\phi _1)$ is simply a morphism from $T_{X_0}$ to
$T_{X_1}$, while for $n=2$, $T(\phi _1,\phi _2; t_1)$ is a homotopy between
$T(\phi _1\phi _2)$ and $T(\phi _2)T(\phi _1)$. For $n\geq 3$ the data
consists of homotopies between the homotopies, etc., fitting together in a
natural way.

We say that $T$ is {\em multi-identity normalized} if $T(1_X, \ldots , 1_X;
t)= 1_{T_X}$. (We will not use this notion very often.) For example, if
$\Xx$ is a category with one element, then a multi-identity normalized
flexible functor to $Top$ is the same as a topological space.

There is a functor $\Xi :M\Xx
\rightarrow \Xx$ obtained by forgetting the $t$'s and composing the
composable sequences $\phi$. It is continuous if the morphism sets of $\Xx$
are given the discrete topology. Furthermore, if $f:X\rightarrow Y$ is a
morphism in $\Xx$ then $\Xi ^{-1}(f)$ is a nonempty contractible set in
$M(X,Y)$ (the contraction sends all $t_i$ to zero by multiplying by
$\lambda \in [0,1]$). In fact, $M\Xx$ is universal for this property.

\begin{proposition}
\mylabel{flex01}
Suppose $\Cc$ is a continuous category (i.e. where the morphism sets are
topological spaces), and suppose $\Phi :\Cc \rightarrow \Xx$ is a
continuous functor (where the morphism sets of $\Xx$ are given the discrete
topology). Suppose that for each morphism $f:X\rightarrow Y$ in $\Xx$, the
space $\Phi ^{-1}(f)$ is nonempty and contractible (resp. weakly
contractible). Then the space of continuous functors $F: M\Xx \rightarrow
\Cc$ such that $\Phi F=\Xi$ is nonempty and contractible (resp. weakly
contractible).
\end{proposition}
{\em Proof:}
Look at the ways of defining $F$, by induction on the dimension of
primitive cells in $M\Xx$ (that is, the hypercubes that make up $M\Xx$ in
the above definition). The space of ways of determining $F$ on each
primitive cell of $M\Xx$ (given the values already specified on the
boundary) is contractible (resp. weakly contractible). The space of ways of
determining $F$ can be seen as the inverse limit of maps where the fibers
are (weakly) contractible, so it is (weakly) contractible. \eop

\subnumero{Functoriality}
If $G:\Xx \rightarrow \Yy$ is a functor, and $T:\Yy \rightarrow \Cc$ is a
flexible functor, then we obtain naturally a flexible functor $G^{\ast}T:\Xx
\rightarrow \Cc$. This satisfies the strict associativity
$H^{\ast}G^{\ast}T=(GH)^{\ast}T$. In fact there is a natural continuous
functor $MG:M\Xx \rightarrow M\Yy$, and $G^{\ast}T=T\circ MG$.

\subnumero{Postnikov tower for continuous categories}

The notion of continuous category is an example of a more general notion
of {\em enriched category} see Kelly \cite{Kelly}.
If $\Zz$ is a category admitting products and an initial object, we can
define the notion of $\Zz$-category $\Cc$: here the ``set of morphisms''
$\Cc (A,B)$ is an object of $\Zz$, and the composition law is a morphism in
$\Zz$; associativity is given by equality of two morphisms in $\Zz$ and the
identities are morphisms from the initial object to the morphism objects
(\cite{Kelly}). A continuous category is just a $Top$-category.

\begin{lemma}
\mylabel{flex02}
Suppose $F:\Zz\rightarrow \Zz '$ is a functor compatible with products. Then we
obtain a functor $\Cc \mapsto F\Cc$ from the category of $\Zz$-categories to the
category of $\Zz '$-categories. If $F$ and $F'$ are two such functors and $\eta
:F\rightarrow F'$ is a natural transformation compatible with products, then
$\eta$ induces a natural transformation $\eta _{\Cc}:F\Cc \rightarrow F'\Cc$
(that is, $\eta _{\Cc}$ is a $\Zz '$-functor from $F\Cc$ to $F'\Cc$).
\end{lemma} {\em Proof:} Keep the same object sets. Put $F\Cc (A,B):= F(\Cc
(A,B))$. From compatibility with products, applying $F$ to the multiplication $$
\Cc (A,B)\times \Cc (B,C)\rightarrow \Cc (A,C) $$
gives a multiplication
$$
F\Cc (A,B)\times F\Cc (B,C)\rightarrow F\Cc (A,C). $$
Given a natural transformation $\eta$, put $\eta _{\Cc}(A)=A$ and define $$
\eta _{\Cc , A,B}:= \eta _{\Cc (A,B)}:F\Cc (A,B)\rightarrow F'\Cc (A,B). $$
Via the isomorphisms $F\Cc (A,B)\times F\Cc (B,C)\cong F(\Cc (A,B)\times
\Cc (B,C))$ and the same for $F'$, we have
$$
\eta _{\Cc (A,B)}\times \eta _{\Cc (B,C)}= \eta _{\Cc (A,B)\times \Cc (B,C)}
$$
(definition of compatibility of $\eta$ with products), and from the
naturality of $\eta$ with respect to the morphism
$$
\Cc (A,B)\times \Cc (B,C)\rightarrow \Cc (A,C) $$
we get a commutative diagram
$$
\begin{array}{ccc}
F(\Cc (A,B)\times \Cc (B,C)) & \rightarrow & F\Cc (A,C) \\ \downarrow &&
\downarrow \\
F'(\Cc (A,B)\times \Cc (B,C)) & \rightarrow & F'\Cc (A,C) , \end{array}
$$
which when composed with the previous statement gives the commutative diagram $$
\begin{array}{ccc}
F\Cc (A,B)\times F\Cc (B,C) & \rightarrow & F\Cc (A,C) \\ \downarrow &&
\downarrow \\
F'\Cc (A,B)\times F'\Cc (B,C) & \rightarrow & F'\Cc (A,C) \end{array}
$$
needed to give $\eta _{\Cc}$ a structure of functor. \eop

Now we use the Postnikov decomposition in $Top$ to get a Postnikov decomposition
for a continuous category $\Cc$. We need to make the Postnikov decomposition
functorial and compatible with products; for this it is necessary to pass
through the simplicial point of view.

We will define a functor $\tau _{\leq n}:Top\rightarrow Top$ which is compatible
with products, and which induces the
coskeleton on homotopy (that is, truncates homotopy groups $\pi _i$ for $i>n$).
This gives a functor $\Cc \rightarrow \tau_{\leq n}\Cc$ from the category of
spacial categories to itself. We will have natural transformations $\tau _{\leq
n}\rightarrow \tau _{\leq n-1}$ and an auxiliary functor $\tau _{\infty}$ with
natural transformations
$$
\Cc \leftarrow \tau _{\infty} \Cc \rightarrow \tau
_{\leq n},
$$
all compatible with products. The first is a weak equivalence and
the second induces isomorphisms on homotopy groups $\pi _i$ for $i\leq n$.
Furthermore, the composition
$$
\tau _{\infty} \Cc \rightarrow \tau _{\leq n}
\rightarrow \tau _{\leq n-1}
$$
is equal to the projection for $n-1$.

In order to define the $\tau_{\leq n}$ in a way compatible with products, we
pass through the category of simplicial sets $SplSet$. Let $Sing :Top
\rightarrow SplSet$ denote the singular complex functor; let
$Re: SplSet \rightarrow
Top$ denote the realization functor; and let $cosk_n:SplSet \rightarrow
SplSet$ denote the $n$th simplicial coskeleton functor. Each of these
functors is
compatible with products. For $Sing$ and $cosk_n$ this is clear (cf Artin-Mazur
\cite{Artin-Mazur} for the coskeleton; it is adjoint to the skeleton, and the
adjunction formula implies compatibility with products). For $Re$, this is the
theorem about dividing a product of simplices into simplices.

For Kan complexes, the simplicial coskeleton functor $cosk _n$ truncates
homotopy groups $\pi _i$ for $i\geq n$. This is not true for complexes which are
not Kan. However, applying $Sing$ to a space automatically give a Kan complex.
Put
$$
\tau_{\leq n}:= Re \circ cosk_{n+1} \circ Sing .
$$
This has the
required effect on homotopy groups. We have a natural transformation $cosk
_{n+1}\rightarrow cosk _n$ compatible with products and the initial object,
which gives rise to a natural transformation $\tau _{\leq n}\rightarrow \tau
_{\leq {n-1}}$.
Put
$$
\tau _{\infty}:=Re\circ Sing .
$$
We have natural transformations $S\rightarrow cosk_{n+1}S$, and hence
$$
Re \circ Sing (X) \rightarrow Re \circ cosk_{n+1} \circ Sing
$$
giving
$$
\tau _{\infty} \rightarrow \tau _{\leq n}.
$$
On the other hand, we have a natural transformation (weak equivalence)
$$
Re(Sing (X)) \rightarrow X,
$$
giving the weak equivalence $\tau _{\infty}X\rightarrow X$. Finally,
note (by the adjunction formula \cite{Artin-Mazur}) that the projection
$S\rightarrow cosk _nS$ is equal to the composition
$$
S\rightarrow cosk_{n+1}S\rightarrow cosk _nS.
$$
This gives the same result for the $\tau$.

Applying Lemma \ref{flex02}, we obtain a sequence of functors $\Cc \mapsto \tau
_{\leq n}\circ \Cc$ from the category of topological categories to itself:
$$
Hom _{\tau _{\leq n}\circ \Cc}(X,Y) := \tau _{\leq n} Hom _{\Cc}(X,Y).
$$
We have natural transformations
$$
\Cc \leftarrow \tau _{\infty} \circ \Cc \rightarrow \tau _{\leq n} \circ \Cc
$$
as well as
$$
\tau _{\leq n}\circ \Cc \rightarrow \tau _{\leq n-1}\circ \Cc
$$
giving a commutative triangle with the morphisms
$$
\tau _{\infty} \circ \Cc \rightarrow \tau _{\leq n}\circ  \Cc
\;\; \; \mbox{and}\;\;\;   \tau _{\infty} \circ
\Cc \rightarrow \tau _{\leq n-1} \circ \Cc .
$$
This collection of functors and natural transformations is the {\em Postnikov
tower} for $\Cc$.

{\em Remark:} What we obtain as $\tau _{\leq n} \circ \Cc$ in the above
discussion is equivalent to the truncation $\tau _{\leq n+1} \Cc$ of
\cite{Tamsamani} \cite{nCat} \cite{descente} etc. (this is because here we are
applying the truncation operation to the morphism spaces in what is already a
$1$-category enriched over $Top$).

For example we obtain a Postnikov tower for $Top$ itself. The stage
$\tau _{\leq 1} \circ Top=\tau _{\leq 2}(Top)$ occured in Gabriel-Zisman
\cite{GabrielZisman}.

\subnumero{Functors from free continuous categories}

The ``obstruction theory'' that follows (in the next several sections) is due to
Cooke \cite{Cooke} and  Dwyer-Kan-Smith \cite{DwyerKanSmith}.
See also Schw\"{a}nzl and Vogt \cite{SchwanzlVogt}, Baues and Wirsching
\cite{BauesWirsching}, Edwards and Hastings
\cite{EdwardsHastings}.

We go to a general situation. Let $\Nn$ be a continuous category and
suppose
that $\Ff \rightarrow \Gg$ is a morphism of continuous categories. Suppose
$G_0\in Hom (\Nn , \Gg )$. We can form the {\em fiber square} $$
\begin{array}{ccc} \Pp &\rightarrow & \Ff \\
\downarrow && \downarrow \\
\Nn & \rightarrow & \Gg
\end{array}
$$
where the horizontal map on the bottom is $G_0$ and where $$
Ob (\Pp )= Ob (\Nn )\times _{Ob (\Gg )}Ob (\Ff ), $$
and for two objects $(X,A)$ and $(X',A')$ in $Ob (\Pp )$ (mapping to $Y$
and $Y'$ in $Ob (\Gg )$ we have
$$
\Pp ((X,A),(X',A')):= \Nn (X, X') \times ^{Path}_{\Gg (Y,Y')}\Ff (A,A'). $$
Here the path-space fiber product of $f:R\rightarrow S$ and $g:T\rightarrow
S$ is defined by
$$
R\times ^{Path}_ST := \{ (r,t,\gamma ):\; r\in R ,\; t\in T,\; \gamma \in
Path _S(f(r),g(t)) \} .
$$
{\em Notation:} We can denote this path-space fiber product spacial
category by $$
\Pp = \Nn \times ^{Path}_{\Gg}\Ff .
$$

\begin{lemma}
\mylabel{flex03}
In the above situation, the homotopy fiber of the map $$
Hom (\Nn , \Ff )\rightarrow Hom (\Nn , \Gg ) $$
over $G_0$ is equal to the space of sections $\Nn \rightarrow \Pp$ of the
fibration $\Pp \rightarrow \Nn$, where $\Pp = \Nn \times ^{Path}_{\Gg}\Ff$
is defined as above, compatibly with the map $$
Sect (\Pp /\Nn )\rightarrow Hom (\Nn , \Ff ). $$
\end{lemma}
{\em Proof:}
The homotopy fiber of $Hom (\Nn , \Ff )\rightarrow Hom (\Nn , \Gg )$ is
equal to the space of pairs $(F_1,\{ G_t\} )$ where $F_1:\Nn \rightarrow
\Ff $ and $G_t$ is a continuous family of functors from $\Nn $ to $\Gg$
indexed by $t\in [0,1]$, with $G_1$ the projection of $F_1$ and $G_0$ as
given. On the other hand, a section of $\Pp /\Nn$ associates to each $X\in
Ob (\Nn )$ an $A\in Ob (\Ff )$ projecting to the same $Y\in Ob (\Gg )$, and
to each $r\in \Nn (X,X')$ a triple $(r, t, \gamma )$ such that $t\in \Ff
(A,A')$ and $\gamma$ is a path of elements in $\Gg (Y,Y')$ joining the
images of $r$ and $t$. The associations $X\mapsto Y$ and $r\mapsto \gamma
(t)$ provide a continuous family of functors $G_t$ (and vice-versa), while
the associations $X\mapsto A$ and $r\mapsto t$ provide a functor $F_1$. The
family $G_t$ has $G_0$ the given functor and $G_1$ the projection of $F_1$.
Thus, the homotopy fiber of $Hom (\Nn , \Ff )\rightarrow Hom (\Nn , \Gg )$
over $G_0$ is equal to $Sect (\Pp /\Nn )$.
\eop

Suppose that the morphism $\Ff \rightarrow \Gg$ induces an isomorphism on
the set of objects. Then $\Pp \rightarrow \Nn$ also induces an isomorphism
on the set of objects, and we can denote objects of $\Nn$ and $\Pp$ by the
same letter $X$ (and their images in $\Ff$ or $\Gg$ by $Y$). The map $\Pp
(X,X')\rightarrow \Nn (X,X')$ is a fibration whose fiber over $r\in \Nn
(X,X' )$ is equal to the homotopy fiber of $\Ff (Y,Y')\rightarrow \Gg
(Y,Y')$ over $G_0(r)$.

{\em Definition:}
We say that a continuous category $\Nn$ is {\em CW-free} if it is an
increasing union of subcategories $\Nn _i$ indexed by an ordinal $I$, such
that putting $\Nn _{i-1}:=\bigcup _{j<i}\Nn _j$ we have that $\Nn _i$ is
obtained from $\Nn _{i-1}$ by attaching a cell $\alpha _i\subset \Nn
_{i}(X_i,Y_i)$ and taking the free category over $\Nn _{i-1}$ and $\alpha
_i$. (Precisions ???) We can suppose that the whole boundary $\partial
\alpha _i$ is contained in $\Nn _{i-1}$ (by treating one skeleton at a
time, except for the stuff generated by smaller skeleta ...).

This definition is basically the same as the corresponding notion for
simplicially enriched categories underlying the approach of Dwyer-Kan to
homotopy coherence cf \cite{DwyerKan1} \cite{DwyerKan2} \cite{DwyerKan3}
\cite{DwyerKanEquivs}.

Suppose $\Mm \rightarrow \Xx$ and $\Nn \rightarrow \Xx$ are functors to a
discrete category $\Xx$, inducing isomorphisms on the set of objects, and
suppose $\Nn \rightarrow \Mm$ is a continuous functor compatible with the
projection to $\Xx$.
Suppose that for each morphism $\phi$ in $\Xx$ the inverse image $\Mm (\phi
)$ is weakly contractible. Suppose that $\Nn$ is CW-free. Suppose $\Pp
\rightarrow \Mm$ is an isomorphism on the set of objects and induces
fibrations on the morphism spaces. Let $Sect (\Nn , \Pp /\Mm )$ denote the
space of morphisms from $\Nn$ to $\Pp$ yielding the same composition $\Nn
\rightarrow \Mm$, and let $Hom _{\Xx}(\Nn ,\Pp )$ denote the space of
morphisms commuting with the projections to $\Xx $.

\begin{lemma}
\mylabel{flex04}
In the above situation, the map
$$
Sect (\Nn, \Pp /\Mm )\rightarrow Hom
_{\Xx}(\Nn ,\Pp )
$$
is a weak homotopy equivalence.
\end{lemma}
{\em Proof:}
We can form the fiber product $\Qq = \Pp \times _\Xx \Mm$ (since $\Xx$ is
discrete, this is the same as the path fiber product). The map $\Qq
\rightarrow \Mm$ is a fibration, and we have $$
Sect (\Nn ,\Qq /\Mm ) = Hom _{\Xx }(\Nn , \Pp ). $$
There is a morphism $\Pp \rightarrow \Qq$ of fibrations over $\Mm$ sending
$r$ to $(r, \varphi (r))$ where $\varphi (r)$ is the image of $r$ in $\Mm$.
This induces the map
$$
Sect (\Nn, \Pp /\Mm )\rightarrow Sect (\Nn ,\Qq /\Nn )=Hom _{\Xx}(\Nn ,\Pp ).
$$
Thus it suffices to prove that $\Pp \rightarrow \Qq$ induces an equivalence
on the space of sections. Note first of all that for any $\phi
:X\rightarrow X'$ in $\Xx $ we have that $\Pp (\phi )\rightarrow \Qq (\phi
) = \Pp (\phi )\times \Mm (\phi )$ is a weak equivalence, and hence a weak
equivalence in each fiber. Then we have a
\newline
{\em Sublemma:} If $\Nn$ is CW-free and $\Pp \rightarrow \Qq$ is a morphism
of fibrations over $\Mm$ inducing a weak equivalence in each fiber, then
$Sect (\Nn ,\Pp /\Mm )\rightarrow Sect (\Nn ,\Qq /\Mm )$ is a weak
equivalence.
\newline
The sublemma evidently will complete the proof. We prove the sublemma by
transfinite induction on the sequence of subcategories involved in the
definition of CW-free. For each $j$ let $\Pp _j\rightarrow \Nn _j$ (resp.
$\Qq_j\rightarrow \Nn _j $) be the fibration pulled back from $\Mm$ to
$\Nn$ and restricted to $\Nn _j$. Let $i$ be the first index where $Sect
(\Pp _i/\Nn _i)\rightarrow Sect (\Qq _i/\Nn _i )$ is not a weak
equivalence. Note that $Sect (\Pp _{i-1}/\Nn _{i-1})$ is a projective limit
of a sequence of fibrations (we shall see below that the maps are
fibrations); and projective limit of fibrations commutes with taking
homotopy groups, so the morphism $$
Sect (\Pp _{i-1}/\Nn _{i-1})\rightarrow Sect (\Qq _{i-1}/\Nn _{i-1}) $$
is a weak equivalence. But now the fact that $\Nn _i$ is freely generated
over $\Nn _{i-1}$ by the cell $\alpha _i$ implies that the map $$
Sect (\Pp _i/\Nn _i )\rightarrow Sect (\Pp _{i-1}/\Nn _{i-1}) $$
is a fibration, and the fiber over $\sigma _{i-1}$ is the space of sections
of the pullback of $\Pp $ over the cell $\alpha_i$, compatible with $\sigma
_i$ on the boundary $\partial (\alpha _i )$. The map between the fibers for
$\Pp $ and $\Qq$ is equal to the map induced by the map of fibrations
pulled back to $\alpha$. Since this map of fibrations is a weak equivalence
on the fiber, the map between the space of sections is a weak equivalence.
Now we have a morphism of fibrations $$
\begin{array}{ccc}
Sect (\Pp _i/\Nn _i )&\rightarrow &Sect (\Pp _{i-1}/\Nn _{i-1})\\
\downarrow & & \downarrow \\
Sect (\Pp _i/\Nn _i )&\rightarrow &Sect (\Pp _{i-1}/\Nn _{i-1}) \end{array}
$$
which is a weak equivalence on the base and on the fiber. From the long
exact homotopy sequence and the $5$-lemma, it is a weak equivalence on the
total space, contradicting our inductive hypothesis. This proves the
sublemma, and hence the lemma. \eop

Suppose $\Nn$ is CW-free, and that we have a functor $\Nn \rightarrow \Xx$
where $\Xx$ is a discrete category, an isomorphism on objects, and with
weakly contractible $\Nn (\phi )$ for each morphism $\phi$ in $\Xx$.
Suppose $\Nn ' \subset \Nn$ is a CW-free cellular subcategory. Suppose
$\Cc$ is a continuous category. We use the Postnikov tower of continuous
categories $\tau _{\leq n}\Cc$, to analyze $Hom (\Nn , \Cc )\rightarrow Hom
(\Nn ' , \Cc )$.

First of all, we note that $Hom (\Nn , \Cc )\rightarrow Hom (\Nn ', \Cc )$
is a fibration. This is because $\Nn$ can be obtained from $\Nn '$ by a
(transfinite) sequence of additions of cells; thus $\Nn '$ appears as one
of the $\Nn _j$ in the above argument. We have seen that $Hom (\Nn ,\Cc
)\rightarrow Hom (\Nn _j , \Cc )$ is a fibration.

We have a diagram
$$
\begin{array}{ccc}
Hom (\Nn ,\tau _{\infty} \Cc ) &\rightarrow & Hom (\Nn ,\Cc ) \\ \downarrow
&& \downarrow \\
Hom (\Nn ',\tau _{\infty} \Cc ) &\rightarrow & Hom (\Nn ',\Cc ) \end{array}
$$
where the vertical arrows are fibrations and the horizontal arrows are weak
equivalences (because the homotopy groups can be calculated by looking at
an inductive process of adding cells to obtain $\Nn$ or $\Nn '$, and $\tau
_{\infty}\Cc \rightarrow \Cc$ is a weak equivalence). In this situation the
fibers of the vertical fibrations are weakly equivalent. We would like to
analyse the fiber over an element $G\in Hom (\Nn ', \tau _{\infty} \Cc )$.
Let $Hom (\Nn , \tau _{\infty}\Cc ; G)$ denote the space of morphisms
restricting to $F$ on $\Nn '$. The composition with $\tau_{\infty} \Cc
\rightarrow \tau _{\leq n}\Cc $ induces elements which we also denote by
$G\in Hom (\Nn , \tau _{\leq n} \Cc ; G)$.

Note that all of the morphisms in the Postnikov tower for $\Cc$ induce
isomorphisms on the set of objects.

The homotopy fiber of
$$
Hom (\Nn ,\tau _{\leq n}\Cc ;G)\rightarrow Hom (\Nn ,\tau _{\leq n-1}\Cc ;G) $$
over $F$
is equal to the space of sections of the fibration $\Pp _{n}\rightarrow
\Nn$ given by path fiber product
$$
\Pp _{n}:= \Nn \times ^{Path}_{\tau _{\leq n-1}\Cc }\tau _{\leq n}\Cc , $$
and equal to the given section $G$ over $\Nn '$.

By the previous lemma, we have
$$
Sect (\Pp _n/\Nn )=Hom _{\Xx} (\Nn ,\Pp _{n}), $$
and
$$
Sect (\Nn ', \Pp _n/\Nn ) = Hom _{\Xx} (\Nn ',\Pp _{n}). $$
Furthermore, the square made up of these two maps and the restrictions,
commutes. Thus the fiber of
$$
Sect (\Pp _n/\Nn )\rightarrow Sect (\Nn ', \Pp _n /\Nn ) $$
is equal to the fiber of
$$
Hom _{\Xx} (\Nn ,\Pp _{n})\rightarrow Hom _{\Xx} (\Nn ',\Pp _{n}) $$
over the corresponding element (which we also denote by $G$).

On the other hand, the fiber of $\Pp _{n}(X,X')$ over $r\in \Nn (X,X' )$ is
equal to the homotopy fiber of $\tau _{\leq n}\Cc (FX,FX')\rightarrow \tau
_{\leq n-1}\Cc (FX,FX' )$ over $F(r)$. This homotopy fiber is a $K(G(\phi
), n)$ where $G(\phi )$ is a group depending only on the morphism $\phi$ in
$\Xx$ which is the image of $r$. Since $\Nn (\phi )$ is weakly
contractible, $\Pp _{n}(\phi )$ is weakly a $K(G(\phi ), n)$. We have
obtained the following situation:

{\em The homotopy fiber of
$$
Hom (\Nn ,\tau _{\leq n}\Cc ;G)\rightarrow Hom (\Nn ,\tau _{\leq n-1}\Cc ;G) $$
over $F$
is equal to the fiber of
$$
Hom _{\Xx}(\Nn ,\Pp )\rightarrow Hom _{\Xx} (\Nn ',\Pp ) $$
over $G$,
where $\Pp \rightarrow \Xx$
is a morphism of spacial categories such that $\Pp (\phi )$ is weakly a
$K(G(\phi ),n)$ for each $\phi \in \Xx$. }

(We say that $\Pp$ is an {\em Eilenberg-MacLane category} over $\Xx$ in the
above situation.)

In the next subsection we will investigate conditions for the map $$
Hom _{\Xx}(\Nn ,\Pp )\rightarrow Hom _{\Xx} (\Nn ',\Pp ) $$
to be a weak equivalence; note that it is a weak equivalence if and only if
the fiber is nonempty and weakly contractible.

\subnumero{Functors to Eilenberg-MacLane categories}

Suppose that $\Nn$ is a CW-free continuous category with a functor to a
discrete category $\Xx$; and suppose that $\Pp \rightarrow \Xx$ is an
Eilenberg-MacLane category over $\Xx$. Suppose for now that $\Pp (\phi )=
K(G(\phi ), n)$ for $n\geq 2$. The groups $G(\phi )$ are then well defined,
and form a system (see below).

Let $\Sigma _k(\phi )$ denote the set of primitive (cubic) cells in degree
$k$ in $\Nn (\phi )$. For $\sigma \in \Sigma _k(\phi )$ can write $$
\partial \sigma = \sum _i(\prod _j\alpha _i^j(\sigma ))\beta _i (\sigma )
(\prod _k\gamma _i^k(\sigma )) +\sum _l \delta _l(\sigma ) $$
where $\beta _i(\sigma )\in \Sigma _{k-1}(\psi _i(\sigma ))$ and $\alpha
^j_i(\sigma )\in \Sigma _0(\xi _i^j(\sigma ))$ and $\gamma ^k_i(\sigma )\in
\Sigma _0(\zeta ^k_i(\sigma ))$, whereas the $\delta _l(\sigma )$ are
decomposable into products of at least two positive-degree pieces.

Let $\Nn ^k$ denote the subcategory generated by all cells in degrees $\leq
k$. We have a fibration
$$
Hom _{\Xx}(\Nn ^k , \Pp )\rightarrow Hom _{\Xx}(\Nn ^{k-1},\Pp ), $$
and the fiber $Fib _k(F_{k-1})$ over $F_{k-1}$ is the space of ways of
mapping each cell $\sigma \in \Sigma _k(\phi )$ into $\Pp (\phi )$
coinciding with the map $F_{k-1}$ on the boundary.
If $k\geq n+2$ then the space of ways of mapping $\sigma$ into $\Pp (\phi
)$ relative to the boundary, is weakly contractible. If $k\leq n$ then
there exists
a way of mapping each $\sigma$ coinciding with the given map on the
boundary (since $\pi _{k-1}(F(\phi ))$ is trivial in that case), and the
space of such maps for a given $\sigma$ is a $K(G(\phi ), n-k)$ (it is a
principal homogeneous space over $G(\phi )$ for $n=k$). If $k=n+1$ the
space is nonempty and weakly contractible if a certain obstruction in
$G(\phi )$ vanishes, otherwise it is empty.

We obtain:
$$
Fib _{n+1}(F_{n})=\left\{
\begin{array}{ll}
\emptyset & obs (F_n )\neq 0 \,\, in \,\, \prod _{\sigma \in \Sigma
_{n+1}(\phi )}G(\phi )\\
\ast & obs (F_n ) =0 \,\, in \,\, \prod _{\sigma \in \Sigma _{n+1}(\phi
)}G(\phi )
\end{array}
\right.
$$
and
$$
Fib _k(F_{k-1})= \prod _{\sigma \in \Sigma _{k}(\phi )}K(G(\phi ), n-k) =
K(\prod _{\sigma \in \Sigma _{k}(\phi )}G(\phi ), n-k)
$$
for $k\leq n$.

In such a situation (a sort of reverse Postnikov tower of fibrations, where
the homotopy degrees decrease in the tower instead of increasing), the
homotopy groups of the total space are calculated by a complex whose
elements are the homotopy groups of the components. We treat this generally
for the moment (note that the discussion below is also a special case of
the spectral sequence for towers discussed in \S 5).

Suppose we have a sequence of fibrations $U_k\rightarrow U_{k-1}$, with
nonempty weakly contractible fiber for $k\geq n-2$. The limit $U$ is then
weakly equal to $U_{n+1}$.
Suppose $u\in U_n$. Suppose that the fiber of $U_k\rightarrow U_{k-1}$ over
the image of $u$, is a $K(G_k,n-k)$ (and here, we can fix the groups $G_k$
because we have the basepoint $u$ for the fiber). The long exact homotopy
sequence (dropping the basepoint, which will always be the image of $u$)
gives $$
\pi _i (U_k)=\pi _i (U_{k-1})
$$
unless $i,i-1 = n-k$, that is unless $i=n-k$ or $i=n-k+1$. We obtain $\pi
_i(U_n)= \pi _i (U_{n-i+1})$, and $\pi _i(U_{n-i-1})=0$. We have
$$
0\rightarrow \pi _i (U_{n-i+1})\rightarrow \pi _i(U_{n-i})\rightarrow
G_{n-i+1} $$
and
$$
\pi _{i+1}(U_{n-i-1})\rightarrow
G_{n-i}\rightarrow \pi _i(U_{n-i})\rightarrow \pi _i(U_{n-i-1})=0 . $$
The compositions
$$
G_{n-i}\rightarrow \pi _i(U_{n-i})\rightarrow G_{n-i+1} $$
give the morphisms for a complex, and from the two exact sequences above,
we have that $\pi _i (U_n)$ is the cohomology of this complex. The complex
is $$
G_{\cdot}(u)=G_0\rightarrow G_1\rightarrow \ldots \rightarrow
G_{n-1}\rightarrow G_n,
$$
and
$$
\pi _i(U,u)= H^{n-i}(G_{\cdot}(u)).
$$
If, furthermore, the obstruction for the existence of $\tilde{u}$ in
$U_{n+1}$ mapping to $u\in U_n$ is an element of a group $G_{n+1}$ then we
can think of $G_{n+1}$ as being the next term in the complex; this is not
quite rigorous because the complex depends on the choice of $u$.

The part of the complex
$$
G_0\rightarrow \ldots \rightarrow G_{n-1} $$
depends only on the choice of $u\in U_{n-1}$ and $U_{n-1}$ is contractible,
so this part of the complex is fixed. In general, $G_{n}$ could be only a
set, and $G_{n-1}$ a nonabelian group. The obstruction is a boolean
variable. The group $G_{n-1}$ acts on $G_n$, preserving the boolean value
of the obstruction.

Return now to our previous situation, and suppose that $n\geq 2$. Then the
groups $G(\phi )$ are well defined and abelian, and for two composable
morphisms $\phi$ and $\psi$ we have
$$
\begin{array}{ccc}
G(\phi )&\rightarrow &G(\phi \psi )\\
a&\mapsto & a\psi
\end{array}
$$
and
$$
\begin{array}{ccc}
G(\psi )&\rightarrow &G(\phi \psi )\\
b&\mapsto & \phi b
\end{array}
$$
giving
$$
\begin{array}{ccc}
G(\phi )\times G(\psi )&\rightarrow &G(\phi \psi )\\ (a,b)&\mapsto & a\psi
+ \phi b.
\end{array}
$$
The actions $a\mapsto a\psi$ and $b\mapsto \phi b$ are compatible with the
multiplication in $\Xx$ and commute.

We obtain a complex
$CP_{\cdot}$ of abelian groups defined by $$
PC_k=\prod _{\sigma \in \Sigma
_{k}(\phi )}G(\phi ).
$$
The end of this complex goes into degree $n+1$ (and past, even), and the
action of $PC_{n-1}$ on $PC_n$ is by translation via the morphism
$PC_{n-1}\rightarrow PC_n$, and the boolean variable is determined by the
morphism $PC_n\rightarrow PC_{n+1}$.

The differentials are given (in terms of the formula for the boundary of a
cell $\sigma$) by $$
(dg)_{\sigma} =
$$
$$
\sum _i(\prod _j\xi _i^j(\sigma ))g_{\beta _i (\sigma )} (\prod _k\zeta
^k_i(\sigma )).
$$
(With appropriate signs.)
For the one-dimensional cells, note that the boundary is a sum of products
of zero dimensional cells. Corresponding to a certain product here, there
is included above one term for each term in the product.

The proofs of the statements claimed in the two previous paragraphs are
left to the reader.

If $\Nn '\subset \Nn$ is a cellular subcategory with the same properties,
then there is a morphism of complexes
$$
PC_{\cdot}(\Nn , G(\cdot ))\rightarrow PC_{\cdot}(\Nn ',G(\cdot )). $$
If this induces an isomorphism on cohomology, then the map $$
Hom _{\Xx}(\Nn ,\Pp )\rightarrow Hom _{\Xx }(\Nn ',\Pp ) $$
is a weak equivalence.

Leaving aside the cases $n=0$ and $n=1$ for the moment, if $\Nn '\subset
\Nn $ induces isomorphisms on cohomology of the complexes $PC_{\cdot}$ for
all systems of abelian groups $G(\phi )$ as above, then it induces weak
equivalences of the homotopy fibers of $$
Hom (\Nn ,\tau _{\leq n}\Cc )\rightarrow Hom (\Nn , \tau _{\leq n-1}\Cc ), $$
and so by induction (except we haven't yet treated the initial cases), we
get that for any $n$,
$$
Hom (\Nn , \tau _{\leq n}\Cc )\rightarrow Hom (\Nn ', \tau _{\leq n}\Cc ) $$
is a weak equivalence.
{\em If, furthermore, $\Cc$ is $n_0$-truncated}, then the fact that a
morphism of continuous categories inducing weak equivalences on the $Hom$
spaces, induces a weak equivalence on morphisms from $\Nn$ (which can be
seen by looking at the inductive construction of $\Nn$ obtained by adding
primitive cells one after the other), implies that $$ Hom (\Nn , \Cc
)\leftarrow Hom (\Nn , \tau _{\infty}\Cc )\rightarrow Hom (\Nn , \tau
_{\leq n}\Cc )
$$
are weak equivalences, and hence in our situation $$
Hom (\Nn ,\Cc )\rightarrow Hom (\Nn ',\Cc ) $$
is a weak equivalence.

Suppose now that $\Nn$ is a cellular subcategory of $\Mm$, the free spacial
category over $\Xx$. Thus $\Nn$ is determined by the subsets $$
\Sigma _k (\Nn )\subset \Sigma _k(\Mm )
$$
where $\Sigma _k(\Mm )$ is the set of all composable $k+1$-tuples of
morphisms in $\Xx$. These subsets are subject to some conditions to insure
that $\Nn$ is generated by these primitive cells.

We can denote a cell by its composable sequence $(\phi _0,\ldots , \phi
_k)$. The primitive cells in its boundary are
$$
(\phi _0,\ldots , \widehat{\phi _i},\ldots , \phi _k), $$
$$
\phi _0 (\phi _1,\ldots , \phi _k),
$$
and
$$
(\phi _0,\ldots , \phi _{k-1})\phi _k .
$$
In this case,
$$
PC_k(\Nn , G(\cdot ))=
$$
$$
\prod _{(\phi _0,\ldots , \phi _k)\in \Sigma _k(\Nn )}G(\phi _0\cdots \phi
_k). $$
An element will be denoted by a function $$
g(\phi _0,\ldots , \phi _k),
$$
and the differential becomes
\begin{eqnarray*}
(dg)(\phi _0,\ldots , \phi _{k+1}) &=& \phi _0g(\phi _1,\ldots ,\phi
_{k+1})\\ && + (-1)^{k+1}g(\phi _0,\ldots , \phi _k)\phi _{k+1} \\ && +
\sum _{i=1}^{k}(-1)^ig(\phi _0,\ldots , \widehat{\phi _i},\ldots , \phi
_{k+1}).
\end{eqnarray*}
We will interpret this in a more abstract way.

For each morphism $\phi$ in
$\Xx$, let $\Lambda_{\Nn }( \phi )$ denote the simplicial complex whose
simplices in degree $k$ are expressions of the form $\alpha (\phi _1,\ldots
, \phi _k)\beta$ such that $\alpha \phi _1\cdots \phi _k\beta = \phi$, with
$(\phi _1,\ldots , \phi _k)$ a primitive cell in $\Nn$ and $\alpha$ and
$\beta$ morphisms in $\Xx$. The boundary simplices are given by forgetting
elements $\phi _i$; if $\phi _1$ or $\phi _k$ are forgotten then they are
put into $\alpha$ or $\beta$. If $u,\phi , v$ is a composable triple, then
we obtain a morphism of simplicial complexes $$
\alpha (\phi _1,\ldots , \phi _k)\beta \mapsto u\alpha (\phi _1,\ldots ,
\phi _k)\beta v. $$
Let $Fl(\Xx )$ denote the category of morphisms of $\Xx$, where a morphism
from $\psi $ to $\phi$ is an expression $\psi = u\phi v$. Then $$
\phi \mapsto \Lambda _{\Nn} (\phi )
$$
is a contravariant functor from $Fl (\Xx )$ to the category of simplicial
sets. On the other hand, our system of groups $G(\cdot )$ is just a
contravariant functor from $Fl (\Xx )$ to the category of abelian groups.

We have
$$
PC_{\cdot}(\Nn ,G(\cdot ))= Tot
(Hom_{\Xx }(\Lambda _{\Nn},G)))(\pm 1).
$$
Here $Hom_{\Xx}(\Lambda _{\Nn},G))$ is the simplicial abelian group
obtained by taking the space of morphisms of presheaves in each degree. The
operation $Tot$ is just the operation of combining the face maps into a
differential to obtain the total complex. This identification comes about
because an element of $PC_{\cdot}(\Nn ,G(\cdot ))$ is a function $g(\phi
_0,\ldots , \phi _k)$ but which may be considered as a function $g(\alpha
(\phi _0,\ldots , \phi _k)\beta )$ on $\Lambda _{\Nn}$ satisfying the rule
$$ g(u\alpha (\phi _0,\ldots , \phi _k)\beta v)= ug(\alpha (\phi _0,\ldots
, \phi _k)\beta )v $$
This rule gives
$$
g(\alpha (\phi _0,\ldots , \phi _k)\beta )= \alpha g(\phi _0,\ldots , \phi
_k)\beta . $$
Note the shift by one, between $\Lambda _{\Nn}$ and $PC_{\cdot}$ (I never
know which direction for these shifts, thus indicated by $\pm 1$!!!).

\begin{lemma}
\mylabel{flex05}
Suppose $\Nn'$ is a cellular subcategory of $\Nn '$ with the required
properties discussed above. If, for each $\phi \in Fl (\Xx )$, the morphism
(inclusion) of simplicial complexes
$$
\Lambda _{\Nn '}(\phi )\rightarrow \Lambda _{\Nn}(\phi ) $$
induces an isomorphism on cohomology with integral coefficients, then for
any system of abelian groups $G(\cdot )$, $$
PC_{\cdot}(\Nn , G(\cdot )) \rightarrow
PC_{\cdot}(\Nn ', G(\cdot ))
$$
is an isomorphism on cohomology, and consequently for any $\Pp \rightarrow
\Xx$ with Eilenberg-Maclane fibers $K(G(\phi ), n)$ for $n\geq 2$, we have
that $$
Hom _{\Xx }(\Nn ,\Pp )\rightarrow Hom _{\Xx }(\Nn ' , \Pp ) $$
is a weak equivalence.
\end{lemma}
{\em Proof:}
The set of morphisms between presheaves may be interpreted as the global
sections of a presheaf of morphisms,
$$
Hom _{\Xx} (U,V)=\Gamma (Hom (U,V)).
$$
For any $G$ we can choose a
resolution by presheaves of the
form $\prod _{\eta}Ind(H^i_{\eta},\eta )$ where each $H^i_{\eta}$ is an
injective group and
$$
Ind (H^i_{\eta},\eta )(\psi ):=
\prod _{(\alpha , \beta ): \eta \rightarrow \psi}H^i_{\eta}. $$
Applying the functors $Tot (Hom _{\Xx} (\Lambda _{\Nn},\cdot ))$ or $Tot
(Hom _{\Xx} (\Lambda _{\Nn '},\cdot ))$, we obtain a morphism of double
complexes corresponding to the inclusion of $\Nn '$ in $\Nn$. In the
columns, we have the morphisms of the total complexes $$
Tot (Hom _{\Xx} (\Lambda _{\Nn}, \prod _{\eta}Ind (H_{\eta},\eta ))) \rightarrow
$$
$$
Tot (Hom _{\Xx} (\Lambda _{\Nn '}, \prod _{\eta}Ind (H_{\eta},\eta ))). $$
But
$$
Hom _{\Xx} (\Lambda _{\Nn }, Ind (H,\eta )) =
$$
$$
Hom (\Lambda _{\Nn }(\eta ),H)
$$
and the same for $\Nn '$. Our morphism in the $i$th column is therefore a
product of morphisms of complexes of the form $$
Hom (\Lambda _{\Nn }(\eta ),H)\rightarrow Hom (\Lambda _{\Nn '}(\eta ),H);
$$
by hypothesis these induce isomorphisms on cohomology. Thus the product is
an isomorphism on cohomology so in the cohomology of the columns of our
double complexes, the morphism induces an isomorphism in cohomology.
Therefore the morphism between double complexes is an isomorphism in
cohomology. On the other hand, take the cohomology in the rows. In the
$k$th row we have the functor $Hom (\Lambda _{\Nn ,k},\cdot )$ applied to
the resolution of $G$. We claim that this functor is exact in the argument
$G$.

The functor $Hom _{\Xx} (\Lambda _{\Nn ,k},\cdot )$ decomposes as a product
over the elementary cells denoted $(\phi _1,\ldots , \phi _k)\in \Nn$ (with
the composition denoted $\phi := \phi _1\cdots \phi _k$): $$
Hom _{\Xx} (\Lambda _{\Nn ,k},G)= \prod _{(\phi_1,\ldots , \phi _k)\in \Nn}
G(\phi ).
$$
Each factor in the product is exact in $G$, so the product is exact in $G$.

Now in the
cohomology of the row we obtain only $Hom (\Lambda _{\Nn ,k},G )$ in the
first column. From this follows that the cohomology of the double complex
is the same as the cohomology of the complex for the presheaf $G$; and the
same being true for $\Nn '$, we can conclude that we obtain an isomorphism
in cohomology of the complexes for the presheaf $G$ (that is, the statement
of the lemma---note that the second part has been discussed previously).
\eop

\subnumero{Specifying flexible functors on subcategories}

We apply the above analysis to a more concrete situation, which arises in
practice. Suppose $\{ \Xx _{\alpha}\} _{\alpha \in A}$ is a family of
subcategories of a category $\Xx$. Suppose that $$
Ob (\Xx ) = \bigcup _{\alpha \in A} Ob (\Xx _{\alpha}). $$
Suppose that $T^{\alpha}:\Xx _{\alpha}\rightarrow \Cc $ are flexible
functors to a continuous category $\Cc$, which {\em agree exactly} on the
intersections $\Xx _{\alpha \beta}:=\Xx _{\alpha}\cap \Xx _{\beta}$ (the
intersection is the subcategory obtained by taking the intersections of the
object and morphism sets). We would like to investigate the space of
flexible functors $T:\Xx \rightarrow \Cc$ which restrict to the
$T^{\alpha}$ on the $\Xx _{\alpha}$. Note that this space is well defined,
since the objects $T_X$ are already fixed for all $X\in \Xx$ (due to the
fact that any $X$ appears in at least one $\Xx _{\alpha}$). The space of
such $T$ is thus a subset of a product of spaces of maps between the fixed
objects.

Put $\Xx _{\alpha _1,\ldots , \alpha _n}$ be the intersection of $\Xx
_{\alpha _1},\ldots , \Xx _{\alpha _n}$. Let
$\Nn _{\alpha _1,\ldots , \alpha _n}$ denote the free cellular category on
$\Xx _{\alpha _1,\ldots , \alpha _n}$ (considered as a subcategory of $\Mm
= M\Xx$). Let $\Nn$ be the subcategory of $\Mm$ generated by the $\Nn
_{\alpha}$.

The set of primitive cells of $\Nn$ is the union of the sets of primitive
cells of the $\Nn _{\alpha}$.
More precisely,
let $\Sigma ^{\alpha_1,\ldots , \alpha _n}(\phi )$ denote the complex of
primitive cells in $\Nn _{\alpha _1,\ldots , \alpha _n}(\phi )$ for $\phi
\in \Xx _{\alpha _1,\ldots , \alpha _n}$.
We form a complex of presheaves of abelian groups on $Fl (\Xx )$, $$
\zz\Lambda ^{\alpha_1,\ldots , \alpha _n}) $$
whose value on a morphism $\psi$ in $\Xx$ is freely generated by the
symbols $f^{\ast}(\sigma )$ where $f:\psi \rightarrow \phi $ is a morphism
in $Fl(\Xx )$ (for an arrow $\phi$ in $\Xx _{\alpha_1,\ldots , \alpha _n}$)
and $\sigma \in \Sigma ^{\alpha_1,\ldots , \alpha _n}(\phi )$.

There is a long exact sequence of complexes of presheaves of abelian groups, $$
\ldots \rightarrow
\bigoplus _{\alpha _1,\ldots ,\alpha _n} \zz\Lambda ^{\alpha_1,\ldots ,
\alpha _n}) \rightarrow $$
$$
\ldots \rightarrow \bigoplus _{\alpha} \zz\Lambda ^{\alpha} \rightarrow \zz
\Lambda _{\Nn} ,
$$
where $\zz \Lambda _{\Nn}$ denotes the complex of sheaves of free abelian
groups generated by the $\Lambda _{\Nn}$.

(Again the proof is left to the reader---the last term is because the
primitive cells of $\Nn$ are the union of the primitive cells of the $\Nn
_{\alpha}$ ...).

Our complex $PC_{\cdot}(\Nn , G(\cdot ))$ is just $Hom _{Fl (\Xx )}(\zz
\Lambda _{\Nn}, G)$, which is the same as
$$
Ext _{Fl (\Xx )}(\zz \Lambda
_{\Nn}, G)
$$
since the complex $\zz \Lambda _{\Nn}$ is co-induced. From the
above long exact sequence, we get a spectral sequence $$
E_1^{i,n} (???)=\bigoplus _{\alpha _1,\ldots , \alpha _n} Ext ^i(\zz
\Lambda ^{\alpha _1,\ldots , \alpha _n}, G) \Rightarrow Ext ^{i+n}(\zz
\Lambda
_{\Nn}, G).
$$
Note finally that
$$
Hom _{Fl (\Xx )}(\zz \Lambda ^{\alpha _1,\ldots , \alpha _n},G) $$
$$
=
Hom _{Fl (\Xx _{\alpha _1,\ldots , \alpha _n})}(\zz \Lambda _{M \Xx
_{\alpha _1,\ldots , \alpha _n}},G|_{Fl(\Xx _{\alpha _1,\ldots , \alpha
_n})})
$$
and recall that $\zz \Lambda _{M
\Xx _{\alpha _1,\ldots , \alpha _n}}$ associates to each arrow $\phi$ the
free abelian group complex associated to the nerve of the category of
objects inside $\phi$.

We treat the case where there is only one index $\alpha$. Then our spectral
sequence is trivial, becoming just
$$
Ext ^{i}_{Fl (\Xx )}(\zz \Lambda
_{\Nn}, G)= Ext ^i_{Fl(\Xx )}(\zz \Lambda _{M\Xx _{\alpha}},G|_{Fl (\Xx
_{\alpha})}).
$$

For an arrow $\phi$ in $\Xx$, we can define a category $Ins _{\Xx
_{\alpha}/\Xx}(\phi )$ to be the category whose objects are objects of $\Xx
_{\alpha}$ lying inside $\phi$ in the sense of $\Xx $, and whose morphisms
are those coming from morphisms in $\Xx _{\alpha}$ between the objects. The
criterion
of Lemma \ref{flex05} is that if, for any morphism $\phi$ of $\Xx$, the
inclusion
of nerves $$
N\, Ins _{\Xx _{\alpha}/\Xx}(\phi )\subset N\, Ins_{\Xx } (\phi ) $$
induces an isomorphism on cohomology with abelian group coefficients, then
the fiber of
$$
Hom _{\Xx }(\Mm ,\Pp )\rightarrow Hom _{\Xx}(\Nn ,\Pp ). $$
is weakly contractible.

\subnumero{The case $n=0$}

For any continuous category $\Cc$, the
construction above with the functor $\pi _0:Top \rightarrow Set$ gives a
discrete category $\pi _0 \Cc$ whose sets of morphisms are the sets of path
components of the morphism sets of $\Cc$. There is a natural functor $\Cc
\rightarrow \pi _0 \Cc$ inducing an isomorphism on object sets. For a
morphism $\phi $ in $\pi _0 \Cc$, the inverse image $\Cc (\phi )$ is just
the path-component corresponding to $\phi$ in the morphisms of $\Cc$.

We say that a map $A\rightarrow B$ between topological spaces is {\em
weakly continuous} if, for any CW complex $Z$ and continuous map
$Z\rightarrow A$, the composed map $Z\rightarrow B$ is continuous. A weakly
continuous map induces maps on homotopy groups. We say that it is a weak
equivalence if it induces isomorphisms on homotopy groups.

\begin{lemma}
\mylabel{flex06}
Suppose $\Cc$ is a $0$-truncated continuous category (that is the path
components of the morphism spaces are weakly contractible). Let $\Mm$ be a
CW-free continuous category. Then the functor $\pi _0$ induces a map $$ Hom
(\Mm ,\Cc )\rightarrow Hom (\pi _0 \Mm , \pi _0 \Cc ), $$ the $Hom$ on the
left denoting the space of continuous functors. This map is a weakly
continuous weak equivalence if the $Hom$ on the right is given the discrete
topology.
\end{lemma}
{\em Proof:}
It is clearly weakly continuous. We have to show that the fibers are
nonempty and weakly contractible. We do this by transfinite induction on
the generating cells of $\Mm$, as usual. Fix an element $$
F\in Hom (\Mm
,\Cc )\rightarrow Hom (\pi _0 \Mm , \pi _0 \Cc ), $$
and for any subcategory $\Nn\subset \Mm$ let $Hom ^F(\Nn ,\Cc )$ denote the
space of morphisms $G$ such that for $\phi\in \Nn$, $G(\phi )$ is in the
path component which is the image by $F$ of the path component of $\phi$.
Express $\Mm $ as an increasing union of subcategories $\Mm _j$ such that
$\Mm _j$ is generated over $\Mm _{<j}$ by a single cell $\alpha _j$. Then
$$ Hom ^F(\Mm _j, \Cc )\rightarrow Hom ^F(\Mm _{<j}, \Cc ) $$
is a fibration with weakly contractible fiber. Furthermore, $Hom ^F(\Mm
_{<j}, \Cc )$ is an inverse limit of fibrations with weakly contractible
fiber, over any $Hom ^F(\Mm _{i}, \Cc )$ for $i<j$; thus $$
Hom ^F(\Mm _j, \Cc )\rightarrow Hom ^F(\Mm _{i}, \Cc ) $$
is a fibration with weakly contractible fiber for any $i<j$. Finally note
that for the first value of $i$, $ Hom ^F(\Mm _{i}, \Cc )$ is weakly
contractible. Thus $Hom ^F(\Mm ,\Cc )$ is weakly contractible as desired.
\eop

\begin{corollary}
\mylabel{flex07}
Suppose $\Cc$ is $0$-truncated and $\Nn \subset \Mm$ is a cellular
subcategory. The restriction morphism
$$
Hom (\Mm ,\Cc )\rightarrow Hom (\Nn ,\Cc ) $$
has fibers whose path-components are weakly contractible. The
path-components in the fiber over $F$ correspond to the functors $\pi _0\Mm
\rightarrow \pi _0 \Cc $ extending $\pi _0 F$. In particular, if there
exists a unique such extension, then the fiber of the restriction morphism
is weakly contractible. \end{corollary}
{\em Proof:} Note that the restriction morphism is a fibration, by our
standard arguments. To find the weak type of the fiber, we can compare the
weak type of the two spaces.
\eop

Suppose $\Mm$ is the free continuous category over $\Xx$ and $\Nn$ is the
cellular subcategory generated by a family of subcategories $\Xx
_{\alpha}\subset \Xx$. Then $\pi _0\Nn $ is the {\em free category
generated by the $\Xx _{\alpha}$, amalgamated over $\Xx _{\alpha \beta}$}.
We explain this more precisely. The object set of $\pi _0\Nn$ is the union
of the object sets of $\Xx _{\alpha}$. In order to determine $\pi _0 \Nn$,
it suffices to consider the $0$-cells and $1$-cells. The $0$-cells of $\Nn$
are formal products of composable sequences of morphisms each coming from
one of the $\Xx _{\alpha}$. There is a $1$-cell linking $\phi _1\cdots \phi
_n$ to $\phi _1\cdots (\phi _i\phi _{i+1})\cdots \phi _n$ whenever $\phi
_i$ and $\phi _{i+1}$ are in the same subcategory $\Xx _{\alpha}$. The
morphisms in $\pi _0\Nn$ are equivalence classes of formal products of
morphisms coming from the $\Xx _{\alpha}$, modulo the equivalence relation
generated by the $1$-cells described above.

\begin{lemma}
\mylabel{flex08}
In the above situation, if $\Aa$ is any category, a functor $F:\pi _0\Nn
\rightarrow \Aa$ is the same thing as a family of functors $F_{\alpha}:\Xx
_{\alpha}\rightarrow \Aa$ such that $F_{\alpha}|_{\Xx _{\alpha \beta}}=
F_{\beta}|_{\Xx _{\alpha \beta}}$.
\end{lemma}
{\em Proof:}
Given a family of such functors, and given a composable sequence $\phi
_1\cdots \phi _n$, let $\phi _i\in \Xx _{\alpha _i}$. We put
$$
F(\phi _1\cdots \phi _n):=F_{\alpha _1}(\phi _1)\cdots F_{\alpha _n}(\phi
_n). $$
This is compatible with the equivalence relation given by the $1$-cells of
$\Nn$, so it defines a functor $\pi _0\Nn \rightarrow \Aa$. On the other
hand, note that we have functors $i_{\alpha}:\Xx _{\alpha}\rightarrow \pi
_0 \Nn$ such that $i_{\alpha}|_{\Xx _{\alpha \beta}}= i_{\beta}|_{\Xx
_{\alpha \beta}}$. If $F:\pi _0\Nn \rightarrow \Aa$, set $F_{\alpha}:=
F\circ i_{\alpha}$. These two constructions are inverses. \eop

\begin{corollary}
\mylabel{flex09}
In the above situation, suppose $\Cc$ is $0$-truncated. The path components
of the fiber of
$$
Hom (\Mm ,\Cc )\rightarrow Hom (\Nn ,\Cc ) $$
over an element $F$ on the right corresponding to $\{ F_{\alpha}:\Xx
_{\alpha}\rightarrow \pi _0\Cc \} $, are weakly connected, and correspond
to the functors $G:\Xx \rightarrow \pi _0 \Cc$ such that $G|_{\Xx
_{\alpha}}=F_{\alpha}$.
\end{corollary}
{\em Proof:}
Put together Corollary \ref{flex07} and the previous lemma. \eop

\subnumero{The case $n=1$}

Now suppose $\Mm$ is the free continuous category over $\Xx$, and that we
have a functor $\rho :\Pp \rightarrow \Xx$ of continuous categories such
that $\rho$ induces an isomorphism on the set of objects, and for each
morphism $\phi \in \Xx$, the inverse image $\Pp (\phi )$ is a connected
$K(\pi ,1)$ space.

We do not assume that there is a privileged choice of base point in $\Pp
(\phi )$. We investigate this question first.

\begin{lemma}
\mylabel{flex10}
Suppose $\Nn\subset \Mm$ is a cellular subcategory, and suppose that there
exists a continuous functor $q:\Nn \rightarrow \Pp$ over $\Xx$. Suppose
that for each $\phi \in \Xx$, the inverse image $\Nn (\phi )$ is connected,
and the map $q(\phi ):\Nn (\phi )\rightarrow \Pp (\phi )$ is homotopic to a
constant map.
Then there exists a functor $p:\Mm \rightarrow \Phi $ extending $q$. \end{lemma}
{\em Proof:}
Let $\Mm _1$ be the subcategory of $\Mm$ generated over $\Nn$ by the
primitive $0$ and $1$-cells of $\Mm$ which are not in $\Nn$. Let $\Mm _2$
be the subcategory generated over $\Nn$ by the primitive cells of dimension
$\leq 2$. Since each $\Nn (\phi )$ is connected, we can choose a functor
$r:\Mm _1\rightarrow \Nn$ equal to the identity on $\Nn$. We obtain a
functor $$
p_1:= qr: \Mm _1\rightarrow \Pp .
$$
On the other hand, if $\alpha$ is a primitive $2$-cell of $\Mm$, then
$p_1|_{\partial \alpha }$ is homotopic to a constant map, from the second
hypothesis. Therefore $p_1|_{\partial \alpha }$ extends to a map defined on
$\alpha$. Choosing one such extension for each primitive $2$-cell, we
obtain a functor $p_2:\Mm _2\rightarrow \Pp$ extending $p_1$. Now we can do
the same for all other primitive cells, since the higher homotopy groups of
the $\Pp (\phi )$ vanish. We obtain the desired $p$. \eop

{\em Remark:}
Note that the hypotheses of the lemma are automatically satisfied if the
$\Nn (\phi )$ are connected and simply connected.

We now suppose given a functor $p:\Mm \rightarrow \Pp$. For each $\phi \in
\Xx$, put
$$
G(\phi ):= \pi _1(\Pp (\phi ), p(\phi )). $$
For a composable pair $\phi , \psi $ we obtain a path $p(\phi , \psi )$
going from $p(\phi \psi )$ to $p(\psi )p(\phi )$ in $\Pp (\phi \psi )$. On
the other hand, the multiplication map from the structure of $\Pp$, $$
m:\Pp (\psi )\times \Pp (\phi )\rightarrow \Pp (\phi \psi ) $$
gives a map
$$
G(\psi )\times G(\phi )\rightarrow \pi _1(\Pp (\phi \psi ), p(\psi )p(\phi
)). $$
Conjugating with the path $p(\phi ,\psi )$ we obtain 	a map (of groups) $$
\mu :G(\psi )\times G(\phi )\rightarrow G(\phi \psi ). $$
The precise formula, with the convention that a composition of paths
$\alpha \beta$ indicates $\alpha$ on the interval $[0,1/2]$ then $\beta$ on
$[1/2,1]$, is
$$
\mu (\alpha , \beta )= p(\phi , \psi ) m(\alpha , \beta ) p(\phi ,\psi
)^{-1}. $$
We
construct a category $\Gg$ together with a functor $\varpi : \Gg
\rightarrow \Xx$ such that $\varpi$ is an isomorphism on the set of
objects, and such that $\Gg (\phi )$ is identified with $G(\phi )$. The
composition in $\Gg $ is defined by the multiplication $\mu$. To verify
that this defines a category, we must check that the composition defined by
$\mu$ is associative. Suppose $\alpha \in G(\phi )$, $\beta \in G(\psi )$,
and $\gamma \in G(\tau )$ with $\phi , \psi , \tau$ a composable sequence
in $\Xx$. Then $$
\mu (\alpha , \mu (\beta , \gamma )) =
$$
$$
p(\phi , \psi \tau ) m(\alpha , \mu (\beta , \gamma )) p(\phi ,\psi \tau )^{-1}=
$$
$$
p(\phi , \psi \tau ) m(\alpha ,
p(\psi , \tau ) m( \beta ,\gamma )
p(\psi ,\tau )^{-1}) p(\phi ,\psi \tau )^{-1} $$
$$
p(\phi , \psi \tau ) m({p(\phi )} ,
p(\psi , \tau ))m(\alpha , m( \beta ,\gamma )) m(*_{p(\phi )} , p(\psi
,\tau )^{-1}) p(\phi ,\psi \tau )^{-1} $$
$$
A m(\alpha , m(\beta , \gamma ))A^{-1}
$$
where
$$
A=
p(\phi , \psi \tau ) m({p(\phi )},
p(\psi , \tau )).
$$
Similarly,
$$
\mu (\mu (\alpha , \beta , \gamma ) =
$$
$$
B m(m(\alpha , \beta ), \gamma )B^{-1}
$$
where
$$
B=
p(\phi \psi , \tau ) m(
p(\phi , \psi ),{p(\tau )} ).
$$
Since the composition $m$ of $\Pp$ is associative, we have $$
m (\alpha , m (\beta , \gamma )) =
m (m (\alpha , \beta , \gamma ) ,
$$
so it suffices to show that $A=B$. But the homotopy $p(\phi , \psi , \tau
)$ is exactly a map of a square to $\Pp (\phi \psi \tau )$ such that two of
the sides go to $A$ and two go to $B$. Thus $A$ and $B$ are homotopic,
giving the associativity of the composition $\mu$.

The structure of group on $G(\phi )$ is compatible with the multiplication,
so the category $\Gg$ is a group in the category of categories over $\Xx$.
In other words, we have a product
$$
\Gg \times _{\Xx }\Gg \rightarrow \Gg
$$
which is associative in the sense that the diagram $$
\begin{array}{ccc}
\Gg \times _{\Xx }\Gg \times _{\Xx }\Gg&\rightarrow &\Gg \times _{\Xx }\Gg
\\ \downarrow &&\downarrow \\
\Gg \times _{\Xx }\Gg &\rightarrow &\Gg
\end{array}
$$
commutes; and there is a morphism $i:\Gg \rightarrow \Gg$ over $\Xx$ giving
the inverse, as well as a section $\Xx \rightarrow \Gg$ giving the
identity.

A {\em $\Gg$-torsor} $\Jj$ over $\Xx$ is a category $\Jj$ with a functor
$\Jj \rightarrow \Xx$ inducing an isomorphism on the set of objects, and a
functor $$
\Gg \times _{\Xx} \Jj \rightarrow \Jj
$$
such that the diagram
$$
\begin{array}{ccc}
\Gg \times _{\Xx }\Gg \times _{\Xx }\Gg&\rightarrow &\Gg \times _{\Xx }\Gg
\\ \downarrow &&\downarrow \\
\Gg \times _{\Xx }\Gg &\rightarrow &\Gg
\end{array}
$$
commutes, and such that $\Gg (\phi )$ acts simply transitively on $\Jj
(\phi )$ for all morphisms $\phi $ in $\Xx$.
A {\em trivialization} of a $\Gg$-torsor $\Jj$ is a functor $u:\Gg
\rightarrow \Jj$ compatible with the left action of $\Gg$ on itself. This
is necessarily an isomorphism. Note that a trivialization is the same thing
as a section $\Xx \rightarrow \Jj$: the section corresponding to $u$ is
$\phi \mapsto u(1_{\phi })$.

If $\Jj$ is a $\Gg$-torsor, we define another relative group $Ad (\Jj
)\rightarrow \Xx$ as follows. An element of $Ad (\Jj )(\phi )$ is an
expression of the form $j^{-1}gj$ for $j\in \Jj (\phi )$ and $g\in \Gg
(\phi )$. Such an expression is equivalent to another $(j')^{-1}g'j'$ if
$j'=hj$ and $g'=hgh^{-1}$ for $h\in \Gg (\phi )$. The multiplication of the
group law is given by formally multiplying after putting the two elements
into a form with the same element $j$. On the other hand, the
multiplication in the category is given by the formula
$$
(j^{-1}gj)\cdot ((j')^{-1}g'j'):= (j\cdot j')^{-1}(g\cdot g')(j\cdot j'). $$
This is evidently associative, and it is compatible with the above
equivalence relation (and the structure of group) since $$
(hj)\cdot (h'j')= (h\cdot h')(j\cdot j'). $$

Note that $Ad (\Jj )(\phi )$ is equal to the group of automorphisms of $\Jj
(\phi )$ considered as a left $\Gg (\phi )$-set, with the automorphisms
acting on the right by the rule
$$
j (j^{-1}gj)= gj.
$$

\begin{lemma}
\mylabel{flex11}
Suppose $\Nn \subset \Mm$ is a cellular subcategory. Let $q$ denote the
restriction of $p$ to $\Nn$. The space of continuous functors $\Mm
\rightarrow \Pp$ over $\Xx$, restricting to $q$ on $\Nn$, is a union of
path-components $T_{(\Jj ,t)}$ parametrized by pairs consisting of a
$\Gg$-torsor $\Jj$ over $\Xx$ and a trivialization $t$ of the pullback of
$\Jj$ to $\pi _0 \Nn$. The path component $T_{(\Jj ,t)}$ has trivial
homotopy groups in degrees $\geq 2$, and its fundamental group is
isomorphic to the group of sections $\Xx \rightarrow Ad (\Jj )$ pulling
back to the identity section on $\pi _0 \Nn$ (which is also the group of
automorphisms of the pair $(\Jj ,t )$). \end{lemma}
{\em Proof:}
Suppose $f:\Mm \rightarrow \Pp$ is a continuous functor. Put $\Jj (\phi )$
equal to the space of paths from $p(\phi )$ to $f(\phi )$, with composition
of arrows defined as above. We get a category $\Jj $ over $\Xx$ with
structure of left $\Gg$-torsor. This is clearly an invariant of the path
component containing $f$, and if $f$ restricts to $q$ (the restriction of
$p$) on $\Nn$ then we obtain a trivialization of $\Jj$ over $\pi _0\Nn$.
Note that $Ad (\Jj )$ is now the same relative group as would have been
obtained by starting with the basepoint $f$. We claim that the
trivializations of $(\Jj , t)$ trivial on $\Nn$ are in one to one
correspondence with the homotopy classes of homotopies from $f$ to $p$.
Given such a homotopy, we obtain a trivialization. Conversely, given a
trivialization corresponding to a section $s:\Xx \cong \Jj$ (equal to the
identity over $\pi _0\Nn$), use the path $s(\phi )$ from $p(\phi )$ to
$f(\phi )$ to make a homotopy between $p$ and $f$. Such a homotopy will
exist and be unique (up to homotopy), if the two paths $$ p(\phi , \psi
)m(s(\phi ), s(\psi ))
$$
and
$$
s(\phi \psi )f(\phi , \psi )
$$
from $p(\phi \cdot \psi )$ to $f(\phi )\cdot f(\psi )$ are homotopic. This
is the case because $s$ is a section (in view of the category
multiplication in $\Jj$, defined analogously to that of $\Gg$ but using the
paths $p(\phi ,\psi )$ and $f(\phi , \psi )$). The existence and uniqueness
statements from here are due to the fact that maps of two-cells into the
$P(\phi \psi )$ are unique up to homotopy, and exist if and only if the
boundary is a trivial path (and maps of $n$-cells exist and are unique for
$n\geq 3$). This statement shows that the fundamental group of the path
component containing $p$ is equal to the space of sections of $\Gg$
(trivial on $\Nn$); and that any functor $f$ giving a pair $(\Jj , t)$
isomorphic to $(\Gg , e)$, lies in the same path component. In particular,
since the relative group associated to any $f$ is $Ad (\Jj )$, the
fundamental group of the path component containing $f$ is the space of
sections of $Ad (\Jj )$ trivial on $\Nn$. To finish, it suffices to see
that $f$ and $f'$ are in the same path component if their torsors $(\Jj
,t)$ and $(\Jj ',t')$ are isomorphic. For this it suffices to note that
there is a construction for the torsor $(\Jj '_f,t'_f)$ of $f'$ from the
point of view of $f$, and the same for the torsor $(\Jj _f,t_f)$ of $f$
from the point of view of $f$ (and this construction depends only on the
isomorphism class of the original torsor, so $(\Jj '_f,t'_f)\cong (\Jj
_f,t_f)$); but as $f$ is in its own path component, $(\Jj _f,t_f)$ is
trivial; thus $(\Jj '_f,t'_f)$ is also trivial so the above proof shows
that $f'$ is in the path component of $f$. \eop

\begin{corollary}
\mylabel{flex12}
Suppose that for any relative group $\Gg$ over $\Xx$ and any $\Gg$-torsor
$\Jj$, we have
$$
Hom _{\Xx }(\Xx ,\Jj )\rightarrow Hom _{\Xx}(\pi _0 \Nn , \Jj ). $$
Then for any relative degree $1$ Eilenberg-MacLane category $\Pp$ over
$\Xx$ with a section $p:\Mm \rightarrow \Pp$, the fiber of $$
Hom (\Mm , \Pp )\rightarrow Hom (\Nn , \Pp ) $$
over the restriction $q$ of $p$ to $\Nn$, is weakly contractible.
\end{corollary}
{\em Proof:}
The fiber is path-connected, since by hypothesis any pair $(\Jj ,t )$ is
trivial; and the fundamental group is trivial since the space of
trivializations restricting to $t$ has exactly one element. \eop

Finally, we return to the case where $\Nn$ is
generated by a family of subcategories $\Xx _{\alpha}\subset \Xx$. We don't
answer the question of the existence of a section $p$ restricting to a given one
$q$ any further than in Lemma \ref{flex11} above. We consider the criterion of
the corollary. If $\Jj$ is a $\Gg$-torsor over $\Xx$, the space of sections $Hom
_{\Xx}(\pi _0 \Nn , \Jj )$ is equal to the space of families of sections
$$ \{
u_{\alpha}:\Xx _{\alpha} \rightarrow \Jj \;\mbox{s.t.} \;\; u_{\alpha}|_{\Xx
_{\alpha \beta}} =u_{\beta}|_{\Xx _{\alpha \beta}} \} .
$$
The criterion is
therefore that this space of sections be equal to the space of sections $u:\Xx
\rightarrow \Jj$.

\newpage

\setcounter{section}{1}

\numero{Morphisms}

We define the space of morphisms from one flexible functor $T^0$ to another one
$T^1$. This leads to a ``Segal category'' (cf Segal \cite{SegalTopology},
Dwyer-Kan-Smith \cite{DwyerKanSmith}) of
flexible functors. This is probably closely related to the simplicial category
$Coh(\Xx , Top)$ defined by Cordier and Porter \cite{CordierPorter} although we
don't prove this relationship.  It is also very similar to the original approach
of Vogt \cite{Vogt1} where he obtained a {\em restricted Kan simplicial set}. In
the present paper we topologize the morphism sets so we obtain a simplicial
space
rather than simplicial set, and the Segal condition is analogue to the
restricted Kan condition.  See \cite{descente} for a somewhat different
approach. As we shall see later (in a reworking of Vogt's original theorem
\cite{Vogt1}), we get the same answer as the simplicial category of fibrant
simplicial presheaves \cite{Jardine1} \cite{JardineBool}.

We denote by $I$ the category with two objects, denoted $0$ and $1$, and
one morphism $1\rightarrow 0$ (aside from the identities). We denote by
$I^n$ the cartesian product, and $I^{(n)}$ the symmetric product. We denote
the objects of $I^{(n)}$ by $0,1,\ldots , n$; there is exactly one morphism
$i\rightarrow j$ when $i\geq j$ (and none if $i<j$). We denote by
$\overline{I}$ the category with objects $0$ and $1$, and morphisms
$0\rightarrow 1$ and $1\rightarrow 0$ (whose compositions are the
identities). Again, $\overline{I}^n$ is the cartesian product and
$\overline{I}^{(n)}$ is the symmetric product with objects denoted
$0,\ldots , n$ (and one isomorphism between each pair of objects). We
consider $I$, $I^n$, $I^{(n)}$ as subcategories respectively of
$\overline{I}$, $\overline{I}^n$, $\overline{I}^{(n)}$.

Suppose $T^0$ and $T^1$ are flexible functors from $\Xx$ to a continuous
category $\Cc$. A {\em morphism} from $T^0$ to $T^1$ is a flexible functor
$F$ from $\Xx \times I$ to $\Cc$ such that $F|_i = T^i$ for $i=0,1$.
This is
the same definition as used by Vogt \cite{Vogt1}.
Note that there is no obvious
natural way of composing morphisms; we will treat this in more detail
during this
section. However, there are natural notions of {\em composable sequences} of
morphisms of flexible functors. A {\em composable sequence of length $n$} of
morphisms of flexible functors is a morphism
$$
F: \Xx \times I^{(n)}\rightarrow \Cc
$$
(again this definition was originally made by Vogt \cite{Vogt1}).
Note that for any $0\leq i\leq j\leq n$ we have an inclusion
$$
\kappa (n; i,j): I^{(j-i)}\rightarrow I^{(n)} $$
sending the object $k$ to the object $k+i$. If $F$ is a composable sequence
of length $n$ then for any $0\leq i\leq j\leq n$ we obtain the composable
sequence $F_{ij}:=\kappa (n; i,j)^{\ast}F$ of length $j-i$. In particular,
we obtain flexible functors $F_i:\Xx \rightarrow \Cc$ and morphisms $F_{i\,
i+1}$ from $F_i$ to $F_{i+1}$.

An {\em equivalence} from $T^0$ to $T^1$ is a flexible functor $F$ from
$\Xx \times \overline{I}$ to $\Cc$ such that $F|_i = T^i$ for $i=0,1$. Note
that an equivalence induces morphisms from $T^0$ to $T^1$ and from $T^1$ to
$T^0$.

If $T^0$, \ldots , $T^n$ are fixed flexible functors, then the set of
composable sequences $F$ of length $n$ with $F_i=T^i$ forms a topological
space which we denote by $Hom (T^0,\ldots , T^n)$. If $\Mm$ is the free
continuous category on $\Xx \times I^{(n)}$ and if $\Nn$ is the cellular
subcategory generated by the subcategories $\Xx \times \{ i\}$, then $Hom
(T^0,\ldots , T^n)$ is the fiber of $Hom (\Mm ,\Cc )$ over the element of
$Hom (\Nn ,\Cc )$ corresponding to $T^0,\ldots , T^n$.

Similarly, the set of
equivalences between $T^0$ and $T^1$ forms a topological space which
we denote $Equiv(T^0, T^1)$. If $\Mm$ is the free continuous category on $\Xx
\times \overline{I}$ and $\Nn$ the subcategory generated by the $\Xx \times \{
i\}$, then $Equiv(T^0,
T^1)$ is the fiber of $Hom (\Mm , \Cc )$ over the point of $Hom (\Nn ,\Cc
)$ corresponding to $T^0,T^1$. We can also define a space of composable
$n$-tuples of equivalences $Equiv (T^0,\ldots , T^n)$ as above.

The following theorem is almost the same as what Vogt did \cite{Vogt1}; he
didn't topologize the spaces of composable $n$-uples, and he obtained the
restricted Kan condition for the resulting simplicial set.  See Dwyer-Kan-Smith
\cite{DwyerKanSmith}.

\begin{theorem}
\mylabel{morph01}
Suppose $T^0, \ldots , T^n$ are flexible functors from $\Xx$ to a
continuous category $\Cc$. Suppose $G_{i}\in Hom (T^{i-1}, T^{i})$ for
$i=1,\ldots , n$. Then the space of composable $n$-tuples $F\in Hom
(T^0,\ldots , T^n)$ such that $F_{i-1,i}=G_{i}$ is weakly contractible.
\end{theorem}
{\em Proof:}
Let $\Mm$ be the free continuous category on $\Yy :=\Xx \times I^{(n)}$ and
let $\Nn$ be the cellular subcategory generated by the subcategories $\Yy
_i:= \Xx \times (\kappa (n; i-1,i)(I))$. We apply the criteria developed in
the first chapter.

{\em First claim:} For any $\phi$ in $\Yy$, $\Nn (\phi )$ is contractible.

We give the proof in the case $n=2$, for now. The points of $\Nn (\phi )$
may be written in the form $(W, \epsilon , t)$ where
$W= (\phi _1,\ldots ,\phi _k)$ is a word of composable morphisms in $\Xx$,
with composition equal to $\phi$; $\epsilon = (\epsilon _0,\ldots ,
\epsilon _k)$ is an increasing sequence of elements of $\{ 0,\ldots , n\}$;
and $t=(t_1,\ldots ,t_{k-1})$ is a point in $[0,1]^{k-1}$. These points are
subject to the usual relations for $t_i=0$. Define a map $$
F:\Nn
(\phi )\times [0,3]\rightarrow \Nn (\phi ) $$
as follows. For $s\in [0,1]$, put
$$ F(W,\epsilon , t; s)= (W,\epsilon , t') $$ where $t'_i=t_i$ if $\epsilon
_i\geq 1$, and $t'_i=(1-s)t_i$ if $\epsilon _i=0$. Put
$$
F'(W,\epsilon , t;1)=(W,\epsilon ', t')
$$
where $t'$ is as defined above for $F(W,\epsilon ,t;1)$, and $\epsilon
'_i=max (1,\epsilon _i)$ for all $i$ except $i=0$ (where $\epsilon _0=0$).
Note that $F(W,\epsilon ,t;1)\sim F'(W,\epsilon ,t;1)$ because of the rules
for changing the $\epsilon _i$. Put
$$
F'' (W,\epsilon ,t;1)=(W'',\epsilon '', t'') $$
where
$$
W'' = I\cdot W
$$
(this denotes concatenation, putting $I$ at the front of $W$); $$
\epsilon ''_0=0,
$$
$$
\epsilon ''_i=max (1,\epsilon _{i-1}), \; \; i>0; $$
and
$$
t''_i=t'_{i-1},
$$
but with $t''_1:=0$. Note again that $F'' (W,\epsilon ,t;1)\sim
F(W,\epsilon ,t;1)$, because with $t''_0=t''_{k''-1}=0$, the addition of
$I$ at the beginning of $W''$ doesn't change anything. For $s\in [1,2]$,
put $$
F (W,\epsilon ,t;s)= (W'',\epsilon '', t^3) $$
where $t^3_1= s-1$; $t^3_i=t''_i$ if $\epsilon _i=2$; and $t^3_i=
(2-s)t''_i$ if $\epsilon _i =1$, $i>1$. Finally, for $s\in [2,3]$, put $$
F(W,\epsilon , t;s)= (W'', \epsilon '' , t^4) $$
where $t^4_1=1$; $t^4_i=0$ if $\epsilon _i=1$ but $i>1$; and
$t^4_i=(3-s)t^3_i$ if $\epsilon _i=2$. Note that $F$ is continuous on $\Nn
(\phi )\times [0,3]$
(it is compatible with the relations corresponding to $t_i=0$ because we
have only changed the $t_i$ by multiplying by a function of $s$); and that
$F(W,\epsilon ,t;3)\sim (W^0,\epsilon ^0,t^0)$ where $$
W^0=(I,\phi _1,\phi _2,\ldots , \phi _n ), $$
$$
\epsilon ^0=(0,1,2),
$$
and
$$
t^0=(1).
$$
This gives the desired contraction in the case $n=2$.

In general, by using the same type of contraction as the $F$ for $s\in
[0,1]$ we can contract $\Nn (\phi )$ to the subspace consisting of those
points for which $\epsilon _1=1$, $u_1=I$ and $t_1=0$. Then contract to the
subspace consisting of points with $\epsilon _1=1$, $\phi _1=I$, $t_1=1$,
and $t_i=0$ for any other $i$ with $\epsilon _i=1$. Then change this to a
point with $\epsilon _1=1$, $u_1=I$, $t_1=1$, and $\epsilon _2=2$; add an
identity, so that $\epsilon _1=1$, $\phi _1=I$, $t_1=1$, $\epsilon _2=2$,
$\phi _2=I$, $t_2=0$; then change to $t_2=1$ and $t_i=0$ for any other
$t_i$ with $\epsilon _i=2$; and so on.

\medskip

From this claim, note that for any continuous category $\Cc$ and any
functor $u:\Nn \rightarrow \Cc$, there exists an extension to $v:\Mm
\rightarrow \Cc$ restricting to $u$ on $\Nn$ (by Proposition
\ref{flex01}). In particular, this solves the existence problems for
degrees $0$ and $1$ in the Postnikov induction.

Note also that the claim gives $\pi _0 \Nn = \Yy$.

We now treat the case of degree $0$. If $\Cc$ is a $0$-truncated homotopy
category, then the fiber of $Hom (\Mm ,\Cc )\rightarrow Hom (\Nn , \Cc )$
is weakly contractible, since $\pi _0 \Nn =\pi _0 \Mm =\Yy$ (see Corollary
\ref{flex07}).

For the case of degree $1$, note that for any relative group
$\Gg\rightarrow \Yy$ and any $\Gg$-torsor $\Jj$, the set of trivializations
of $\Jj$ over $\pi _0\Nn$ is the same as the set of trivializations over
$\Yy$, since $\pi _0\Nn =\Yy$. By Corollary \ref{flex12}, for any functor
$\Pp \rightarrow \Yy$ where the $\Pp (\phi )$ are Eilenberg-MacLane spaces
of degree $1$, the fiber of $$
Hom _{\Yy}(\Mm , \Pp )\rightarrow Hom _{\Yy}(\Nn , \Pp ) $$
is weakly contractible (it is nonempty as a corollary of the claim above).

Finally, we consider the case of a functor $\Pp \rightarrow \Yy$ whose
fibers $\Pp (\phi )$ are Eilenberg-MacLane spaces of degree $m\geq 2$. We
would like to see that the fiber of
$$
Hom _{\Yy}(\Mm , \Pp )\rightarrow Hom _{\Yy}(\Nn , \Pp ) $$
is weakly contractible. By the criterion of Lemma \ref{flex05}, it suffices
to see that the subcomplex $\Lambda _{\Nn}(\phi )$ has the same cohomology
with abelian group coefficients, as $\Lambda _{\Mm }(\phi
)=N(Ins_{\Yy}(\phi ))$. For this, we may assume that $\phi = (\psi ,
0\leftarrow n)$ where $\psi$ is a morphism in $\Xx$ (for otherwise it
amounts to the same statement for a smaller $n$). Then
$$
Ins _{\Yy}(\phi )=Ins _{\Xx}(\psi )\times Ins _{I^{(n)}}(0\leftarrow n). $$
But $Ins _{I^{(n)}}(0\leftarrow n)=I^{(n)}$. Let $N_{0,n}$ denote the nerve
of $I^{(n)}$ and for any $i$ let $N_{i-1,i}$ denote the nerve of $\kappa
(n; i-1,i)(I)\subset I^{(n)}$. Then
$$
N(Ins _{\Yy}(\phi ))=N(Ins _{\Xx}(\psi ))\times N_{0,n} $$
and
$$
\Lambda _{\Nn }(\phi )=\bigcup _{i=1}^nN(Ins _{\Xx}(\psi ))\times N_{i-1,i} . $$
$$
=N(Ins _{\Xx}(\psi ))\times \left( \bigcup _{i=1}^n N_{i-1,i}\right) . $$
Notice, however, that
$$
\left( \bigcup _{i=1}^n N_{i-1,i}\right) \subset N_{0,n} $$
induces an isomorphism on cohomology (both sides are contractible). By the
Kunneth formula, the inclusion $\Lambda _{\Nn}(\phi )\subset \Lambda
_{\Mm}(\phi )$ induces an isomorphism on cohomology with integer
coefficients, hence with any abelian group coefficients. By Lemma
\ref{flex05}, the fibers of $$
Hom _{\Yy }(\Mm , \Pp )\rightarrow
Hom _{\Yy }(\Nn ,\Pp )
$$
are weakly contractible.

Now if $\Cc$ is any continuous category, then by the Postnikov tower for
$\Cc$ explained in the previous chapter, the fibers of $$
Hom (\Mm , \Cc )\rightarrow Hom (\Nn ,\Cc ) $$
are weakly contractible. This completes the proof of the theorem. \eop

For $n=2$ this theorem permits us to speak of the {\em composition} of two
morphisms of flexible functors $G_1:T^0\rightarrow T^1$ and
$G_2:T^1\rightarrow T^2$. In effect, the space of composable pairs $F$
restricting to $(G_1,G_2)$ is weakly contractible, so there is not too much
ambiguity in speaking of ``the'' composable pair defined by $(G_1,G_2)$.
For any such composable pair $F$, we obtain a morphism $T^0\rightarrow T^2$
(from the third inclusion $I\rightarrow I^{(2)}$), which we think of as
``the'' composition.

More generally, the above theorem says that we obtain a ``Segal category''
\cite{SegalTopology} \cite{DwyerKanSmith} \cite{effective} of flexible functors
$Flex(\Xx , Top )$, by setting the object-set $Flex(\Xx ,Top)_0$ to be the set
of flexible functors and setting
$$
Flex(\Xx ,Top)_{n/}(T^0,\ldots , T^n):= Hom(T^0,\ldots , T^n)
$$
(in the notation of \cite{effective}). One expects that
$Flex(\Xx ,Top)$ is equivalent to the Segal category obtained by taking the
nerve of Cordier's and Porter's $Coh(\Xx , Top)$ or Jardine's simplicial
category of simplicial presheaves on $\Xx$. In fact, the arguments in
\cite{CordierPorter}, \cite{DHK}, \cite{Jardine1}, the present paper, and
\cite{nCat} \cite{descente} probably serve to prove this (but I haven't
verified the details).

In spite of the ambiguity in the actual choice of composition in the
general case, there is a canonical choice in the case of the composition of
a morphism of flexible functors with a morphism of functors, or vice-versa.
These canonical compositions may be viewed as obtained by simply
transporting structure of the morphism of flexible functors, by the
morphism of functors on one side or the other. In terms of formulas, if $F$
is a morphism of flexible functors and $G$ a morphism of functors, we can
set $$
(FG)(\phi , \epsilon , t)=
F(\phi , \epsilon , t)\circ G_{X_0}
$$
when $\epsilon _0=0$ (here $X_0$ is the
target of the first morphism $\phi _1$ in the word $\phi$). Similarly, if
$F$ is a morphism of functors and $G$ a morphism of flexible functors then
we can set
$$
(FG)(\phi , \epsilon , t)=
F_{X_n}\circ G(\phi , \epsilon , t)
$$
for $\epsilon _n=1$,
where $X_n$ is the source of the last $\phi _n$ in $\phi$.

The following theorem is a refinement of a theorem of Vogt \cite{Vogt1} (see
also Cordier-Porter \cite{CoPo2} for the generalization to coherent diagrams in
any simplicial category) stating that a morphism of diagrams which is level-wise
a homotopy equivalence, has an inverse. We prove here that a suitably defined
``space of inverses'' is contractible.

One can remark that our space of
inverses is somewhat similar to the space of ``strong homotopy inverses''
mentionned by Cordier-Porter in \cite{CoPo2}, however their notion involved
only a first-order homotopy; our strong homotopy inverses are actually coherent
diagrams indexed by the category $\overline{I}$, and thus include all
higher-order homotopies.

A similar
later result is Theorem 2.5.1 of \cite{limits}.

\begin{theorem}
\mylabel{morph02}
Suppose $T^0$ and $T^1$ are flexible functors from $\Xx$ to $\Cc$ and $G$
is a morphism from $T^0$ to $T^1$. The space of equivalences $F$ from $T^0$
to $T^1$ restricting to the morphism $G$ is nonempty if and only if the
morphisms $G(1_X,0\rightarrow 1 )$ are invertible in $\pi _0 \Cc $ for all
$X\in \Xx$. If the space of $F$ is nonempty, then it is weakly
contractible.
\end{theorem}
{\em Proof:}
Let $\Yy = \Xx \times \overline{I}$ and let $\Yy _1=\Xx \times I\subset
\Yy$. Let $\Mm$ be the free continuous category on $\Yy$ and let $\Nn$ be
the cellular subcategory defined by $\Yy _1$ (it is the free continuous
category on $\Yy _1$). The morphism $G$ is a functor from $\Nn$ to $\Cc$,
and the space of $F$ restricting to $G$ is the fiber of $Hom (\Mm , \Cc )$
over $G$.

It is clear
that if $F$ exists then the morphisms $G(1_X,0\rightarrow 1 )$ are
invertible in $\pi _0 \Cc $. Conversely, if $\Cc$ is $0$-truncated and if
the morphisms $G(1_X,0\rightarrow 1 )$ are invertible in $\pi _0\Cc$, then
there exists $F:\Mm \rightarrow \Cc$ restricting to $G$ on $\Nn$ (and the
space of such $F$ is weakly contractible by Lemma \ref{flex06}). To see
this, note that to give $F$ is the same as to give $\pi _0 F: \Yy
\rightarrow \pi _0\Cc$. Put $$
a_X:= (1_X, 0\rightarrow 1 )
$$
and
$$
b_X:= (1_X, 1\rightarrow 0).
$$
A morphism in $\Yy$
is either a morphism of $\Yy _1$ or else can be written as $b_Xf$ where $f$
is a morphism in $\Xx \times \{ 1\}$. If $f$ is a morphism in $\Yy _1$, put
$F(f)=G(f)$. Our hypothesis says that $G(a_X)$ is invertible; note that the
inverse is unique. We put $F(b_Xf):= G(a_X)^{-1}G(f)$. We have to check
that this defines a functor. A morphism from $(X,0)$ to $(Y,1)$ may be
written as $ga_X$; then
$$
F(ga_X)F(b_Xf)=G(g)G(a_X)G(a_X)^{-1}G(f)=G(gf)=F(gf)=F(ga_Xb_Xf). $$
If $g$ is a morphism from $(X,0)$ to $(Y,0)$, then we have $$
a_Yg= g'a_X
$$
for $g':(X,0)\rightarrow (Y,0)$ corresponding to the same morphism in
$\Xx$. This gives
$$
gb_Xf=b_Yg'f,
$$
so
$$
F(gb_Xf)= G(a_Y)^{-1}G(g'f),
$$
while
$$
G(g)=G(a_Y)^{-1}G(g')G(a_X).
$$
Thus
$$
F(g)F(b_Xf)= G(g)G(a_X)^{-1}G(f)=G(a_Y)^{-1}G(g')G(f) = F(gb_Xf). $$
The verifications for multiplication on the other side are similar.

Suppose now that $\Pp \rightarrow \Yy$ is a relative Eilenberg-MacLane
category of degree $1$, such that there exists a section $q$ over $\Nn$. We
claim that there exists a section $p:\Mm \rightarrow \Pp$. For each object
$X$, choose $p(b_X)\in \Pp (b_X)$. Put $p(a_X)=q(a_X)$. Choose paths $u_X$
from $p(1_{X,0})=e_{X,0}$ to $p(b_X)p(a_X)$. Then we obtain a path
$p(a_X)u_X$ from $p(a_X)$ to $p(a_X)p(b_X)p(a_X)$. Choose paths $v_X$ from
$p(1_{X,1})=e_{X,0}$ to $p(a_X)p(b_X)$ such that the path $v_Xp(a_X)$ from
$p(a_X)$ to
$p(a_X)p(b_X)p(a_X)$ is homotopic to $p(a_X)u_X$. For making this choice,
note that multiplication by $p(a_X)$ is a homotopy equivalence from $\Pp
(f)$ to $\Pp (a_Xf)$ and similarly for $p(b_X)$ (since $\Pp _{X,1}$ is
path-connected, so multiplication by $p(a_X)p(b_X)$ is homotopy equivalent
to multiplication by $p(a_Xb_X)=e_{X,1}$ and vice-versa). Now for morphisms
$f$ in $\Xx \times \{ 0\}$, put
$$
p(a_Y^{\epsilon}fb_X^{\delta }):=p(a_Y)^{\epsilon}q(f)p(b_X)^{\delta } $$
for $\epsilon$ and $\delta$ equal to $0$ or $1$. Put
$$
p(a_Z^{\epsilon }fb_Y,a_Ygb_X^{\delta }):= $$
$$
p(a_Z)^{\epsilon} [q(f,g)\ast (q(f)u_Yq(g))]p(b_X)^{\delta}. $$
Put
$$
p(a_Z^{\epsilon }f,gb_X^{\delta }):=
$$
$$
p(a_Z)^{\epsilon} q(f,g)p(b_X)^{\delta}. $$
(here the $\ast$ denotes composition of paths, while juxtaposition denotes
composition in $\Pp$). We have to check that a certain diagram commutes up
to homotopy, in the case of three morphisms. We treat the case where the
three
morphisms are
$$
a_Z^{\epsilon }fb_Y,\;\; a_Ygb_X,\;\; a_Xhb_W^{\delta} . $$
The other cases are similar.

Two sides of the diagram compose as
$$
p(a_Z^{\epsilon }fgb_X, a_Xhb_W^{\delta})\ast $$
$$
[p(a_Z^{\epsilon }fb_Y,a_Ygb_X)p(a_Xhb_W^{\delta})] $$
$$
=[p(a_Z)^{\epsilon}[q(fg,h)\ast (q(fg)u_Xq(h))]p(b_W)^{\delta}] $$
$$
\ast
p(a_Z)^{\epsilon} [q(f,g)\ast (q(f)u_Yq(g))]p(b_X)p(a_Xhb_W^{\delta}) $$
$$
= [p(a_Z)^{\epsilon}[q(fg,h)\ast (q(fg)u_Xq(h))]p(b_W)^{\delta}] $$
$$
\ast
p(a_Z)^{\epsilon} [q(f,g)\ast (q(f)u_Yq(g))]p(b_X)p(a_X)q(h)p(b_W)^{\delta}. $$
We may remove the outer $p(a_Z)^{\epsilon}$ and $p(b_W)^{\delta}$ from the
rest of the calculation (putting them back in at the end), since in general
$$
(u\gamma v)\ast (u\gamma ' v)= u(\gamma \ast \gamma ')v. $$
Our composition becomes
$$
q(fg,h)\ast (q(fg)u_Xq(h))\ast [[q(f,g)\ast (q(f)u_Yq(g))]p(b_X)p(a_X)q(h)] $$
$$
=
q(fg,h)\ast
(q(fg)u_Xq(h)) \ast (q(f,g)p(b_X)p(a_X)q(h))\ast (q(f)u_Yq(g)p(b_X)p(a_X)q(h)).
$$
Note that the operation $\ast$ is associative up to homotopy. The square
$$
q(f,g)u_X q(h)
$$
is a homotopy between
$$
(q(fg)u_Xq(h)) \ast (q(f,g)p(b_X)p(a_X)q(h)) $$
and
$$
(q(f,g)q(h))\ast (q(f)q(g)u_X q(h)).
$$
Thus we may make this replacement above. Our expression becomes $$
q(fg,h)\ast
(q(f,g)q(h))\ast (q(f)q(g)u_X q(h))\ast
(q(f)u_Yq(g)p(b_X)p(a_X)q(h)).
$$
Similarly, the square
$$
q(f)u_Yq(g)u_Xq(h)
$$
is a homotopy between
$$
(q(f)q(g)u_X q(h))\ast
(q(f)u_Yq(g)p(b_X)p(a_X)q(h))
$$
and
$$
(q(f)u_Y q(g)q(h))\ast (q(f)p(b_Y)p(a_Y)q(g)u_Xq(h)), $$
and making this replacement our expression becomes
$$
q(fg,h)\ast
(q(f,g)q(h))\ast (q(f)u_Y q(g)q(h))\ast (q(f)p(b_Y)p(a_Y)q(g)u_Xq(h)).
$$
Next, $q(f,g,h)$ is a homotopy between $q(fg,h)\ast (q(f,g)q(h))$ and
$q(f,gh)\ast (q(f)q(g,h))$; the expression becomes
$$
q(f,gh)\ast (q(f)q(g,h))
\ast (q(f)u_Y q(g)q(h))\ast (q(f)p(b_Y)p(a_Y)q(g)u_Xq(h)).
$$
Finally, with the homotopy given by the square
$$
q(f)u_Y q(g,h)
$$
as before, and after
putting in the outer elements dropped above, our expression becomes
$$
p(a_Z)^{\epsilon}[q(f,gh)\ast (q(f)u_Yq(gh))
$$
$$
\ast (q(f)p(b_Y)p(a_Y)q(g,h))
\ast (q(f)p(b_Y)p(a_Y)q(g)u_Xq(h))]p(b_W)^{\delta}
$$
$$
=
p(a_Z)^{\epsilon}[q(f,gh)\ast (q(f)u_Yq(gh))]p(b_W)^{\delta}
$$
$$
\ast
p(a_Z)^{\epsilon}
[(q(f)p(b_Y)p(a_Y)q(g,h)) \ast (q(f)p(b_Y)p(a_Y)q(g)u_Xq(h))]p(b_W)^{\delta}
$$
$$
=
p(a_Z^{\epsilon}fb_Y, a_Yghb_W^{\delta})
$$
$$
\ast p(a_Z^{\epsilon}fb_Y)p(a_Ygb_X, a_X hb_W^{\delta}), $$
which is the composition of the other two sides of the diagram.
Verification of the other diagrams is similar.

In fact, the definition of $p$ needs to be modified slightly, because of
the condition $p(1_{X,1})=e_{X,1}$. We enforce this condition by decree. We
need to change the definition of the paths to take this into account. Paths
$p(u,v)$ where one of $u$ or $v$ is the identity, are defined as in the
following example:
$$
p(e_{X,1}, a_X f):= p(a_X)q(1_X, f).
$$
Suppose $f:Y\rightarrow X$ and
$g:X\rightarrow Y$ such that $fg=e_X$. Then we put $$
p(a_Xfb_Y,a_Ygb_X):= v_X\ast p(a_X)[q(f,g)\ast (q(f)u_Yq(g))]p(b_X) $$
and
$$
p(a_Xf,gb_X):= v_X\ast (p(a_X)q(f,g)p(b_X)). $$
Finally, one should check the commutativity of the above diagram in the
cases where one of the morphisms is the identity. These are verified as
before, inserting $v_X$ or $v_Y$ where necessary and using the condition
established at the start.

We treat the case where the three
morphisms are
$$
a_Xfb_Y,\;\; a_Ygb_X,\;\; a_Xh ,
$$
with $fg=e_X$ (but no other composition equal to the identity). Two sides
of the diagram compose as
$$
p(e_{X,1}, a_Xh)\ast
$$
$$
[p(a_Xfb_Y,a_Ygb_X)p(a_Xh)]
$$
$$
=(p(a_X)q(e_X,h))\ast
[v_X\ast p(a_X)[q(f,g)\ast (q(f)u_Yq(g))]p(b_X)]p(a_Xh) $$
$$
= (p(a_X)q(e_X,h))\ast
(v_Xp(a_X )q(h))\ast [p(a_X) [q(f,g)\ast (q(f)u_Yq(g))]p(b_X)p(a_X)q(h)]. $$
The property $v_Xp(a_X)=p(a_X)u_X$ allows us to take the first factor of
$p(a_X)$ outside to give
$$
p(a_X)[q(fg,h))\ast(q(fg)u_Xq(h)) \ast
(q(f,g)p(b_X)p(a_X)q(h))\ast
(q(f)u_Yq(g)p(b_X)p(a_X)q(h))].
$$
Note that we have inserted $q(fg)=1_{X,0}$ and replaced $e_X$ by $fg$. This
expression is the same as in the previous argument (but note that there, we
had already taken the front $p(a_X)$ out temporarily). The rest of the
calculation is the same. Note that in our example, no other compositions
are the identity so the other composition of two sides of the square is the
same as above.
(In the case $fg=e_X$ and $gh=e_Y$ (in which case $f=h$!) the same
calculation is done at the start and at the end; they meet at the same
formula in the middle.)

We have now verified that if there exists a section $q:\Nn \rightarrow \Pp$
then there exists a section $p:\Mm \rightarrow \Pp$. The section $p$ we
have constructed does not restrict to $q$ on $\Nn$. In fact, we have only
used the existence of $q$ over the cellular subcategory $\Nn '$ generated
by $\Xx \times \{ 0\}$, and $p$ is equal to $q$ over $\Nn '$.

We now show that every section over $\Nn$ extends to a section over $\Mm$.
For this, we use Lemma \ref{flex11}. Let $\Gg \rightarrow \Yy$ be the
relative group corresponding to $\Pp$ with section $p$. Let $\Gg _1$ denote
the restriction of $\Gg$ to $\Yy _1$. We have to show that every $\Gg
_1$-torsor $\Jj_1$ extends to a $\Gg$-torsor $\Jj$ over $\Xx$. Suppose
$\Jj_1\rightarrow \Yy _1$ is a $\Gg _1$-torsor. For each $X\in \Xx$, choose
an element $\alpha _X\in \Jj _1(a_X)$. Put $\Jj (u)=\Jj _1(u)$ for
morphisms $u$ in $\Yy _1$. For a morphism $fb_X$ in $\Yy$ (with
$f:(X,0)\rightarrow (Y,0)$), put $$
\Jj _1 (f)\cong \Jj (fb_X)
$$
with this isomorphism written as $j\mapsto j\beta _X$. Introduce the
convention $\alpha _X\beta _X= e_{X,1}$ and $\beta _X\alpha _X = e_{X,0}$.
The
action of $\Gg (fb_X)$ given as follows. An element of $\Gg (fb_X)$ can be
written uniquely as $u1_{b_X}$ for $u\in \Gg (f)$, where $1_{b_X}$ is the
identity of $\Gg (b_X)$. Then put
$$
(u1_{b_X})\ast (j\beta _X):=
(u\ast j)\beta _X.
$$
Define the horizontal compositions as follows: $$
k(j\beta _X):= (kj)\beta _X
$$
for $k\in \Jj _1(g)$ with $g: (Y,0)\rightarrow (Z,0)$; $$
(\alpha _Zk)(j\beta _X):= l
$$
for the element $l:(X,1)\rightarrow (Z,1)$ in $\Jj _1$ such that $$
l\alpha _X = \alpha _Z kj
$$
(note that composition with $\alpha _X$ gives an isomorphism $$
Hom _{\Jj _1}((X, 1),(Z,1))\cong Hom _{\Jj _1}((X,0),(Z,1)) $$
because if $u_0$ is an element on the left, then an arbitrary element on
the right can be written uniquely as
$$
(w1_{a_X})\ast (u_0\alpha _X)
$$
for $w1_{a_X}\in \Gg $, and this in turn is equal to $(w\ast u_0)\alpha _X$); $$
(j\beta _X)(\alpha _X k):= jk
$$
for $k: (W,0)\rightarrow (X,0)$ in $\Jj _1$; and finally $$
(j\beta _X)(\alpha _X k\beta _W):= jk,
$$
noting that any morphism $(W,1)\rightarrow (X,1)$ in $\Jj _1$ can be
written as $\alpha _X k\beta _W$ by an argument similar to that
parenthetized above. The composition defined in this way is associative and
compatible with the action of $\Gg$, so we obtain a $\Gg$-torsor $\Jj$
restricting to $\Jj _1$. This translates to the statement that any section
$q$ over $\Nn$ extends to a section $p$ over $\Mm$. To show that the path
component containing $p$ is uniquely determined by the path component
containing $q$, and that the fundamental groups of these path components
are isomorphic, we have to show that for a $\Gg$-torsors $\Jj $ over $\Yy$,
the set of sections $t$ of $\Jj$ over $\Yy$ is the same as the set of
sections $t_1$ over $\Yy _1$ (then see a previous argument of changing
points of view, to see that the set of isomorphisms between two torsors is
the same over $\Yy$ and $\Yy _1$). Fix a section $t_1$ over $\Yy_1$; we
have to show that there exists a unique section $t$ restricting to $t_1$.
The uniqueness follows from the fact that any morphism in $\Yy$ (but not in
$\Yy _1$) can be written as $fb_X$; if $t$ is a section restricting to
$t_1$ on $\Yy _1$ then $$
t(fb_X)= t(f)t(b_X);
$$
but $t(a_X)t(b_X)=e_{X,1}$ and this uniquely determines $t(b_X)$; another
would be of the form $u\ast t(b_X)$ for $u\in \Gg (b_X)$, and then $$
e_{X,1}=t(a_X)(u\ast t(b_X))=(1_{a_X}u)\ast (t(a_X)t(b_X))= 1_{a_X}u $$
so $1_{b_X}=1_{b_X}1_{a_X}u=u$.
Similarly, for existence of $t$ we define $t(b_X)$ to be the unique element
such that $t(a_X)t(b_X)=e_{X,1}$, and we put $$
t(fb_X):= t_1(f)t(b_X).
$$
This is a section---the formulas can be verified by using formulas for $t_1$,
and the fact that certain formulas have unique solutions in $\Jj$.

We have now completed the treatment of the degree $1$ case: the morphism $$
Hom _{\Yy}(\Mm , \Pp )\rightarrow Hom _{\Yy } (\Nn ,\Pp ) $$
is a weak equivalence, so its fibers are nonempty and weakly contractible.

We turn to the criterion established in Lemma \ref{flex05} for degree
$m\geq 2$. We need to show that for any morphism $\phi$ in $\Yy$, the
inclusion $$
\Lambda _{\Nn}(\phi )\subset \Lambda _{\Mm} (\phi ) $$
induces an isomorphism on cohomology. But $\Lambda _{\Mm} (\phi )$ is the
nerve of the category of objects inside $\phi$. If we write $\phi = (\psi ,
c)$ for $\psi $ a morphism in $\Xx$ and $c$ a morphism in $\overline{I}$,
then $Ins _{\Yy }(\phi )$ is homotopic to $Ins _{\Xx}(\psi )$ (it is the
product of this with the nerve of $Ins _{\overline{I}}(c)\cong
\overline{I}$ which is contractible). On the other hand, $\Lambda _{\Nn}
(\phi )$ is the nerve of the category $Ins _{\Yy _1/\Yy }(\phi )$ of
objects of $\Yy_1$ which are $\Yy$-inside $\phi$. This category splits as
a product $$
Ins _{\Yy _1/\Yy }(\phi )= Ins _{\Xx}(\psi )\times Ins _{I/\overline{I}}(c). $$
The nerve is therefore the product of nerves. Since there is exactly one
morphism in $\overline{I}$ between any two objects, $$
Ins _{I/\overline{I}}(c)\cong I
$$
and its nerve is contractible. Thus $\Lambda _{\Nn}(\phi )$ is also
homotopic to $Ins _{\Xx}(\psi )$. The inclusion is compatible with these
homotopies, so it is a homotopy equivalence. By Lemma \ref{flex05}, the
maps $$
Hom _{\Yy}(\Mm , \Pp )\rightarrow Hom _{\Yy}(\Nn ,\Pp ) $$
are weak equivalences for any relative Eilenberg-MacLane category $\Pp
\rightarrow \Xx$ of degree $m\geq 2$. By Postnikov induction, for any
continuous category $\Cc$, the fiber of $$
Hom _{\Yy}(\Mm , \Cc )\rightarrow Hom _{\Yy}(\Nn ,\Cc ) $$
over an element $F$ as in the theorem, is nonempty and weakly contractible.
This completes the proof of the theorem. \eop

\newpage

\setcounter{section}{2}

\numero{A canonical homotopy equivalence}

In this section we show how to ``strictify'' a flexible functor. The first
result of this type was the strictification result for fibered categories in
SGA 1 \cite{SGA1}: any fibered category can be replaced by an equivalent split
one. This result essentially concerns homotopy-coherent diagrams of $K(\pi ,
1)$'s. The next occurrence of such a result was in Segal's paper
\cite{SegalTopology} where it is mentionned very rapidly in an appendix. While
this is not altogether evident, Segal's method is actually the same as the
method of SGA 1 (and it is this method which we shall present below). The next
occurrence in the literature is the result of Dwyer-Kan, in the context of their
``free simplicial category'' approach to the notion of homotopy coherence
\cite{DwyerKanEquivs}.  As noted above, Cordier has shown that the Dwyer-Kan
approach to homotopy-coherence is the same as the Segal-Leitch-Vogt approach
\cite{Cordier1} \cite{Cordier2}, so in a certain sense Segal's result subsumes
that of Dwyer-Kan. See also Cordier-Porter \cite{CoPo2} where Segal's
strictification result is taken up with more details.

Our treatment, which is based on the idea of SGA 1, is of course the
same as these previous ones (so there is nothing new in the present section).
Nonetheless, it will be helpful to treat explicitly the necessary homotopy
equivalences involved.

Suppose $\Xx$ is a category and $T:\Xx \rightarrow Top$ is a flexible
functor. We will define a flexible functor $K'T:\Xx \rightarrow Top$ which
is in fact a functor.

For each object $X\in \Xx$, let
$(\Xx /X)^{\ast}$ denote the category obtained by adding an initial object
denoted $\ast$ to $\Xx/X$. Define $(K'T)_X$ to be the space of
contravariant flexible
functors $P:(\Xx /X)^{\ast}\rightarrow Top$ such that $P|_{\Xx}=T$ and
$P|_{\{ \ast \} }=\ast$ is the unique flexible functor with values in the
subcategory of $Top$ consisting of the space with a single point (choose
one such space once and for all!). More concretely, $(K'T)_X$ is the space
of ways of specifying a function
$$
P(\phi _1,\ldots , \phi _k; t_0,\ldots , t_{k-1}) \in T_{X_k} $$
for a composable sequence $\phi$ in $\Xx/X$, satisfying: the usual
relations for $t_i=0$ (for $t_0=0$ we drop off $\phi _1$---note that there
is really another morphism, the morphism to the initial object, which
should be considered as at the beginning of the word $\phi$); and the
relation of composition for $t_i=1$ but where the composition is of the map
$T(\phi '',t'')$ applied to the point $P(\phi ' , t')$ (here $\phi ',t'$
are the parts before the $t_i=1$ and $\phi '',t''$ the parts after).

We will now define a homotopy equivalence $M$ between $K'T$ and $T$.

First, some notation for homotopy equivalences. Recall that a homotopy
equivalence between $S$ and $T$ is a flexible functor $M:\overline{I}\times
\Xx \rightarrow Top$, whose restrictions to $\{ 0\} \times \Xx $ and $\{
1\} \times \Xx $ are $S$ and $T$, respectively. We write this out more
explicitly. A composable sequence of maps in $ \overline{I}\times \Xx$
consists of a composable sequence $\phi _1,\ldots , \phi _n$ in $\Xx$,
together with a specification for each $X_i$ of either $0$ or $1$. We make
this specification with a supplementary list $\epsilon _1,\ldots ,\epsilon
_{n-1}$ with $\epsilon _i=0\; \mbox{or} \; 1$, as well as the specification
of $\epsilon _0$ and $\epsilon _{n}$ (also equal to $0$ or $1$). If
$\epsilon _i=0$, it means that the object $0$ is specified to go with
$X_i$, whereas if $\epsilon _i=1$ it means that object $1$ is specified to
go with $X_i$. A homotopy equivalence $M$ then consists of the
specification of
$$
M^{\epsilon _0 ,\epsilon _n} (\phi _1,\ldots , \phi _n; \epsilon _1,\ldots
, \epsilon _{n-1})(t_1,\ldots , t_{n-1})
$$
where
$$
M^{0,0}: S_{X_0}\rightarrow S_{X_n};
$$
$$
M^{0,1}: S_{X_0}\rightarrow T_{X_n};
$$
$$
M^{1,0}: T_{X_0}\rightarrow S_{X_n};
$$
and
$$
M^{1,1}: T_{X_0}\rightarrow T_{X_n}.
$$
This is subject to the following conditions. First of all, if $\epsilon
_0=\ldots =\epsilon _n =0$ then $M$ agrees with $S$, and if $\epsilon
_0=\ldots =\epsilon _n =1$ then $M$ agrees with $T$. If $t_i=0$ then
$$
M^{\epsilon _0 ,\epsilon _n} (\phi _1,\ldots , \phi _n; \epsilon _1,\ldots
, \epsilon _{n-1})(t_1,\ldots , t_{n-1})
$$
$$
=
M^{\epsilon _0 ,\epsilon _n} (\phi _1,\ldots ,\phi _{i}\phi _{i+1} ,\ldots
, \phi _n; \epsilon _1,\ldots ,\widehat{\epsilon _i},\ldots , \epsilon
_{n-1})(t_1,\ldots ,\widehat{t_i},\ldots , t_{n-1}). $$
And if $t_i=1$ then
$$
M^{\epsilon _0 ,\epsilon _n} (\phi _1,\ldots , \phi _n; \epsilon _1,\ldots
, \epsilon _{n-1})(t_1,\ldots , t_{n-1})
$$
$$
=
M^{\epsilon _i ,\epsilon _n} (\phi _{i+1},\ldots , \phi _n; \epsilon
_{i+1},\ldots , \epsilon _{n-1})(t_{i+1},\ldots , t_{n-1}) $$
$$
\circ
M^{\epsilon _0 ,\epsilon _i} (\phi _1,\ldots , \phi _{i}; \epsilon
_1,\ldots , \epsilon _{i-1})(t_1,\ldots , t_{i-1}).
$$

Now we would like to define a homotopy equivalence between $K'T$ and $T$.
Note that an element of $K'T_{X_n}$ consists of a function $P(\xi _0, \xi
_1,\ldots , \xi _k)(s_1,\ldots , s_k)$, satisfying the usual properties
with respect to $T$, where the arguments $\xi _i : Y_i\rightarrow Y_{i-1}$
form a composable sequence in $\Xx$ and $\xi _0:Y_0\rightarrow X_n$ (thus,
$\xi _0$ gives $\xi _1,\ldots , \xi _k$ a structure of composable sequence
in $\Xx /X_n$). Note that $$
K'T (\phi _1,\ldots , \phi _n)(t_1,\ldots , t_{n-1})(P)(\xi _0,\ldots , \xi
_k)(s_1,\ldots , s_k)
$$
$$
:= P(\phi _1\cdots \phi _n\xi _0, \xi _1,\ldots , \xi _k)(s_1,\ldots , s_k).
$$

\small

We define our homotopy equivalence $M$ as follows, with $P$ denoting a
point of $K'T_{X_0}$ and $p$ denoting a point of $T_{X_0}$. $$
M^{0,0}(\phi _1,\ldots , \phi _n; \epsilon _1,\ldots , \epsilon
_{n-1})(t_1,\ldots , t_{n-1})(P)(\xi _0,\ldots , \xi _k )(s_1 ,\ldots ,
s_k) $$
$$
:=
P(\phi _1,\ldots , \phi _{n-1} , \phi _n\xi _0 , \xi _1,\ldots , \xi
_k)(\epsilon _1t_1,\ldots , \epsilon _{n-1}t_{n-1} , s_1,\ldots , s_k ). $$
$$
M^{0,1}(\phi _1,\ldots , \phi _n; \epsilon _1,\ldots , \epsilon
_{n-1})(t_1,\ldots , t_{n-1})(P)
$$
$$
:=
P(\phi _1,\ldots , \phi _{n-1} , \phi _n)(\epsilon _1t_1,\ldots , \epsilon
_{n-1}t_{n-1} ).
$$
$$
M^{1,0}(\phi _1,\ldots , \phi _n; \epsilon _1,\ldots , \epsilon
_{n-1})(t_1,\ldots , t_{n-1})(p)(\xi _0,\ldots , \xi _k )(s_1 ,\ldots ,
s_k) $$
$$
:=
T(\phi _1,\ldots , \phi _{n-1} , \phi _n\xi _0 , \xi _1,\ldots , \xi
_k)(\epsilon _1t_1,\ldots , \epsilon _{n-1}t_{n-1} , s_1,\ldots , s_k )(p).
$$
$$
M^{1,1}(\phi _1,\ldots , \phi _n; \epsilon _1,\ldots , \epsilon
_{n-1})(t_1,\ldots , t_{n-1})(p)
$$
$$
:=
T(\phi _1,\ldots , \phi _{n-1} , \phi _n)(\epsilon _1t_1,\ldots , \epsilon
_{n-1}t_{n-1} )(p).
$$
We have to check the properties of the points in $K'T$ resulting from $M$,
in the first and third cases; and in all cases we have to check the
properties of $M$. Note that the required properties for $s_i=0$ or for
$t_i=0$ are easily verified by looking at the formulas (note that with our
new notation for $P$, the relationship between the places of the $\phi _i$
and the $t_i$ is the same as for $T$, $M$, etc.). The required properties
for $s_i=1$ are similarly easily verified.

{\em The case of $M^{0,0}$, with $t_i=1$ and $\epsilon _i=0$:} In this case
we have
$$
M^{0,0}(\phi _1,\ldots , \phi _n; \epsilon _1,\ldots , \epsilon
_{n-1})(t_1,\ldots , t_{n-1})(P)(\xi _0,\ldots , \xi _k )(s_1 ,\ldots ,
s_k) $$
$$
=
P(\phi _1,\ldots , \phi _{n-1} , \phi _n\xi _0 , \xi _1,\ldots , \xi
_k)(\epsilon _1t_1,\ldots , \epsilon _{n-1}t_{n-1} , s_1,\ldots , s_k ) $$
$$
=
P(\phi _1,\ldots , \phi _{i-1}, \phi _i\phi _{i+1} , \ldots , \phi _{n-1} ,
\phi _n\xi _0 , \xi _1,\ldots , \xi _k)(\epsilon _1t_1,\ldots ,
\widehat{\epsilon _{i}t_{i}}, \ldots , \epsilon _{n-1}t_{n-1} , s_1,\ldots
, s_k ). $$
On the other hand,
$$
M^{0,0} (\phi _{i+1},\ldots , \phi _n; \epsilon _{i+1},\ldots , \epsilon
_{n-1})(t_{i+1},\ldots , t_{n-1}) $$
$$
\circ
M^{0,0} (\phi _1,\ldots , \phi _{i}; \epsilon _1,\ldots , \epsilon
_{i-1})(t_1,\ldots , t_{i-1})
(P)(\xi _0,\ldots , \xi _k )(s_1 ,\ldots , s_k) $$
$$
= M^{0,0} (\phi _{i+1},\ldots , \phi _n; \epsilon _{i+1},\ldots , \epsilon
_{n-1})(t_{i+1},\ldots , t_{n-1}) $$
$$
(M^{0,0} (\phi _1,\ldots , \phi _{i}; \epsilon _1,\ldots , \epsilon
_{i-1})(t_1,\ldots , t_{i-1})
(P))(\xi _0,\ldots , \xi _k )(s_1 ,\ldots , s_k) $$
$$
=
M^{0,0} (\phi _1,\ldots , \phi _{i}; \epsilon _1,\ldots , \epsilon
_{i-1})(t_1,\ldots , t_{i-1})
(P)
$$
$$
(\phi _{i+1},\ldots , \phi _{n-1} , \phi _n\xi _0 , \xi _1,\ldots , \xi
_k)(\epsilon _{i+1}t_{i+1},\ldots , \epsilon _{n-1}t_{n-1} , s_1,\ldots ,
s_k ) $$
$$
=
P(\phi _1,\ldots , \phi _{i-1}, \phi _i\phi _{i+1} , \ldots , \phi _{n-1} ,
\phi _n\xi _0 , \xi _1,\ldots , \xi _k)
$$
$$
(\epsilon _1t_1,\ldots , \epsilon _{i-1}t_{i-1}, \epsilon
_{i+1}t_{i+1},\ldots , \epsilon _{n-1}t_{n-1} , s_1,\ldots , s_k ) $$
$$
=
P(\phi _1,\ldots , \phi _{i-1}, \phi _i\phi _{i+1} , \ldots , \phi _{n-1} ,
\phi _n\xi _0 , \xi _1,\ldots , \xi _k)
$$
$$
(\epsilon _1t_1,\ldots , \widehat{\epsilon _{i}t_{i}}, \ldots , \epsilon
_{n-1}t_{n-1} , s_1,\ldots , s_k ) $$
as desired.

{\em The case of $M^{0,0}$, with $t_i=1$ and $\epsilon _i=1$:} In this case
we have
$$
M^{0,0}(\phi _1,\ldots , \phi _n; \epsilon _1,\ldots , \epsilon
_{n-1})(t_1,\ldots , t_{n-1})(P)(\xi _0,\ldots , \xi _k )(s_1 ,\ldots ,
s_k) $$
$$
=
P(\phi _1,\ldots , \phi _{n-1} , \phi _n\xi _0 , \xi _1,\ldots , \xi
_k)(\epsilon _1t_1,\ldots , \epsilon _{n-1}t_{n-1} , s_1,\ldots , s_k ) $$
$$
=
T(\phi _{i+1} ,
\ldots , \phi _{n-1} , \phi _n\xi _0 , \xi _1,\ldots , \xi _k)(\epsilon
_{i+1}t_{i+1},
\ldots , \epsilon _{n-1}t_{n-1} , s_1,\ldots , s_k ) $$
$$
(P(\phi _1,\ldots , \phi _{i})(\epsilon _1t_1,\ldots ,\epsilon
_{i-1}t_{i-1})). $$
On the other hand,
$$
M^{1,0} (\phi _{i+1},\ldots , \phi _n; \epsilon _{i+1},\ldots , \epsilon
_{n-1})(t_{i+1},\ldots , t_{n-1}) $$
$$
\circ
M^{0,1} (\phi _1,\ldots , \phi _{i}; \epsilon _1,\ldots , \epsilon
_{i-1})(t_1,\ldots , t_{i-1})
$$
$$
(P)(\xi _0,\ldots , \xi _k )(s_1 ,\ldots , s_k) $$
$$
= M^{1,0} (\phi _{i+1},\ldots , \phi _n; \epsilon _{i+1},\ldots , \epsilon
_{n-1})(t_{i+1},\ldots , t_{n-1}) $$
$$
(M^{0,1} (\phi _1,\ldots , \phi _{i}; \epsilon _1,\ldots , \epsilon
_{i-1})(t_1,\ldots , t_{i-1})
(P))(\xi _0,\ldots , \xi _k )(s_1 ,\ldots , s_k) $$
$$
=
T(\phi _{i+1},\ldots , \phi _{n-1} , \phi _n\xi _0 , \xi _1,\ldots , \xi
_k)(\epsilon _{i+1}t_{i+1},\ldots , \epsilon _{n-1}t_{n-1} , s_1,\ldots ,
s_k ) $$
$$
(M^{0,1} (\phi _1,\ldots , \phi _{i}; \epsilon _1,\ldots , \epsilon
_{i-1})(t_1,\ldots , t_{i-1})
(P))
$$
$$
=
T(\phi _{i+1},\ldots , \phi _{n-1} , \phi _n\xi _0 , \xi _1,\ldots , \xi
_k)(\epsilon _{i+1}t_{i+1},\ldots , \epsilon _{n-1}t_{n-1} , s_1,\ldots ,
s_k ) $$
$$
(P(\phi _1,\ldots , \phi _{i})(\epsilon _1t_1,\ldots , \epsilon _{i-1}t_{i-1}))
$$
as desired.

{\em The case of $M^{0,1}$, with $t_i=1$ and $\epsilon _i=0$:} In this case
we have
$$
M^{0,1}(\phi _1,\ldots , \phi _n; \epsilon _1,\ldots , \epsilon
_{n-1})(t_1,\ldots , t_{n-1})(P)
$$
$$
=
P(\phi _1,\ldots , \phi _{n-1} , \phi _n)(\epsilon _1t_1,\ldots , \epsilon
_{n-1}t_{n-1} ) $$
$$
=
P(\phi _1,\ldots , \phi _{i-1}, \phi _i\phi _{i+1} , \ldots , \phi _{n-1} ,
\phi _n)(\epsilon _1t_1,\ldots , \widehat{\epsilon _{i}t_{i}}, \ldots ,
\epsilon _{n-1}t_{n-1} ). $$
On the other hand,
$$
M^{0,1} (\phi _{i+1},\ldots , \phi _n; \epsilon _{i+1},\ldots , \epsilon
_{n-1})(t_{i+1},\ldots , t_{n-1}) $$
$$
\circ
M^{0,0} (\phi _1,\ldots , \phi _{i}; \epsilon _1,\ldots , \epsilon
_{i-1})(t_1,\ldots , t_{i-1})
(P)
$$
$$
= M^{0,1} (\phi _{i+1},\ldots , \phi _n; \epsilon _{i+1},\ldots , \epsilon
_{n-1})(t_{i+1},\ldots , t_{n-1}) $$
$$
(M^{0,0} (\phi _1,\ldots , \phi _{i}; \epsilon _1,\ldots , \epsilon
_{i-1})(t_1,\ldots , t_{i-1})
(P))
$$
$$
=
M^{0,0} (\phi _1,\ldots , \phi _{i}; \epsilon _1,\ldots , \epsilon
_{i-1})(t_1,\ldots , t_{i-1})
(P)
$$
$$
(\phi _{i+1},\ldots , \phi _{n-1} , \phi _n)(\epsilon _{i+1}t_{i+1},\ldots
, \epsilon _{n-1}t_{n-1} ) $$
$$
=
P(\phi _1,\ldots , \phi _{i-1}, \phi _i\phi _{i+1} , \ldots , \phi _{n-1} ,
\phi _n)
$$
$$
(\epsilon _1t_1,\ldots , \epsilon _{i-1}t_{i-1}, \epsilon
_{i+1}t_{i+1},\ldots , \epsilon _{n-1}t_{n-1} ) $$
$$
=
P(\phi _1,\ldots , \phi _{i-1}, \phi _i\phi _{i+1} , \ldots , \phi _{n-1} ,
\phi _n)(\epsilon _1t_1,\ldots , \widehat{\epsilon _{i}t_{i}}, \ldots ,
\epsilon _{n-1}t_{n-1}) $$
as desired.

{\em The case of $M^{0,1}$, with $t_i=1$ and $\epsilon _i=1$:} In this case
we have
$$
M^{0,1}(\phi _1,\ldots , \phi _n; \epsilon _1,\ldots , \epsilon
_{n-1})(t_1,\ldots , t_{n-1})(P)
$$
$$
=
P(\phi _1,\ldots , \phi _n )(\epsilon _1t_1,\ldots , \epsilon _{n-1}t_{n-1}) $$
$$
=
T(\phi _{i+1} ,
\ldots , \phi _{n-1} , \phi _n )(\epsilon _{i+1}t_{i+1}, \ldots , \epsilon
_{n-1}t_{n-1} )
$$
$$
(P(\phi _1,\ldots , \phi _{i})(\epsilon _1t_1,\ldots ,\epsilon
_{i-1}t_{i-1})). $$
On the other hand,
$$
M^{1,1} (\phi _{i+1},\ldots , \phi _n; \epsilon _{i+1},\ldots , \epsilon
_{n-1})(t_{i+1},\ldots , t_{n-1}) $$
$$
\circ
M^{0,1} (\phi _1,\ldots , \phi _{i}; \epsilon _1,\ldots , \epsilon
_{i-1})(t_1,\ldots , t_{i-1})
(P)
$$
$$
= M^{1,1} (\phi _{i+1},\ldots , \phi _n; \epsilon _{i+1},\ldots , \epsilon
_{n-1})(t_{i+1},\ldots , t_{n-1}) $$
$$
(M^{0,1} (\phi _1,\ldots , \phi _{i}; \epsilon _1,\ldots , \epsilon
_{i-1})(t_1,\ldots , t_{i-1})
(P))
$$
$$
=
T(\phi _{i+1},\ldots , \phi _{n} )(\epsilon _{i+1}t_{i+1},\ldots , \epsilon
_{n-1}t_{n-1} )
$$
$$
(M^{0,1} (\phi _1,\ldots , \phi _{i};
\epsilon _1,\ldots , \epsilon _{i-1})(t_1,\ldots , t_{i-1}) (P))
$$
$$
=
T(\phi _{i+1},\ldots , \phi _{n})(\epsilon _{i+1}t_{i+1},\ldots , \epsilon
_{n-1}t_{n-1})
$$
$$
(P(\phi _1,\ldots , \phi _{i})(\epsilon
_1t_1,\ldots , \epsilon _{i-1}t_{i-1}))
$$
as desired.

{\em The case of $M^{1,0}$, with $t_i=1$ and $\epsilon _i=0$:} In this case
we have
$$
M^{1,0}(\phi _1,\ldots , \phi _n; \epsilon _1,\ldots , \epsilon
_{n-1})(t_1,\ldots , t_{n-1})(p)(\xi _0,\ldots , \xi _k )(s_1 ,\ldots ,
s_k) $$
$$
=
T(\phi _1,\ldots , \phi _{n-1} , \phi _n\xi _0 , \xi _1,\ldots , \xi
_k)(\epsilon _1t_1,\ldots , \epsilon _{n-1}t_{n-1} , s_1,\ldots , s_k )(p)
$$
$$
=
T(\phi _1,\ldots , \phi _{i-1}, \phi _i\phi _{i+1} , \ldots , \phi _{n-1} ,
\phi _n\xi _0 , \xi _1,\ldots , \xi _k)(\epsilon _1t_1,\ldots ,
\widehat{\epsilon _{i}t_{i}}, \ldots , \epsilon _{n-1}t_{n-1} , s_1,\ldots
, s_k )(p). $$
On the other hand,
$$
M^{0,0} (\phi _{i+1},\ldots , \phi _n; \epsilon _{i+1},\ldots , \epsilon
_{n-1})(t_{i+1},\ldots , t_{n-1}) $$
$$
\circ
M^{1,0} (\phi _1,\ldots , \phi _{i}; \epsilon _1,\ldots , \epsilon
_{i-1})(t_1,\ldots , t_{i-1})
(p)(\xi _0,\ldots , \xi _k )(s_1 ,\ldots , s_k) $$
$$
= M^{0,0} (\phi _{i+1},\ldots , \phi _n; \epsilon _{i+1},\ldots , \epsilon
_{n-1})(t_{i+1},\ldots , t_{n-1}) $$
$$
(M^{1,0} (\phi _1,\ldots , \phi _{i}; \epsilon _1,\ldots , \epsilon
_{i-1})(t_1,\ldots , t_{i-1})
(p))(\xi _0,\ldots , \xi _k )(s_1 ,\ldots , s_k) $$
$$
=
M^{1,0} (\phi _1,\ldots , \phi _{i}; \epsilon _1,\ldots , \epsilon
_{i-1})(t_1,\ldots , t_{i-1})
(p)
$$
$$
(\phi _{i+1},\ldots , \phi _{n-1} , \phi _n\xi _0 , \xi _1,\ldots , \xi
_k)(\epsilon _{i+1}t_{i+1},\ldots , \epsilon _{n-1}t_{n-1} , s_1,\ldots ,
s_k ) $$
$$
=
T(\phi _1,\ldots , \phi _{i-1}, \phi _i\phi _{i+1} , \ldots , \phi _{n-1} ,
\phi _n\xi _0 , \xi _1,\ldots , \xi _k)
$$
$$
(\epsilon _1t_1,\ldots , \epsilon _{i-1}t_{i-1}, \epsilon
_{i+1}t_{i+1},\ldots , \epsilon _{n-1}t_{n-1} , s_1,\ldots , s_k )(p) $$
$$
=
T(\phi _1,\ldots , \phi _{i-1}, \phi _i\phi _{i+1} , \ldots , \phi _{n-1} ,
\phi _n\xi _0 , \xi _1,\ldots , \xi _k)
$$
$$
(\epsilon _1t_1,\ldots , \widehat{\epsilon _{i}t_{i}}, \ldots , \epsilon
_{n-1}t_{n-1} , s_1,\ldots , s_k )(p) $$
as desired.

{\em The case of $M^{1,0}$, with $t_i=1$ and $\epsilon _i=1$:} In this case
we have
$$
M^{1,0}(\phi _1,\ldots , \phi _n; \epsilon _1,\ldots , \epsilon
_{n-1})(t_1,\ldots , t_{n-1})(p)(\xi _0,\ldots , \xi _k )(s_1 ,\ldots ,
s_k) $$
$$
=
T(\phi _1,\ldots , \phi _{n-1} , \phi _n\xi _0 , \xi _1,\ldots , \xi
_k)(\epsilon _1t_1,\ldots , \epsilon _{n-1}t_{n-1} , s_1,\ldots , s_k )(p)
$$
$$
=
T(\phi _{i+1} ,
\ldots , \phi _{n-1} , \phi _n\xi _0 , \xi _1,\ldots , \xi _k)(\epsilon
_{i+1}t_{i+1},
\ldots , \epsilon _{n-1}t_{n-1} , s_1,\ldots , s_k ) $$
$$
(T(\phi _1,\ldots , \phi _{i})(\epsilon _1t_1,\ldots ,\epsilon
_{i-1}t_{i-1})(p)).
$$
On the other hand,
$$
M^{1,0} (\phi _{i+1},\ldots , \phi _n; \epsilon _{i+1},\ldots , \epsilon
_{n-1})(t_{i+1},\ldots , t_{n-1}) $$
$$
\circ
M^{1,1} (\phi _1,\ldots , \phi _{i}; \epsilon _1,\ldots , \epsilon
_{i-1})(t_1,\ldots , t_{i-1})
(p)(\xi _0,\ldots , \xi _k )(s_1 ,\ldots , s_k) $$
$$
= M^{1,0} (\phi _{i+1},\ldots , \phi _n; \epsilon _{i+1},\ldots , \epsilon
_{n-1})(t_{i+1},\ldots , t_{n-1}) $$
$$
(M^{1,1} (\phi _1,\ldots , \phi _{i}; \epsilon _1,\ldots , \epsilon
_{i-1})(t_1,\ldots , t_{i-1})
(p))(\xi _0,\ldots , \xi _k )(s_1 ,\ldots , s_k) $$
$$
=
T(\phi _{i+1},\ldots , \phi _{n-1} , \phi _n\xi _0 , \xi _1,\ldots , \xi
_k)(\epsilon _{i+1}t_{i+1},\ldots , \epsilon _{n-1}t_{n-1} , s_1,\ldots ,
s_k ) $$
$$
(M^{1,1} (\phi _1,\ldots , \phi _{i}; \epsilon _1,\ldots , \epsilon
_{i-1})(t_1,\ldots , t_{i-1})
(p))
$$
$$
=
T(\phi _{i+1},\ldots , \phi _{n-1} , \phi _n\xi _0 , \xi _1,\ldots , \xi
_k)(\epsilon _{i+1}t_{i+1},\ldots , \epsilon _{n-1}t_{n-1} , s_1,\ldots ,
s_k ) $$
$$
(T(\phi _1,\ldots , \phi _{i})(\epsilon _1t_1,\ldots , \epsilon
_{i-1}t_{i-1})(p))
$$
as desired.

{\em The case of $M^{1,1}$, with $t_i=1$ and $\epsilon _i=0$:} In this case
we have
$$
M^{1,1}(\phi _1,\ldots , \phi _n; \epsilon _1,\ldots , \epsilon
_{n-1})(t_1,\ldots , t_{n-1})(p)
$$
$$
=
T(\phi _1,\ldots , \phi _{n-1} , \phi _n)(\epsilon _1t_1,\ldots , \epsilon
_{n-1}t_{n-1} )(p) $$
$$
=
T(\phi _1,\ldots , \phi _{i-1}, \phi _i\phi _{i+1} , \ldots , \phi _{n-1} ,
\phi _n)(\epsilon _1t_1,\ldots , \widehat{\epsilon _{i}t_{i}}, \ldots ,
\epsilon _{n-1}t_{n-1} )(p). $$
On the other hand,
$$
M^{0,1} (\phi _{i+1},\ldots , \phi _n; \epsilon _{i+1},\ldots , \epsilon
_{n-1})(t_{i+1},\ldots , t_{n-1}) $$
$$
\circ
M^{1,0} (\phi _1,\ldots , \phi _{i}; \epsilon _1,\ldots , \epsilon
_{i-1})(t_1,\ldots , t_{i-1})
(p)
$$
$$
= M^{0,1} (\phi _{i+1},\ldots , \phi _n; \epsilon _{i+1},\ldots , \epsilon
_{n-1})(t_{i+1},\ldots , t_{n-1}) $$
$$
(M^{1,0} (\phi _1,\ldots , \phi _{i}; \epsilon _1,\ldots , \epsilon
_{i-1})(t_1,\ldots , t_{i-1})
(p))
$$
$$
=
M^{1,0} (\phi _1,\ldots , \phi _{i}; \epsilon _1,\ldots , \epsilon
_{i-1})(t_1,\ldots , t_{i-1})
(p)
$$
$$
(\phi _{i+1},\ldots , \phi _{n-1} , \phi _n)(\epsilon _{i+1}t_{i+1},\ldots
, \epsilon _{n-1}t_{n-1} )
$$
$$
=
T(\phi _1,\ldots , \phi _{i-1}, \phi _i\phi _{i+1} , \ldots , \phi _{n-1} ,
\phi _n)
$$
$$
(\epsilon _1t_1,\ldots , \epsilon _{i-1}t_{i-1}, \epsilon
_{i+1}t_{i+1},\ldots , \epsilon _{n-1}t_{n-1} )(p) $$
$$
=
T(\phi _1,\ldots , \phi _{i-1}, \phi _i\phi _{i+1} , \ldots , \phi _{n-1} ,
\phi _n)
$$
$$
(\epsilon _1t_1,\ldots , \widehat{\epsilon _{i}t_{i}}, \ldots , \epsilon
_{n-1}t_{n-1})(p) $$
as desired.

{\em The case of $M^{1,1}$, with $t_i=1$ and $\epsilon _i=1$:} In this case
we have
$$
M^{0,1}(\phi _1,\ldots , \phi _n; \epsilon _1,\ldots , \epsilon
_{n-1})(t_1,\ldots , t_{n-1})(p)
$$
$$
=
T(\phi _1,\ldots , \phi _n )(\epsilon _1t_1,\ldots , \epsilon
_{n-1}t_{n-1})(p) $$
$$
=
T(\phi _{i+1} ,
\ldots , \phi _{n-1} , \phi _n )(\epsilon _{i+1}t_{i+1}, \ldots , \epsilon
_{n-1}t_{n-1} )
$$
$$
(T(\phi _1,\ldots , \phi _{i})(\epsilon _1t_1,\ldots ,\epsilon
_{i-1}t_{i-1})(p)).
$$
On the other hand,
$$
M^{1,1} (\phi _{i+1},\ldots , \phi _n; \epsilon _{i+1},\ldots , \epsilon
_{n-1})(t_{i+1},\ldots , t_{n-1}) $$
$$
\circ
M^{1,1} (\phi _1,\ldots , \phi _{i}; \epsilon _1,\ldots , \epsilon
_{i-1})(t_1,\ldots , t_{i-1})
(p)
$$
$$
= M^{1,1} (\phi _{i+1},\ldots , \phi _n; \epsilon _{i+1},\ldots , \epsilon
_{n-1})(t_{i+1},\ldots , t_{n-1}) $$
$$
(M^{1,1} (\phi _1,\ldots , \phi _{i}; \epsilon _1,\ldots , \epsilon
_{i-1})(t_1,\ldots , t_{i-1})
(p))
$$
$$
=
T(\phi _{i+1},\ldots , \phi _{n} )(\epsilon _{i+1}t_{i+1},\ldots , \epsilon
_{n-1}t_{n-1} )
$$
$$
(M^{1,1} (\phi _1,\ldots , \phi _{i};
\epsilon _1,\ldots , \epsilon _{i-1})(t_1,\ldots , t_{i-1}) (p))
$$
$$
=
T(\phi _{i+1},\ldots , \phi _{n})(\epsilon _{i+1}t_{i+1},\ldots , \epsilon
_{n-1}t_{n-1})
$$
$$
(T(\phi _1,\ldots , \phi _{i})(\epsilon
_1t_1,\ldots , \epsilon _{i-1}t_{i-1})(p)) $$
as desired.

Finally we have to check that $M$ restricts to $K'T$ and $T$. At $\epsilon
=0$ we have
$$
M^{0,0}(\phi _1,\ldots , \phi _n; 0,\ldots , 0)(t_1,\ldots ,t_{n-1}(P)(\xi
_0,\ldots , \xi _k)(s_1,\ldots , s_k)
$$
$$
=
P(\phi _1, \ldots , \phi _n\xi _0, \ldots , \xi _k)(0,\ldots , 0,s_1,\ldots
, s_k)
$$
$$
=
P(\phi _1 \cdots \phi _n\xi _0, \ldots , \xi _k)(s_1,\ldots , s_k)
$$
$$
= K'T(\phi _1,\ldots , \phi _n)(t_1,\ldots , t_{n-1})(P) (\xi _0, \ldots ,
\xi _k)(s_1,\ldots ,
s_k)
$$
as required. At $\epsilon =1$ we have
$$
M^{1,1}(\phi _1,\ldots , \phi _n; 1,\ldots , 1)(t_1,\ldots , t_{n-1})(p) $$
$$
=
T(\phi _1,\ldots , \phi _n)(t_1,\ldots , t_{n-1})(p) $$
as required.

\normalsize

This completes the verification that $M$ provides a homotopy equivalence
between $K'T$ and $T$. We obtain the conclusion that every flexible functor
is homotopy equivalent to a functor.

\subnumero{Equivalence between $K'(G^{\ast}T)$ and $G^{\ast}(K'T)$} Suppose
$G:\Yy \rightarrow \Xx $ is a functor, and $T:\Xx \rightarrow Top$ is a
flexible functor. We obtain two functors, $K'(G^{\ast}T)$ and
$G^{\ast}(K'T)$, from $\Yy $ to $Top$. We would like to see that they are
equivalent, together with their natural morphisms from $T$.

There is a natural morphism of functors $\Psi : G^{\ast}(K'T) \rightarrow
K'(G^{\ast}T)$ defined by
$$
\Psi _S (P)(\xi _0,\ldots , \xi _n)(s_1,\ldots , s_n) $$
$$
:=
P(G(\xi _0),\ldots , G(\xi _n))(s_1,\ldots , s_n) $$
for $S\in \Cc$, $P\in K'T_{G(S)}$, $\xi _0:Y_0\rightarrow S$, and $\xi
_i:Y_i\rightarrow Y_{i-1}$ a composable sequence in $\Yy$ (giving $G(\xi
_0):G(Y_0)\rightarrow G(S)$, and $G(\xi _i):G(Y_i)\rightarrow G(Y_{i-1})$ a
composable sequence in $\Xx$). This gives a morphism of flexible functors
in a trivial way. We have equivalences $M$ between $K'(G^{\ast}T)$ and
$G^{\ast}T$, and $G^{\ast}N$ between $G^{\ast}(K'T)$ and $G^{\ast}T$. We
would like to define a flexible functor $\Xi:\Cc \times I_{u\leftarrow v}
\times I_{a\leftrightarrow b}\rightarrow Top$ such that $\Xi _a$ is the
identity of $G^{\ast}T$, and $\Xi _b$ is the morphism $\Psi$ defined above;
and $\Xi ^u$ is the homotopy equivalence $G^{\ast}N$, while $\Xi ^v$ is the
homotopy equivalence $M$.

Define maps
$$
\Xi ^{\epsilon _0, \epsilon _n}_{\delta _0, \delta _n}(\phi , \epsilon ',
\delta '
)(t)
$$
where $\phi = \phi _1,\ldots , \phi _n$, $\epsilon '= \epsilon _1,\ldots ,
\epsilon _{n-1}$, $\delta '= \delta _1,\ldots , \delta _{n-1}$ and $t=t_1,\ldots
, t_{n-1}$. The indices $\epsilon _i$ and $\delta _i$ indicate the functor as
follows: $$ \epsilon =0, \, \delta =0\, : \;\; G^{\ast}(K'T) $$ $$
\epsilon =0, \, \delta =1\, : \;\; K'(G^{\ast}T) $$
$$
\epsilon =1, \, \delta =0\, : \;\; G^{\ast}T $$
$$
\epsilon =1, \, \delta =1\, : \;\; G^{\ast}T. $$
There is the constraint on $\delta$ that $\delta _i\geq \delta _{i-1}$ (in
other words, $\delta$ consists of a string of $0$'s followed by a string of
$1$'s). In particular, there is no map of the form $\Xi ^{\epsilon
_0,\epsilon _n}_{1,0}$.

In defining the maps $\Xi$, use the same formulas as above, but noting that
the arguments are a composable sequence in $\Yy$; apply $G$ if necessary to
obtain a composable sequence in $\Xx$. In the case of morphisms finishing
at $G^{\ast}(K'T)$, there are some arguments consisting of a composable
sequence in $\Xx$. In this case, the condition on $\delta$ implies that the
input cannot consist of a point $P\in K'(G^{\ast}T)_{X_0}$---which would
take as argument only a composable sequence in $\Yy$. The definitions and
verifications otherwise follow the same pattern as for $M$ defined above.

Thus the natural morphism $\Xi$ is equivalent to the identity morphism via
the equivalences $M$ and $G^{\ast}N$.

{\em Remark:} The identity morphism $I:T\rightarrow T$ is defined by $$
I^{\epsilon _0, \epsilon _n}(\phi _1,\ldots , \phi _n; \epsilon _1,\ldots ,
\epsilon _{n-1})(t_1,\ldots ,t_{n-1})
$$
$$
:=T(\phi _1,\ldots , \phi _n)(t_1,\ldots , t_{n-1}). $$

\subnumero{Philosophical remarks}
The results of this section show that, in a certain sense, it wasn't
necessary to consider flexible functors, since a flexible functor can
always be replaced by an equivalent functor. However, the notion of
flexible functor was used in defining the notion of morphism of flexible
functors. For two given functors, the space of morphisms of functors will
probably not be the same as the space of morphisms of flexible
functors---it seems clear that the latter is the right choice. We could, of
course, say that we will only look at functors, but look at morphisms
between them which are morphisms of flexible functors. By the results of
this section, this will give an equivalent theory. In fact, we often need
to use the fact that a flexible functor can be replaced by a functor, in
the arguments to follow.

The best justification for
introducing the notion of flexible functor is that it leads the way in
finding out the right properties, definitions, and arguments: for example,
the notion of morphism of flexible functor is the obvious one which goes
with the notion of flexible functor, but it would not necessarily have been
obvious if we started out looking only at functors.

\newpage

\setcounter{section}{3}

\numero{Descent conditions and sheaves}

In this section we give a treatment of descent and ``higher sheaves''
which is alternative to the closed-model-category treatment started by K.
Brown \cite{BrownThesis} an mostly due to Joyal \cite{Joyal}, Jardine
\cite{Jardine1} \cite{Jardine2}, and Thomason \cite{Thomason1}. The end
result is
basically the same as in those references.

\subnumero{$\infty$-truncation}
Recall that $\tau _{\infty}: Top \rightarrow Top$ denotes the functor
sending a topological space $A$ to the realization of the simplicial set of
singular simplices of $A$.
There is a natural morphism $\tau _{\infty}A\rightarrow A$, inducing an
isomorphism on homotopy groups (in other words, a weak equivalence). Recall also
that $\tau _{\infty}$ is compatible with products (since the realization of the
product of two simplicial sets is the product of the realizations).

If $T$ is a functor from $\Xx$ to $Top$ then we obtain a functor $\tau
_{\infty}T$. By the results of the previous section, we can replace a
flexible functor $T$ by the functor $KT$ and then take $\tau _{\infty}KT$.

We often use a subscript $\infty$ to denote the result of applying the
construction $\tau _{\infty}$. For example, $$
Mor _{\infty}(S,T):= \tau _{\infty}Mor (S,T). $$

The interest of this construction is that if $A$ is a topological space,
then $\tau _{\infty}A$ provides a CW-complex with the same weak homotopy
type as $A$, in a canonical way.

\subnumero{Global sections}
Let $\Xx ^{\ast}$ denote the category obtained by adding an initial object
denoted $\ast$ to $\Xx$. Suppose $T:\Xx \rightarrow Top $ is a
contravariant flexible functor (we sometimes call this a {\em flexible
presheaf}). We define a first version of the {\em space of global sections}
$\Gamma '(\Xx ,T)$ to be the space of contravariant flexible
functors $T':\Xx ^{\ast}\rightarrow Top$ such that $T'|_{\Xx}=T$ and
$T'|_{\{ \ast \} }=\ast$ is the unique flexible functor with values in the
subcategory of $Top$ consisting of the space with a single point.

Put $\Gamma '_{\infty}(\Xx , T):= \tau _{\infty} \Gamma '(X,T)$.

This space is the ``homotopy-limit'' or {\em holim} \cite{BousfieldKan}
\cite{Vogt1} \cite{Vogt2} of the homotopy-coherent diagram $T$.

Let $\ast _{\Xx}$ denote the unique flexible functor from $\Xx$ into the
subcategory of $Top$ consisting of the space with a single point. There is
a functor $\rho :\Xx \times I \rightarrow \Xx ^{\ast}$ sending $\Xx \times
\{ 0\}$ isomorphically to $\Xx$ and $\Xx \times \{ 1\}$ to $\ast$. If $\eta
\in \Gamma '
(\Xx ,T)$ then $\rho ^{\ast}\eta$ is a morphism of flexible functors from
$\ast _{\Xx}$ to $T$.

\begin{lemma}
\mylabel{desc01}
The map $\eta \mapsto \rho ^{\ast}\eta$ induces a weak equivalence $$
\rho ^{\ast}: \Gamma '(\Xx ,T)\rightarrow Mor (\ast _{\Xx }, T) $$
to the space of morphisms of flexible functors from $\ast _{\Xx }$ to $T$.
In particular, it induces an equivalence $\Gamma '_{\infty}(\Xx
,T)\rightarrow Mor _{\infty}(\ast _{\Xx }, T)$.
\end{lemma}
{\em Proof:}
Let $\Yy$ denote the category $\Xx \times I \cup \{ e\} $ where $e$ is a
final object with no morphisms to other objects. Let $\Xx _0$ and $\Xx _1$
denote the two subcategories of $\Xx \times I$ isomorphic to $\Xx$, and let
$\Xx ' = \Xx _0\cup \{ e\}$ and $\Xx '' = \Xx _1\cup \{ e\}$. Note that
$\Xx '\cong \Xx '' \cong \Xx \cup \{ e\}$. Let $$
Flex ^R (\Yy ' , Top )
$$
denote the
space of flexible functors from a subcategory $\Yy ' \subset \Yy$
satisfying the following restrictions: that on $\Yy ' \cap \Xx _1$, it
should be equal to $T$; and that elsewhere it should be equal to $\ast$. We
consider the following maps
$$
Flex ^R (\Xx '' , Top ) \stackrel{A}{\leftarrow} Flex ^R (\Yy , Top )
\stackrel{B}{\rightarrow} Flex ^R ( \Xx \times I , Top ). $$
We will show that they are weak equivalences. The pullback via the morphism
$\Yy \rightarrow \Xx ''$ which sends $\Xx '$ to $e$, is a half-inverse to
$A$, so if we show that $A$ is a weak equivalence then this pullback will
be too; and the morphism in question is the composition of $B$ with this
pullback. Thus this will complete the proof.

We treat the map $A$ first. Suppose that $F''\in Flex ^R(\Xx '', Top)$. We
would like to show that the space of elements in $Flex ^R (\Yy , Top )$
restricting to $F''$ is weakly contractible. Let $\Yy ^o$ denote the
subcategory of $\Yy$ which has the same objects, but whose set of morphisms
is the union of the sets of morphisms of $\Xx '$ and $\Xx ''$. We have a
flexible functor $F^o:\Yy ^o \rightarrow Top$ given by $F''$ on $\Xx ''$
and $\ast$ on $\Xx '$, and we are looking at the space of flexible functors
$F : \Yy \rightarrow Top$ restricting to $F^o$. For this we apply the
theory developed in Section 1. We treat the case $n\geq 2$ first. Let
$\Lambda $ denote the presheaves on
$Fl (\Yy )$ which was called $\zz \Lambda$ in Section 1, and let $\Lambda
^o$ denote the analogously defined complex of presheaves on $Fl (\Yy ^o )$.
Let $p: Fl (\Yy
) \rightarrow Fl (\Yy ^o)$ denote the natural morphism. Let $G$ denote the
presheaf of groups on $Fl (\Yy )$ which arises in the problem. Due to the
fact that our topological spaces are points over $\Xx '$, the presheaf $G$
has the following property. Let $s:Fl (\Yy )\rightarrow \Yy $ denote the
functor taking an arrow to its source. There is a presheaf of groups $H$ on
$\Yy$ and a morphism $s^{\ast}H \rightarrow G$ such that the kernel and
cokernel are supported on the category $Fl (\Xx _1)$. Furthermore, $H$ is
supported on $\Xx _1$ (i.e. its value is the trivial group on all other
objects).

We have to show that the morphism
$$
Hom ( \Lambda , G) \rightarrow Hom (
\Lambda ^o , p^{\ast}G)
$$
is a quasi-isomorphism. Note that the complexes $\Lambda $ and $\Lambda ^o$
are already complexes of acyclic objects, so these $Hom$s are also $RHom$s.
If $f$ is a functor from a
category $\Aa$ to another category $\Bb$ then there is a functor
$f_{\dagger}$ from presheaves (of abelian groups) on $\Aa$ to presheaves on
$\Bb$ defined by $$
f_{\dagger} M (X) = \frac{\bigoplus _{X\rightarrow f(Y), \, Y\in \Aa}
M(Y)}{\sim}
$$
where $\sim$ is the equivalence relation generated by the condition that
$m_{Y'}$ (corresponding to $X\rightarrow f(Y')$) is equivalent to its
restriction
from $f(Y')$ to $f(Y)$, for $X\rightarrow f(Y)$, whenever there is a
morphism $Y\rightarrow Y'$. This is an adjoint to the pullback: $$
Hom _{\Bb }( f_{\dagger}M, N)= Hom_{\Aa}(M,f^{\ast}N). $$
Now back to our proof.
Let $C$ be the cone on the map $p_{\dagger}\Lambda ^o\rightarrow \Lambda$
and let $D$ be the cone on the map $s^{\ast}H\rightarrow G$. To prove the
desired quasiisomorphism, it suffices to show that $Hom (C,D)$ and $Hom (C,
s^{\ast}H)$ have vanishing cohomology (this implies that $Hom (C,G)$ does,
giving the desired quasiisomorphism).
We treat $Hom (C,D)$ first.
Since $C$ is a complex of
$Hom$-acyclic objects, the functor $Hom (C,\cdot )$ is exact. Thus it
suffices to prove that $Hom (C, K_i)$ has vanishing cohomology for $K_i$
equal to the cohomology presheaves of $D$ (that is, the kernel and cokernel
of $s^{\ast}H\rightarrow G$). These are supported on $Fl (\Xx _1)$.
Let $q:Fl (\Xx _1 )\rightarrow Fl (\Yy )$ be the inclusion. For any pair of
presheaves $J,K$ on $Fl (\Yy )$ such that $K$ is supported on $Fl (\Xx
_1)$, one can check that $$
Hom (J,K) = Hom ( q^{\ast} J, q^{\ast}K). $$
A similar statement holds with $\Yy$ replaced by $\Yy ^o$ (and denoting the
corresponding map by $q_o$).
To get the acyclicity of $Hom (C, K_i)$ it suffices to show that $$
Hom (\Lambda , K_i)\rightarrow Hom (\Lambda ^o , p^{\ast}K_i ) $$
is a quasiisomorphism. For this, we apply the above criterion. Since
$q_o^{\ast}p^{\ast}K_i = q^{\ast}K_i$, it suffices to show that $$
Hom (q^{\ast}\Lambda , q^{\ast}K_i)\rightarrow Hom (q_o^{\ast}\Lambda ^o,
q^{\ast}K_i)
$$
is a quasiisomorphism, and for this it suffices to show that
$q_o^{\ast}\Lambda ^o\rightarrow q^{\ast}\Lambda$ is aquasiisomorphism of
complexes of presheaves on $Fl (\Xx _1)$. In other words, we have to show
that if $u\in Fl (\Xx _1)$ then $\Lambda ^o (u)\rightarrow \Lambda (u)$ is
a quasiisomorphism. But the first is the total complex of the nerve of the
category of objects of $\Yy ^o$ inside $u$, while the second is the total
complex of the nerve of the category of objects of $\Yy$ inside $u$. These
categories are the same, since any factorization of $u$ in $\Yy$ is in fact
a factorization in $\Yy ^o$ (once a morphism goes out of $\Xx _1$ it can't
get back in). This completes the proof that $Hom (C,D)$ has vanishing
cohomology.

Now we treat $Hom (C, s^{\ast}H)$. For this it suffices to show that $$
Hom (s_{\dagger}\Lambda , H ) \rightarrow Hom (s^o_{\dagger}\Lambda ^o , H) $$
is a quasiisomorphism, where $s^o:Fl (\Yy ^o)\rightarrow \Yy $ is the
natural source map. Use the formula
$$
\Lambda = \bigoplus _{(\varphi _1,\ldots , \varphi _k)}(I_{\varphi _1\cdots
\varphi _k})_{\dagger}\zz
$$
where $(\varphi _1,\ldots , \varphi _k)$ is a composable sequence and $$
I_{\varphi _1\cdots \varphi _k}: \ast \rightarrow Fl (\Yy ) $$
is the inclusion of a one point category sending it to the composed
morphism $\varphi _1\cdots \varphi _k$. There is a similar formula for
$\Lambda ^o$. From these and the compatibility of the operation
lower-dagger with composition (following from the adjoint property and this
compatibility for pullbacks), we find that
$(s_{\dagger}\Lambda )(X)$ is the total complex of the nerve of the
category of objects in $\Yy$ below $X$, while $(s^o_{\dagger}\Lambda
^o)(X)$ is the total complex of the nerve of the category of objects of
$\Yy ^o$ which are below $X$ in the category $\Yy$. These two nerves are
both contractible, because $\Yy$ has a final object shared with $\Yy ^o$
(giving a final object of the two categories of objects below $X$). Thus
the morphism between the complexes is a quasiisomorphism. This completes
the proof that $Hom (C, s^{\ast}H)$ has vanishing cohomology, so this
completes the verification of the desired quasiisomorphism for the morphism
$A$ in the case $n\geq 2$.

The case $n=0$ is easy and is left to the reader. For the case of $n=1$,
the first question is from extendability of the section of a relative
category $\Pp \rightarrow \Yy$, from $\Yy ^o$ to $\Yy$. This works because
any $\Pp$ which can arise is a pullback of one over $\Xx ''$, by the
projection $\Yy \rightarrow \Xx ''$. Any section over $\Yy ^o$ will also be
a pullback, so just extend it to the pullback section over $\Yy$. Now
suppose we have a relative group $\Gg \rightarrow \Yy$. Note that $\pi _0
(\Nn )= \Yy ^o$ in this case, since $\Nn$ is generated by the subcategory
$\Yy ^0$. The path components of the space of sections of the corresponding
$\Pp$ are parametrized by pairs $(\Jj , t)$ where $\Jj$ is a $\Gg$-torsor
over $\Yy$ and $t$ is a trivialization over $\Yy ^o$. But the
trivialization of $\Jj $ over $\Xx ''$ already gives an element of $\Jj
(\phi )$ for each arrow $\phi $ going from $\Xx ''$ to $\Xx '$ (note that
$\Gg (\psi )$ has only one element if $\psi$ is an arrow in $\Xx '$, so
$\Jj (\psi )$ also has only one element, and composition with this element
gives an isomorphism between $\Jj (\phi )$ and $\Jj (\psi \phi )$---thus
the $\Jj (\phi )$ depend only on the source of $\phi $ if $\phi$ has target
in $\Xx '$). These arrows are the only ones which need to be added to $\Yy
^o$ to get $\Yy$, so this gives a trivialization of $\Jj$ over $\Yy$. Thus
there is only one path component (corresponding to the trivial torsor). The
fundamental group is the space of sections of $\Gg$ which are trivial on
$\Yy ^o$. By the same consideration as above, it is easy to see that any
such section must be trivial. Thus the fundamental group of our path
component is trivial (and it is $1$-truncated in any case), so the space of
sections is weakly contractible. This completes the case $n=1$ and since we
have treated $n\geq 2$ above, we conclude that $A$ is a weak equivalence.

We turn now to $B$.
In this case, let $\Nn$ denote the subcategory of the full free category
$\Mm$ over $\Yy$, consisting of cells from $\Xx \times I$ or from $\Xx '$.
We claim, first of all, that for any $\phi$, $\Nn (\phi )$ is contractible.
We only need to treat the case where $\phi$ is a map from an object $X\in
\Xx _1$ to the final object $e$. The points in $\Nn$ correspond to pairs
$(\psi , t)$ where $\psi = (\psi _1, \ldots , \psi _k)$ is a composable
sequence with composition $\phi$, and $t$ is a collection of parameters
such that $t_i=1$ for at least one of the objects $X_i$ which is in $\Xx
_0$. To contract this space (keeping at least one such $t$ and keeping the
relations for $t_j=0$), first scale down to $0$ all of the $t$'s
corresponding to objects in $\Xx _1$. Then replace the composable sequence
with an equivalent one (making $k$ smaller), such that $X_{k-1}$ is in $\Xx
_0$. Then change to an equivalent sequence with $X'_k$ added between
$X_{k-1}$ and $X_k$, where $X'_k$ is the same object as $X_k$ (the source
of $\phi$) but considered as being in $\Xx _0$. The last map is the
identity cross the map $1\rightarrow 0$, and the parameter corresponding to
$X'_k$ is $0$. Then scale this parameter up to $1$, then scale the rest of
the parameters down to zero. The result is in a normal form, giving the
contraction.

With this contraction, the cases $n=0$ and $n=1$ are trivial since $\pi _0
(\Nn )=\Yy$. To treat the cases $n\geq 2$ we use the theory of section 1,
particularly the part at the end for the case where $\Nn$ is generated by
several subcategories (in this case, two). From the exact sequence given
there, in order to show that the complex $\Lambda _{\Nn}$ is contractible,
it suffices to show that the complexes corresponding to the subcategories
$\Xx '$, $\Xx \times I$, and $\Xx _0$ are contractible. For each of these
categories $\Yy '$, the presheaf of complexes associates to a morphism
$\phi$ the nerve of the category of objects of $\Yy ' $ which are inside
$\phi$ with respect to $\Yy$. We can assume that $\phi$ goes from an object
$(X,1)$ in $\Xx _1$ to the final object $e$. Then the category of objects
inside $\phi$ is the same as the category of objects of $\Yy '$ which are
under $(X,1)$ with respect to $\Yy$. If $\Yy ' = \Xx \times I$ then this
category has an initial object $(X,1)$; if $\Yy ' = \Xx '$ then it has a
final object $e$. If $\Yy '=\Xx _0$ then it also has an initial object,
namely $(X,0)$
(compare with the proof of contractibility of $\Nn (\phi )$ given above).
In all three cases, the nerve is contractible as required
This completes the proof that $B$ is a weak equivalence. \eop

In view of this lemma, and to make future arguments more convenient, we
define $$
\Gamma (\Xx ,T):= Mor (\ast _{\Xx}, T),
$$
and
$$
\Gamma _{\infty}(\Xx ,T):= Mor _{\infty}(\ast _{\Xx}, T). $$
From the above lemma (and the fact that weak equivalences between CW
complexes are equivalences), the map $\eta \mapsto \rho ^{\ast}\eta$
induces a homotopy equivalence $$
\Gamma '_{\infty}(\Xx ,T)\rightarrow \Gamma _{\infty}(\Xx ,T). $$
This has the advantage of making it easy to define the functoriality with
respect to morphisms of flexible functors. If $T^0$ and $T^1$ are flexible
presheaves then we obtain a weak equivalence $$
a:Mor (\ast _{\Xx } , T^0, T^1)\rightarrow \Gamma (\Xx , T^0)\times Mor
(T^0, T^1),
$$
giving an equivalence
$$
a_{\infty}: Mor_{\infty} (\ast _{\Xx } , T^0, T^1)\rightarrow \Gamma
_{\infty}(\Xx , T^0)\times Mor_{\infty} (T^0, T^1). $$
Taking the ``third side of the triangle'' gives a morphism
$$
b:Mor (\ast _{\Xx } , T^0, T^1)\rightarrow Mor (\ast _{\Xx } , T^0)= \Gamma
(\Xx ,
T^1)
$$
and hence
$$
b_{\infty}: Mor_{\infty} (\ast _{\Xx } , T^0, T^1)\rightarrow \Gamma
_{\infty}(\Xx , T^1).
$$
Let $\Pp$ be the space of flexible functors $\overline{I}\rightarrow Top $
extending the functor $I\rightarrow Top$ which corresponds to the morphism
$a_{\infty}$ considered above. Taking the arrow in the other direction in
$\overline{I}$ (for each point of $\Pp$) and composing with $b_{\infty}$
yields $$
\Pp \times Mor_{\infty} (T^0, T^1)\times \Gamma _{\infty}(\Xx , T^0)
\rightarrow \Gamma
_{\infty}(\Xx , T^1),
$$
hence
$$
\Pp \times Mor_{\infty} (T^0, T^1)\rightarrow Hom (\Gamma
_{\infty}(\Xx , T^0), \Gamma
_{\infty}(\Xx , T^1)).
$$
The fact that $\Pp$ is a weakly contractible space (Theorem \ref{morph02})
means that we obtain a ``very well defined'' map from $Mor_{\infty} (T^0,
T^1)$ to $Hom (\Gamma _{\infty}(\Xx , T^0), \Gamma
_{\infty}(\Xx , T^1))$.

We also note a stricter functoriality with respect to the variable $\Xx$.
Namely, if $F:\Yy \rightarrow \Xx$ is a functor, then we obtain a morphism
$$
F^{\ast}:\Gamma (\Xx , T)\rightarrow \Gamma (\Yy ,F^{\ast}T). $$
Similarly, we obtain
$$
F^{\ast}:\Gamma _{\infty}(\Xx , T)\rightarrow \Gamma _{\infty}(\Yy
,F^{\ast}T). $$
These satisfy $(FG)^{\ast}=G^{\ast}F^{\ast}$.

We say that a flexible functor $T$ is {\em multi-identity normalized} if
$T(1_X, \ldots , 1_X; t)=1_{T_X}$ for all $t$.

Suppose $\Xx$ is a category with a final object $e$. Then there is a
functor $p:\Xx \times I\rightarrow \Xx$ which sends $\Xx \times \{ 0\}$ to
$\Xx$ by the identity, and which sends $\Xx \times \{ 1\} $ to $e$. If $T$
is a (multi-identity normalized) contravariant flexible functor from $\Xx$
to $Top$, then $p^{\ast} T$ is a morphism (of flexible functors) from the
constant functor $\underline{T_e}$ (whose value on any object $X\in \Xx$ is
$T_e$) to $T$. There is a canonical map $T_e\rightarrow Mor (\ast _{\Xx},
\underline{T_e})$ sending $x$ to the morphism of functors which sends
$\ast$ to $x$ (considered as a morphism of flexible functors). The fiber of
$Mor (\ast _X , \underline{T_e}, T)$ over $p^{\ast}T\in Mor
(\underline{T_e}, T)$ maps to $Mor (\ast _{\Xx }, \underline{T_e})$ by a
weak equivalence. There is a canonical weakly contractible space
parametrizing inverses to this weak equivalence (that is, maps from
$Mor_{\infty}(\ast
_{\Xx }, \underline{T_e})$ to the fiber). Thus we obtain a very well
defined map from $\tau _{\infty}T_e$ into this fiber, and by projection we
obtain a map
$$
\tau _{\infty}T_e\rightarrow \Gamma _{\infty}(\Xx , T). $$
In fact, there is a canonical way to get a map $T_e\rightarrow \Gamma (\Xx
, T)$: to a point $x\in T_e$, associate the morphism of flexible functors
obtained by looking at the subspace of $p^{\ast}T$ given by $\{ x\} $ in
the values $T_e$ over objects of $\Xx \times \{ 1\}$, and equal to the
values over objects of $\Xx \times \{ O\}$. This gives an element of
$\Gamma (\Xx ,T)$, which varies continuously with $x$. One can see that
this canonical map yields, under $\tau _{\infty}$, a representative for the
map displayed above.

Note that if $\Yy \subset \Xx$ is a subcategory, we obtain by the second
functoriality above $\Gamma (\Xx ,T)\rightarrow \Gamma (\Yy , T|_{\Yy})$,
so composing we get a canonical morphism $T_e\rightarrow \Gamma (\Yy , T)$.
Similarly we obtain $\tau _{\infty}T_e\rightarrow \Gamma _{\infty}(\Yy ,
T)$

Finally, we note that the functorialities in the two arguments are
compatible, up to homotopy.

Suppose $R$, $S$ and $T$ are flexible functors from $\Zz$ to a continuous
category $\Cc$, and suppose $p:\Yy \rightarrow \Zz$ is a functor. Then the
functor $\Yy \times I^{(2)}\rightarrow \Zz \times I^{(2)}$ induces a
pullback morphism
$$
p^{\ast} : Mor _{\Zz ,\infty}(R,S,T)\rightarrow Mor _{\Yy
,\infty}(p^{\ast}R, p^{\ast}S,
p^{\ast}T).
$$
This is compatible with the maps to
$Mor _{\Zz ,\infty}(R,S)$ {\em etc.\ } and the similar pullbacks for those
spaces of morphisms. Thus $p^{\ast}$ commutes with the (weak ??)
equivalence
$$
Mor _{\Zz ,\infty}(R,S,T)\rightarrow Mor _{\Zz ,\infty}(R,S)\times Mor
_{\Zz ,\infty}(S,T). $$
Inverting this equivalence commutes with $p^{\ast}$ up to homotopy, so we
obtain a diagram
$$
\begin{array}{ccc}
Mor _{\Zz ,\infty}(R,S)\times Mor _{\Zz ,\infty}(S,T) & \rightarrow & Mor
_{\Zz ,\infty}(R,T)\\ \downarrow && \downarrow \\
Mor _{\Yy ,\infty}(p^{\ast}R,p^{\ast}S)\times Mor _{\Yy
,\infty}(p^{\ast}S,p^{\ast}T)
& \rightarrow & Mor _{\Yy ,\infty}(p^{\ast}R,p^{\ast}T) \end{array}
$$
which commutes up to homotopy.

Applying this to the case $\Cc = Top$ and $R=\ast$ we find a diagram $$
\begin{array}{ccc}
\Gamma _{\infty}(\Zz , S)\times Mor _{\Zz ,\infty}(S,T) & \rightarrow &
\Gamma _{\infty}(\Zz , S)\\ \downarrow && \downarrow \\ \Gamma
_{\infty}(\Yy , p^{\ast}S)\times Mor _{\Yy ,\infty}(p^{\ast}S,p^{\ast}T) &
\rightarrow & \Gamma _{\infty}(\Yy , p^{\ast}T) \end{array}
$$
commuting up to homotopy.

\subnumero{The definition of flexible sheaf}
Suppose $\Xx$ is a {\em site}, that is a category with products and fibered
products, and a notion of {\em covering family} $\{ U_{\alpha}\rightarrow X\}$
of any object $X$ (satisfying certain axioms). Recall that a {\em sieve} over
$X$ is a subcategory $\Bb \subset \Xx /X$ of the category of objects of $\Xx$
over $X$, such that if $Y\rightarrow X$ is in  $\Bb$ and if $Z\rightarrow Y$ is
a morphism in $\Xx$, then the composition $Z\rightarrow X$ is in $\Bb$; and such
that $\Bb$ contains a covering family of $X$.

We say that a contravariant flexible functor $T:\Xx \rightarrow Top$ is a
{\em flexible sheaf} (resp. {\em weak flexible sheaf}) if for any object
$X\in \Xx$ and any sieve $\Bb \subset \Xx /X$, the morphism
$$
T_X\rightarrow \Gamma (\Bb , T|_{\Bb})
$$
is a homotopy equivalence (resp. weak homotopy equivalence). Note that $X$
is a final object of $\Xx /X$---the morphism in question here is the one
defined above.

Note that $T$ is a weak flexible sheaf if and only if, for any object $X\in
\Xx$ and any sieve $\Bb \subset \Xx /X$, the morphism $$
\tau _{\infty}T_X\rightarrow \Gamma _{\infty}(\Bb , T|_{\Bb}) $$
is a homotopy equivalence.

It turns out (cf \cite{descente}) that the notion of being a
flexible sheaf corresponds (via the above construction $K$ for example)
with the notion of a simplicial presheaf being {\em flasque with respect to the
objects of the site $\Xx$} in the terminology of Jardine \cite{Jardine1}.
This property is also often known as the ``Mayer-Vietoris property'' cf
Brown-Gersten \cite{BrownGersten}, and also as the property of ``cohomological
descent'' (Thomason \cite{Thomason1} \cite{Wiebel}).

\subnumero{Direct limits}
Suppose $\Jj$ is a category and $F:\Jj \rightarrow Top$ is a contravariant
functor. The {\em direct limit} (see Bousfield-Kan \cite{BousfieldKan},
Vogt \cite{Vogt1} \cite{Vogt2}) $\lim _{\rightarrow }F$ is the space obtained by
glueing together $\Delta ^n\times F(T_0)$ for each composable sequence
$$
T_0\stackrel{\phi _1}{\leftarrow}T_1\ldots T_{n-1}\stackrel{\phi
_n}{\leftarrow}T_n $$
in $\Jj$. The glueing is obtained by identifying the face $(\partial _i
\Delta ^n )\times F(T_0)$ with the cell corresponding to the composable
sequence obtained by replacing the pair $\phi _{i-1} , \phi _i$ by the
composition $\phi _{i-1} \phi _i$. For $i=0$, the face $(\partial _0 \Delta
^n )\times F(T_0)$
is glued to the simplex $\Delta ^{n-1}\times F(T_1)$ via the map $F(\phi _1 )$.

Suppose $\lambda :\Yy \rightarrow \Xx$ is a functor. Recall that $\lambda$
is {\em cartesian} if, for every morphism
$\phi:X'\rightarrow X$ in $\Xx$ and every $Y\in \lambda ^{-1}(X)$, there is
a universal pair $(Y',\psi )$ where $Y'\in \lambda ^{-1}(X')$ and $\psi \in
\lambda ^{-1}(\phi )$. We say that $\lambda$ is {\em supercartesian} if the
pullback $(Y',\psi )$ is unique (and not just unique up to isomorphism, as
is a consequence of the cartesian property). Note that a supercartesian
functor is necessarily split cartesian, that is there is a contravariant
functor $\Xx \rightarrow Cat$ defined by $X\mapsto \lambda ^{-1}(X)$; and
the splitting is unique.

Suppose $\lambda :\Yy \rightarrow \Xx$ is supercartesian, and suppose
$F:\Yy \rightarrow Top$ is a contravariant functor. Then we can define the
{\em relative direct limit}
$$
\lim _{\rightarrow }(\Yy /\Xx , F),
$$
a contravariant functor from $\Xx$ to $Top$, as follows. For $X\in \Xx$, put $$
\lim _{\rightarrow }(\Yy /\Xx , F)(X):= \lim _{\rightarrow} (F|_{\lambda
^{-1}(X)}.
$$
If $\phi :X'\rightarrow X$ is a morphism in $\Xx$, we obtain a functor
$\phi ^{\ast}: \lambda ^{-1}(X)\rightarrow \lambda ^{-1}(X')$ from the
supercartesian property. Furthermore, the fact that $F$ is a functor gives
a natural morphism
$$
F(Y)\rightarrow F(\phi ^{\ast}Y).
$$
Using these morphisms, we obtain a morphism from the set of composable
sequences in $\lambda ^{-1}(X)$ to the set of composable sequences in
$\lambda ^{-1}(X')$ sending $\{ Y_i , \phi _i\}$ to $\phi ^{\ast}Y_i , \phi
^{\ast}\phi _i \}$, and also morphisms of spaces $$
\Delta ^{n}\times
F(Y_0)\rightarrow \Delta ^n\times F(\phi ^{\ast}Y_0). $$
These preserve the
data used in glueing the cells together to form the direct limit, so we
obtain a morphism of direct limits
$$
\lim _{\rightarrow} (F|_{\lambda
^{-1}(X)}\rightarrow \lim _{\rightarrow} (F|_{\lambda ^{-1}(X')}.
$$
This gives the morphism of functoriality $\lim _{\rightarrow }(\Yy /\Xx ,
F)(\phi )$, which satisfies the required associativity.

\subnumero{Construction of $HT$}

Suppose $T$ is a contravariant flexible functor from a site $\Xx$ to $Top$.
If $\Bb$ is a sieve over $X$ then we obtain, for each object
$f:Y\rightarrow X$ in $\Xx /X$, a sieve $\Bb _Y$ in $\Xx /Y$ (the set of
objects $U\rightarrow Y$ such that the composition $U\rightarrow
Y\rightarrow X$ is in $\Bb$). We have a natural functor (composition with
$Y\rightarrow X$) going from $\Bb _Y$ to $\Bb$, which gives the pullback
morphism $$
f^{\ast}:\Gamma (\Bb,T|_{\Bb })\rightarrow \Gamma (\Bb _Y,T|_{\Bb _Y}). $$
If $g:Z\rightarrow Y$ is a morphism, we obtain a functor $\Bb _Z\rightarrow
\Bb _Y$, hence the same pullback for $g$ is
$$
g^{\ast}:\Gamma (\Bb _Y,T|_{\Bb _Y})\rightarrow \Gamma (\Bb _Z,T|_{\Bb _Z}). $$
This system of pullbacks is strictly associative,
$g^{\ast}f^{\ast}=(fg)^{\ast}$.

Let ${\bf Sv}(\Xx )$ be the category whose objects are pairs $(X,\Bb )$
with $X\in \Xx$ and $\Bb \subset \Xx /X$ a sieve. The morphisms from
$(X,\Bb )$ to $(X' ,\Bb ')$ are the morphisms $f:X\rightarrow X'$ such that
$\Bb \subset f^{\ast}\Bb '$. We obtain a functor $$
\Gamma ^{\bf Sv}(T):(X,\Bb )\mapsto \Gamma (\Bb ,T|_{\Bb} ) $$
from ${\bf Sv}(\Xx )$ to $Top$. Put
$$
HT := \lim _{\rightarrow } ({\bf Sv}(\Xx ) /\Xx , \Gamma ^{\bf Sv}(T)). $$
This is a contravariant functor from $\Xx$ to $Top$ which associates to
each $X\in \Xx$ the direct limit
$$
(HT)_X=\lim _{\rightarrow , \Bb } \Gamma (\Bb ,T|_{\Bb }) $$
over sieves $\Bb$ of $\Xx /X$.

We also obtain functors
$$
\Gamma ^{\bf Sv}_{\infty}(T):(X,\Bb )\mapsto \Gamma _{\infty}(\Bb ,T|_{\Bb} ) $$
from ${\bf Sv}(\Xx )$ to $Top$, and we can put $$
H_{\infty}T := \lim _{\rightarrow } ({\bf Sv}(\Xx ) /\Xx , \Gamma ^{\bf
Sv}_{\infty}(T)).
$$
This is a contravariant functor from $\Xx$ to $Top$ which associates to
each $X\in \Xx$ the direct limit
$$
(H_{\infty}T)_X=\lim _{\rightarrow , \Bb } \Gamma _{\infty}(\Bb ,T|_{\Bb }) $$
over sieves $\Bb$ of $\Xx /X$.

{\em Caution:} The functor $\tau _{\infty}$ does not commute with direct
limits, so $H_{\infty}T$ is not equal to $\tau _{\infty}(HT)$. However, the
following lemma implies that there is a very well defined homotopy
equivalence between them.

\begin{lemma}
\mylabel{desc21}
There is a commutative diagram
$$
\begin{array}{ccc}
\tau _{\infty}H_{\infty}T &\rightarrow & \tau _{\infty}(HT)\\
\downarrow&& \downarrow \\
H_{\infty}T& \rightarrow &HT
\end{array}
$$
where the top and left arrows are homotopy equivalences and the bottom and
right arrows are weak equivalences over each object. \end{lemma}
{\em Proof:}
This is a general statement about the way $\tau _{\infty}$ commutes with
homotopy direct limits (we could have stated the same lemma for any direct
limit of topological spaces but we have given here the form we will use).
To see that the bottom arrow is a weak equivalence, use the fact that the
homotopy groups of a direct limit are the direct limit of the homotopy
groups. The facts that the bottom and right arrows are weak equivalences,
and that the left arrow is an equivalence (since $H_{\infty}T$ is already
of CW type), imply that the top arrow is a weak equivalence, hence an
equivalence since the two sides are of CW type.
\eop

Recall that the functor $K'T:\Xx \rightarrow Top$ is defined by $K'T(X):=
\Gamma '(\Xx /X,T|_{\Xx /X})$ (recall the definition of $\Gamma '$ from
above). Put $KT(X):= \Gamma (\Xx
/X,T|_{\Xx /X})$. The morphism $\Gamma '(\Xx /X,T|_{\Xx /X})\rightarrow
\Gamma (\Xx
/X,T|_{\Xx /X})$ is functorial in the first argument, so we obtain a
morphism of functors $K'T\rightarrow KT$. This is a homotopy equivalence,
by Lemma \ref{desc01}. As $\Xx /X$ is a sieve over $X$, $KT(X)$ appears
naturally as a subspace of $HT(X)$. This gives an inclusion of functors
$KT\hookrightarrow HT$.

In the appendix to this section, we define a morphism of flexible functors
$M:T\rightarrow K'T$ (in fact,a homotopy equivalence). By composition we
get a morphism
$$
T\rightarrow HT
$$
(note that there is a canonical choice for the composition of a morphism of
flexible functors and a morphism of functors---cf the remark below ???).

If $T$ is a functor, then the morphism $T\rightarrow KT$ is a morphism of
functors, so we obtain a morphism of functors $T\rightarrow HT$. We also
obtain a morphism of functors $\tau _{\infty}T\rightarrow H_{\infty}T$.

\subnumero{The main proposition}

As explained above, the spaces $\Gamma (\Bb _Y,T|_{\Bb _Y})$ vary
functorially in $Y\in \Xx /X$. We obtain a functor $H_{\Bb}T:\Xx
/X\rightarrow Top$ defined by
$$
H_{\Bb }T_Y:= \Gamma (\Bb _Y,T|_{\Bb _Y}). $$
The morphism
$$
\Gamma (\Xx /Y, T|_{\Xx /Y})\rightarrow \Gamma (\Bb _Y,T|_{\Bb _Y}) $$
is functorial in $Y\rightarrow X$, so it gives a morphism of functors $$
KT|_{\Xx /X}\rightarrow H_{\Bb}T.
$$
We can compose this with the morphism of flexible functors $T\rightarrow
KT$ to get a morphism of flexible functors
$$
T|_{\Xx /X}\rightarrow H_{\Bb} T.
$$

Similarly, we obtain $H_{\Bb ,\infty}T:\Xx /X\rightarrow Top$ defined by
$$
H_{\Bb ,\infty}T_Y:= \Gamma (\Bb _Y,T|_{\Bb _Y}). $$

\begin{lemma}
\mylabel{desc02}
The two maps $\Gamma (\Bb ,T|_{\Bb})\rightarrow \Gamma (\Bb , (H_{\Bb
}T)|_{\Bb})$, one
defined by the composition
$$
\Gamma (\Bb , T|_{\Bb})= (H_{\Bb }T)_X\rightarrow \Gamma (\Bb , (H_{\Bb
}T)|_{\Bb})
$$
and the other induced (as described in ???) by the morphism $F:T|_{\Xx
/X}\rightarrow
H_{\Bb}T$, are homotopic.
\end{lemma}
{\em Proof:}
This statement is conserved by functoriality in terms of equivalences of
flexible functors $T$. Since any $T$ is equivalent to a functor, we may
assume from now on that $T$ is a functor.

We describe explicitly the
points of $\Gamma ((\Bb , (H_{\Bb }T)|_{\Bb})$ (we leave it to the reader
to see that this description follows from the definitions). Such a point
$U$ consists of the following data: for each composable sequence $$ (\phi ,
\eta , \psi )= (\phi _1,\ldots , \phi _a , \eta , \psi _1, \ldots , \psi
_b) $$
of morphisms in $\Bb$ (with $\phi _i : X_i\rightarrow X_{i-1}$, $\psi _j :
Y_j\rightarrow Y_{j-1}$, and $\eta : Y_0\rightarrow X_a$); each $$
(\epsilon , \xi )= (0=\epsilon _0, \ldots , \epsilon _a =1, 0=\xi _0,
\ldots , \xi _b=1)
$$
with $\epsilon _i , \xi _j \in \{ 0 , 1\}$ and these being increasing
sequences; and each
$$
(s , t )= (s_1,\ldots , s_{a-1}, t_1,\ldots , t_{b-1}); $$
we have a point
$$
U(\phi , \eta , \psi ; \epsilon , \xi ;s,t)\in T_{Y_b}. $$
This is subject to the following axioms: the standard ones for $s_i=0$ or
$t_j=0$; if $s_i=1$ with $\epsilon _i=0$ then $$
U(\phi , \eta , \psi ; \epsilon , \xi ;s,t)= U(\phi ', \eta , \psi ;
\epsilon ', \xi ;s',t) $$
where the prime denotes throwing away everything below the object $X_i$; if
$s_i=1$ with $\epsilon _i=1$ then
$$
U(\phi , \eta , \psi ; \epsilon , \xi ;s,t)= U(\phi ', \eta ', \psi ;
\epsilon ', \xi ;s',t) $$
where $\phi ' , \epsilon ' , s'$ have everything thrown away above the
object $X_i$, and $$
\eta ' = \phi _{i+1}\cdots \phi _a \eta ; $$
and similar cases for $\psi , \xi , t$: if $t_i=1$ with $\xi _i=1$ then $$
U(\phi , \eta , \psi ; \epsilon , \xi ;s,t)= (\psi _{i+1}\cdots \psi
_b)^{\ast}U(\phi , \eta , \psi ' ; \epsilon , \xi ';s,t') $$
where the prime denotes throwing away everything above the object $X_i$,
and $(\psi _{i+1}\cdots \psi _b)^{\ast}$ is the restriction for the functor
$T$; and finally, if $s_i=1$ with $\xi _i=0$ then $$
U(\phi , \eta , \psi ; \epsilon , \xi ;s,t)= U(\phi , \eta ', \psi ' ;
\epsilon , \xi ';s,t') $$
where $\psi ' , \xi ' , t'$ have everything thrown away below the object
$X_i$, and
$$
\eta ' = \eta \psi _1 \cdots \psi _j.
$$

Recall that a point $V$ of $\Gamma (\Bb , T|_{\Bb})$ can be described in a
similar way. It is a function $V(\phi ; \epsilon ; t)$ for a composable
sequence $\phi$ in $\Bb$, a sequence of indices $\epsilon _i\in \{ 0,1\}$,
and a sequence of $t_i\in [0,1]$. It satisfies similar axioms (which we
don't repeat here).

The images of $V$ by the two maps in the lemma are the points $U_0$ and
$U_1$ given by
$$
U_0 (\phi , \eta , \psi ; \epsilon , \xi ; s,t)= (\eta \psi _1\cdots \psi
_b)^{\ast}V(\phi ; \epsilon ; s)
$$
and
$$
U_1 (\phi , \eta , \psi ; \epsilon , \xi ; s,t)= V(\psi ; \xi ; t). $$
To prove the lemma, it suffices to construct a continuous function $$
U_r(\phi , \eta , \psi ; \epsilon , \xi ; s,t) $$
for $r\in [0,1]$,
satisfying the properties given above with respect to the inner variables,
and giving $U_0$ and $U_1$ for $r=0$ and $r=1$ respectively.

Put
$$
U_{1/2}(\phi , \eta , \psi ; \epsilon , \xi ; s,t):= $$
$$
V(\phi _1, \ldots , \phi _i, \zeta , \psi _{j+1} , \ldots , \psi _b ; \beta
; s_1,\ldots , s_{i-1}, t_{j+1},\ldots , t_{b-1}) $$
where $i$ is the largest with $\epsilon _i=0$ and $j$ the smallest with
$\xi _j=1$; the map in the middle is
$$
\zeta = \phi _{i+1}\cdots \phi _a \eta \psi _1\cdots \psi _j ; $$
and the index $\beta $ is the one which makes $\zeta$ the map between index
zero and index one.

Let $\tilde{\epsilon}$ denote the index for the composable sequence $(\phi
, \eta , \psi )$ obtained by taking $\epsilon$ for the part $\phi$ and
adding on $1$'s in the remaining places. For $r\in [0,1/2]$ we will put $$
U_r(\phi , \eta , \psi ; \epsilon , \xi ; s , t):= V(( \phi , \eta , \psi);
\tilde{\epsilon}; \sigma (\epsilon , \xi , r, s, t)) $$
for a function
$$
\sigma (\epsilon , \xi ,r, s , t)\in Q_{\tilde{\epsilon}}([0,1]^{a+b} ) $$
satisfying the properties necessary to give a homotopy between $U_0$ and
$U_{1/2}$. Here $Q_{\tilde{\epsilon}}([0,1]^{a+b} )$ is the quotient of
$[0,1]^{a+b} $ obtained by identifying points $\sigma $ and $\sigma '$ if
they formally act the same way in the argument of $V$, that is if they
agree between the last $1$ with $\tilde{\epsilon}_i=0$ and the first one
with $\tilde{\epsilon}_i =1$. The function $V(\alpha
,\tilde{\epsilon},\sigma )$ factors through a function on the quotient
$Q([0,1]^{a+b} )$. We will denote points in the quotient by the same
$a+b$-tuples as for points in $[0,1]^{a+b}$. Note also that if $s,t$ is a
point in the argument of $U$ then for $x,y\in [0,1]$, the concatenation
$(s, x,y,t)$ is a point in $[0,1]^{a,b}$.

Note that we
obtain $U_0$ by putting
$$
\sigma (\epsilon , \xi , 0, s,t)= (s, 1,1, t), $$
whereas we obtain $U_{1/2}$ by putting
$$
\sigma (\epsilon , \xi , 1/2, s,t)= (s', 0,0, t'), $$
where $s'_j= (1-\epsilon _j)s_j$ and $t'_j= \xi _j t_j$. The properties
which we need are:
\newline
(1)\, If $s_i=0$ then $\sigma _i =0$ (up to equivalence in the quotient).
\newline
(2)\, If $t_j=0$ then $\sigma _{a+2+j}=0$ (again, up to equivalence in the
quotient).
\newline
(3)\, If $s_i=1$ for $\epsilon _i=0$ then (as points in the quotient) $$
\sigma = (\ldots , 1 , \sigma (\epsilon ' , \xi , r , s' ,t)) $$
where
$\epsilon '$ and $s'$ are obtained by throwing away everything up to and
including $\epsilon _i$ and $s_i$.
\newline
(4)\, If $s_i= 1$ for $\epsilon _i=1$ then (as points in the quotient) $$
\sigma = (\sigma '(\epsilon ' , \xi , r , s' ,t)) $$
where $\epsilon ' , s'$ are obtained by throwing away everything from
$\epsilon _i , s_i$ on, and $\sigma '(\ldots )$ is obtained from $\sigma
(\ldots )$ by adding a sequence of zeros from the place $\sigma _i$ to the
place $\sigma _a$. \newline
(5)\, If $t_i= 1$ for $\xi _i=0$ then (as points in the quotient) $$
\sigma = (\sigma '(\epsilon , \xi ', r , s ,t')) $$
where $\xi ' , t'$ are obtained by throwing away everything up to and
including $\xi _i , t_i$, and $\sigma '(\ldots )$ is obtained from $\sigma
(\ldots )$ by adding a sequence of zeros from the place $\sigma _{a+1}$ to
the place $\sigma _{a+1+i}$.
\newline
(6)\, If $t_i=1$ for $\xi _i=1$ then (as points in the quotient) $$
\sigma = ( \sigma (\epsilon , \xi ', r , s ,t'), 1 , \ldots ) $$
where
$\xi '$ and $t'$ are obtained by throwing away everything from $\xi _i,
t_i$ on.

These conditions (including the expressions for $r=0$ and $r=1/2$) give
boundary conditions for the function $\sigma (\epsilon , \xi , r , s ,t)$
(they are defined inductively in terms of the choices of $\sigma$ for
smaller $(a+b)$---we make our choices of $\sigma$ by induction). The
boundary conditions are compatible (it is for this compatibility that we
have to take the quotient $Q_{\tilde{\epsilon}}$), giving a continuous
function defined on the boundary of the space of $r,s,t$.
To show that such a function $\sigma$ exists, it suffices to see that the
target space
$Q_{\tilde{\epsilon}}([0,1]^{a+b})$ is contractible.

We can define a contraction
$$
F_q: [0,1]\times
([0,1]^{a+b})\rightarrow ([0,1]^{a+b})
$$
which
preserves any coordinate which is equal to $1$, by setting $(x,\sigma )=
\sigma '$ with $\sigma '_i = 1- (1-x)(1-\sigma _i)$. This passes to the
quotient to give a map (retraction to the point $(1,\ldots , 1)$) $$
F: [0,1]\times Q_{\tilde{\epsilon}}([0,1]^{a+b}) \rightarrow
Q_{\tilde{\epsilon}}([0,1]^{a+b}). $$
Thus $Q_{\tilde{\epsilon}}([0,1]^{a+b})$ is contractible. This allows us to
construct the homotopy between $U_0$ and $U_{1/2}$.

Next, we need a homotopy between $U_{1/2}$ and $U_1$. For this, put
$\tilde{\xi}:= (0,\ldots , 0 , \xi )$, and for $r\in [1/2 , 1]$, set $$
U_r(\phi , \eta , \psi ; \epsilon , \xi ; s , t):= V(( \phi , \eta , \psi);
\tilde{\xi}; \tau (\epsilon , \xi , r, s, t)) $$
for a function
$$
\tau (\epsilon , \xi ,r, s , t)\in Q_{\tilde{\xi}}([0,1]^{a+b} ) $$
satisfying the properties necessary to give a homotopy between $U_0$ and
$U_{1/2}$. Here $Q_{\tilde{\xi}}([0,1]^{a+b} )$ is again the quotient of
$[0,1]^{a+b} $ obtained by identifying points $\sigma $ and $\sigma '$ if
they formally act the same way in the argument of $V$, that is if they
agree between the last $1$ with $\tilde{\xi}_i=0$ and the first one with
$\tilde{\xi}_i =1$.

As before, the properties required of $\tau$ (analogous to those listed
above for $\sigma$) may be interpreted as boundary conditions given by a
continuous function from the boundary of the space of $r,s,t$ into the
space $Q_{\tilde{\xi}}([0,1]^{a+b} )$ (again, this function is defined by
our inductive choice of $\tau$ for smaller values of $a+b$). Since
$Q_{\tilde{\xi}}([0,1]^{a+b} )$ is contractible, there exists a choice of
$\tau$ as required. This completes the proof of the lemma. \eop

\begin{lemma}
\mylabel{desc03}
Suppose $F:S\rightarrow T$ is a morphism of flexible functors over a
category $\Yy$ with initial object $e$. Then the diagram $$
\begin{array}{ccc}
S_e & \rightarrow & \Gamma (\Yy ,S)\\
\downarrow & & \downarrow\\
T_e&\rightarrow &
\Gamma (\Yy , T)
\end{array}
$$
commutes up to homotopy.
\end{lemma}
{\em Proof:}
The statement is preserved by homotopy equivalences of the morphism $F$,
and $F$ is equivalent to a morphism of functors. Then we may take as the
morphism $\Gamma (\Yy , F)$ the morphism induced in the usual way. The
diagram in question becomes commutative.
\eop

\begin{corollary}
\mylabel{desc04}
Suppose $F:S\rightarrow T$ is a morphism of flexible functors from the site
$\Xx$ to $Top$. Suppose $X\in \Xx$ and $\Bb$ is a sieve over $X$. Then the
diagram $$
\begin{array}{ccc}
S_X & \rightarrow & \Gamma (\Bb ,S|_{\Bb})\\ \downarrow & & \downarrow\\
T_X&\rightarrow &
\Gamma (\Bb , T|_{\Bb})
\end{array}
$$
commutes up to homotopy.
\end{corollary}
{\em Proof:}
By the lemma, the diagram
$$
\begin{array}{ccc}
S_X & \rightarrow & \Gamma (\Xx /X ,S|_{\Xx /X})\\ \downarrow & & \downarrow\\
T_X&\rightarrow &
\Gamma (\Xx /X , T|_{\Xx /X})
\end{array}
$$
commutes up to homotopy. On the other hand, the diagram $$
\begin{array}{ccc}
\Gamma (\Xx /X ,S|_{\Xx /X}) & \rightarrow & \Gamma (\Bb ,S|_{\Bb})\\
\downarrow & & \downarrow\\
\Gamma (\Xx /X ,T|_{\Xx /X})&\rightarrow & \Gamma (\Bb , T|_{\Bb})
\end{array}
$$
commutes up to homotopy, by the compatibility discussed at the beginning of
the chapter. The morphism $S_X \rightarrow \Gamma (\Bb , S|_{\Bb})$ is
obtained by composing the top rows of the two diagrams, and the same for
$T$ with the bottom rows. Putting together the two homotopy-commutative
diagrams, we get the corollary.
\eop


\begin{proposition}
\mylabel{desc05}
Suppose $\Xx$ is a site and $T:\Xx \rightarrow Top$ is a contravariant
flexible functor. Then $T$ is a flexible sheaf if and only if the morphism
$T\rightarrow HT$ is a homotopy equivalence.
\end{proposition}
{\em Proof:}
Suppose that $H$ is a flexible sheaf. Then, for each $X\in \Xx$ and each
sieve $\Bb$ over $X$, the morphism
$$
T_X\rightarrow \Gamma (\Bb ,T|_{\Bb })
$$
is a homotopy equivalence. But this implies that if $\Bb '\subset \Bb$ is a
smaller sieve, then
$$
\Gamma (\Bb ,T|_{\Bb })
\rightarrow \Gamma (\Bb ',T|_{\Bb '})
$$
is a homotopy equivalence. In other words, all of the maps in the direct limit
$$
(HT)_X= \lim _{\rightarrow , \Bb } \Gamma (\Bb , T|_{\Bb}) $$
are homotopy equivalences. Note that the direct limit is over a directed
set. In that case, the direct limit is homotopy equivalent to the first
element, which is $(KT)_X$. As $T_X\rightarrow KT_X$ is a homotopy
equivalence, we conclude that $T_X \rightarrow (HT)_X$ is a homotopy
equivalence. Since this is true for all $X$, by Theorem \ref{morph02} the
morphism of flexible functors $T\rightarrow HT$ extends to a homotopy
equivalence (note also that the space of such extensions is weakly
contractible). This provides the first half of the proposition.

Suppose that $T\rightarrow HT$ extends to a homotopy equivalence, so there
is a morphism of homotopy functors $HT\rightarrow T$ providing a homotopy
inverse at each object (that is, the morphism provided by the homotopy
equivalence). We obtain a diagram
$$
\begin{array}{ccc}
\Gamma (\Bb , T|_{\Bb}) & \leftarrow & T_X \\ \uparrow & & \uparrow \\
\Gamma ( \Bb , (HT)|_{\Bb}) & \leftarrow & (HT)_X \\ \uparrow & & \uparrow \\
\Gamma (\Bb , (H_{\Bb }T)|_{\Bb}) & \leftarrow & (H_{\Bb}T)_X \\ \uparrow &
\nearrow & \uparrow \\
\Gamma (\Bb , T|_{\Bb} ) & \leftarrow & T_X \end{array}
$$
where we claim that all squares and triangles are homotopy commutative.

(NB---we need to discuss the notion of homotopy between two morphisms, and
homotopy commutative diagrams, in the section on morphisms...???)

The first and second squares commute by the previous corollary.

The first triangle commutes by Lemma \ref{desc02}.

The second triangle trivial---note that the diagonal arrow is equality.

The compositions along the sides are homotopic to the identity.

For the right side, this follows from the fact that the composition of the
bottom two arrows is homotopic to the natural inclusion $T_X\rightarrow
HT_X$ (second of the above lemmas), and the fact that the other part of the
homotopy equivalence provides a homotopy inverse.

On the left side, we first claim that the composition of morphisms of
flexible functors $T|_{\Xx /X}\rightarrow H_{\Bb}T\rightarrow (HT)|_{\Xx
/X}$ is homotopic to the natural morphism $T|_{\Xx /X}\rightarrow
(HT)|_{\Xx /X}$. To verify this, note that both morphisms factor through
$T|_{\Xx /X}\rightarrow (KT)|_{\Xx /X}$, so it suffices to verify that the
composition of morphisms of functors $$
(KT)|_{\Xx /X} \rightarrow H_{\Bb}T\rightarrow (HT)|_{\Xx /X}
$$
is homotopic to the morphism $(KT)|_{\Xx /X} \rightarrow (HT)|_{\Xx /X}$.
Suppose $Y\in \Xx /X$.
Then $(KT)_Y = \Gamma (\Xx /Y , T_{\Xx /Y})$ whereas $$
(HT)_Y= \lim _{\rightarrow ,\Bb '} \Gamma (\Bb , T_{\Bb '}). $$
Since the sieves $\Xx /Y$ and $\Bb _Y$ appear in the direct limit, there
is, corresponding to the inclusion of sieves $\Bb _Y\subset \Xx /Y$, a
morphism
$$ \Gamma (\Xx /Y , T_{\Xx /Y})\times [0,1]\rightarrow \lim _{\rightarrow
,\Bb '} \Gamma (\Bb , T_{\Bb '}) $$
sending $\Gamma (\Xx /Y , T_{\Xx /Y})\times \{ 0\}$ in by the natural
inclusion, and at $\{ 1\}$ giving the pullback map $\Gamma (\Xx /Y , T_{\Xx
/Y})\rightarrow \Gamma (\Bb _Y , T_{\Bb _Y})$. Note that if $Z\rightarrow
Y$ is a morphism, then the pullback of $\Bb _Y$ to
$Z$ is $\Bb _Z$. Thus this collection of morphisms gives a morphism of functors
$$
(KT)|_{\Xx /X} \times \underline{[0,1]} \rightarrow (HT)|_{\Xx /X}, $$
restricting to the
natural morphism on $(KT)|_{\Xx /X} \times \{ 0\}$. It restricts to the
composition with $ H_{\Bb}T\rightarrow
(HT)|_{\Xx /X}$ on
$(KT)|_{\Xx /X} \times \{ 1\}$, because $H_{\Bb}T$ associates to each $Y\in
\Xx /X$ the subspace $\Gamma (\Bb _Y, T_{\Bb _Y})$ of the direct limit
$(HT)_Y$. This morphism of functors gives a homotopy between the natural
map $(KT)|_{\Xx /X} \rightarrow (HT)|_{\Xx /X}$ and the composition $$
(KT)|_{\Xx /X} \rightarrow H_{\Bb}T\rightarrow (HT)|_{\Xx /X},
$$
giving the claimed statement.

It follows that the composition of the morphisms of flexible functors (over
$\Xx /X$) occuring on the left side is homotopic to the identity (by the
usual properties of compositions). We obtain the same thing after
restricting to $\Bb$. To conclude, note that if the composition of a
sequence of morphisms of homotopy functors is homotopic to a morphism $f$,
then the composition of the induced morphisms on $\Gamma (\Bb , \cdot )$ is
homotopic to the morphism induced by $f$. (This follows from the definition
of $\Gamma$ as a space of morphisms, and the properties of composition of
morphisms.)

Finally, it follows from the form of this diagram and the homotopy
commutativity that the morphism $\Gamma (\Bb , T|_{\Bb})\rightarrow T_X$
obtained by composing the diagonal arrow and the two upper right arrows, is
a homotopy inverse to the map $T_X\rightarrow \Gamma (\Bb , T|_{\Bb})$.
This is true for every $X$ and every sieve, so $T$ is a flexible sheaf.
\eop

The following lemma is needed for the next proposition.
\begin{lemma}
\mylabel{desc05a}
Suppose $A$ is a contravariant functor on a category $\Yy$. Then the morphism $$
\Gamma (\Yy , \tau _{\infty} A)\rightarrow \Gamma (\Yy , A ) $$
is a weak equivalence.
\end{lemma}
{\em Proof:}
More generally, if $B \rightarrow A$ is a morphism of functors inducing
weak equivalences over each object of $\Yy$, then the induced morphism on
the spaces of global sections is a weak equivalence. This is because the
space of global sections is defined in terms of maps from CW complexes into
the spaces. \eop

\begin{proposition}
\mylabel{desc05b}
Suppose $T:\Xx \rightarrow Top$ is a contravariant functor of CW type. Then
$T$ is a weak flexible sheaf if and only if the morphism $\tau _{\infty} T
\rightarrow H_{\infty}T$ is a homotopy equivalence. \end{proposition}
{\em Proof:}
The first half of the proof goes through as before. For the second half,
suppose $\tau _{\infty} T \rightarrow H_{\infty}T$ is a homotopy
equivalence. Look at the diagram $$
\begin{array}{ccccc}
\Gamma _{\infty}(\Bb , \tau _{\infty}T|_{\Bb}) & \rightarrow & \Gamma (\Bb
, \tau _{\infty}T|_{\Bb}) & \leftarrow & \tau _{\infty}T_X \\ \uparrow &&
\uparrow & & \uparrow \\
\Gamma _{\infty}( \Bb , (H _{\infty}T)|_{\Bb}) & \rightarrow & \Gamma ( \Bb
, (H _{\infty}T)|_{\Bb}) & \leftarrow & (H_{\infty}T)_X \\ \uparrow &&
\uparrow & & \uparrow \\
\Gamma _{\infty}(\Bb , (H_{\Bb ,\infty}T)|_{\Bb}) & \rightarrow & \Gamma
(\Bb , (H_{\Bb ,\infty}T)|_{\Bb}) & \leftarrow & (H_{\Bb,\infty }T)_X \\
\uparrow && \uparrow & & \uparrow \\
\Gamma _{\infty}(\Bb , \tau _{\infty}T|_{\Bb} ) & \rightarrow & \Gamma (\Bb
, \tau _{\infty}T|_{\Bb} ) & \leftarrow & \tau _{\infty}T_X \end{array}
$$
similar to the one in the previous proof. Here the diagonal arrow on the
bottom is
$$
\Gamma _{\infty}(\Bb , \tau _{\infty}T|_{\Bb} ) \rightarrow (H_{\Bb,\infty
}T)_X .
$$
Since the horizontal arrows on the left are weak equivalences, we can lift
the leftward-pointing horizontal arrows on the right side to obtain a
diagram $$
\begin{array}{ccc}
\Gamma _{\infty}(\Bb , \tau _{\infty}T|_{\Bb}) & \leftarrow & \tau _{\infty}T_X
\\ \uparrow & & \uparrow \\
\Gamma _{\infty}( \Bb , (H _{\infty}T)|_{\Bb}) & \leftarrow & (H_{\infty}T)_X \\
\uparrow & & \uparrow \\
\Gamma _{\infty}(\Bb , (H_{\Bb ,\infty}T)|_{\Bb}) & \leftarrow &
(H_{\Bb,\infty }T)_X \\
\uparrow & \nearrow & \uparrow \\
\Gamma _{\infty}(\Bb , \tau _{\infty}T|_{\Bb} ) & \leftarrow & \tau
_{\infty}T_X .
\end{array}
$$
One can check that it is commutative by using the commutativity of the
diagram in the proof of the previous proposition. In order to do this, one
needs to know that the morphism
$$
\Gamma (\Bb , H_{\Bb , \infty} T|_{\Bb })\rightarrow \Gamma (\Bb ,
H_{\Bb}T|_{\Bb})
$$
is a weak equivalence. This follows from the previous lemma. After we know
that the diagram is commutative, the same argument as for the previous
proposition, implies that the morphism
$$
\tau _{\infty}T\rightarrow \Gamma _{\infty}(\Bb , \tau
_{\infty}T|_{\Bb})\;\;\;\;\;\;\;\;\;\;\; (\dagger ) $$
is a
homotopy equivalence. Since $T$ is of CW type, the morphism (of functors)
$\tau _{\infty}T\rightarrow T$ is a homotopy equivalence. We get an
exactly-commutative diagram $$
\begin{array}{ccc}
T_X&\rightarrow & \Gamma (\Bb , T|_{\Bb})\\ \uparrow & & \uparrow \\
\tau _{\infty}T_X & \rightarrow & \Gamma (\Bb , \tau _{\infty}T|_{\Bb})
\end{array}
$$
where the vertical arrows are homotopy equivalences. Note that our morphism
$(\dagger )$ is just a lifting to a diagonal map from the upper left to the
lower right, in the diagram obtained by applying $\tau _{\infty}$ to this
diagram. The fact that $(\dagger )$ is an equivalence implies then that
$$
\tau _{\infty} T_X \rightarrow \Gamma _{\infty}(\Bb , T|_{\Bb }) $$
is an equivalence. This is the necessary statement. \eop

\subnumero{Covering families}
Suppose $\Xx $ is a site, and $\Uu = \{ U_{\alpha}\}$ is a covering family
of an object $X$. If $T$ is a contravariant flexible functor from $\Xx$ to
$Top$ we define a space $\Gamma (\Uu , T)$ homotopy equivalent to $\Gamma
(\Bb ^{\Uu}, T)$ where $\Bb ^{\Uu}$ is the sieve defined by $\Uu$
(consisting of all morphisms $V\rightarrow X$ which factor through some
$U_{\alpha}$ in $\Uu$).

We denote fiber products of elements in $\Uu$ by
$$
U_{\alpha _0\cdots \alpha _n}:= U_{\alpha _0}\times _X \ldots \times _X
U_{\alpha _n}.
$$

The following definitions and considerations are due to Boardman and Vogt
\cite{BoardmanVogt} see also the Stasheff polytopes \cite{Stasheff}, and
\cite{StasheffRecent}.

A {\em tree diagram of width $n$} is a tree embedded in the plane, with
$n+2$ strands hanging down from a certain point (the ``top''). The strands
are infinitely long. The edges between the vertices or between the vertices
and the top are assigned lengths in $(0,\infty ]$, and when the length
approaches $0$ the two vertices should come together. Choose the embedding
so that the length of an edge is equal to the vertical distance between the
vertices. Constrain the lengths so that the vertical distance between two
vertices is $\leq 1$ (even if the vertices are not on the same edge). Call
the bands between horizontal lines which contain vertices, the {\em
intervals} of $\tau$; the lengths of edges should really be thought of as
lengths of these intervals.

If $\alpha = \alpha _0,\ldots , \alpha _n$ is an index and $\tau$ is a tree
diagram of width $n$, then we assign $\alpha _i$ to the region in the plane
between the $i$th and $i+1$st strand. For each interval $\epsilon$, we
obtain an index $\alpha (\tau , \epsilon )$ corresponding to a horizontal
slice in the band $\epsilon$ (take all the $\alpha _j$ which appear in
regions which intersect the slice).

Cutting at an interval $\epsilon$, we obtain a tree diagram $\tau '$ which
is the part above $\epsilon$ (the strands in $\epsilon$ are made into the
infinite ends of the new diagram), and a {\em disconnected tree diagram}
$\tau ''$ corresponding to the part below $\epsilon$ (with the strands in
$\epsilon$ becoming the various top points). If $\epsilon _1,\ldots ,
\epsilon _k$ are the bands below $\epsilon := \epsilon _0$ (with $\epsilon
_k$ being the bottom) then we obtain a sequence of fiber products with
projection maps $$
U_{\alpha (\tau '', \epsilon _{i-1})}\stackrel{\phi _i}{\leftarrow}
U_{\alpha (\tau '', \epsilon _{i})},
$$
and the lengths of $\epsilon _i$ provide a sequence $t_1,\ldots , t_{k-1}$.
Note that the first space is $U_{\alpha '}$ where $\alpha '$ is the bottom
index of $\tau '$, and the last is $U_{\alpha}$. Our flexible functor thus
gives a map
$$
T(\alpha '' , \tau ''):= T(\phi _1,\ldots , \phi _k; t_1,\ldots , t_{k-1}):
T_{U_{\alpha '}}\rightarrow T_{U_{\alpha}} . $$

The space $\Gamma (\Uu , T)$ is defined to be the space of all collections
of functions $$
Q(\alpha _0,\ldots , \alpha _n; \tau )\in T_{U_{\alpha _0\cdots \alpha _n}}
$$
where $\tau$ is a tree diagram of width $n$. $Q$ is assumed to be
continuous in the variable $\tau$ and to satisfy the following axiom: if an
interval $\epsilon$ in the tree diagram $\tau$ has width
$1$, then let $\tau '$ be the tree diagram truncated at $\epsilon$ and let
$\tau ''$ be the part below $\epsilon$. Let $\alpha '$ be the index
corresponding to the bottom of $\tau '$ and let $\alpha ''$ be the
indexation of $\tau ''$. Then we
require that that
$$
T(\alpha '', \tau '')(Q(\alpha ', \tau ')) = Q(\alpha , \tau ). $$

If $\Uu = \{ U_{\alpha}\} _{\alpha \in J}$ is a covering family in $\Xx
/X$, let $\Jj$ denote the category of all formal products of elements of
$J$. In other words, the objects of $\Jj$ are of the form
$B_{\alpha_1,\ldots , \alpha _k}$ for $\alpha_1,\ldots , \alpha _k\in J$,
with morphisms corresponding to projections and diagonal maps. The covering
$\Uu$ gives a functor $$
\xi _{\Uu} : \Jj \rightarrow \Xx /X
$$
defined by
$$
\xi _{\Uu}( B_{\alpha_1,\ldots , \alpha _k}):= U_{\alpha _1,\ldots , \alpha
_j}. $$
(Note: the category $\Jj$ is defined to be the universal one such that we
get $\xi _{\Uu}$ for any covering indexed by $J$ in any site.) The functor
$\xi _{\Uu}$ actually has image in the sieve $\Bb ^{\Uu}$ determined by
$\Uu$.

\begin{lemma}
\mylabel{desc09a}
We have
$$
\Gamma (\Uu , T)=\Gamma '(\Jj , \xi _{\Uu}^{\ast}T). $$
\end{lemma}
{\em Proof:}
This follows immediately from the definition of the space of sections
$\Gamma '(\Jj , \xi _{\Uu}^{\ast}T)$.
\eop

The boundary of a tree diagram consists exactly of those trees with an edge
of length one.

\begin{lemma}
\mylabel{desc06a}
{\rm (Boardman and Vogt)}
The space of tree diagrams of width $n$ is topologically an $n$-disk, with
boundary topologically an $n-1$-sphere.
\end{lemma}
{\em Proof:} \cite{BoardmanVogt}.
One can see how this works by checking the cases of small values of $n$
(there is a subdivision of the $n$-disk whose cells correspond to the cells
in the space of tree diagrams of width $n$).
\eop

\begin{lemma}
\mylabel{desc06}
If $\Uu$ is a finite covering family and $T$ is of CW-type, then $\Gamma
(\Uu , T)$ is of CW-type.
\end{lemma}
{\em Proof:}
If $T$ is $m$-truncated, then the space of choices of $Q$ on the cells of
dimension $n>m$ is contractible (that is we get a sequence of fibrations
with contractible fiber; the inverse limit has contractible fiber). Hence
$\Gamma (\Uu , T)$ is the space of ways of specifying $Q$ for tree diagrams
of width $\leq m$. Since $\Uu$ is finite, there are a finite number of
cells. The space of ways of mapping each cell, given the boundary, is of CW
type; as there are a finite number of cells, the space is of CW type. \eop

\begin{lemma}
\mylabel{desc09}
The pullback morphism
$$
\xi _{\Uu}^{\ast} : \Gamma '(\Bb ^{\Uu}, T|_{\Bb ^{\Uu}})\rightarrow \Gamma
'(\Jj , \xi _{\Uu}^{\ast}T) = H(\Uu , T)
$$
is a homotopy equivalence. The composition $T_X\rightarrow \Gamma '(\Bb
^{\Uu },T)\rightarrow H(\Uu , T)$ is equal to the usual map.
\end{lemma}
{\em Proof:}
The composition is equal to the usual map $T_X\rightarrow \Gamma '(\Jj ,
\xi _{\Uu}^{\ast}T)$, and this corresponds to the usual map into $H(\Uu ,
T)$ via the equality of Lemma \ref{desc09a}. We need to provide the homotopy
inverse $H(\Uu , T)\rightarrow \Gamma (\Bb ^{\Uu}, T)$.

We start with some observations similar to those in Section 1. For $Y\in
\Bb ^{\Uu }$, let
$N_{\Bb ^{\Uu }}(Y)$ be the space of pairs $(\phi , t)$ where $\phi = \phi
_1,\ldots , \phi _k)$ is a composable sequence of maps $\phi
_{j}:X_{j-1}\leftarrow X_{j}$ with $X_k=Y$, and $t= (t_0,\ldots ,
t_{k-1})\in [0,1]^k$, subject to the standard identifications for $t_i=0$.
This collection of spaces is endowed with maps $$
N_{\Bb ^{\Uu }}(Y)\times M_{\Bb
^{\Uu }}(\phi) \rightarrow N_{\Bb ^{\Uu }}(Z) $$
whenever $\phi : Z\rightarrow Y$, compatible with the composition in
$M_{\Bb ^{\Uu }}$ (this action of $M$ on $N$ is obtained, as was the case
for the composition in $M$, by inserting a $t_j=1$ in the appropriate
place). The space of sections $\Gamma ' (\Bb ^{\Uu }, T)$ may be identified
with the space of ways of mapping each $N_{\Bb ^{\Uu }}(Y)$ into $T_Y$
compatibly with the action of $M_{\Bb ^{\Uu }}$ (via the functor $T$).

Now for any $Y\in \Bb ^{\Uu }$,
let $N'(Y)\subset N_{\Bb ^{\Uu }}(Y)$ denote the subset consisting of pairs
$(\phi ,t)$ where there is a $j_0$ such that $t_{j_0}=1$ and for $j\leq
j_0$ the space $X_j$ is one of the $U_{\alpha _1\ldots \alpha _m}$ and the
morphism $\phi _{j}:X_{j-1}\leftarrow X_{j}$ is a standard projection. We
have a morphism $H(\Uu , T)\rightarrow \Gamma (\Bb ^{\Uu},N'; T)$ (in fact,
this is an isomorphism). We claim that $N'(Y)$ is contractible. It will
follow from Proposition \ref{flex01} that there is a natural contractible
space parametrizing morphisms $\Gamma (\Bb ^{\Uu},N'; T)\rightarrow \Gamma
(\Bb ^{\Uu}, T)$, so we obtain a natural contractible space of morphisms
from $H(\Uu , T)$ to $\Gamma (\Bb ^{\Uu}, T)$.

The map $\Gamma (\Bb ^{\Uu},T)\rightarrow H(\Uu , T)=\Gamma (\Bb ^{\Uu},
N',T)$ is induced by the inclusions $N'(Y)\subset N(Y)$. Our system of maps
$N(Y)\rightarrow N'(Y)$, composed with this system of inclusions, gives a
system of maps $N(Y)\rightarrow N(Y)$ compatible with the action; this
system of maps is homotopic to the
identity (by applying the same type of argument as in Proposition
\ref{flex01}, using the fact that the $N(Y)$ are contractible), so the
composition $$
\Gamma (\Bb ^{\Uu},T)\rightarrow \Gamma (\Bb ^{\Uu}, N',T)\rightarrow
\Gamma (\Bb ^{\Uu},T)
$$
is homotopic to the identity. On the other hand, the spaces $N'(Y)$ are
unions of standard cells of $N(Y)$, and this collection of cells is closed
under the operation of taking the boundary of a standard cell. Hence, in
the inductive construction of Proposition \ref{flex01}, we may use the
identity in $N'(Y)$ when it is available as a solution to the extension
problem for a given cell; thus we may assume that the system of morphisms
$N(Y)\rightarrow N'(Y)$ is the identity on $N'(Y)$. Thus the composition $$
\Gamma (\Bb ^{\Uu},
N',T)\rightarrow \Gamma (\Bb ^{\Uu},T)\rightarrow \Gamma (\Bb ^{\Uu},
N',T)
$$
is the identity. This proves that the map we construct is a homotopy
inverse to the standard one.

We construct the contration of $N'(Y)$. Fix a morphism $f:Y\rightarrow
U_{\beta}$ (this exists because $Y\in \Bb ^{\Uu}$). Suppose $(\phi , t)\in
N'(Y)$ with $\phi = (\phi _1,\ldots , \phi _n)$. Define a sequence of
points $(\psi ^s, t^s)$ in $N'(Y)$ inductively on $s\in [0, 1]$ as follows.
Set $m=t_1+\ldots +t_{n-1}$; this depends continuously on $(\phi ,t)$, and
$m\geq 1$ since by hypothesis there is at least one $t_i=1$. Let $X_i^s$
denote the sequence of spaces involved. For $s=s_{i+1}:=(t_{i+1}+\ldots +
t_{n-1})/m$, we will have $X^s_j= X_{j-1}\times U_{\beta}$ for $j\geq i+1$
and $X^s_j=X_{j}$ for $j\leq i-1$. We will have the corresponding values of
$t_j$; but also $t_i=0$ and $X^s_i$ is either equal to $X_i$ or to
$X_i\times U_{\beta}$ (it doesn't matter which, because of the equivalence
relation and the fact $t_i=0$). The morphisms are given by the $\phi$ in
the first factor, and the identity in the second factor; except at the end
where $X^s_n=Y$ for all $s$ and the morphism to $X^s_{n-1}$ is given by the
$\phi _n$ in the first factor, and the morphism $f$ in the second factor if
necessary. As $s$ passes from $s_i$ to $s_i +t_i/2m$, the value of $t^s_i$
increases linearly from $0$ to $t_i$ (with derivative $2m$); the rest of
the $t^s_j$ stay the same. As $s$ passes from $s_i+t_i/2m$ to
$s_{i+1}=s_i+t_i/m$, the value of $t^s_{i+1}$ decreases linearly from $t_i$
to $0$. This defines the homotopy on the interval $[s_{i+1},s_{i}]$;
continue in the same way for the interval $[s_{i},s_{i-1}]$ and so on (note
that $s_n=0$ and the $s_i$ are decreasing since we go backwards; this is
because we have numbered the composable sequence in the direction opposite
to that of the arrows). This process depends continuously on $(\phi , t)$
modulo the equivalence relations, since if a $t_i=0$ then $s_i=s_{i+1}$ and
nothing is changed in the interval. Grouping these together, we obtain a
homotopy indexed by $s\in [0,1]$.
Note
furthermore that we have preserved the conditions for inclusion in $N'(Y)$,
because our changes consist of adding on a factor of $U_{\beta}$, preserving
the condition on the form of the first part of the series; and we have
taken care not to destroy the condition $t_{j_0}=1$ (the $t_{j_0}$ can be
changed only after replacing it by another one which would take the value
one in this case).

At the end, we have $X_j^1= X_j\times U_{\beta}$ for $j=0,\ldots , n-1$.
This is equivalent to the sequence of length $n+1$ where we add on
$X_{-1}=U_{\beta}$ and put $\phi _{0}=pr_2$ and $t_{-1}=0$. Then, for $s\in
[1,2]$, increase $t_{-1}$ linearly from $0$ to $1$. Finally, for $s\in
[2,3]$, multiply the $t_j$ by $3-s$ for $j\geq 0$; the point rests in
$N'(Y)$ because of $t_{-1}=1$. At $s=3$ the point is equivalent to the
point $(f,1)$ (with $X_{-1}=U_{\beta}$, $X_0=Y$, and $t_{-1}=1$). This
completes our contraction of $N'(Y)$ to the point $(f,1)$.

For the last statement of the lemma, note that the map $T_X\rightarrow
H(\Uu , T)$ is induced by the fact that the $U_{\alpha _1\ldots \alpha _m}$
are in $\Bb ^{\Uu}$. Thus the composition in question is equal to this map.
\eop

\begin{corollary}
\mylabel{desc12}
A flexible functor $T$ is a flexible sheaf (resp. weak flexible sheaf) if
and only if, for all coverings $\Uu$ of objects $X\in \Xx$, the morphism
$T_X\rightarrow H(\Uu , \Xx)$ is a homotopy equivalence (resp. weak
equivalence). \end{corollary}
{\em Proof:}
This follows from the previous lemma and the fact that the $\Bb ^{\Uu}$ are
cofinal in the set of all sieves $\Bb$ (note that if the condition of the
corollary is satisfied for all coverings, the map from $T_X$ to the
telescope $(HT)_X$ over
all sieves is a (weak) homotopy equivalence because the map to the
telescope over a cofinal subset is a (weak) homotopy equivalence). \eop

We say that the site $\Xx$ is {\em quasi-compact} if any sieve $\Bb$
contains a finite covering family $\{ U_{\alpha}\}_{\alpha = 1,\ldots ,
k}$. Recall that a homotopy functor $T$ {\em of CW type} if each $T_X$ is
homotopy equivalent to a CW-complex. We say that a homotopy functor $T$ of
CW type is {\em $n$-truncated} if for every point $x\in T_X$ we have $\pi
_i(T_X,x)=0$ for $i\geq n+1$.

\begin{corollary}
\mylabel{desc08}
Suppose $\Xx$ is a quasi-compact site and $T$ is an $n$-truncated flexible
functor of CW-type. Then $HT$ is an $n$-truncated flexible functor of
CW-type. In particular, $HT$ and $H_{\infty}T$ are homotopy equivalent.
\end{corollary}
{\em Proof:}
We have
$$
(HT)_X= \lim _{\rightarrow ,\Bb} \Gamma (\Bb , T|_{\Bb}), $$
but this limit is homotopy equivalent to the limit over the cofinal
subsystem of sieves of the form $\Bb ^{\Uu}$, where the $\Uu$ are finite
covering families (this is cofinal since $\Xx$ is quasi-compact). By the
previous lemma, $\Gamma (\Bb ^{\Uu},T)\sim \Gamma (\Uu , T)$ and by Lemma
\ref{desc06}, this is of CW-type. A direct limit of spaces of CW-type is
again of
CW-type, so $(HT)_X$ is of CW-type. Note that $HT$ is also $n$-truncated.
One can see that each $\Gamma (\Bb , T)$ is $n$-truncated, by looking at
the fibrations for the sequence of primitive cells---these are fibrations
where the fibers are all mapping spaces into $n$-truncated spaces, hence
the fibers are $n$-truncated; and the limit of the sequence of fibrations
is $n$-truncated. Finally, the fact that $HT$ is of CW type implies that it
is equivalent to $\tau _{\infty}HT$, and by Lemma \ref{desc21} this is
equivalent
to $H_{\infty}T$.
\eop

{\em Remark:} If $\Xx$ is not quasi-compact, then $HT$ will not, in
general, be of CW type. For example, suppose $\Xx$ is an infinite discrete
topological space. If $\Sigma$ denotes the constant presheaf of sets whose
value on each element is a set with two elements, then $\Gamma (\Xx ,
\Sigma )$ is the space of dyadic numbers with its Tychonoff topology, which
is not a CW complex. If $T$ is a flexible functor of CW type and $HT$ is
not of CW type, then $T$ cannot be a flexible sheaf; it could still be a
weak flexible sheaf, though.

\begin{corollary}
\mylabel{desc22}
If $\Xx$ is a quasi-compact site and $T$ is a truncated flexible functor,
then $T$ is a flexible sheaf if and only if it is a weak flexible sheaf.
\end{corollary}
{\em Proof:}
Combine
Corollary \ref{desc08} and the characterizations of Propositions \ref{desc05}
and  \ref{desc05b}.
\eop

\subnumero{The associated flexible sheaf on a quasi-compact site}

We suppose, for the rest of this section, that $\Xx$ is quasi-compact. We
will construct the flexible sheaf associated to a truncated flexible
functor $T$ of CW type. In the next section we will give a different
construction of the weak flexible sheaf associated to a truncated flexible
functor on any site (not necessarily quasi-compact). Thus, for
quasi-compact sites we will have given two different constructions.

It should be noted that in the ``closed model category'' approach
(Jardine \cite{Jardine1}), one obtains the associated flexible sheaf quite
easily
by taking the fibrant replacement (which is constructed using the small object
argument).

{\em Remark:} Suppose $T$ is a functor. Then the morphism $F:T\rightarrow
HT$ is a morphism of functors, since it is the composition $T\rightarrow
KT\hookrightarrow HT$.

Suppose $f:S\rightarrow T$ is a morphism of functors. Then $f$ induces a
morphism $Hf:HS\rightarrow HT$ of functors. This is given by an induced
morphism
$$
\Gamma f:\Gamma (\Bb ,S|_{\Bb})\rightarrow \Gamma (\Bb ,T|_{\Bb}) $$
for each sieve $\Bb$ over $Y\in \Xx$. This induced morphism is defined by $$
(\Gamma f (P))(\phi _1,\ldots , \phi _n, t_0,\ldots , t_{n-1}):=
f_{X_n}(P(\phi _1,\ldots , \phi _n, t_0,\ldots , t_{n-1})) $$
for $P\in \Gamma (\Bb ,S)$. Note that if $t_i=1$ then we have $$
(\Gamma f (P))(\phi , t) :=f_{X_n}(P(\phi ,t)) $$
$$
= f_{X_n}(S(\phi '')(P(\phi ', t')))
$$
$$
= T( \phi '' )f_{X_n}(P(\phi ', t'))
$$
$$
=T(\phi '', t'')(\Gamma f (P))(\phi ' , t'), $$
where $\phi ' , t'$ refer to the first parts of the sequences, and $\phi
'', t''$ to the second parts. This shows that $(\Gamma f (P))$ satisfies
the required conditions (the conditions for $t_i=0$ being obvious). The
morphisms $\Gamma f$ are compatible with restriction to smaller sieves, so
we get a morphism of direct limits $(Hf)_Y(HS)_Y\rightarrow (HT)_Y$.

Suppose $Z\rightarrow Y$ is a morphism in $\Xx$. If $\Bb$ is a sieve over
$Y$ then we obtain a sieve $\Bb _Z$ over $Z$ consisting of all the
morphisms to $Z$ which, when composed with the morphism to $Y$, are in
$\Bb$. There is an obvious functor $\Bb _Z\rightarrow \Bb$ and the pullback
of $S|_{\Bb}$ or $T|_{\Bb}$ by this functor is $S|_{\Bb_Z}$ or
$T|_{\Bb_Z}$. We obtain the pullback morphisms
$$
\Gamma (\Bb , S|_{\Bb})\rightarrow \Gamma (\Bb _Z, |_{\Bb_Z}) $$
and similarly for $T$; and the morphisms $\Gamma f$ are compatible with
these pullbacks. As the restriction from $Y$ to $Z$ for the functor $HT$ is
given by these pullback morphisms, the morphisms $(Hf)_Y$ fit together to
give a morphism of functors $Hf$.

\begin{lemma}
\mylabel{desc14}
Suppose $T$ is a functor. Then we obtain a natural morphism of functors
$i_T:T\rightarrow HT$, and similarly $i_{HT}:HT\rightarrow H(HT)$. On the
other hand, by the above remarks we have an induced morphism
$H(i_T):HT\rightarrow H(HT)$. For any $X\in \Xx$, the two morphisms
$i_{HT,X}$ and $H(i_T)_X$ from $(HT)_X$ to $(H(HT))_X$ are homotopic.
\end{lemma}
{\em Proof:}
This is more or less the same statement as Lemma \ref{desc02}. \eop

Suppose $G_i$ is a sequence of functors, with morphisms $g_i:G_i\rightarrow
G_{i+1}$. Suppose that the $G_i$ are all $n$-truncated and of CW-type. Then
we obtain a sequence of functors $HG_i$ with morphisms $Hg_i$.

We claim that the natural morphism $Tel (HG_i)\rightarrow H(Tel (G_i))$ is
a homotopy equivalence. Note first of all that both sides are $n$-truncated
of CW-type. Thus it suffices to prove that it is a weak equivalence at each
$X\in \Xx$. Note then that both sides may be considered as telescopes over
all sieves $\Bb $ over $X$, homotopy equivalent to telescopes over sieves
$\Bb ^{\Uu }$. It suffices to prove that for any covering $\Uu$ of $X$, the
morphism
$$
Tel (H(\Uu , G_i ))\rightarrow H(\Uu , Tel (G_i)) $$
induces isomorphisms of homotopy groups. But since the $G_i$ are
$n$-truncated, we may replace the $H(\Uu , G_i)$ by $H^{\leq n}(\Uu ,G_i)$
(which has the same definition as $H(\Uu , G_i)$ but looking only at tree
diagrams of width $\leq n$), and the same for $Tel (G_i)$.

Thus it suffices to prove that the morphism $$
Tel (H^{\leq n}(\Uu , G_i ))\rightarrow H^{\leq n}(\Uu , Tel (G_i)) $$
induces an isomorphism on homotopy groups. On the left, the homotopy groups
are just the direct limit of the $\pi _i H^{\leq n}(\Uu , G_i )$, while on
the right the
$\pi _i H^{\leq n}(\Uu , Tel (G_i))$
may be viewed as the inverse limit over a finite diagram where the elements
are spaces of maps from a finite CW complex into the component spaces of
$Tel (G_i)$ (and the maps in the diagram are fibrations). But the set of
homotopy classes of maps of a compact space into a mapping telescope is
just the direct limit of the space of homotopy classes of maps from that
space into each component; thus on the right we obtain the same direct
limit as on the left, so we get an isomorphism. We obtain the desired
isomorphism of homotopy groups, proving that $Tel (HG_i)\rightarrow H(Tel
(G_i))$ is a homotopy equivalence.

We apply this to the sequence of functors defined inductively by $G_1=HT$
and $G_{i+1}=H(G_i)$ with $g_i:=i_{G_i}$ the natural inclusion. Note that,
since the $G_i$ are functors, the $g_i$ are morphisms of functors. Put $$
FT:= Tel (G_i) = Tel (HT\rightarrow H(HT)\rightarrow H(H(HT))\rightarrow
\ldots ).
$$
We would like to show that the natural morphism $FT\rightarrow H(FT)$ is a
homotopy equivalence. Note first of all that $FT$ is $n$-truncated and of
CW-type. Thus it suffices to prove that $(FT)_X\rightarrow (H(FT))_X$
induces an isomorphism of homotopy groups for all $X\in \Xx$.

For simplicity, we ignore base points in the following discussion.

Note that
$$
\pi _j (FT) = \lim _{\rightarrow }\pi _j (G_i). $$
We will consider the map $\pi _j (G_i)\rightarrow \pi _j(H(FT))$ induced by
the map $G_i\rightarrow FT$. The composed map $$
G_i\rightarrow H(FT )=H(Tel (G_i ))
$$
is equal to the composition
$$
G_i \rightarrow H(G_i) \rightarrow Tel (H(G_k))\rightarrow H(Tel (G_k)). $$
Note that in the $Tel (H(G_k))$ the maps are the $H(i_{G_k})$. The composed
map $G_i\rightarrow G_{i+1} \rightarrow Tel (H(G_k))$ is homotopic to the
map from $G_i$, by the nature of the telescope construction.

In fact, the construction $H$ is natural on functors, and the natural
inclusion $i$ is natural, that is if $f:S\rightarrow T$ is a morphism of
functors, then $i_t\circ f = H(f)\circ i_S$. Thus we obtain a morphism of
telescopes
$$
Tel (G_i)\rightarrow Tel (H(G_k)).
$$
The composition into $H(Tel (G_k))$ is equal to the natural inclusion
$i_{Tel (G_i)}$. Thus (in view of the above claim) it suffices to prove
that the morphism
$$
Tel (G_i)\rightarrow Tel (H(G_k)).
$$
induces an isomorphism of homotopy groups. Equivalently, we have to show that
$$
\lim _{\rightarrow } \pi _i (G_k) \rightarrow \lim _{\rightarrow} \pi _i
(H(G_k))
$$
is an isomorphism.

We will show that, for any $k$, if $\eta \in \pi _i (HG_k)$ then there is
$\eta '\in \pi _i (G_{k+1})$ such that the image of $\eta$ in $\pi _i
(HG_{k+1})$ (by the morphism $Hg_k$) is equal to the image of $\eta '$ by
the morphism $i_{G_{k+1}}$. And we will show that if $\zeta \in \pi _i
(G_k)$ such that $i_{G_k}(\zeta )=0$ in $\pi _i (HG_k)$, then $g_k (\zeta
)=0$ in $\pi _i(G_{k+1})$. From these two statements, it follows that the
morphism of direct limits is an isomorphism. To prove the two statements,
recall that $G_{k+1}=HG_k$ and $g_k= i_{G_k}$. The second statement follows
immediately, since the two maps in question are the same. For the first
statement, let $\eta '$ be the same element as $\eta$; from Lemma
\ref{desc14}, we get $Hg_k(\eta ):=H(i_{G_k})(\eta )= i_{HG_k}(\eta )=:
i_{G_{k+1}}(\eta ')$.

We have proved that the morphism $i_{FT}:FT\rightarrow H(FT)$ (which is a
morphism of functors), is a homotopy equivalence, under the assumptions
that $T$ is an $n$-truncated homotopy functor of CW type, and $\Xx$ is
quasi-compact. By
the characterization of Proposition \ref{desc05}, $FT$ is a flexible sheaf.
We call $FT$ the {\em flexible sheaf associated to $T$}. We have the
following universal property.

\begin{theorem}
\mylabel{desc15}
Suppose $\Xx$ is quasi-compact and $T$ is an $n$-truncated homotopy functor
of CW type. Then for any homotopy sheaf $R$ and any morphism
$A:T\rightarrow R$, the space of ways of completing this to a homotopy
commutative triangle $T\rightarrow FT \rightarrow R$ is nonempty and weakly
contractible. \end{theorem}

The proof is the same as for the case which will be treated
in the next section, so we refer there for the proof.

\newpage

\setcounter{page}{84}

\setcounter{section}{4}

\numero{Homotopy group sheaves}

Suppose $T$ is a
flexible functor from $\Xx$ to $Top$, and suppose $i\geq 1$. Suppose $t\in
\Gamma
(\Xx , T)$. Then we define
$$
\pi ^{\rm pre}_i(T, t)
$$
to be the presheaf which associates to each object $Y\in \Xx$ the group
$\pi _i (T_Y, t_Y)$. Note that the higher homotopies in the definitions of
$T$ and $t$ yield natural restriction morphisms $$
\pi _i (T_Y,t_Y)\rightarrow \pi _i (T_Z,t_Z) $$
whenever $Z\rightarrow Y$ is a morphism in $\Xx$. These are strictly
compatible with composition, giving a contravariant functor from $\Xx$ to
the category of groups (abelian groups if $i\geq 2$). This gives the
desired presheaf. Similarly, for $i=0$ define a presheaf of sets $\pi ^{\rm
pre}_0(T)$ on $\Xx$ (which doesn't depend on a choice of basepoint).

Suppose $\Xx $ is a site. Then
let $\pi _i(T,t)$ be the sheaf associated to the presheaf $\pi ^{\rm
pre}_i(T,t)$ (dropping the basepoint if $i=0$).

Suppose $T\rightarrow S$ is a morphism of flexible functors, and $\ast
\rightarrow T\rightarrow S$ is a diagram of flexible functors (that is, a
flexible functor from $\Xx \times I^{(2)}$ to $Top$ restricting to the
given ones, and to the given morphism $T\rightarrow S$). Let $t$ denote the
point $\ast \rightarrow T$ and let $s$ denote the point $\ast \rightarrow
S$. Then we obtain a morphism
of presheaves
$$
\pi ^{\rm pre}(T,t)\rightarrow \pi ^{\rm pre}_i(S,s), $$
sheafifying to give a morphism of sheaves if $\Xx$ is a site. Again, there
is no need for the stuff concerning the basepoints if $i=0$.

Recall the following definition from \cite{Illusie1} (although he didn't use
the same terminology!): a morphism $T\rightarrow S$ of flexible functors on
a site $\Xx$ is an {\em Illusie weak equivalence} if it induces
isomorphisms on all homotopy group sheaves. We will show below that if $S$
and $T$ are (truncated?) weak flexible sheaves of CW-type, then an Illusie
weak equivalence is in fact an equivalence.

\begin{lemma}
\mylabel{bigPi01}
Suppose $T$ is a flexible functor on a site $\Xx$. Then the morphism
$T\rightarrow HT$ is an Illusie weak equivalence. \end{lemma}
{\em Proof:}
???
\eop

\subnumero{Towers}
A {\em tower of topological spaces} is a sequence of spaces $A_i$ indexed
by $i\in \nn$, together with maps $\tau _i : A_i\rightarrow A_{i-1}$. This
may be considered as a contravariant functor $A$ from the category $\Nn$
(consisting of one object for each natural number and one morphism from $i$
to $j$ whenever $i\leq j$) to $Top$. The {\em limit} of the tower is the
space $\Gamma (\Nn , A)$. Note that if the maps $\tau _i$ are fibrations,
then this limit is homotopic to the usual projective limit. We say that the
tower is {\em vertically truncated} if the maps $\tau _i$ are weak
equivalences for $i\gg 0$. We say that it is {\em horizontally truncated}
if each space $A_i$ is truncated, and {\em uniformly horizontally
truncated} if there is an $n$ such that each $A_i$ is $n$-truncated. Note
that vertically truncated plus horizontally truncated implies uniformly
horizontally truncated; in this case we say $A$ is completely truncated.

The {\em fibers} of the tower are
$$
F_i(\eta ):= \; \mbox{fiber of} \; \tau _i: A_i\rightarrow A_{i-1} \;
\mbox{over} \; \eta \in A_{i-1}.
$$
If $\eta \in \Gamma (\Nn , A)$ then we obtain the fibers $F_i(\eta )$ for
all $i$ (using as basepoint the projection of $\eta$ in $A_{i-1}$).

\subnumero{Nonabelian spectral sequences} The usual notion of spectral
sequence involves abelian component groups. This can be generalized
slightly as follows. (I assume that this materiel is well-known but don't give
a specific reference.)

A {\em nonabelian complex ending at $i=n$} (perhaps it would be better to
replace this notion by the notion of ``crossed complex'' ?) is a sequence of
groups $G_i$ for $i< n$ with maps $d_i: G_i\rightarrow G_{i+1}$, a pointed set
$(G_n,p)$ with an action $d _{n-1}$ of $G_{n-1}$ on $G_n$, and a boolean
function
$d_n:G_n\rightarrow \{ 0,1 \}$, subject to the following axioms. The group
$G_i$ is abelian for $i\leq n-2$. The compositions $d_id_{i-1}$ are zero
for $i\leq n-2$, $d_{n-2}(G_{n-2})$ is a normal subgroup of the stabilizer
of $p$ under the action $d _{n-1}$, the action preserves the boolean
function, and $d_n(p)=0$. The {\em cohomology} $H^i$ of such a complex is
the kernel of $d_i$ modulo the image of $d_{i-1}$ for $i\leq n-2$;
$H^{n-1}$ it is the stabilizer of $p$ modulo the image of $d_{n-2}$, and
$H^n$ is the quotient of $d_n^{-1}(0)\subset G_n$ by the action of
$G_{n-1}$.

A {\em nonabelian spectral sequence in $p+q\leq n$} consists of the
following data, subject to the axioms given afterwards. The data consist of
families of nonabelian complexes
$$
\{ E_r^{p,q}, d_r: E^{p,q}_r \rightarrow E^{p+r, q+1-r}_r\} $$
for integers $r\geq 2$, such that the complexes $E^{p,q}_r$ end at $p+q=n$;
together with isomorphisms
$$
H^{p,q}_r \cong E_{r+1}^{p,q}
$$
between the cohomology of the complexes $E_r$ and the next terms $E_{r+1}$.

Note that for $p+q=n-1$ the cohomology groups are nonabelian, and for
$p+q=n$ they are sets, pointed by the image of the basepoint in the
previous set.

We say that a nonabelian spectral sequence {\em abuts to a collection $\{
K^i \} _{i\leq n}$},
if $K^i$ is a set for $i=n$, a group for $i\leq n-1$, abelian for $i\leq
n-2$, if for each $i$ there is a cofiltration (that is, a collection of
quotients) $\{ K^i\rightarrow F_jK^i\}$ such that at each $(p,q)$ the
sequence $E^{p,q}_r$ eventually stabilizes, and the stabilized value is $$
E^{p,q}_{\infty}={\rm ker}(F_qK^{p+q}\rightarrow F_{q+1}K^{p+q}). $$
For
$p+q=n-1$, the quitients are quotient groups, and for $p+q=n$, the
quotients are pointed sets. In the last case, the kernel denotes the set of
elements mapping to the base point.
When the spectral sequence starts with a certain term $E_2$ and abuts to
$K^i$ we commonly write
$$
E_2^{p,q}\Rightarrow K^{p+q}.
$$

\begin{proposition}
\mylabel{pi01}
Suppose $A$ is a vertically truncated tower of spaces, and suppose $\eta
\in \Gamma (\Nn , A)$. Then there is a (nonabelian) spectral sequence in
$p+q\leq 0$ with
$$
E_2^{p,q}= \pi _{-p-q}(F_{-q},\eta )\Rightarrow \pi _{-p-q}(\Gamma (\Nn ,
A),\eta ). $$
The filtration on $\pi _{-p-q}(\Gamma (\Nn , A))$ is $$
F_q\pi _{-p-q}(\Gamma (\Nn , A),\eta ):= \;\; \mbox{the image in}\;\; \pi
_{-p-q}(A_{-q},\eta ).
$$
\end{proposition}
{\em Proof:}
Left to the reader.(??? ref ???)
\eop

{\em Remark:} The spectral sequence degenerates in the case of a tower with
only one nontrivial fiber, or in the case of a Postnikov tower. A fibration
may be considered as a tower with two nontrivial fibers, and here the
spectral sequence gives the long exact sequence of homotopy. In the case of
a tower where
$F_i$ a $K(G, n-i)$ space, then the $E_2$ term is nonzero for $p=2q-n$,
which form a line along the $d_2$ differential. The $d_2$ differential
gives a complex, and the cohomology of the complex is the end of the
spectral sequence. This gives the result discussed in \S 1.

\newpage

\setcounter{section}{5}

\numero{Fibrations and fiber products}

Suppose $f:A\rightarrow C$ and $g:B\rightarrow C$ are two maps of
topological spaces. We define the {\em pathwise fiber product} $A\times
^{\rm path}_C B$ to be the space of triples $(a,b, \gamma )$ where $a\in
A$, $b\in B$, and $\gamma$ is a path in $C$ from $f(a)$ to $g(b)$. We
obtain a diagram $$
\begin{array}{ccc}
A\times ^{\rm path}_C
B&\rightarrow &A \\
\downarrow & & \downarrow \\
B&\rightarrow & C
\end{array}
$$
with a homotopy of commutativity (given by the paths $\gamma$).

Suppose $R$, $S$, and $T$ are three contravariant flexible functors from a
category $\Xx$ to $Top$. Suppose we have morphisms $F:R\rightarrow T$ and
$G: S\rightarrow T$. We say that a diagram (that is a flexible functor from
$\Xx \times I^2$ to $Top$)
$$
\begin{array}{ccc}
Q&\rightarrow &R \\
\downarrow & & \downarrow \\
S&\rightarrow & T
\end{array}
$$
is {\em path-cartesian} if for each object $X\in \Xx$, the diagram $$
\begin{array}{ccc}
Q_X&\rightarrow &R_X \\
\downarrow & & \downarrow \\
S_X&\rightarrow & T_X
\end{array}
$$
together with its homotopy of commutativity (given by the structure of
diagram above) is homotopic to the diagram plus homotopy $$
\begin{array}{ccc}
R_X\times ^{\rm path}_{T_X}S_X&\rightarrow &R_X \\ \downarrow & & \downarrow \\
S_X&\rightarrow & T_X \, .
\end{array}
$$

{\em Notation:} For brevity, we sometimes denote a commutative square as
above by
$$
Q\rightarrow R\times S\stackrel{\textstyle\rightarrow}{\rightarrow}T. $$

\begin{lemma}
\mylabel{fibr001}
If $Q\rightarrow R\times S\stackrel{\textstyle\rightarrow}{\rightarrow}T$
is a path-cartesian diagram of flexible functors, and if $U$ is a flexible
functor, then this induces a path-cartesian diagram of spaces $$
Mor _{\infty}(U,Q)\rightarrow Mor _{\infty}(U,R) \times Mor
_{\infty}(U,S)\stackrel{\textstyle\rightarrow}{\rightarrow} Mor
_{\infty}(U,T).
$$
\end{lemma}
{\em Idea of proof:}
The statement of the lemma amounts to saying that the morphism from the
space of diagrams $$
U\rightarrow Q\rightarrow R\times S
\stackrel{\textstyle\rightarrow}{\rightarrow}T $$
(that is, flexible functors on $(I^2)^{\ast}\times \Xx$ where the asterisk
denotes the category with a final object---corresponding to the
$U$---added), to the space of diagrams
$$
U\rightarrow R\times S\stackrel{\textstyle\rightarrow}{\rightarrow}T, $$
is a weak equivalence. This is shown by the same type of argument as
outlined in \S 2, using the Postnikov tower for $Top$. In the following
paragraph we indicate the idea for the higher part (not for the first two
stages). I have not checked the details.

At each stage, we look at the space of flexible functors from
$(I^2)^{\ast}\times \Xx$ into $\tau _{\leq n}Top$, restricting to the given
ones on the subcategories isomorphic to $(I^2)\times \Xx$ corresponding to
the given diagrams (with $Q$ and $U$). However, the space of morphisms is
not weakly contractible at each stage in the Postnikov tower; instead, this
tower yields a spectral sequence by the discussion at the beginning of the
next section (note that there is no circularity!). In each row of the
spectral sequence is the cohomology of the mapping cone of the map on
complexes $PC_{\cdot}$ corresponding to the inclusion $\Nn \subset \Mm$ of
continuous semi-categories involved in the argument ($\Mm$ is the one for
$(I^2)^{\ast}\times \Xx$ and $\Nn$ corresponds to the union of the two
subcategories in question). In turn, the cohomology in the row can be
calculated by a spectral sequence: it becomes the spectral sequence for the
presheaf (\v{C}ech) cohomology on $Fl(\Xx )$ of the complex of presheaves
$$ (\phi :X_0\leftarrow X_n) \mapsto
$$
$$
\pi _n(Hom (U_{X_0}, Q_{X_n}))\rightarrow \pi _n(Hom (U_{X_0}, (R\times
S)_{X_n})) \rightarrow
\pi _n(Hom (U_{X_0}, T_{X_n})),
$$
which is exact in the middle. Denote by $F_n$ the kernel, and note that the
cokernel is then $F_{n-1}$. The differential $d_3$ comes from the class in
$Ext ^3_{Fl(\Xx )}(F_{n-1}, F_n)$---a class which can be seen as coming
from the long exact sequence of homotopy. Somewhat surprisingly, this class
vanishes: this is a consequence of the fact that we have maps of flexible
sheaves $Q\rightarrow R\times S$, and not just the long exact sequence.
Thus $d_3$ vanishes and the spectral sequence degenerates. The result is
that the cohomology of the complex $PC_{\cdot}$ for the $n$th row in the
first spectral sequence, is an extension of the cohomology of $F_n$ on
$Fl(\Xx
)$, by the cohomology of $F_{n-1}$ on $Fl(\Xx )$ shifted by two. Now back
to the first spectral sequence: the differential $d_2$ is an isomorphism
between the quotient equal to the cohomology of $F_{n-1}$ in the $n-1$st
row, and the subobject equal to the cohomology of $F_{n-1}$ in the $n$th
row. Thus the differential $d_2$ kills off everything, so we obtain a weak
equivalence in the limit of the tower.
\eop

\begin{corollary}
\mylabel{fibr01}
If $Q\rightarrow R\times S\stackrel{\textstyle\rightarrow}{\rightarrow}T$
and $Q'\rightarrow
R\times S\stackrel{\textstyle\rightarrow}{\rightarrow}T$ are two
path-cartesian diagrams, then the space of homotopies between them (i.e.
functors from $\Xx \times I^2 \times \overline{I}$ to $Top$) is nonempty
and weakly contractible. \end{corollary}
{\em Proof:}
This follows from the previous lemma by a composition-type argument (one
may have to reprove that the space of certain compositions of diagrams is
weakly contractible, but this will be the same as usual). \eop

{\em Remark:} If $R$, $S$, and $T$ are {\em functors} (and the morphisms
$F$ and $G$, morphisms of functors) then the pathwise fiber product
$R\times ^{\rm path}_TS$ is also a functor, and the diagram (NB this is a
flexible diagram because of the path ...) is path-cartesian. On the other
hand, we can always replace the diagram $R\times S
\stackrel{\textstyle\rightarrow}{\rightarrow}T$ by a homotopic diagram of
functors. The above construction gives a path-cartesian diagram. Thus, for
any pair of morphisms of flexible functors $F:R\rightarrow T$ and
$G:S\rightarrow T$ there exists a path-cartesian diagram. The above lemma
says that it is essentially unique.

In a slight abuse of notation, we call the flexible functor $Q$ the {\em
pathwise fiber product of $R$ and $S$ over $T$.}

There is probably a canonical choice for this pathwise fiber product even
in the flexible case---we leave this to the reader.

\subnumero{Fibrations}
A {\em fibration diagram} is a path-cartesian diagram of the form $$
\begin{array}{ccc}
R&\rightarrow &S \\
\downarrow & & \downarrow \\
\ast &\rightarrow & T\, .
\end{array}
$$
We say that $R$ is the {\em fiber} of the morphism $S\rightarrow T$ over
the point $x\in \Gamma (\Xx ,T)$ corresponding to the morphism $\ast
\rightarrow T$. From the above discussion, if $F:S\rightarrow T$ is any
morphism of flexible functors, and if $x\in \Gamma (\Xx , T)$, then there
exists a fibration diagram (which we denote more briefly by $R\rightarrow
S\rightarrow T$ with the point understood). Any two such diagrams are
homotopic by a very well defined homotopy. Thus we may speak of the fiber
of a morphism over a point.

\subnumero{Path and loop spaces}
We can define path spaces and loop spaces. If $x,y\in \Gamma (\Xx , T)$ we
define the {\em path space} $P^{x,y}T$ to be the pathwise fiber product of
$\ast$ and $\ast$ over $T$, via the morphisms $x$ and $y$. We have a
path-cartesian diagram
$$
\begin{array}{ccl}
P^{x,y}T&\rightarrow &\ast \\
\downarrow & & \, \downarrow y\\
\ast &\stackrel{x}{\rightarrow }& T\, .
\end{array}
$$
If $x=y$ then we call this the {\em loop space} $\Omega ^xT:= P^{x,x}T$.

\begin{lemma}
\mylabel{fibr02}
If
$$
\begin{array}{ccc}
\Gamma (\Xx ,Q)&\rightarrow &\Gamma (\Xx ,R) \\ \downarrow & & \, \downarrow \\
\Gamma (\Xx ,S)&\rightarrow & \Gamma (\Xx ,T) \end{array}
$$
is a path-cartesian diagram of flexible functors on a category $\Xx$, then
the induced diagram
$$
\begin{array}{ccc}
\Gamma (\Xx ,Q)&\rightarrow &\Gamma (\Xx ,R) \\ \downarrow & & \, \downarrow \\
\Gamma (\Xx ,S)&\rightarrow & \Gamma (\Xx ,T) \end{array}
$$
is a path-cartesian diagram of spaces. If $\Xx$ is a site then the diagram $$
\begin{array}{ccc}
HQ&\rightarrow &HR \\
\downarrow & & \, \downarrow \\
HS&\rightarrow & HT
\end{array}
$$
is a path-cartesian diagram of flexible functors; and the same for $H_{\infty}$.
\end{lemma}
{\em Proof:}
The statement for global sections follows directly from Lemma \ref{fibr001}
applied with $U= \ast _{\Xx}$. The second statement follows by taking the
direct limit over sieves. \eop

\begin{corollary}
\mylabel{fibr03}
Suppose
$$
\begin{array}{ccc}
Q&\rightarrow &R \\
\downarrow & & \, \downarrow \\
S&\rightarrow & T
\end{array}
$$
is a path-cartesian diagram of flexible functors on a site $\Xx$. If $R$,
$S$ and $T$ are weak flexible sheaves, then so is $Q$.
\end{corollary}
{\em Proof:}
The morphism $Q\rightarrow HQ$ fits into a morphism of path-cartesian
diagrams, where by hypothesis the other three morphisms ($R\rightarrow HR$
etc.) are weak equivalences on each object. This implies that
$Q_X\rightarrow (HQ)_X$ is a weak equivalence for each object $X$. \eop

In particular, if $T$ is a weak flexible sheaf and $x, y\in \Gamma (\Xx , T)$
then the path space $P^{x,y}T$ is a weak flexible sheaf. In the next subsection,
we give a sort of converse to this. For the moment, note the following
consequence of the lemma (which works for loop spaces, too, of course).

\begin{corollary}
\mylabel{fibr04}
Suppose $T$ is a flexible functor on a site $\Xx$ and $x,y\in \Gamma (\Xx
,T)$. Then $P^{x,y}(HT)$ is weakly equivalent to $H(P^{x,y}T)$. Similarly,
if $x,y\in \Gamma (\Xx ,\tau _{\infty}T)$ and if $x',y'$ denote their
images in $\Gamma (X,T)$ then
$P^{x,y}(H_{\infty}T)$ is weakly equivalent to $H_{\infty}(P^{x',y'}T)$.
\end{corollary}
{\em Proof:}
By the lemma, $H(P^{x,y}T)$ fits into a path-cartesian diagram with $\ast$,
$\ast$ and $HT$.
\eop

\subnumero{Flexible functors whose path spaces are weak sheaves} Suppose
that $T$ is a flexible functor from a site $\Xx$ to $Top$. We say that {\em
the path spaces of $T$ are weak sheaves} if for any object $X\in \Xx$ and
any pair of points $x,y\in T_X$, the path space $P^{x,y}(T|_{\Xx /X})$ is a
weak flexible sheaf on the site $\Xx /X$.

\begin{lemma}
\mylabel{fibr05}
Suppose that the path spaces of $T$ are weak sheaves. Then for each object
$X\in \Xx$, the morphism $\pi _0 (T_X) \rightarrow \pi _0 ((HT)_X)$ is
injective. For each object $X\in \Xx$ and each point $x\in T_X$, the
morphisms $\pi _i (T_X , x)\rightarrow \pi _i ((HT)_X, x)$ are isomorphisms
for $i\geq 1$.

Similarly, for each $X\in
\Xx$ and nested pair of sieves $\Bb \subset \Bb '$ over $X$, the morphism
$\pi _0 (\Gamma (\Bb ', T|_{\Bb '})) \rightarrow \pi _0 (\Gamma (\Bb ,
T|_{\Bb}))$ is injective; and for $x\in \Gamma (\Bb ', T|_{\Bb '})$, the
morphisms $\pi _i (\Gamma (\Bb ', T|_{\Bb '}), x)\rightarrow \pi _i (\Gamma
(\Bb , T|_{\Bb}), x)$ are isomorphisms for $i\geq 1$.
\end{lemma}
{\em Proof:}
The fact that $P^{x,x}$ is a weak sheaf means that the map
$$
\pi _i (T_X,x)= \pi _{i-1}((P^{x,x}T)_X,x)\rightarrow
$$
$$
\pi
_{i-1}((HP^{x,x}T)_X,x)= \pi _{i-1}((P^{x,x}HT)_X,x)=\pi _i ((HT)_X,x)
$$
is an isomorphism. For $x,y\in T_X$, the fact that $P^{x,y}$ is a weak
sheaf means that
$$
\pi _0(P^{x,y}T_X)\rightarrow \pi _0((HP^{x,y}T)_X)=\pi _0(P^{x,y}(HT)_X) $$
is an isomorphism; but this space is nonempty if and only if $x$ and $y$
are in the same path component (of $T_X$ on the left, of $(HT)_X$ on the
right). This implies that $\pi _0(T_X)\rightarrow \pi _0((HT)_X)$ is
injective. The statements in the second paragraph are the same. \eop

Note that the second paragraph applies in particular to the case $\Bb '=
\Xx /X$; then we can replace $\Gamma (\Bb ', T|_{\Bb '})$ by $T_X$ in the
statement.

\begin{lemma}
\mylabel{fibr06}
Suppose that the path spaces of $T$ are weak sheaves. Then the path spaces
of $HT$ are weak sheaves.
\end{lemma}
{\em Proof:}
Suppose $a,b\in (HT)_X$. By the direct limit construction of $(HT)_X$, we
may assume that there is a sieve $\Bb _0$ over $X$ and sections $a',b'\in
\Gamma (\Bb _0 , T)$ yielding the points $a,b$ in $(HT)_X$. Using Lemma
\ref{desc02}, we may suppose that the images of $a,b$ in $\Gamma (\Bb _0,
HT)$ are the same as the images of $a',b'$ obtained from the morphism
$T\rightarrow HT$. For any $Y\in \Bb _0$ we obtain a sequence of weak
equivalences $$
(P^{a',b'}T)_Y\sim H(P^{a'_Y,b'_Y}T|_{\Xx /Y})_Y\sim $$
$$
(P^{a'_Y,b'_Y}
H(T|_{\Xx /Y}))_Y= (P^{a'_Y,b'_Y}HT)_Y\sim (P^{a,b}HT)_Y. $$
The first is since the path spaces of $T$ are weak sheaves; the second by
Corollary \ref{fibr04}; the third since the construction $HT$ over the
object $Y$ depends only on $T|_{\Xx /Y}$; and the last using Lemma
\ref{desc02} as mentioned above. Composing these, we obtain a weak
equivalence $$
P^{a',b'}T\sim P^{a,b}HT
$$
of flexible functors over $\Bb _0$.
Now suppose
$\Bb \subset \Bb _0$. Then
$$
P^{a',b'}\Gamma (\Bb , T)= \Gamma (\Bb , P^{a',b'}T)\sim \Gamma (\Bb ,
P^{a,b}HT), $$
and taking the direct limit (over sieves $\Bb \subset \Bb _0$) we obtain $$
(P^{a,b}HT)_X \sim \lim _{\rightarrow } P^{a',b'}\Gamma (\Bb , T) \sim \lim
_{\rightarrow} \Gamma (\Bb ,
P^{a,b}HT)\sim H(P^{a,b}HT)_X.
$$
For the first and last, note that the construction $H$ may be obtained (up
to equivalence) by taking the limit over any cofinal set of sieves, in
particular over the sieves contained in a given one. Thus the morphism
$$
P^{a,b}HT\rightarrow H(P^{a,b}HT)
$$
is a weak equivalence. This implies that $P^{a,b}HT$ is a weak sheaf. \eop

{\em Remark:} The proof above is much easier if the points $a,b$ come from
$T_X$, but this is not necessarily the case.

\begin{theorem}
\mylabel{fibr07}
Suppose that $T$ is a contravariant flexible functor on a site $\Xx$, such
that the path spaces of $T$ are weak sheaves. Then $HT$ is a weak flexible
sheaf.
\end{theorem}
{\em Proof:}
By replacing $\Xx$ by $\Xx /X$ we may assume that $\Xx$ has a final object
$X$ and we prove that for any sieve $\Bb$ over $X$, the morphism
$(HT)_X\rightarrow \Gamma (\Bb , HT)$ is a weak equivalence. By Lemma
\ref{fibr06}, the path spaces of $HT$ are sheaves; then by Lemma
\ref{fibr05}, to prove the weak equivalence it suffices to prove that
$(HT)_X\rightarrow \Gamma (\Bb , HT)$ induces a surjection on $\pi _0$.

Suppose $y\in \Gamma (\Bb , HT)$.
Then for each $U\in \Bb$ there is a sieve $\Bb '_U$ over $U$ such that for
$V\in \Bb '_U$, $y_U|_{V}$ is in a path component of $(HT)_V$ which
contains a point of $T_V$. This follows from the construction of $HT$:
$(HT)_U$ is the direct limit of the $\Gamma (\Bb ', T|_{\Bb '})$, so any
point $y_U$ occurs in the space corresponding to one of the sieves $\Bb
'_U$. By one of the axioms for a site, the union of all of the $\Bb '_U$
forms a covering family of $X$, generating a sieve $\Dd$. If $V\in \Dd$
then $y_V$ is in a path component of $(HT)_V$ containing a point $z_V$ of
$T_V$. We can now construct a point $z\in \Gamma (\Dd , T)$ mapping to the
path component of $y$ in $\Gamma (\Dd , HT)$. To do this, proceed by
transfinite induction on the primitive cells (of dimension $\geq 1$) in the
continuous category $\Mm$ associated to $\Xx$. Each time we add a cell, we
can map it into the space $T_V$ involved, in a way homotopic to the
original map into $(HT)_V$ plus the previous homotopies on the boundary
pieces, because the inclusion from the union of components of $T_V$
containing the boundary into the path component of $(HT)_V$ containing the
original cell is a weak equivalence (note that the boundary is nonempty
since the cells are of positive dimension, which in turn is due to the fact
that we have already chosen the zero-dimensional cells $z_V$).

We have shown that the image of our point $y$ in $\Gamma (\Dd , HT)$ is in
a path component coming from a point in $\Gamma (\Dd , T)$, hence in
particular from a point $z$ in $(HT)_X$. (Here there is also the question
of the homotopy between two maps from $\Gamma (\Dd ,T)$ to $\Gamma (\Dd ,
HT)$, the one factoring through $(HT)_X$ and the other coming from the
morphism $T\rightarrow HT$---this is essentially the main lemma above the
main proposition in the section on descent.)

Finally, note that by Lemma \ref{fibr05}, $\Gamma (\Bb , HT )$ is weakly
equivalent to a union of some path components of $\Gamma (\Dd , HT)$. Thus,
if there is a path from $y$ to $z$ in $\Gamma (\Dd , HT)$ then there is a
path from $y$ to $z$ in $\Gamma (\Bb , HT)$. This completes the proof of
surjectivity on $\pi _0$. We obtain the theorem.
\eop

\begin{corollary}
\mylabel{fibr08}
Suppose $T$ is an $n$-truncated contravariant flexible functor from a site
$\Xx$ to $Top$. Define inductively $H^k T:= H(H^{k-1}T)$. Then $H^{n+2}T$
is a weak flexible sheaf. If $\Xx$ is quasi-compact and $T$ is of CW-type
(i.e. takes values in $Top ^{CW}$, the subcategory of spaces homotopy
equivalent to CW complexes) then $H^{n+2}T$ is a flexible sheaf.
\end{corollary}
{\em Proof:}
First we treat the case $n=0$. If $T$ is a $0$-truncated flexible functor,
it is a flexible sheaf if and only if $\pi _0 (T)$ is a sheaf. And $HT$ is
again a $0$-truncated flexible functor with
$$
\pi _0(HT)_X=\lim _{\rightarrow , \Bb
}\Gamma (\Bb , \pi _0 (T)).
$$
In other words, the operation $H$ reduces, on the level of $\pi _0$, to the
usual operation used to construct the sheaf associated to a presheaf of
sets. As is well known, applying this operation twice gives a sheaf. Thus
$\pi _0(H^2T)$ is a sheaf of sets, so $H^2T$ is a flexible sheaf.

For $n>0$ we proceed by induction on $n$. Now, note that the operation $H$
commutes with taking path spaces. The path spaces $P^{x,y}T$ are
$n-1$-truncated, so by the inductive hypothesis $P^{x,y}(H^{n+1}T)\sim
H^{n+1}(P^{x,y}T)$ are weak flexible sheaves. The theorem implies that
$H^{n+2}T$
is a weak flexible sheaf.

For the second statement, note that if $\Xx$ is quasi-compact and $T$ is of
CW-type, then $HT$ is also of CW-type. A weak flexible sheaf of CW-type is
a flexible sheaf. \eop

{\em Remark:} Since $H_{\infty}^{n+2}T$ and $H^{n+2}$ are weakly
equivalent, under the hypotheses of the previous corollary,
$H_{\infty}^{n+2}T$ is a flexible sheaf of CW-type.

\begin{corollary}
Put $F_{\infty}T:= \lim _{\rightarrow , m} H_{\infty}^mT$ (with the
morphisms in the direct system being the same as those used at the end of
the previous section). If $T$ is an $n$-truncated flexible functor of
CW-type then $F_{\infty}T$ is a flexible sheaf of CW type, with a standard
morphism $T\rightarrow F_{\infty}T$. \end{corollary} {\em Proof:}
The morphisms from $H^{n+2}_{\infty}T$ on are equivalences, so
$F_{\infty}T\sim H_{\infty}^{n+2}T$. The standard morphism $\tau
_{\infty}T\rightarrow F_{\infty}T$ yields a morphism from $T$, if $T$ is of
CW type. \eop

{\em Remark:} Here, as below,
we allow ourselves to pick an inverse when the inverse is ``very well
defined'', that is the set of choices is parametrized by a weakly
contractible space. We ignore what needs to be said to account for the fact
that the choices made in various places aren't compatible, but only up to
very well defined homotopy.

\begin{theorem}
\mylabel{fibr09}
Suppose $T$ is an $n$-truncated contravariant flexible functor of CW type
and $S$ is a weak flexible sheaf on a site $\Xx$. Then the space of
composable pairs of morphisms $T\rightarrow F_{\infty}T\rightarrow S$ such
that the morphism $T\rightarrow F_{\infty}T$ is the standard one, maps to
$Mor (T,S)$ by a fibration with weakly contractible fiber. \end{theorem}
{\em Proof:}
The operation $T\mapsto \Gamma _{\infty}(\Bb , T|_{\Bb})$ is functorial
with respect to morphisms of the flexible functor $T$; thus $T\mapsto
H_{\infty}T$ and even $T\mapsto F_{\infty}T$ are functorial in $T$; and the
morphism
$T\rightarrow F_{\infty}T$ is a natural
transformation (all of this being in a ``very well defined'' sense). We show
that the space of homotopy classes of composable pairs $T\rightarrow
F_{\infty}T\rightarrow S$ restricting to the given morphism $T\rightarrow
S$, is nonempty and has only one path component. Note that we may assume
that $S$ is of CW type. To see that the set of diagrams in question is
nonempty, apply the functoriality of $F_{\infty}$ to the morphism
$T\rightarrow S$: we get a commutative diagram $$
\begin{array}{ccc}
T&\rightarrow & F_{\infty}T \\
\downarrow && \downarrow \\
S&\rightarrow & F_{\infty}S,
\end{array}
$$
where the bottom arrow is an equivalence (since $S$ is a weak sheaf).
Composing the right vertical arrow with this inverse of the bottom
equivalence gives the desired diagram. To see that there is only one path
component, suppose that we have a diagram $T\rightarrow
F_{\infty}T\rightarrow S$. Apply $F_{\infty}$ to get a morphism of
commutative triangles (we write the triangles vertically): $$
\begin{array}{ccc}
T&\rightarrow & F_{\infty}T \\
\downarrow && \downarrow \\
F_{\infty}T&\rightarrow & F_{\infty}(F_{\infty}T) \\ \downarrow && \downarrow \\
S&\rightarrow & F_{\infty}S .
\end{array}
$$
Here the middle and bottom horizontal arrows are equivalences, as is the
upper right vertical arrow. The composition of the middle horizontal
equivalence and that of the upper right is homotopic to the identity of
$F_{\infty}T$ (this is essentially Lemma \ref{desc02} again). The
composition of the two vertical arrows on theright is the morphism
obtained from $T\rightarrow S$ by functoriality. The commutativity of the
diagram (up to homotopy) implies that the diagram on the left is homotopic
to the diagram constructed above, so there is only one path component.

The same argument works if we have a morphism of weak flexible sheaves
$S\rightarrow S'$ and want to find a diagram $T\rightarrow
F_{\infty}T\rightarrow S$ having composition with this morphism equal to a
certain fixed diagram for $T\rightarrow F_{\infty}T\rightarrow S'$. Apply
this to the diagonal morphism $S\rightarrow S\times S$, and note that the
space of diagrams $T\rightarrow S\rightarrow S\times S$ restricting to a
morphism $(\eta , \xi ):T\rightarrow S\times S$ is the path space $P^{\eta
, \xi}Mor (T,S)$. We find that all the higher path spaces of the space of
diagrams as desired in the lemma, have exactly one connected component.
This implies that the space of diagrams in question is weakly contractible. \eop

{\em Remark:}
If $T$ is an $n$-truncated flexible functor on a quasi-compact site, then
$FT= F_{\infty}T$.
The above theorem implies Theorem \ref{desc15}. Alternatively, note that
the same proof works for Theorem \ref{desc15}.

\newpage

\setcounter{section}{6}

\numero{An analogue of Whitehead's theorem}

In this section, we discuss some properties of the $\Delta$-categories of
flexible functors and flexible sheaves, and the relationship with what is
known in Illusie's thesis and the work of Jardine.

The following notion is now called ``Segal category'' in \cite{effective} and
\cite{descente} for example. This notion (which is of course implicit in
the work of Segal) first appeared explicitly in the paper of Dwyer-Kan-Smith
\cite{DwyerKanSmith}. They prove that the homotopy categories of Segal
categories and simplicial categories (both localized by inverting the
``equivalences'') are equivalent.

We keep here the original notation of {\tt v1}.

A {\em $\Delta$-category} is a functor $M$ from the category $\Delta$ of
standard simplices, to the category of topological spaces, such that the
space $M_0$ is discrete, and
for any $n$ the morphism
$$
M_n \rightarrow M_1 \times _{M_0}\ldots \times _{M_0}M_1 $$
is a weak equivalence. This is modelled on Segal's machinery for loop
spaces \cite{Adams} \cite{SegalTopology}. We say that the set $M_0$ is the {\em
set of objects}, and for any objects $X_0,\ldots , X_k$ the fiber of $M_k$ over
the point $(X_0,\ldots , X_k)$ (via the map associating to a $k$-simplex the
$k+1$-tuple of its vertices), is the space of composable $k$-tuples of morphisms
between the spaces $X_0,\ldots , X_k$. The condition above means that the space
of composable $k$-tuples maps to the space of $k$-tuples of morphisms, by a
weak equivalence. We write $M_1(X_0,X_1)$ for this space of morphisms, in
case $k=1$.

A {\em $\Delta$-functor} between two $\Delta$-categories $M$ and $N$, is
just a morphism
of flexible functors from $M$ to $N$.

If $\Xx$ is a category, then we obtain a $\Delta$-category in the obvious
way (it is just the simplicial nerve of $\Xx$ with all sets given the
discrete topology). Similarly we can do this for a continuous category
$\Cc$.

{\em Note ({\tt v2}:} the following statement from {\tt v1} should be rephrased,
it is basically the statement we are conjecturing above; it should concern the
Segal category of flex functors $Flex(\Xx , Top)$...???

\begin{lemma}
The space of $\Delta$-functors from $\Xx$ to a $\Delta$-category $M$
associated to a continuous category $\Cc$ is ``the same'' (?) as the space
of flexible functors from $\Xx$ to $\Cc$.
\end{lemma}
{\em Proof:}
???
\eop

If $M$ is a $\Delta$-category then we can define the homotopy category $\pi
_0 M$ to be the category whose nerve is the simplicial set $\pi _0 (M)$.
The following theorem is a generalization of Vogt's theorem \cite{Vogt1}
\cite{CoPo2}: Vogt's theorem is the same statement without a Grothendieck
topologi (i.e. in the case where the site $\Xx$ has the coarse topology). The
present statement shows that the ``flexible sheaf'' condition is exactly what
is needed to make Vogt's theorem work for obtaining the localization of
simplicial presheaves by Illusie weak equivalence.

\begin{theorem}
\label{genvogt}
There is a natural equivalence of categories between $\pi _0 Flex(\Xx )$ (where
here $Flex(\Xx )$ is the $\Delta$-category of truncated weak flexible sheaves on
a site $\Xx$) and the full subcategory of Illusie's derived category of the
category of simplicial presheaves on
$\Xx$, consisting of the truncated ones.
\end{theorem}
{\em Proof:}
To each simplicial presheaf $P$ on $\Xx$ we can take the flexible sheaf
associated to the realization, $FRP$. This provides a functor from the
category of simplicial presheaves to $\pi _0 Flex(\Xx )$. We have to show
that it factors through the derived category---to show this it suffices to
show that it takes Illusie's weak equivalences to homotopy equivalences in
$Flex(\Xx )$. Suppose that $F:P\rightarrow Q$ is an Illusie weak
equivalence, that is it induces isomorphisms on the homotopy group sheaves.
This implies that the morphism between associated sheaves of realizations
satisfies the same property; so we may assume that $P\rightarrow Q$ is a
morphism of functors which are weak sheaves and of the form weak sheaf
associated to the realization of simplicial presheaves. We have to show that it
is a homotopy equivalence. Suppose $q\in Q_X$, and let $F$ be the fiber of the
morphism $P|_{\Xx /X} \rightarrow Q|_{\Xx /X}$ over the point $q$. We will show
that $F$ is weakly contractible (that is, its value on every object is weakly
contractible). First of all, the fact that the morphism of sheaves $\pi _0
(P)\rightarrow \pi _0 (Q)$ is surjective means that there is a sieve $\Bb$
over $X$ such that for all $U\in \Bb$, $F_U$ is nonempty. If $p\in F_U$,
and if $F|_{\Xx /U}$ is $n$-truncated, then the homotopy group presheaf
$\pi _n(F|_{\Xx /U},p)$ is already a sheaf (by an argument similar to that of
\ref{fibr08}). By hypothesis, this sheaf is zero (from the long exact homotopy
sequence and the isomorphism of homotopy sheaves of $P$ and $Q$ ). Thus $F|_{\Xx
/U}$ is $n-1$-truncated. By induction, it is $0$-truncated. Hence $F_U$ is
weakly
contractible. This being true for all elements of the sieve $\Bb$, we can
easily construct a point in $\Gamma (\Bb , F|_{\Bb})\sim F_X$. The same
argument then gives $F_X$ weakly contractible. Thus $P\rightarrow Q$ is a
weak equivalence on each fiber. But the realization of a simplicial
presheaf is of CW-type, and the operation of taking the associated weak
sheaf preserves CW type, so $P$ and $Q$ are of CW type. Thus the morphism
is objectwise a homotopy equivalence. We have seen (Theorem \ref{morph02}) that
this implies that it extends to an equivalence of flexible functors, that is it
is invertible in $\pi _0 Flex (\Xx )$.

We get a factorization to a functor from Illusie's derived category to $\pi
_0 Flex (\Xx )$. To see that it is essentially surjective, note that any
flexible functor can be replaced by an equivalent functor. We will see
below that there exist functors $T$ (the CW-free ones) such that for any
other functor $S$, the space of morphisms of functors from $T$ to $S$ is
equal to the space of morphisms of flexible functors. Furthermore, any
flexible functor is equivalent to a CW-free one. From this, if
$P\rightarrow Q$ is a morphism of flexible functors, we can choose a
CW-free functor $T$ with an equivalence of flexible functors $T\rightarrow
P$; if $P$ was a functor then this can be chosen to be a morphism of
functors, and it is then a quasi-isomorphism in Illusie's sense (say on the
singular simplicial presheaves associated). The composed morphism
$T\rightarrow Q$ can, if $Q$ was a functor, be replaced by a homotopic
morphism of functors, giving a morphism of singular simplicial presheaves.
The diagram $P\leftarrow T\rightarrow Q$ is exactly a morphism in Illusie's
derived category, and this maps to the class of our original morphism in
$\pi _0Flex (\Xx )$. Finally we have to show injectivity. It suffices to
consider the case of two morphisms $f_1,f_2:P\rightarrow Q$ where $P$ is
CW-free. They are equal in the derived category if there is a
quasiisomorphism $g:T\rightarrow P$ such that $f_1g=f_2g$. Call this map
$h$. There is an inverse $g':P\rightarrow T$, a priori a morphism of
flexible functors but since $P$ is CW free we may replace it by a morphism
of functors. Then $f_1$ and $f_2$ are both homotopic (in the space of
morphisms of flexible functors) to the morphism $hg'$. By the result about
morphisms from CW-free functors to be given below, this means that they are
homotopic as morphisms of functors. This shows that $f_1=f_2$ in the
homotopy category.

One last remark:  the localization of the category of simplicial presheaves by
Illusie weak equivalence, is the same as the localization of the homotopy
category of simplicial presheaves by weak equivalence. This can be seen from
\cite{Quillen} using the existence of  a closed model structure
\cite{Jardine1}.  (Admittedly it would be better to have a proof here which
doesn't rely on the existence of a closed model structure...???)
\eop

\subnumero{CW-free functors}

The ``CW-free'' functors we discuss below are the elementary cofibrations
in the closed model structure on the category of simplicial diagrams defined by
Bousfield-Kan \cite{BousfieldKan} (using an argument due to Quillen
\cite{Quillen}) and later taken up by Hirschhorn \cite{Hirschhorn}. In view of
these origins, in \cite{descente} we call these the elementary cofibrations for
the ``HBKQ model category structure''.

Recall that we say that a functor $T:\Xx \rightarrow Top$ is obtained from
$R:\Xx \rightarrow Top$ by a {\em free addition of an $n$-cell} if there is
an object $Y\in \Xx$ and a decomposition $$
T_X= R_X \cup _{\sim} \bigcup _{\phi :X\rightarrow Y}B^n_{\phi} $$
where the equivalence relation comes from an attaching map $\alpha
:S^{n-1}\rightarrow R_Y$ and where $T(\psi ): B^n_{\phi}\cong B^n_{ \phi
\psi}$ is the identity on the ball. The attaching map for $B^n_{\phi}$ is
$R(\phi )\alpha :S^{n-1}\rightarrow R_X$.

We say that a functor $T$ is CW-free if there is an increasing union of
closed subfunctors $T_i$ indexed by a well ordered set $I$, such that for
any $i$, $T_i$ is obtained from $T_{<i}$ by free addition of an $n_i$-cell
for some $n_i$. We may proceed by skeleta, and assume that the $n_i$ is an
increasing function of $i$ if we wish.

\begin{lemma}
Suppose $T$ is a CW-free functor from a category $\Xx$ to $Top$. Suppose
$S$ is a functor. Then the space of morphisms of functors $Mor ^{\rm
strict}(T,S)$ maps to the space of morphisms of flexible functors $Mor
(T,S)$ by a weak equivalence.
\end{lemma}
{\em Proof:}
We may proceed by induction on the sequence of subspaces used to define
$T$. It suffices to show that the map on relative morphism spaces $$
Mor ^{\rm strict} (B^{n+1}_{\Xx /Y},S^{n}_{\Xx /Y},S; f) \rightarrow Mor
(B^{n+1}_{\Xx /Y},S^{n}_{\Xx /Y},S; f) $$
is a weak equivalence. We know in both cases that this is the fiber of the
morphism space from $T_i$ to $S$, over an element in the morphism space
from $T_{<i}$ to $S$ (and that the map of restricting morphisms is a
fibration, and so on). The first question is, if the second space is
nonempty then is the first space nonempty? This can be answered by asking a
question about whether two points are in the same connected component in
the space of morphisms from $S^n_{\Xx /Y}$ to $S$ sending the basepoint
(north pole, say) to a point $\eta$ in $S_Y$. Given that both spaces are
nonempty, we can complete the map from the ball to a map from the sphere by
adding a given map on the other hemisphere. The space of morphisms on the
ball relative the boundary is homotopic to the space of morphisms from the
sphere relative the base point. We are reduced to showing that the map of
morphism spaces for pointed morphisms $$
Mor ^{\rm strict} ((S^{n}_{\Xx /Y},o),S; \eta ) \rightarrow Mor
((S^{n}_{\Xx /Y},o),S; \eta )
$$
is a weak equivalence. But this is equivalent to the map of loop spaces $$
Mor ^{\rm strict} (\ast _{\Xx /Y}, (\Omega ^{\eta})^nS) \rightarrow Mor
(\ast _{\Xx /Y}, (\Omega ^{\eta})^nS ). $$
These are (exactly)
equal to the morphism spaces from $\ast$ to the loop space, over the
category $\Xx /Y$.
The first is clearly equal to $(\Omega ^{\eta})^nS_Y$, while the second is
equal to the space of sections $\Gamma (\Xx /Y, (\Omega ^{\eta})^nS |_{\Xx
/Y})$. By Lemma \ref{desc01} and the discussion in \S 3, the morphism from the
space over the object $Y$ to the space of sections over $\Xx /Y$ is a weak
equivalence. This completes the proof of the lemma. \eop

Note that this lemma takes care of what we needed for the previous proof.
It gives, as a corollary, an analogue of Whitehead's theorem in our case.

\begin{theorem}
Suppose $T$ is a truncated CW-free functor, and suppose $S$ is a truncated
CW-free
functor which is a weak flexible sheaf. Suppose that $f:S\rightarrow T$ is
a morphism of functors which induces isomorphisms on all homotopy group
sheaves. Then there is a morphism of functors
$g:T\rightarrow S$ such that $gf$ is homotopic to the identity. If, furthermore,
$T$ is a weak flexible sheaf and $S$ is of CW type then the other composition
$fg$ is homotopic to the identity. \end{theorem} {\em Proof:} The composed
morphism $\tau_{\infty}S\rightarrow F_{\infty}T$ to the associated sheaf is a
weak equivalence in the sense of Illusie, between sheaves of CW type. By the
same argument as above, it is a weak equivalence on each point, so it is a
homotopy equivalence on each point, and so it extends to an equivalence of
flexible sheaves. The resulting inverse map from $T$ to $S$ can be replaced by a
map of functors $g$, by the previous lemma, and the composition $gf$ is
homotopic to the identity as a map of flexible functors. Applying the previous
lemma to the space of maps from $S$ to itself, we get that $gf$ is homotopic to
the identity as a map of functors. If $T$ is already a flexible sheaf then there
is no need to look at $F_{\infty}T$, and we get also that $fg$ is homotopic to
the identity (first as a map of flexible functors then as a map of functors,
again by the previous lemma). \eop

\newpage

\setcounter{section}{7}

\numero{Mapping sheaves}

Suppose $S$ and $T$ are contravariant flexible functors from $\Xx$ to a
continuous category $\Cc$. We define the functor $Map (S,T):\Xx \rightarrow
Top$ by
$$
Map (S,T)(X):= Mor _{\Xx /X}(S|_{\Xx /X}, T|_{\Xx /X}). $$
This is functorial in $X$: if $f:Y\rightarrow X$ then we have a functor
$\Xx /Y\rightarrow \Xx /X$ and the pullback for $\Xx /Y\times I \rightarrow
\Xx /X\times I$ gives a map
$$
f^{\ast}: Mor _{\Xx /X}(S|_{\Xx /X}, T|_{\Xx /X})\rightarrow Mor _{\Xx
/Y}(S|_{\Xx /Y}, T|_{\Xx /Y}), $$
satisfying the strict associativity $f^{\ast} g^{\ast} = (gf)^{\ast}$.

For example, $Map (\ast _{\Xx}, T)$ is the functor we called $KT$ in
previous sections.

Suppose $A\subset B$ is an inclusion of contravariant functors from $\Xx$
to $Top$, and suppose $T$ is a flexible functor from $\Xx$ to $Top$.
Suppose $G\in Mor _{\Xx /X} (B|_{\Xx /X},T|_{\Xx /X})$. We can define the
{\em restriction of $G$ to $A$} as the composition of $G$ with the
inclusion $A|_{\Xx /X}\rightarrow B|_{\Xx /X}$. Note that this composition
can be defined canonically, since the inclusion is a morphism of functors.
If $F\in Mor _{\Xx }(A,T)$, define the relative mapping space
$$
Map (B, T;
F)(X):= \{ G\in Mor _{\Xx /X }(B, T) \; \mbox{s.t.} \; G|_A = F\} . $$

Suppose $B$ is a contravariant functor from $\Xx$ to $Top$. Define the {\em
suspension $\Sigma B$} to be the contravariant functor assigning to $X\in
\Xx$ the suspension of $B_X$. If $A\subset B$ is an inclusion, then we get
an inclusion $\Sigma A \subset \Sigma B$.

If $\eta \in \Gamma (\Xx , T)$ and $S$ is any contravariant functor, denote
by $\eta _S$ the composition of $\eta \in Mor (\ast , T)$ with
$S\rightarrow \ast$.

\begin{lemma}
\mylabel{mapp01}
If $A\subset B$ is an inclusion of functors and $\eta \in \Gamma (\Xx , T)$
then $$
Map (\Sigma B , T; \eta _{\Sigma A}) \sim Map (S , \Omega ^{\eta} T;
c^{\eta}_A) $$
where $c^{\eta}\in \Gamma (\Xx , \Omega ^{\eta} T)$ is the section
corresponding to the constant path.
\end{lemma}
{\em Proof:}
???
(we need a more concrete definition of the loop space ... ???) \eop

Let $\Yy = \Xx/X$ with its functor $\rho$ to $\Xx$. For any $n$, define a
contravariant functor $S^n_{\Yy}$ on $\Xx$ which assigns to each object
$U\in \Xx$ a disjoint union of $n$-spheres with one for each object of
$\Yy$ over $U$. In effect, there is one sphere for each morphism
$U\rightarrow X$. Note that the $0$-sphere $S^0$ is a disjoint union of two
points. Let $A^n_{\Yy} \subset S^n_{\Yy}$ denote the union of eastern
hemispheres in the spheres. (These are homotopic for any $n$.)

\begin{corollary}
\mylabel{mapp02}
If $\eta \in \Gamma (\Xx , T)$ then
$$
Mor _{\Yy} (S^n_{\Yy} , T; \eta _{A^n_{\Yy}}) \sim \Gamma (\Yy , \rho
^{\ast} ((\Omega ^{\eta})^nT))).
$$
\end{corollary}
{\em Proof:}
Note first of all that
$$
Map (S^0_{\Yy}, T; \eta _{A^n_{\Yy}}) = Map _{\Yy}(\ast , \rho ^{\ast} T)
\sim \rho ^{\ast} T.
$$
Applying the lemma $n$ times (and noting that $(S^n_{\Yy}, A^n_{\Yy})=
\Sigma ^n (S^0_{\Yy}, A^0_{\Yy})$) we get
$$
Map (S^n_{\Yy}, T; \eta _{A^n_{\Yy}})\sim Map (S^0_{\Yy}, (\Omega
^{\eta})^nT; \eta _{A^0_{\Yy}}).
$$
Taking global sections gives the desired statement. \eop

Now suppose that the pair $A\subset B$ is equivalent by excision to the
pair $A^n_{\Yy} \subset S^n_{\Yy}$ for some $\Yy = \Xx /Y$. Then we get $$
Mor (B, T; \eta _A) \sim \Gamma (\Yy , \rho ^{\ast}(\Omega ^{\eta})^nT). $$

\begin{proposition}
\mylabel{mapp03}
Suppose that $A\subset B$ is also equivalent by excision to the pair
$S^{n-1}_{\Yy} \subset D^n_{\Yy}$ where $D^n_{\Yy}$ denotes the union of
one $n$-disk
for each object of $\Yy$. Suppose that $T$ is a flexible sheaf; suppose
$X\in \Xx$ and $F\in Mor _{\Xx /X}(A, T)$; then for a sieve $\Bb$ in $\Xx
/X$, the map
$$
Mor _{\Xx /X}(B|_{\Xx /X},T|_{\Xx /X}; F|_{\Xx /X})\rightarrow Mor
_{\Bb}(B|_{\Bb}, T|_{\Bb}; F|_{\Bb}) $$
is a weak equivalence.
\end{proposition}
{\em Proof:}
Note that the restriction $D^n_{\Yy}|_{\Xx /X}$ is the same as $D^n_{\Yy
'}$ where $\Yy ' = \Xx /X\times Y \rightarrow \Xx /X$ (and the same for the
$S^{n-1}_{\cdot}$). Thus we may replace $\Xx $ by $\Xx /X$, so it suffices
to give the proof in the case where $X$ is a final object in $\Xx$. With
this reduction, our map is equivalent by excision to $$
Mor _{\Xx} (D^n_{\Yy}, T ; F|_{S^{n-1}_{\Yy}})\rightarrow Mor _{\Bb}
(D^n_{\Yy}|_{\Bb}, T |_{\Bb}; F|_{S^{n-1}_{\Yy}}|_{\Bb}). $$
We claim that this is equivalent to
$$
Mor _{\Yy} (D^n , \rho ^{\ast}T; F|_{S^{n-1}_{\Yy}}) \rightarrow Mor _{\Bb
'} (D^n, T |_{\Bb '}; F|_{S^{n-1}_{\Yy}}|_{\Bb '}), $$
where $\Bb '$ is the sieve over $Y$ induced by $\Bb$. This can be seen by
following precisely the definition of morphism: a morphism consists of the
data, for each composable sequence in $\Xx \times I$, going from $(X_n,1)$
to $(X_0,0)$, of a map from $(D^n_{\Yy})_{X_0}$ to $T_{X_n}$; but since
$(D^n_{\Yy})_{X_0}$ is the union of one disk for each morphism
$X_0\rightarrow Y$, we can think of this as a morphism from $D^n$ to
$T_{X_n}$ for each composable sequence in $\Xx /Y$; note that the $T_{X_n}$
can then be interpreted as $(\rho ^{ast}T)_{X_n}$.

Denote $F|_{S^{n-1}_{\Yy}}$ by $F'$.
Let $\eta \in \Gamma (\Yy ,
\rho ^{\ast}T)$ be the image of a basepoint $\sigma \in S^{n-1}$ by $F'$.
If the right hand side of our morphism is empty then there is nothing to
prove. Thus we may assume that this is not the case, and choose $$
G\in Mor _{\Bb} (D^n, T |_{\Bb '}; F'|_{\Bb '}). $$
The map
$$
Mor _{\Yy}(S^{n-1} , \rho ^{\ast}T ; \eta _{\sigma}) \rightarrow Mor _{\Bb
'}(S^{n-1} , \rho ^{\ast}T|_{\Bb '}; \eta _{\sigma }|_{\Bb '}) $$
is equivalent, via Corollary \ref{mapp02}, to
$$
\Gamma (\Yy , (\Omega ^{\eta})^{n-1}(\rho ^{\ast}T))\rightarrow \Gamma (\Bb ' ,
(\Omega ^{\eta})^{n-1}(\rho ^{\ast}T)|_{\Bb '}), $$ which is a weak equivalence
because $\rho ^{\ast}T$ is a sheaf on the site $\Yy$ (it follows immediately
from the definition of sheaf on $\Xx$ that the restriction to
$\Xx /Y$ is also a sheaf), so the $n$th loop space is a sheaf (Corollary
\ref{fibr03}). But $G$ provides a homotopy between the element $$
F'|_{\Bb '}\in
Mor _{\Bb '}(S^{n-1} , \rho ^{\ast}T|_{\Bb '}; \eta _{\sigma }|_{\Bb '}) $$
and $\eta _{S^{n-1}_{\Bb}}$. Therefore there exists a homotopy $G'$ between
$F'$ and $\eta _{S^{n-1}}$ over $\Yy$. More precisely, there exists an
element $$
G'\in Mor _{\Yy }(D^n, \rho ^{\ast }T; \eta _{\sigma}) $$
mapping to
$$
F'\in Mor _{\Yy }(S^{n-1} , \rho ^{\ast} T; \eta _{\sigma}), $$
or equivalently an element
$$
G' \in Mor _{\Yy} (D^n, \rho ^{\ast}T; F'). $$
Now by excision, we have the equivalence $$
Mor _{\Yy}(D^n, \rho ^{\ast}T; F')\sim Mor _{\Yy}(S^n, \rho ^{\ast}T;
G'_{A^n}) $$
where $G'_{A^n}$ is the element $G'$ considered as a map on the hemisphere
$A^n\subset S^n$. (Note that when we write $D^n$, $S^n$, $A^n$ without
subscripts we mean the constant functors on the category $\Yy$.) The same
holds over the sieve $\Bb '$, and the restriction to $\Bb$ for the left
side of the previous equation is equivalent to the restriction for the
right side. Thus it suffices to show that $$
Mor _{\Yy}(S^n, \rho ^{\ast}T; G'_{A^n})\rightarrow Mor _{\Bb '}(S^n, \rho
^{\ast}T|_{\Bb '}; G'_{A^n}|_{\Bb '}) $$
is a weak equivalence. By Corollary \ref{mapp02}, this is equivalent to $$
\Gamma (\Yy , (\Omega ^{\eta})^n(\rho ^{\ast}T))\rightarrow \Gamma (\Bb ' ,
(\Omega ^{\eta})^n(\rho ^{\ast}T)|_{\Bb '}), $$
which is a weak equivalence since $(\Omega ^{\eta})^n(\rho ^{\ast}T)$ is a
sheaf on $\Yy$. This completes the proof of the proposition. \eop

Suppose $B:\Xx \rightarrow Top$ is a contravariant functor. We say that $B$
is {\em obtained by a sequence of free additions of cells} if there is a
filtration of $B$ by subspaces $A_j$ indexed by a well-ordered set $J$,
such that (with the usual notations) for any $j$, the pair $(A_j , A_{<j})$
is equivalent by excision to some $(D^n_{\Yy}, S^{n-1}_{\Yy})$ for some
$\Yy = \Xx /Y$; and such that the first $A_{0}$ is empty. Note that we can
assume that there is a final object $e$ in $J$ with $B=A_e$ (at this step,
we might not add anything to $A_{<e}$).

\begin{theorem}
\mylabel{mapp04}
If $S$ is a contravariant flexible functor of CW-type on $\Xx$ then $S$ is
homotopy equivalent to a functor $B$ obtained by a sequence of free
additions of cells.
\end{theorem}
{\em Proof:}
We may assume $S$ is a functor. We start with $A_0$ empty and add the cells
by induction on the dimension $n$, so as to get a map $B\rightarrow S$
inducing an isomorphism on homotopy group presheaves. In each dimension,
proceed in two steps. First add cells so as to kill off the kernel of the
map in the homotopy groups of degree $n-1$; then add cells so as to get a
surjection in the homotopy groups of degree $n$. We
add a cell corresponding to category $\Yy = \Xx /Y$ every time we need to
kill an element or add a new element in the homotopy over $Y$. This gives
cells for all objects of $\Yy$, but if there are too many they will be
taken care of at the next step. Each step does not disturb the previous
step. Finally, note that we obtain a map between functors of CW type
inducing isomorphisms on homotopy, so it is a homotopy equivalence on each
object, hence a homotopy equivalence.
\eop

\begin{theorem}
\mylabel{mapp05}
Suppose $S$ is a contravariant flexible functor of CW-type on a site $\Xx$,
and suppose $T$ is a flexible sheaf. Then $Map (S,T)$ is a weak flexible
sheaf. \end{theorem}
{\em Proof:}
By the previous theorem, we may suppose $S=B$ is obtained by a sequence of
free additions of cells. Let $A_j$ be the corresponding sequence of
subfunctors. Suppose $X\in \Xx$ and $\Bb \subset \Xx /X$ is a sieve. We
prove by transfinite induction that the maps $$
Mor _{\Xx /X}(A_{j}|_{\Xx /X}, T|_{\Xx /X}) \rightarrow Mor
_{\Bb}(A_{j}|_{\Bb}, T|_{\Bb})
$$
are weak equivalences.
Let $j$ be the first element where this is not the case. We have the diagram
$$
\begin{array}{ccc}
Mor _{\Xx /X}(A_{<j}|_{\Xx /X}, T|_{\Xx /X}) & \rightarrow & Mor
_{\Bb}(A_{<j}|_{\Bb}, T|_{\Bb}) \\
\parallel & & \parallel \\
\lim _{\leftarrow} Mor _{\Xx /X}(A_{j'}|_{\Xx /X}, T|_{\Xx /X}) &
\rightarrow & \lim _{\leftarrow}Mor _{\Bb}(A_{j'}|_{\Bb}, T|_{\Bb})
\end{array}
$$
where the inverse limits are over $j'<j$. These are inverse limits of a
sequence of fibrations. By the inductive hypothesis, the map on the bottom
is a map of inverse limits induced by an inverse system of weak
equivalences, so it is a weak equivalence. On the other hand, we have the
diagram $$
\begin{array}{ccc}
Mor _{\Xx /X}(A_{j}|_{\Xx /X}, T|_{\Xx /X}) & \rightarrow & Mor
_{\Bb}(A_{j}|_{\Bb}, T|_{\Bb}) \\
\downarrow & & \downarrow \\
Mor _{\Xx /X}(A_{<j}|_{\Xx /X}, T|_{\Xx /X}) & \rightarrow & Mor
_{\Bb}(A_{<j}|_{\Bb}, T|_{\Bb})
\end{array}
$$
where the vertical maps are fibrations. Furthermore, by Proposition \ref{mapp03}
applied to the pair $(A_j , A_{<j})$, the map between the fibers of the
vertical maps is a weak equivalence. As we have seen, the map on the bottom
is a weak equivalence, therefore the map at the top is a weak equivalence,
contradicting our supposition that the inductive hypothesis wasn't true for
some $j$. Therefore it is true for all $j$, including the final object $e$.
We obtain the statement that
$$
Mor _{\Xx /X}(S|_{\Xx /X}, T|_{\Xx /X}) \rightarrow Mor _{\Bb}(S|_{\Bb},
T|_{\Bb})
$$
is a weak equivalence. To complete the proof of the proposition, note that $$
Mor _{\Xx /X}(S|_{\Xx /X}, T|_{\Xx /X})=Map (S,T)(X)\sim \Gamma (\Xx /X ,
Map (S,T)|_{\Xx /X}),
$$
so it
suffices to show that
$$
Mor _{\Bb}(S|_{\Bb}, T|_{\Bb})\sim _{\rm w.e.} \Gamma (\Bb , Map
(S,T)|_{\Bb}). $$
This follows from the following lemma (noting that $Map (S,T)|_{\Bb}= Map
(S|_{\Bb},T|_{\Bb})$). Note that it suffices to treat the case of sieves
given by covering families.
\eop

\begin{lemma}
\mylabel{mapp06}
Suppose $S$ and $T$ are flexible functors on a category $\Xx$ with final
object, with $S$ of CW-type. Suppose that $\Zz$ is a sieve over the final
object in $\Xx$, given by a covering family $\Uu =\{ U_{\alpha}\} _{\alpha
\in J}$. Then the natural morphism $$
Mor _{\Zz}(S,T)\rightarrow \Gamma (\Zz , Map (S,T)) $$
is a weak homotopy equivalence.
\end{lemma}
{\em Proof:}
We may assume that $S$ is a functor. The homotopy groups of both sides may
be interpreted as homotopy classes in the same functors applied to the
pairs $S\times S^n, S\times A^n)$. To prove
injectivity, it suffices to prove surjectivity for the pairs $S\times D^n,
S\times S^{n-1})$. Thus it suffices to prove that for any pair $A\subset B$
of functors, and any $F\in Mor _{\Zz}(A,T)$, the map $$
Mor _{\Zz}(B,T; F)\rightarrow \Gamma (\Zz , Map (B,T; F)) $$
is surjective on $\pi _0$. Since our original $S$ was of CW-type, we may
write $B$ as a transfinite union of cells (starting with $A$); thus it
suffices to prove the above surjectivity in the case where $(B,A)$ is
equivalent under excision to $(D^n_{\Yy},S^{n-1}_{\Yy})$, for a $\Yy =
\Zz/Y$ with $Y\in \Zz$. But the surjectivity here is equivalent to
injectivity in the case of the pair $((S^{n})_{\Yy},A^{n}_{\Yy})$ (here we
have replaced $n-1$ by $n$ for convenience in what follows). Then as
discussed above, $$
Mor _{\Zz}(B,T; F)\sim Mor _{\Yy} (S^n,T; \eta _{A^n}) $$
where $\eta \in \Gamma (\Yy , T|_{\Yy})$ is the point homotopic to the
restriction of $F$ to the $A^n$ over $\Yy$. Similarly, for any $X\in \Zz$
we have
$$
Map (B,T;F)(X)= Mor _{\Zz /X}(B|_{\Zz /X},T|_{\Zz /X}; F|_{\Zz /X})
$$
$$
\sim Mor
_{\Zz /X\times Y}(S^n, T|_{\Zz /X\times Y};\eta _{A^n}|_{\Zz /X\times Y}).
$$
On the other hand, by Corollary \ref{mapp02},
$$
Mor _{\Yy} (S^n,T; \eta _{A^n})\sim [(\Omega ^{\eta})^n(T|_{\Yy})]_Y =
(\Omega ^{\eta _Y})^n(T_Y),
$$
and
$$
Mor _{\Zz /X\times Y}(S^n, T|_{\Zz /X\times Y};\eta _{A^n}|_{\Zz /X\times
Y}) \sim [(\Omega ^{\eta})^n(T|_{\Yy})]_{X\times Y} = (\Omega ^{\eta
_{X\times Y}})^n(T_{X\times Y}). $$
Let $W$ be the flexible functor on $\Zz$ defined by $W_X= T_{X\times Y}$.
Then $\eta_Y\in T_Y$ corresponds to a point $\eta '\in \Gamma (\Xx , W)$,
and our question is to show that $$
(\Omega ^{\eta_Y})^n(T_Y) \rightarrow \Gamma (\Zz , (\Omega ^{\eta '})^nW) $$
is an equivalence.

We claim that the map $T_Y\rightarrow \Gamma (\Zz , W)$ is an equivalence.
From this, the desired equivalence will follow because taking the loop
space commutes with taking global sections.

Here is where we use the hypothesis that $\Zz$ is the sieve associated to a
covering family $\Uu= \{ U_{\alpha }\} _{\alpha \in J}$. Let $\Jj$ denote
the free category with products on the index set $J$. Let $\xi _{\Uu}: \Jj
\rightarrow \Zz$ be the morphism determined by the covering $\Uu$. Recall
from Lemma \ref{desc09a} (applied to the sieve $\Zz$ over the final object
$e$ of
$\Xx$, determined by the covering $\Uu$), that the pullback map $$
\xi _{\Uu}^{\ast}: \Gamma (\Zz ,
W)\rightarrow \Gamma (\Jj , \xi _{\Uu }^{\ast}W) $$
is a homotopy equivalence (note that $W$ comes from a flexible functor on
$\Xx$). On the other hand, our covering family $\Uu$ gives a covering
family $\Uu ' = \{ U_{\alpha }\times _eY\}_{\alpha \in J}$ of $Y$. Since
the index set is still the same, we obtain a functor
$$
\xi _{\Uu '}:\Jj \rightarrow \Xx /Y ,
$$
with the property that the pullback
$$
\xi ^{\ast}_{\Uu '} : \Gamma (\Xx /Y,
T|_{\Xx /Y} )\rightarrow \Gamma (\Jj , \xi _{\Uu '}^{\ast}T) $$
is a homotopy equivalence---again by Lemma \ref{desc09a} applied this time
to $T$.
Note here that the sieve over $Y$ generated by the covering family $\Uu '$
is equal to $\Xx /Y$ since $Y$ is a member of the sieve $\Zz$ generated by
$\Uu$. Recall that
$$
\Gamma (\Xx /Y,
T|_{\Xx /Y} )= T_Y .
$$
On the other hand, for an object $B_{\alpha _1,\ldots , \alpha _k}$ of $\Jj$,
$$
(\xi _{\Uu}^{\ast}W)_{B_{\alpha _1,\ldots , \alpha _k}}= W_{U_{\alpha
_1}\times _e\cdots \times _e U_{\alpha _k}}, $$
whereas
$$
(\xi _{\Uu '}^{\ast}T)_{B_{\alpha _1,\ldots , \alpha _k}}= T_{(U_{\alpha
_1}\times _eY)\times _Y \cdots \times _Y(U_{\alpha _k}\times _eY)}.
$$
Now
$$
(U_{\alpha _1}\times _eY)\times _Y \cdots \times _Y(U_{\alpha _k}\times _eY)=
(U_{\alpha _1}\times _e\cdots \times _e U_{\alpha _k})\times _eY, $$
and in view of the definition of $W$ this means that $$
\xi _{\Uu}^{\ast}W=\xi _{\Uu '}^{\ast}T . $$
From our two applications of Lemma \ref{desc09a}, we obtain the homotopy
equivalence
$$ \Gamma (\Zz ,
W)\sim \Gamma (\Xx /Y,
T|_{\Xx /Y} )\sim T_Y.
$$
This equivalence is homotopic to the map in question (??? verify ???), so
we get the claim and hence the lemma. \eop

{\em Remark:}
For the lemma, it suffices that $\Zz$ admits products and fibered
products---then construct $\Xx$ by adding a final object. Take for covering
family the set of all objects of $\Zz$.

\subnumero{Relative mapping functors}
Suppose $\Xx$ admits products. Then
it satisfies the hypotheses of
Lemma \ref{mapp06} above, so $Mor (S,T)=\Gamma (\Xx , Map (S,T))$.

Suppose $R\rightarrow T$ and
$S\rightarrow T$ are morphisms of flexible functors on a category $\Xx$. Let $$
Mor _T(R,S)
$$
denote the space of diagrams
$$
\begin{array}{ccc}
R &\rightarrow & S \\
&\searrow & \downarrow \\
&&T\; .
\end{array}
$$
Let $Sect (T,R):= Mor _T(T,R)$.
Note that $Sect (S, S\times _TR)\sim Mor _T(S, R)$. If $V\rightarrow T$ is
a morphism then there is a natural map $$
Mor _T(R,S) \rightarrow Mor _V(R\times _TV, S\times _VT). $$
(Properties ???)

The following theorem establishes the existence of {\em relative mapping
spaces.}
\begin{theorem}
There exists a flexible functor over $T$, $$
U\rightarrow T,
$$
together with an element
$$
\xi \in Mor _U(R\times _TU, S\times _TU) $$
such that for any morphism $V\rightarrow T$, then the map induced by $\xi$
$$
Mor _T(V, U )\rightarrow Mor _V(R\times _TV, S\times _TV) $$
is a weak equivalence.
Consequently, the pair $(U,\xi )$ is unique up to very well defined equivalence.
We denote this pair (or the underlying space $U$) by $Mor (R/T,S/T)$.
\end{theorem}
{\em Proof:}
???
\eop

We can also define a space of equivalences over $T$, denoted $Iso _T(R,S)$,
and an associated relative space of equivalence $$
Equiv (R/T,S/T)
$$
with a similar universal property.

\begin{corollary}
The space $Mor (R/T, S/T)$ is preserved under restriction of the underlying
category. If $\Xx$ admits products, then $Mor (R/\ast , S/\ast )$ is
naturally equivalent to $Map (R,S)$. If
Thus if $\Xx$ admits fiber products and if $t\in T_X$, the fiber of $Mor
(R/T,S/T)|_{\Xx /X}$ over $t$ is equivalent to $Map (R_t, S_t)$ where $R_t$
and $S_t$ are the fibers of the restrictions of $R$ and $S$ to $\Xx /X$,
over $t$.
\end{corollary}
{\em Proof:}
???
\eop

\begin{corollary}
If $\Xx$ is a site, and if $R$, $S$ and $T$ are truncated weak flexible
sheaves, then
$Map (R/T, S/T)$ is a weak flexible sheaf. \end{corollary}
{\em Proof:}
Denote $Map (R/T, S/T)$ again by $U$. We have a map $U\rightarrow T$ and we
know that $T$ is a truncated
weak flexible sheaf and $U_t$ is a truncated weak flexible sheaf for any
$t\in T_X$ (since, by the previous corollary, $U_t\sim Map (R_t, S_t)$ and
we have shown this to be a weak sheaf in Theorem \ref{mapp05}). The degrees of
truncations are uniform in $X$ and $t$. This implies, first of all, that $U$ is
truncated (by the long exact sequence of homotopy).
We claim that this implies that $U$ is a weak sheaf. For this we may assume
that $U$ and $T$ are of CW-type.
There is a map $U\rightarrow H_{\infty}U$ over $T$ (since $T\rightarrow
H_{\infty}T$ is an
equivalence). Furthermore, since the construction $H_{\infty}$ commutes
with taking
fiber products, the induced morphisms on fibers are weak equivalences. This
implies (from the long exact homotopy sequence and the five
lemma---checking the low degree cases by hand) that $U\rightarrow HU$ is a
weak equivalence. Thus $U$ is a weak sheaf.
\eop

\subnumero{Classifying spaces}
Suppose $T$ is a flexible functor on a site $\Xx$. We say that a flexible
sheaf $U$ on $\Xx /X$ is {\em locally equivalent to $T$} if there exists a
sieve $\Bb \subset \Xx /X$ such that for each $Y\in \Bb$, the restriction
$U|_{\Xx /Y}$ is equivalent to $T|_{\Xx /Y}$. We say that a
morphism $R\rightarrow S$ {\em
has fibers locally equivalent to $T$} if for any object $X\in \Xx$ and any
section $\eta \in \Gamma (\Xx /X, S)$, the fiber $Fib (R/S, \eta )$ is
locally equivalent to $T$.

\begin{proposition}
\mylabel{mapp07}
Suppose $T$ is a flexible sheaf on $\Xx$, and suppose $R\rightarrow Z$ is a
morphism of flexible sheaves,
such that
$$
Equiv _Z(T\times Z ,R)\rightarrow Z
$$
is a weak equivalence. Then the fibers of $R\rightarrow Z$ are locally
equivalent to $T$. Furthermore, for any morphism $U\rightarrow S$ with
fibers locally equivalent to $T$, the space of diagrams $$
\begin{array}{ccc}
U & \rightarrow & R \\
\downarrow &&\downarrow \\
S & \rightarrow & Z
\end{array}{ccc}
$$
inducing an equivalence between $U$ and $R\times _Z S$, is weakly
contractible. \end{proposition}
{\em Proof:}
???
\eop

In the situation of the proposition, we say that $Z$ is a {\em classifying
space} for $T$, and that $R\rightarrow Z$ is the universal fibration. Note
that the morphism $T\rightarrow \ast$ corresponds to a basepoint $\zeta \in
\Gamma (\Xx , Z)$. The fiber of $Equiv _Z(T\times Z ,R)\rightarrow Z$ over
$\zeta $ is equivalent to $Equiv (T,T)$. Hence, we have $$
\Omega ^{\zeta} Z \sim Equiv (T,T ).
$$

{\em Remark:} There is a weakly defined multiplication on $End (T,T)$; this
coincides with composition of paths in the loop space. (????).

\subnumero{Universal characterizations and completions} Suppose $\Xx$ is a
site and ${\bf F}$ is a class of flexible functors on $\Xx$ closed under
the operations $T\mapsto Map (U,T)$ for finite CW-complexes $U$ (considered
as constant flexible functors on $\Xx$).

(NB In fact, we would need to do the following for classes which are only
closed under the operations of taking pointed maps of connected CW
complexes...???)

Suppose $\psi :R\rightarrow S$ is a morphism of flexible functors on $\Xx$.
We say that $\psi$ is {\em an equivalence relative to ${\bf F}$} if, for
any $T\in {\bf F}$, the map
$$
Mor (S,T)\rightarrow Mor (R,T)
$$
(obtained by composition with $\psi$) is a weak equivalence.
An equivalent condition is that for any $R\rightarrow T$ with $T\in {\bf F}$,
the space of  diagrams
$$
\begin{array}{ccc}
R&\rightarrow &S\\
&\searrow &\downarrow \\
&& T
\end{array}
$$
is weakly contractible.

Similarly, we
say that $\psi$ is {\em universally an equivalence relative to ${\bf F}$}
if, for any $T\in {\bf F}$, the map
$$
Map (S,T)\rightarrow Map (R,T)
$$
(obtained by composition with $\psi$) is a weak equivalence of flexible
functors. Note that this implies the first condition when $\Xx$ is a site,
since $Mor (S,T)= \Gamma (\Xx , Map (S,T))$ by Lemma \ref{mapp06}, and $\Gamma$
preserves weak equivalences (???).

If $R\rightarrow S$ is an equivalence relative to ${\bf F}$ (resp.
universally an equivalence relative to ${\bf F}$) and if $S\in {\bf F}$
then we say that $S$ {\em is the ${\bf F}$-completion of $R$} (resp. {\em
the universal ${\bf F}$-completion of $R$}).

\begin{lemma}
\mylabel{mapp08}
Suppose $\psi : R\rightarrow S$ is a morphism of flexible functors on
$\Xx$, with the property that for any morphism $T\rightarrow T'$ of
flexible functors in ${\bf F}$, and for any $\eta \in Mor (S,T')$, the map
from the fiber of
$$
Mor (S,T)\rightarrow
Mor (S,T')
$$
over $\eta$, to the fiber of
$$
Mor (R,T)\rightarrow
Mor (R,T')
$$
over the image of $\eta$, induces a surjection on the set of path
components. Then $\psi$ is universally a weak equivalence.

Suppose that for any morphism $T\rightarrow T'$ of flexible functors in
${\bf F}$, any $X\in \Xx$, and any $\eta \in Map (S,T')_X$, the map from
the fiber of
$$
Mor (S|_{\Xx /X},T|_{\Xx /X})\rightarrow Mor (S|_{\Xx /X},T'|_{\Xx /X})
$$
over $\eta$, to the fiber of
$$
Mor (R|_{\Xx /X},T|_{\Xx /X})\rightarrow Mor (R|_{\Xx /X},T'|_{\Xx /X})
$$
over the image of $\eta$, induces a surjection on the set of path
components. Then $\psi$ is universally a weak equivalence. \end{lemma}
{\em Proof:}
???
\eop

We can give an easier criterion if we don't need to worry about the
fundamental group.

\begin{lemma}
\mylabel{mapp09}
Suppose $R$ is a flexible functor such that for any $T\in {\bf F}$, the
fundamental group of every path component of $Mor (R,T)$ (resp. $Map
(S,T)_X$ for all $X\in \Xx$) is abelian and acts trivially on the higher
homotopy groups. Suppose $\psi : R\rightarrow S$ is a morphism of flexible
functors such that for every $T\in {\bf F}$, the morphism $$
Mor (S,T)\rightarrow Mor (R,T)
$$
induces an isomorphism on the set of path components (resp.
$$
Map (S,T)\rightarrow Map (R,T)
$$
induces an isomorphism on the homotopy presheaf $\pi _0^{\rm pre}$). Then
$\psi$ is an
equivalence relative to ${\bf F}$
(resp. universally an equivalence relative to ${\bf F}$).
\end{lemma}
{\em Proof:}
We prove the first case (??? for the universal case). Under the hypotheses
of the lemma, the union over all path components of the $n$th homotopy
groups of $Mor (R,T)$ are just $Mor (R, Map (S^n, T))$. As $Map (S^n, T)$
is also in ${\bf F}$, the principal hypothesis implies that the maps
induced by $Mor (S,T)\rightarrow Mor (R,T)$ on all homotopy groups of all
path components, are surjective. On the other hand, putting $T'= Map (S^1,
T)$ we obtain $Mor (S, T')$ equal to the set of conjugacy classes in the
fundamental groups of path components of $Mor (S,T)$. We obtain the result
that the maps on fundamental groups are injective on the sets of conjugacy
classes, but this implies that the kernels of the morphisms consist of one
conjugacy class, necessarily that of the identity; thus the kernels of the
morphisms on fundamental groups are trivial, that is the morphisms on
fundamental groups are injective and hence isomorphisms. Similarly, the
morphisms on higher homotopy groups of path components are injective on the
spaces of orbits under the fundamental groups. Again, this implies that the
kernel has only one orbit, necessarily the orbit of the trivial element
(which consists only of the trivial element), so the kernel is trivial.
Thus the morphisms on homotopy groups of path components are isomorphisms,
giving a weak equivalence.
\eop

Suppose ${\bf G}= \{ {\bf G}_i\}$ is a collection classes of sheaves (resp.
sheaves of groups for $i=1$ and sheaves of abelian groups for $i\geq 2$) on
$\Xx$, and suppose that ${\bf G}_i \subset {\bf G}_{i-1}$. Let ${\bf F}$ be
the class of flexible functors $T$ such that all homotopy group sheaves
$\varpi _i$ are in ${\bf G}_i$. Suppose that ${\bf F}$ satisfies the
closure hypothesis above (or at least what is necessary to apply the
lemmas). In this case we also use the terminology {\em equivalence relative
to ${\bf G}$} or {\em ${\bf G}$-completion}.

\begin{lemma}
\mylabel{mapp10}
In the situation of the previous paragraph, if $\psi :R\rightarrow S$ is a
morphism of functors inducing an isomorphism $$
\psi ^{\ast}: H^i(S, G)\cong H^i(R,G)
$$
for all $G\in {\bf G}_i$, then $ \psi$ is an equivalence relative to ${\bf
G}$. \end{lemma}
{\em Proof:}
???
\eop

\subnumero{The Leray spectral sequence}
Suppose $G$ is a sheaf of groups (abelian for $n\geq 2$) on $\Xx$. Suppose
$S\rightarrow T$ is a morphism of flexible sheaves, and suppose that
$\varpi _0(T)=\ast$.
Let
$$
Map _T^{\leq j}(S, K(G, n)\times T)
$$
denote the tower obtained by applying the relative Postnikov truncation
(relative to $T$) to $Map _T(S,K(G,n)\times T)$. We obtain a tower of spaces $$
Mor ^{\leq j}_{/T}(S, K(G, n)):=
Sect (T, Map _T^{\leq j}(S, K(G, n)\times T)) $$
whose limit is $Mor (S, K(G, n))$. Suppose $\eta \in Mor (S, K(G,n))$
projects to an element which we denote by $\eta_{j-1}$, of $Mor ^{\leq
j-1}_{/T}(S, K(G, n))$. The fiber $F_j(\eta )$ over $\eta$ is equal to
$Sect (T, Fib _j(\eta ))$ where $Fib _j(\eta )$ is the pullback of $$
Map _T^{\leq j}(S, K(G, n)\times T)\rightarrow Map _T^{\leq j-1}(S, K(G,
n)\times T)
$$
by the section
$\eta _{j-1}:T\rightarrow Map _T^{\leq j-1}(S, K(G, n)\times T)$. The fiber
of $S\rightarrow T$ is locally well defined, that is there is a covering of
$\Xx$ (we assume $\Xx$ has a final object) with sections of $T$ on each
element of the covering, and we can look at the fibers over these sections;
these are locally isomorphic over the fiber products of elements in the
covering. With this understanding, we can say that $Fib _j(\eta
)\rightarrow T$ is a fibration with fiber which is locally of the form
$K(A, j)$ where $A$ is the $n-j$th cohomology of the fiber of $S\rightarrow
T$ with coefficients in $G$.

We get a spectral sequence
$$
E^{p,q}_2= \pi _{-p-q}(Sect (T, Fib _{-q}(\eta )), \eta _j) \Rightarrow \pi
_{-p-q}(Mor (S, K(G,n)). $$
This is the Leray spectral sequence in our case; note that $Fib _j(\eta
)\rightarrow T$ are the analogues of local systems of coefficients on $T$,
and the $\pi _i (Sect (T, Fib _j (\eta )), \eta _j )$ are the analogues of
cohomology of $T$ with coefficients in these local systems. The spectral
sequence converges to cohomology of the upper space $S$.

We introduce some
notation to suggest this analogy. Put
$$
R^i(S/T, G; n):= (Fib _{n-i}(\eta )\rightarrow T) $$
(locally a fibration with fiber $K(A, n-i)$ where $A$ is $H^i$ of the fiber
of $S/T$), and define
$$
H^k(T, R^i(S/T, G; n)):= \pi _{n-i-k}(Sect(T, R^i(S/T, G; n)), \eta ).
$$
Then our spectral sequence (with $-p-q = n-i-k$ and $-q=n-i$, thus $p=k$
and $q= i-n$) becomes
$$
E_2^{k, i-n}= H^k(T, R^i(S/T, G; n))\Rightarrow H^{k+i-n}(S, G). $$
This shows the analogy with the Leray spectral sequence; note that it is
just shifted by decreasing $q$ by $n$.

\begin{corollary}
\mylabel{mapp11}
Suppose
$$
\begin{array}{ccc}
S&\rightarrow &S' \\
\downarrow && \downarrow \\
T &\rightarrow &T'
\end{array}
$$
is a diagram with $T$ and $T'$ connected (in the sheaf-theoretic sense).
Suppose that the morphism on fibers (which is defined locally) induces an
isomorphism on cohomology with coefficients $G$ in a certain class ${\bf
G}$, and that the morphism $T\rightarrow T'$ induces an isomorphism on
cohomology with all local coefficient systems which can arise. Then
$S\rightarrow S'$ induces an isomorphism on cohomology with coefficients
$G$ in ${\bf G}$. \end{corollary}
{\em Proof:}
The morphism induces an isomorphism on the $E_2$ term of the Leray spectral
sequence, so it is an isomorphism on the limit. (??? basepoints ???) \eop

\subnumero{Local systems}
Suppose $T$ is a contravariant flexible functor on $\Xx$. A {\em local
system $L$ over $T$} consists of the data, for each $X\in \Xx$, of a local
system $L_X$ (of abelian groups or other objects) over $T_X$; and for each
$\phi : Y\rightarrow X$, a morphism $L(\phi ): L_X\rightarrow T(\phi
)^{\ast}(L_Y)$; such that for a composable pair
$Z\stackrel{\psi}{\rightarrow}Y \stackrel{\phi}{\rightarrow}X$, the diagram
$$
\begin{array}{ccc}
L_X & \stackrel{L(\phi \psi )}{\rightarrow} & T(\phi \psi )^{\ast}(L_Z) \\
\downarrow {\scriptstyle L(\phi )} & & \downarrow {\scriptstyle \cong }\\
T(\phi )^{\ast}(L(Y)) & \stackrel{T(\phi )^{\ast}(L(\psi ))}{\rightarrow} &
T(\phi )^{\ast}T(\psi )^{\ast} (L_Z)
\end{array}
$$
commutes, where the vertical arrow on the right is the isomorphism given by
the homotopy $T(\phi , \psi )$. Note that if $T$ is a functor, then this
isomorphism is the identity and we can write the condition more easily as
$$
T(\phi )^{\ast}(L(\psi ))\circ L(\phi ) = L( \phi \psi ). $$

Suppose $T$ is a contravariant functor and $L$ is a local system of abelian
groups on $T$. Then we can define the presheaf of complexes
$$
C^{\cdot , {\rm pre}}(T; L)_X := C^{\cdot}_{\rm sing}(T_X; L_X), $$
whose values are the singular cochain complexes on $T_X$ with coefficients
in the local systems $L_X$. Let
$C^{\cdot}(T; L)$ be the sheafification of $C^{\cdot , {\rm pre}}(T; L)$.
We define the {\em cohomology of $T$ with coefficients in $L$} to be the
sheaf hypercohomology $H^{\cdot}(T; L)$ of $\Xx$ with coefficients in
$C^{\cdot}(T; L)$ (that is, take an injective resolution of $C^{\cdot}(T;
L)$, take the global sections over $\Xx$ and take the cohomology).

We can make the complex which calculates this cohomology, canonical, as
follows: let $V$ denote the functor direct limit over all hypercoverings,
of the associated Verdier (hyper-\v{C}ech) complex. Then $VC^{\cdot}(T, L)$
is a complex whose cohomology is $H^{\cdot }(T,L)$. This is compatible with
restrictions, and in fact if $f:S\rightarrow T$ is a morphism of functors
then we obtain a morphism of complexes $VC^{\cdot}(T, L)\rightarrow
VC^{\cdot}(S, f^{\ast}L)$. (???)

??? sheaf condition on $L$ ???

We may think of
$L$ as a flexible functor by looking at its ``espace \' etal\' e'' $L^{\rm
ee}$ with morphism to $T$. Note that here, the higher homotopies for $T$
are lifted to give those of $L^{\rm ee}$, and the morphism is really a
morphism of structures of flexible functor.

We obtain a local system $\Gamma (\Xx , L)$ over $\Gamma (\Xx , T)$,
defined by the condition that its espace \' etal\' e is $\Gamma (\Xx ,
L^{\rm ee})$.

We say that $L$ is a {\em sheafy local system} if for any $X\in \Xx$ and
any sieve $\Bb$ over $X$, the morphism
$$
L_X\rightarrow p_{\Bb}^{\ast} \Gamma (\Bb , L) $$
is an isomorphism, where $p_{\Bb}: T_X\rightarrow \Gamma (X, T)$ is the
natural map.

If $L$ is a local system over $T$, we can define $L'$ by $$
(L')_X:= \lim _{\rightarrow , \Bb } p_{\Bb}^{\ast} \Gamma (\Bb , L) $$
and $L''$ by repeating this operation. Then $L''$ is a sheafy local system,
and we have a map $L\rightarrow L''$. We call this the {\em sheafification
of $L$}.

(Universal property of $L''$??)

\subnumero{Eilenberg-MacLane fibrations} Suppose $T$ is a contravariant
functor which is a flexible sheaf. An {\em Eilenberg-MacLane fibration of
degree $n$ over $T$} is a morphism $\zeta :E\rightarrow T$ of flexible
sheaves such that for any $X\in \Xx$ and any $\eta \in T_X$, the fiber of
$\zeta $ over $\eta$ is an Eilenberg-MacLane sheaf of degree $n$ (that is,
the homotopy group sheaves are trivial except in degree $n$).

For the moment, fix $n\geq 2$. Suppose $\zeta : E\rightarrow T$ is an
Eilenberg-MacLane fibration of degree $n$. Then we define a local system
$L^{\rm pre}$ on $T$ by
setting
$L^{\rm pre}_X$ to be the local system on $T_X$ corresponding to the $n$-th
homology of the fibers of $\zeta _X$. Then let $L$ be the sheafification of
$L^{\rm pre}$.

\begin{proposition}
\mylabel{mapp12}
There is a canonical class $\kappa \in VC^{n+1}(T, L)$ with the property
that for any morphism of functors $f:S\rightarrow T$, the homotopy classes
of sections of the pullback $S\times _T^{\rm path}E$ over $S$ are in
one-to-one correspondence with the elements $u\in VC^{n}(S, f^{\ast}L)$
such that $d(u)= f^{\ast}(\kappa )$.
\end{proposition}
{\em Proof:}
???
\eop

Given a local system $L$ over a flexible functor $T$, we should be able to
associate an Eilenberg-MacLane fibration $K(L, n)\rightarrow T$ with a
section $o$ and associated local system $L$.

Then for any Eilenberg MacLane fibration $\zeta : E\rightarrow T$ of degree
$n$ with local system $L$, we should get a section $\kappa$ of $K(L, n+1)$
and equivalence between $\zeta $ and the relative path-space $$
P^{o,\kappa}_{/T}(K(L, n+1))\rightarrow T. $$


\subnumero{Projection formulas}
Suppose $\Xx$ and $\Jj$ are categories, and suppose $\Sigma$ is a presheaf
of sets over $\Xx \times \Jj$. Suppose $T:\Xx \rightarrow Top$ is a
flexible sheaf. Let $S=lim _{\Xx \times \Jj /\Xx}\Sigma $ be the
contravariant functor from $\Xx$ to $Top$ defined by
$$
S_X:= \lim _{\rightarrow } (\Sigma |_{\{ X\} \times \Jj}). $$

\begin{lemma}
\mylabel{mapp13}
There is a natural (weak ??) homotopy equivalence $$
Mor _{\Xx}(S,T)\rightarrow Mor _{\Xx \times \Jj} (\Sigma , pr_1^{\ast}T). $$
\end{lemma}
{\em Proof:}
???
\eop

Suppose $a:\Yy \rightarrow \Xx$ is a functor, and suppose $S$ and $T$ are
flexible presheaves on $\Yy$. Then we can make a flexible presheaf $$
Map_{\Yy /\Xx}(S,T)
$$
on $\Xx$ whose value on $X$ is the space $$
Mor _{\Yy /X}(S|_{\Yy /X},T|_{\Yy /X}).
$$
Here $\Yy /X$ is the category of pairs $(Y,f)$ where $Y\in \Yy$ and
$f:a(Y)\rightarrow X$ is a morphism in $\Xx$. Note that for $X'\rightarrow
X$ we have a functor $\Yy /X'\rightarrow \Yy /X$, so we obtain a pullback
$$
Mor _{\Yy /X}(S|_{\Yy /X},T|_{\Yy /X})\rightarrow Mor _{\Yy /X'}(S|_{\Yy
/X'},T|_{\Yy /X'}). $$
This is strictly associative, so $Map_{\Yy /\Xx}(S,T)$ becomes a
contravariant functor (which can be considered as a flexible functor by
taking the higher homotopies to be constant).

\begin{lemma}
\mylabel{mapp15}
The morphism given by functoriality for $\Yy /X\rightarrow \Yy$ is a (weak
?) homotopy equivalence
$$
Mor _{\Yy} (S,T)\rightarrow \Gamma (\Xx ,Map_{\Yy /\Xx}(S,T)). $$
\end{lemma}
{\em Proof:}
???
\eop

\begin{lemma}
\mylabel{mapp16}
In the situation of the start of this subsection, there is a natural (weak
??) homotopy equivalence
$$
Mor _{\Xx \times \Jj} (\Sigma , pr_1^{\ast}T)\rightarrow \Gamma (\Jj , Map
_{\Xx\times \Jj /\Jj}(\Sigma ,pr_1^{\ast}T). $$
\end{lemma}
{\em Proof:}
???
\eop

\newpage

\setcounter{section}{8}

\numero{More homotopy group sheaves}

We begin with some lemmas relating the homotopy group sheaves and presheaves.

\begin{lemma}
\mylabel{piA01}
If $T$ is an $n$-truncated weak flexible sheaf on $\Xx$, and $\eta \in
\Gamma (\Xx , T)$, then the presheaf $\pi _n^{\rm pre}(T,\eta )$ is a
sheaf. \end{lemma}
{\em Proof:}
By taking loop spaces, it suffices to prove this when $n=0$. But if $T$ is
a $0$-truncated flexible functor, then for any sieve $\Bb$ over an object
$X$, $\Gamma (\Bb , T)$ is $0$-truncated and $\pi _0\Gamma (\Bb , T)=
\Gamma (\Bb , \pi _0^{\rm pre}(T))$. Thus $T$ is a weak sheaf if and only
if $\pi _0^{\rm pre}(T)$ is a sheaf.
\eop

\begin{corollary}
\mylabel{piA02}
If $T$ is an $n$-truncated weak flexible sheaf, and for $i>m$ the homotopy
group sheaves $\pi _i(T|_{\Xx /X}, x)$ vanish for all $x\in T_X$, then $T$
is $m$-truncated.
\end{corollary}
{\em Proof:} By induction it suffices to prove this when $m=n-1$; but in
that case
$$
\pi _n(T|_{\Xx /X}, x)= \pi _n^{\rm pre}(T|_{\Xx /X}, x)
$$
by the
above lemma, so the vanishing of the homotopy group sheaves implies that
the homotopy groups over each object vanish, and $T$ is $n-1$-truncated.
\eop

\begin{corollary}
\mylabel{piA03}
If $T$ is a truncated weak flexible sheaf such that the homotopy group
sheaves are all trivial, then $T$ is weakly equivalent to $\ast$.
\end{corollary}
{\em Proof:} In this case $T$ is $0$-truncated and $\pi _0(T)=\ast$. The
unique morphism $T\rightarrow \ast$ is a weak equivalence. \eop

\subnumero{Local systems}
Suppose $\Xx$ is a category and $S$ is a flexible functor on $\Xx$. We say
that a morphism $L\rightarrow S$ of flexible functors is {\em relatively
$n$-truncated} if the induced morphisms of homotopy group presheaves $\pi
_i ^{\rm pre}(L|_{\Xx /X},y)\rightarrow \pi _i^{\rm pre}(S|_{\Xx /X},s)$
(where $s$ is the image of $y$ in the sense described in Section ???), are
isomorphisms for $i\geq n+2$ and injective for $i=n+1$. This is equivalent
to saying that for any point $s\in S_X$, the fiber of $L|_{\Xx
/X}\rightarrow S|_{\Xx /X} $ over $s$ is $n$-truncated.

Suppose $\Xx$ is a site and $L$ and $S$ are truncated weak flexible
sheaves. Suppose that the induced morphism on homotopy group sheaves is
injective for $i=n+1$ and an isomorphism for $i\geq n+2$. Then this implies
that the homotopy group sheaves of the fiber of $L|_{\Xx /X}$ over a point
$s\in S_X$ vanish above degree $n$ (by the long exact sequence of homotopy
group sheaves coming from a fibration). This implies that the fiber is
$n$-truncated (because the fiber is also a weak flexible sheaf, by ???).
Therefore the morphism $L\rightarrow S$ is relatively $n$-truncated. Thus
we may measure relative truncation using homotopy group sheaves.

Keep the hypotheses that $\Xx$ is a site and $S$ is a truncated weak
flexible sheaf.
A {\em covering} of $S$ is a relatively $0$-truncated morphism
$L\rightarrow S$ such that $L$ is a truncated weak flexible sheaf (in which
case its degree of truncation is no worse that that of $S$). The condition
that $L$ be a weak flexible sheaf should be emphasized here since we will
not repeat it.

If $L\rightarrow S$ and
$L'\rightarrow S$ are two coverings, the space of morphisms from $L$ to
$L'$ over $S$ is $0$-truncated, so we can be a little bit less careful
about defining things precisely. The fiber product $L\times _SL'$ of two
coverings is again a covering.

A {\em binary operation} on a covering $L$ over $S$ is a morphism $L\times
_SL\rightarrow L$ over $S$ (that is, a diagram including this morphism and
the two morphisms to $S$). To be precise, the choice of path-cartesian
diagram
$$
L\times _SL\rightarrow L\times L \stackrel{\displaystyle
\rightarrow}{\rightarrow} S
$$
is part of the data of a binary operation. But any two choices are
equivalent by a very well defined equivalence, so the definition of the
operation for one choice leads to a definition for every other choice, and
these definitions are all compatible with each other up to very well
defined homotopies, under the equivalences between choices of fiber
products. Suppose $L$ and $L'$ are
coverings of $S$ with binary operations. A morphism $L\rightarrow L'$ is
{\em compatible with the binary operations} if the two resulting morphisms
$L\times _SL\rightarrow L'$ are homotopic. Note that this condition does
not depend on the choice of fiber products or on the choice of the map
$L\times _SL\rightarrow L'\times _SL'$ induced by $L\rightarrow L'$. We say
that two binary operations on the same covering $L$ are {\em equivalent} if
the identity morphism $L\rightarrow L$ is compatible with them (putting one
operation on each of the copies of $L$).

A binary operation on $L$ is {\em associative} if the two resulting maps
$L\times _SL\times _SL\rightarrow L$ are homotopic (in the space of maps
over $S$). Again, this condition does not depend on the choice of fiber
products, and is preserved by equivalence of binary operations. Similarly,
a binary operation is {\em commutative} if the two maps $L\times
_SL\rightarrow L$ are homotopic. Given two operations, we say that they
satisfy the {\em distributive law} if the two morphisms $$
L\times _SL\times _SL\times _SL \rightarrow L $$
corresponding to $(a+b)(c+d)$ and $ac+ad+bc+bd$ are homotopic.

Finally, an associative binary operation {\em admits an identity} if there
is a section $e:S\rightarrow L$ such that the morphism $L\rightarrow L$
corresponding to multiplication by this section, is homotopic to the
identity. If an associative binary operation admits an identity, then the
identity is unique up to very well defined homotopy (by the usual proof of
multiplying two different identities together). Suppose given an
associative binary operation which admits an identity. We say that it {\em
admits an inverse} if there is a morphism $L\rightarrow L$ over $S$ such
that the resulting morphism $L\rightarrow L$ (composition of $a$ with the
supposed inverse of $a$) is homotopic to the composition $L\rightarrow
S\stackrel{e}{\rightarrow}L$.

If $L$ is a covering of $S$ and if
$s\in S_X$ then the {\em fiber} $L_s$ of $L$ over $s$ is defined to be the
set $\pi _0(L_X\times _{S_X}\ast )$. A morphism between coverings induces a
morphism between fibers over any point. If $L\times _SL\rightarrow L$ is a
binary operation on a covering $L$, we obtain a binary operation on each
fiber $L_s$.

\begin{lemma}
\mylabel{bigPi01bis}
Suppose $L$ is a covering of $S$ with a binary operation. Suppose $L$ and
$S$ are of CW-type.
The operation is associative (resp. is commutative, resp. is associative
and admits an identity, resp. is associative and admits an identity and an
inverse)
if and only if, for all $X\in \Xx$ and all point
s $s\in S_X$, the resulting binary operation on $L_s$ is associative (resp. is
commutative, resp. is associative and admits an identity, resp. is associative
and admits an identity and an inverse). Given two operations, they satisfy the
distributive law if and only if, for any $X\in \Xx$ and any $s\in S_X$, the
resulting two operations on $L_s$ satisfy the distributive law. \end{lemma} {\em
Proof:} We analyse the notion of morphism of coverings. Suppose $L\rightarrow S$
and $L'\rightarrow S$ are two coverings, with $L$, $L'$, and $S$ weak flexible
sheaves of CW type. The space of morphisms from $L$ to $L'$ over $S$ is the
space of flexible functors $\Xx \times I^{(2)}\rightarrow Top$ which restrict to
the morphism $L\rightarrow S$ over the subcategory $\Xx \times \{ 0\leftarrow
2\}$ and which restrict to the morphism $L`\rightarrow S$ on the subcategory
$\Xx \times \{ 1\leftarrow 2\}$. There is a subcategory $A\subset I^{(2)}$
consisting of all the objects, and all the morphisms except $0\leftarrow 1$. The
triple $S,L,L'$ with the pair of morphisms $L\rightarrow S$ and $L'\rightarrow
S$ corresponds to a flexible functor $\Xx \times A\rightarrow Top$, and a
morphism of coverings is an extension of this to a flexible functor on $\Xx
\times I^{(2)}$. By the invariance of such things under homotopy equivalences,
we can replace the flexible functor $\Xx \times A\rightarrow Top$ by a functor,
which is homotopy equivalent. The space of ways of completing this new functor
to a flexible functor on $\Xx \times I^{(2)}$ will be the equivalent to the old
one, and this equivalence preserves the morphisms between the fibers (note that
$S$ itself is replaced by an equivalent functor, so there are a few things which
need to be said about fibers over points in the old $S$ versus the new
one---these  can safely be left to the reader). We can now replace $L$ by $\tau
_{\infty}L$ and we can replace $L'$ by the path fibration over $S$, so we may
assume that $L\rightarrow S$ and $L'\rightarrow S$ are morphisms of functors,
with $L_X$ a CW complex for all $X$ and $L'_X\rightarrow S_X$ a fibration for
all $S$. Now the canonical composition of morphisms of flexible functors
$L\rightarrow L'$, with the morphism of functors $L'\rightarrow S$, gives a map
$$ Mor (L,L')\rightarrow Mor (L,S). $$ One can see that it is a fibration, using
the fact that the values of $L$ are CW complexes and the values of
$L'\rightarrow S$ are fibrations. The canonical composition gives an equivalence
between $Mor (L,L')$ and the fiber of $Mor ^2(L,L',S)$ over the given element of
$Mor (L',S)$. Our space of morphisms in question is the fiber of $Mor
^2(L,L',S)$ over the given point in $Mor (L,S)\times Mor (L',S)$ (this map is a
fibration). This fiber can be calculated by first taking the fiber over the
point in $Mor (L',S)$ and then taking the fiber over the point in $Mor (L,S)$.
Via the equivalence given by the canonical composition, this is equivalent to
the homotopy fiber of $Mor (L,L')$ over the given point in $Mor (L,S)$. But
since this map is a fibration, the homotopy fiber is equal to the actual fiber.
Thus the space of morphisms from $L$ to $L'$ over $S$ is equal to the fiber of
$Mor (L,L')$ over the given element of $Mor (L,S)$. We can now make this
explicit. An element here consists, for each composable sequence $\phi _1,
\ldots , \phi _k$ in $\Xx$, each increasing sequence $\epsilon _{\cdot}$ of
numbers $0$ or $1$, and each sequence $t_1,\ldots , t_{k-1}$, of a commutative
diagram $$ \begin{array}{ccc} L_{X_0}&\rightarrow & L'_{X_k} \\ \downarrow &&
\downarrow \\ S_{X_0}&\rightarrow &S_{X_k} \end{array} $$ satisfying the
required properties of continuity in $t$, and the required properties for
$t_i=0$ or $t_i=1$. The space of such diagrams, for a given morphism
$X_0\leftarrow X_k$, is $0$-truncated, because $L$ and $L'$ are $0$-truncated
over $S$. Thus the required family of maps exists (and the space of them is
weakly contractible) for composable sequences of length $k\geq 3$. For $k=2$,
the space of choices is either empty or weakly contractible. For $k=1$ the space
of choices is $0$-truncated. The space of all choices (for all composable
sequences) is equal to the space of choices for all composable sequences of
length one, such that for the composable sequences of length two the required
map exists. In particular, the space of choices is $0$-truncated, and two
choices are homotopic if and only if they are homotopic for each composable
sequence of length one. By analysing this a bit more closely (using the fact
that each morphism in $\Xx \times I$ decomposes as a morphism in $\Xx$ times a
standard horizontal morphism), one sees that the space of choices is equal to
the space of choices of a diagram $$ \begin{array}{ccc} L_X&\rightarrow & L'_X
\\ \downarrow && \downarrow \\ S_X &\rightarrow &S_X \end{array} $$ for each
$X\in \Xx$, subject to a compatibility condition: for each morphism
$X\rightarrow Y$ in $\Xx$, the diagram $$ \begin{array}{ccc} L_X&\rightarrow
&L'_X\\ \downarrow && \downarrow \\ L_Y&\rightarrow &L'_Y \end{array} $$
commutes up to homotopy relative to $S_X$.

It is now easy to see that if two morphisms $f,g:L\rightarrow L'$ of
coverings over $S$ induce the same morphism on all fibers, then they are
homotopic; because the diagrams as above over each $X$ are then homotopic.
This implies the statement of the lemma
for associativity, commutativity, and distributivity.

For the existence of the identity, suppose that an associative binary
operation admits an identity on each point. Then over each $S_X$ we obtain
a covering with associative binary operation admitting an identity in each
fiber; since the identity element in each fiber is unique, it is fixed
under the action of the fundamental group of the path component of $S_X$,
so there is a section $e_X:S_X\rightarrow L_X$. This is the unique section
restricting to the identity elements in the fibers. We have to check that
these
sections satisfy the required compatibility condition. Suppose
$X\rightarrow Y$ is a morphism.
The identity element in the fiber of $L_Y$ over $s\in S_Y$ maps to the
identity element in the fiber of $L_X$ over the image $s_X$, since the map
$L_{Y,s}\rightarrow L_{X,s_X}$ is a map of sets with binary operation. This
allows us to verify commutativity of the required diagram (using the fact
that a morphism of coverings is determined up to homotopy by the morphisms
on the fibers). We obtain the required section $e:S\rightarrow L$.

A similar argument gives the construction of the involution $L\rightarrow
L$ sending an element to its inverse for the operation (again, the
restriction maps on fibers preserve the operations, so they send the
inverses of an element to the inverses of the image of that element). \eop

Similar considerations hold for operations among several coverings (which
will be used in the case of modules over a ring, for example).

A {\em local system of groups} over $S$ is a covering $L\rightarrow S$
together with a
binary operation which is associative, and admits an identity and an
inverse. It is an {\em abelian local system} if, in addition, the operation
is commutative. A {\em local system of rings} is a covering $A\rightarrow
S$ together with two associative binary operations denoted $\cdot$ and $+$,
both admiting identities denoted $1$ and $0$ respectively, such that $+$
admits an inverse, and such that they satisfy the distributive law. If $A$
is a local system of rings, then we can similarly define the notion of a
{\em local system of $A$-modules} $M$.

Suppose $S=\ast$, and $L$ is a covering of $S$. Then $\pi _0 (L)$ is a
sheaf of sets on $\Xx$. Any operations on $L$ as defined above give
corresponding operations on the sheaf $\pi _0(L)$, and conversely, given
any operations on $\pi _0(L)$ we get corresponding (very well defined)
operations on $L$. Similarly for morphisms of coverings and morphisms of
sheaves. By these constructions, we shall think of the theory of coverings
of $\ast$ as being equivalent to the theory of sheaves on $\Xx$.

Suppose $A$ is a local system of rings on $S$, and suppose that $L$, $M$,
and $N$
are local systems of $A$-modules. A map of coverings $L\times M\rightarrow
N$ is {\em $A$-bilinear} if the map on fibers $L_s\times M_s\rightarrow
N_s$ is $A_s$-bilinear for all points $s\in S_X$, all $X\in \Xx$. The {\em
tensor product} $L\otimes _AM$ is a local system of $A$-modules with an
$A$-bilinear map $L\times M\rightarrow L\otimes _A M$ satisfying the
universal property that for any $A$-bilinear map $L\times M\rightarrow N$
there is a factorization $L\otimes _AM\rightarrow N$ which is unique up to
homotopy. This characterizes the tensor product up to an isomorphism which
is unique up to homotopy (and note that the homotopy is very well defined
since the space of maps between coverings is $0$-truncated). The tensor
product may be constructed as above, by taking the tensor product of the
local systems of modules over each space; this satisfies a universal
property in terms of things which are like coverings but which
are not sheaves (and this universal property also provides the uniqueness
needed to make the construction);
then sheafify, and this satisfies the universal property for coverings
(which are, by definition, sheaves).

If $f:S\rightarrow T$ is a morphism of weak flexible sheaves, and if
$L$ is a covering of $T$, then put
$$
f^{\ast}(L)= L\times _TS.
$$
It is a covering of $S$. If $L$ comes with some binary operations
satisfying certain axioms,

\subnumero{Relative truncation}

Suppose $T\rightarrow S$ is a morphism of flexible functors. We define the
{\em relative presheaf truncation}
$$
\tau ^{\rm pre}_{\leq n}(T/S)\rightarrow S $$
by replacing this morphism with a morphism of functors, then applying a
canonical fiberwise truncation over each object. To define a canonical
fiberwise truncation of a morphism $A\rightarrow B$ of topological spaces,
consider the following operation {\em ($k$-killing)}: for each morphism
$S^k\rightarrow A$ and $B^{k+1}\rightarrow B$ such that the composition of
the first with the map $A\rightarrow B$ is the boundary of the second, add
a cell $B^{k+1}$ to $A$ with the given boundary, and with map to $B$ as
given. Note that this is natural with respect to morphisms of the diagram
$A\rightarrow B$. To obtain the truncation $\tau ^{\rm pre}_{\leq n}$,
apply this operation for $k=n$, $k=n+1$, and so forth, then take the union
of all of the spaces.

There is a commutative diagram (that is, a flexible functor from $\Xx
\times I\times I$ to $Top$):
$$
\begin{array}{ccc}
T&\rightarrow &\tau ^{\rm pre}_{\leq n}(T/S)\\ \downarrow &&\downarrow \\
S&= &S.
\end{array}
$$

The relative truncation is relatively $n$-truncated over $S$, and has the
property that
for any morphism $R\rightarrow S$ which is relatively $n$-truncated, the
map from the space of morphisms $\tau ^{\rm pre}_{\leq n}(T/S)\rightarrow
R$ over $S$, to the space of morphisms $T\rightarrow R$ over $S$, is a weak
equivalence. (Proof ???). In particular, the space of morphisms (over $S$,
and with the map from $T$) between any two relative truncations is weakly
contractible and weakly equivalent to the space of equivalences (over $S$
and with the map from $T$) which maps to it.

If $T\rightarrow S$ is a morphism of truncated weak flexible sheaves, define $$
\tau _{\leq n}(T/S)
$$
to be the sheafification of $\tau ^{\rm pre}_{\leq n}(T/S)$. There is again
a morphism $T\rightarrow \tau _{\leq n}(T/S)$ over $S$. This sheafified
relative truncation is again relatively $n$-truncated over $S$ (proof ???)
and satisfies
the same universal property as above, for relatively $n$-truncated
morphisms of weak sheaves $R\rightarrow S$ (this comes from the universal
property of sheafification).

Relative truncation commutes with pathwise fiber products, that is there is
a canonical equivalence (very well defined by the universal properties)
$$
\tau _{\leq n}(T/S)\times _SS' = \tau _{\leq n}(T\times _SS'/S'). $$
In particular, the fiber of the relative truncation over a point $s\in S_X$
is the relative truncation of the fiber of $T$ over $s$.

\subnumero{Relative homotopy groups}

Suppose $T\rightarrow S$ is a morphism of weak flexible sheaves on a site
$\Xx$, and suppose $S\stackrel{\sigma}{\rightarrow}T\rightarrow S$ is a
section.
Put
$$
\Omega ^{\sigma}(T/S):= S\times ^{\rm path}_T S $$
with its map to $S$ (fix the second projection, for example, but note that
both projections are homotopic). We call this the {\em relative loop
space}. There is a map
$$
\Omega ^{\sigma}(T/S)\times ^{\rm path}_S\Omega ^{\sigma}(T/S) \rightarrow
\Omega ^{\sigma}(T/S)
$$
easily obtained from the definition, and which corresponds to concatenation
of loops. If $s\in S_X$ then the fiber of $\Omega ^{\sigma}(T/S)$ over $s$
is the loop space of the fiber of $T$ over $s$, based at $\sigma (s)$. The
operation restricts to the operation of concatenation of loops on this
fiber. In particular, after taking the relative $0$-truncation, we obtain
an associative operation with identity and inverse.

Define
$$
\pi _0(T/S)= \tau _{\leq 0}(T/S)
$$
and for $i\geq 1$ put
$$
\pi _i(T/S, \sigma ):= \tau _{\leq 0} ((\Omega ^{\sigma})^i (T/S)). $$
We call these the {\em homotopy local systems of $T$ over $S$ based at
$\sigma$.}
They are local systems of groups, and abelian groups for $i\geq 2$. (Note
that the commutativity holds by checking it on the fibers.) If $s\in S_X$
then the fiber of the homotopy local system over $s$ is the homotopy group
sheaf of the fiber of $T$ over $s$, that is $$
\pi _i(T/S,\sigma )_s = \pi _i (T|_{\Xx /X}\times _{S|_{\Xx /X}}\ast ,
\sigma (s))
$$
where product is by the morphism $s:\ast \rightarrow S|_{\Xx /X}$.

Formation of the homotopy local systems commutes with fiber products: if
$f:S'\rightarrow S$ and $T\rightarrow S$ are morphisms of truncated weak
flexible sheaves, and if $\sigma :S\rightarrow T$ is a section, then $$
\pi _i ( T\times _SS'/S', \sigma |_{S'})= f^{\ast}\pi _i (T/S, \sigma ).
$$

A morphism $T\rightarrow S$ between truncated weak flexible sheaves is
relatively $n$-truncated if and only if the homotopy local systems $$
\pi _i (T\times _ST/T, 1)
$$
formed using the identity section $T\rightarrow T\times _ST$, vanish for
$i>n$. (proof ???)

\subnumero{Relative Eilenberg-MacLane sheaves}

Suppose $S$ is a truncated weak flexible sheaf of CW type, suppose $G$ is a
local system of groups over $S$, and $n$ a positive integer. If
$T\rightarrow S$ is a morphism, $\sigma:S\rightarrow T$ is a section, and
$\eta : G\rightarrow \pi _n(T/S,\sigma )$ is a morphism of local systems of
groups over $S$,
then we say that {\em $(T/S, \sigma ,\eta )$ is a $K(G/S,n)$} if $T$ is
a truncated weak flexible sheaf of CW type, if $\pi _0 (T/S)=S$, if $\pi _i
(T/S, \sigma )=0$
for $i\neq n$, and if $\eta : G\cong \pi _n (T/S, \sigma )$. If this exists
for $n\geq 2$ then $G$ must be abelian.

We say that $(T/S, \sigma )$ {\em is a relative Eilenberg-MacLane sheaf of
degree $n$ over $S$}
if, when we put $G=\pi _n(T/S, \sigma )$ and $\eta$ equal to the identity,
then $(T/S, \sigma , \psi )$ is a $K(G/S, n)$.

Note that a relative Eilenberg-MacLane sheaf of degree $n$ over $S$ is
relatively $n$-truncated over $S$ (apply Corollary \ref{piA02} to the
fibers).

We make a similar definition for any covering $G$, for $n=0$.

\begin{theorem}
\mylabel{bigPi03}
Suppose $G$ is a local system of groups (abelian for $n\geq 2$) over a
truncated weak flexible sheaf of CW type $S$. Then there exists a triple
$(T/S, \sigma , \eta )$ which is a $K(G/S, n)$. If $(T/S, \sigma , \eta )$
and $T'/S, \sigma ', \eta ')$ are two such, then the space of morphisms
$T\rightarrow T'$ over $S$ and compatible with $\eta $ and $\eta '$ is
weakly contractible; and the space of equivalences over $S$ is also weakly
contractible (and hence every morphism extends to an equivalence).
\end{theorem}

We will prove this theorem in several steps. First of all, we treat the
case where $S=\ast$. We speak of {\em Eilenberg-MacLane sheaves} instead of
relative Eilenberg-MacLane sheaves, in this case.

For any covering $G$ of $S$ ($n=0$), local system $G$ of groups ($n=1$) or
local system $G$ of abelian groups ($n\geq 2$) we can construct a {\em
standard Eilenberg-MacLane sheaf} $K^{\rm st}(\Gg , n)$ as follows:
there is a construction of the Eilenberg-MacLane space $K(G,n)$ which is
functorial in $G$; use this to construct a presheaf $K^{\rm pre}(\Gg , n)$
with a sectoin denoted $o$, whose homotopy group presheaves are trivial
except for $\Gg$ in degree $n$; then put $$
K^{\rm st}(\Gg , n):= F_{\infty}K^{\rm pre}(\Gg , n). $$

(is the next lemma actually contained somewhere else ???)

\begin{lemma}
\mylabel{piA05}
If $S\rightarrow T$ is a morphism of truncated weak flexible sheaves of CW
type and if $t\in \Gamma (\Xx ,T)$ then this can be completed to a
fibration diagram $Q\rightarrow S\rightarrow T$ with basepoint $t$, such
that $Q$ is a weak flexible sheaf. This fibration diagram induces a long
exact sequence of homotopy group sheaves (for a basepoint $q\in \Gamma (\Xx
, Q)$ and its image $s$ in $S$).
\end{lemma}
{\em Proof:}
We can complete the morphism into a fibration diagram $Q\rightarrow
S\rightarrow T$ with $t$ as the basepoint in $T$. Applying the functor
$F_{\infty}$ preserves the property of having a fibration diagram, and it
preserves $S$ and $T$; thus we obtain a fibration diagram with
$F_{\infty}Q$ as the fiber, which is a flexible sheaf. By uniqueness of the
fibration diagram (Lemma \ref{fibr01}), $Q\sim F_{\infty}Q$, so $Q$ was a
flexible sheaf. The fact that the fibration diagram is homotopic to a
fibration over every object gives a long exact sequence of homotopy group
presheaves $$ \ldots \pi _i ^{\rm pre}(Q,q)\rightarrow \pi _i ^{\rm
pre}(S,s)\rightarrow \pi _i ^{\rm pre}(T,t)\rightarrow \pi _{i-1} ^{\rm
pre}(Q,q)\rightarrow \ldots $$
for any $q\in Q_X$ and $s,t$ the images in $S$ and $T$. Sheafifying we get
a long exact sequence of homotopy group sheaves $$
\ldots \pi _i (Q,q)\rightarrow \pi _i (S,s)\rightarrow \pi _i
(T,t)\rightarrow \pi _{i-1} (Q,q)\rightarrow \ldots $$
as desired.
\eop

\begin{lemma}
\mylabel{pi09}
Suppose $(S,s ,\eta )$ and $(T,t ,\zeta )$ are Eilenberg-MacLane sheaves
$K(\Gg ,n)$ and $K(\Gg ', n)$. Suppose we have a morphism $\psi : \Gg
\rightarrow \Gg '$. Then the space
$$
Mor _{\Xx }^{\psi}((S,s , \eta ),(T,t ,
\zeta ))
$$
of pointed morphisms compatible with $\psi $ via the isomorphisms
$\eta , \zeta $ is a nonempty contractible space. \end{lemma}
{\em Proof:}
We first prove this in the case where $S$ is the standard Eilenberg-MacLane
sheaf constructed above. By Theorem \ref{fibr09}, it suffices to prove that
the
space of pointed maps from $(S',s')=(K^{\rm pre}(\Gg , n),o)$ to $(T,t)$ is
nonempty and contractible. But since the values $T_X$ are $n$-truncated,
and the values $S'_Y$ are obtained by attaching cells of dimension $\geq
n$, the spaces of choices of higher homotopies for the morphism
$S'\rightarrow T$ are weakly contractible, and the space of choices in any
case is $0$-truncated. Thus, $Mor _{\Xx }^{\psi}((S,s , \eta ),(T,t , \zeta
))$ is $0$-truncated. Its homotopy group $\pi _0$ is equal to the set of
choices of morphisms $f_{\varphi}: S'_Y\rightarrow T_X$ for each $\varphi :
X\rightarrow Y$ in $\Xx$, subject to the conditions $$
f(\alpha \beta )= f(\beta )S'(\alpha ) = T(\beta )f(\alpha ). $$
The set of choices of a morphism
$S'_Y\rightarrow T_X$ is the same as the set of group homomorphisms $\Gg
_X\rightarrow \Gg '_Y$.
The set of choices of $f$ becomes equal to the set of choices of
$f(1_X):\Gg _X\rightarrow \Gg '_X$ such that $$
f(1_X)S'(\varphi )= T(\varphi )f(1_Y)
$$
for $\varphi : X\rightarrow Y$; this is just $Hom _{\Xx}(\Gg , \Gg ')$.
This completes the case where $S$ is a standard Eilenberg-MacLane sheaf.

Apply this to the case where $\Gg = \Gg '$. This shows that for any
Eilenberg-MacLane sheaf $(T,t)$ of the form $K(\Gg , n)$ there is a
morphism from the standard $(S,s)$ to $(T,t)$ inducing an isomorphism on
the unique nontrivial homotopy group sheaf. Complete this to a fibration
diagram $Q\rightarrow S\rightarrow T$ with $Q$ being the fiber over $t$. By
the long exact sequence for homotopy group sheaves, the homotopy group
sheaves of $Q$ are trivial. By Corollary \ref{piA03}, this implies that $Q$
is weakly contractible. As $S$ and $T$ are of CW-type, this implies that
the morphism $(S,s)\rightarrow (T,t)$ is an equivalence. This shows that
any Eilenberg-MacLane sheaf is equivalent to a standard one. Our argument
above thus applies to any $(S,s)$. \eop

In the situation of the above
lemma, the functor
$$
Y\mapsto Mor ^{\psi}_{\Xx /Y}((S,s,\eta )|_{\Xx /Y},
(T,t, \zeta )|_{\Xx /Y})
$$
is a weak sheaf denoted $Map ^{\psi}((S,s,\eta ),
(T,t, \zeta ))$. By applying the above lemma to each $\Xx /Y$ we find that
it is weakly equivalent to $\ast$.

{\em Proof of Theorem \ref{bigPi03}}
We now turn to the situation of Theorem \ref{bigPi03}. Suppose $(T/S,
\sigma , \eta )$ and $(T'/S, \sigma ', \eta ')$ are respectively $K(G/S,
n)$ and $K(G'/S, n)$. Suppose $\psi : G\rightarrow G'$ is a morphism of
local systems on $S$. Then we can define a relative mapping space
$$
Map ^{\psi} ((T/S, \sigma ', \eta ), (T'/S, \sigma ', \eta ')), $$
with the same properties as in the standard case developed in section ???.
But the fiber of this over any point $t\in S_X$ is the mapping sheaf $Map
^{\psi}((T_s, \sigma (s), \eta ), (T'_s, \sigma '_s, \eta ')$ which we have
seen to be equivalent to $\ast$. This implies that the morphism $$
Map ^{\psi} ((T/S, \sigma ', \eta ), (T'/S, \sigma ', \eta ')), \rightarrow S
$$
is a weak equivalence, and by the defining property of the relative mapping
space, (applied to $S$) this implies that the space of morphisms over $S$
$$
Mor ^{\psi}_S((T, \sigma ,\eta ), (T', \sigma ', \eta ')) $$
is weakly contractible. If $G=G'$ and $\psi $ is the identity, the same
argument goes through to prove
that the space of equivalences over $S$, respecting base point and $\eta$,
is weakly contractible. This gives the second statement of the theorem. To
complete the proof, we just have to note that there is a functorial
construction of a relative Eilenberg-MacLane space over each $S_X$, for any
local system $G_X$ over $S_X$. Piece these together to get the desired
Eilenberg-MacLane presheaf, then sheafify. \eop

Armed with this existence and uniqueness result, we can speak of ``the''
Eilenberg-MacLane space $K(G/S, n)$. Denote by $o$ the basepoint, and by
$\eta$ the isomorphism of homotopy group sheaves.

\subnumero{The abelian category of local systems on $S$} If $S$ is a
truncated weak flexible sheaf of CW type, then the category of pre-local
systems of abelian groups
(that is, local systems in the trivial topology) over $S$ forms an abelian
category in an obvious way. Here we take $\pi _0$ of the morphism spaces as
the morphism sets. The kernels and cokernels are obtained just by taking
kernel and cokernel over each object.

\begin{lemma}
The category of local systems of abelian groups over $S$ forms an abelian
category, and the functor ``sheafification'' from the category of pre-local
systems to the category of local systems is exact. In particular, the
kernel and cokernel of a morphism of local systems are obtained by
sheafifying the pre-local system kernels and cokernels. \end{lemma}
{\em Proof:}
???
\eop

\begin{lemma}
(Suppose the site $\Xx$ has enough points.) The category of local systems
of abelian groups over $S$ has enough injectives. \end{lemma}
{\em Proof:}
???
\eop

Because of this lemma, we obtain derived functors $Ext ^i(L, M)$ of the
functor $L,M\mapsto Hom (L,M)$, with a long exact sequence.

\subnumero{Homology}

Suppose $T\rightarrow S$ is a morphism of truncated weak flexible sheaves
of CW type. We define the {\em homology local systems} $H_i (T/S)$ on $S$
as follows. For each $X$ let
$H^{\rm pre}_i(T/S)_X\rightarrow X$ be the local system of homology groups
of the homotopy fibers of $T_X\rightarrow S_X$. The structure of morphism
of flexible functors on $T$ gives the information necessary to give this
the structure of flexible functor. Then let $H_i (T/S)$ be the weak
flexible sheaf associated to $H^{\rm pre}_i(T/S)$.

Note that if $f:S'\rightarrow S$ is a morphism of truncated weak flexible
sheaves of CW type then
$$
H_i(T\times _SS'/S')=f^{\ast}H_i (T/S).
$$
If $p:S'\rightarrow T$ is a morphism then there is an obvious map of local
systems on $S'$,
$$
Hur_i: \pi _i (T\times _SS'/S', p)\rightarrow H_i (T\times _SS'/S'), $$
which we call the ``Hurewicz map''.
It is compatible with
the previously mentioned base changes.

We can also define homology with local coefficients. If $L$ is a local
system of abelian groups on $S$ and if $S\rightarrow T$ is a morphism then
we obtain local systems $H^i(S/T, L)$ over $T$. These are obtained by
taking the pre-local systems of homology of the fibers with coefficients in
$L$, and sheafifying.
Note that if $L'$ is any pre-local system we can define a pre-local system
of homology $H^{\rm pre}(S/T, L')$.

\begin{lemma}
Suppose $S\rightarrow T$ is a morphism of
truncated weak flexible sheaves of CW type. Suppose $L'$ is a pre-local system
on $S$ and let $L'\rightarrow  L$ be the morphism to its sheafification. Then
the induced morphism $$ H^{\rm pre}_(S/T, L') \rightarrow H^{\rm pre}(S/T ,
L) $$
gives an isomorphism of sheafifications. \end{lemma}
{\em Proof:}
???
\eop

\begin{corollary}
Let $\underline{\zz}$ denote the sheafification of the constant local
system $\zz$. Then there is a natural isomorphism $$
H_i(S/T) \cong H_i(S/T, \underline{\zz}). $$
\end{corollary}
{\em Proof:}
Note that $H_i^{\rm pre}(S/T) = H^{\rm pre}_i(S/T, \zz )$. This induces an
isomorphism between $H_i(S/T)$ and the sheafification of $H^{\rm
pre}_i(S/T, \zz )$. By the lemma, this sheafification is isomorphic to
$H_i(S/T, \underline{\zz})$.
\eop

We have a Serre spectral sequence for homology. Suppose $R\rightarrow
S\rightarrow T$ are morphisms of truncated weak flexible sheaves of CW
type, and $L$ is a local system on $R$.
The Serre spectral sequence over each object $X\in \Xx$ is functorial and
independent of homotopy, so it gives a spectral sequence of pre-local
systems on $T$,
$$
E_2^{i,j}= H^{\rm pre}_i(S/T, H^{\rm pre}_j(R/S ,L))\Rightarrow H^{\rm
pre}_{i+j}(R/T, L).
$$
The sheafification of a spectral sequence remains a spectral sequence (by
Lemma ???,
exact sequences of pre-local systems are turned into exact sequences of
local systems, under sheafification).
The $E_2^{i,j}$-term may be calculated using the previous lemma: first
taking the sheafification of the local system $H^{\rm pre}_j(R/S,L)$ and
then taking the sheafified homology, gives the same thing as sheafifying
the whole term at once. Thus we get a spectral sequence $$
E_2^{i,j}= H_i(S/T, H_j(R/S ,L))\Rightarrow H_{i+j}(R/T, L).
$$
Note that, by the previous corollary, putting $L=\underline{\zz}$ gives a
spectral sequence
$$
E_2^{i,j}= H_i(S/T, H_j(R/S ))\Rightarrow H_{i+j}(R/T).
$$

Suppose $T\rightarrow S$ is a morphism of truncated weak flexible sheaves
of CW type, and suppose that $\pi _0(T/S)$ is trivial (equal to $S$), and
for any morphism $p:S'
\rightarrow T$ and any
$i\leq n$, the relative $\pi _i (T\times _SS'/S', p)$ is trivial (equal to
$S'$). Then we say that $T\rightarrow S$ is {\em relatively $n$-connected.}

\begin{lemma}
Suppose $f:S\rightarrow T$ is a relatively $1$-connected morphism of
truncated weak flexible sheaves of CW type. Then any local system on $S$ is
pulled back from $T$, or more precisely for any local system $L$ on $S$, if
we put $M:= H_0(S/T, L)$ we have a natural isomorphism $L\cong f^{\ast}M$.
Furthermore, $M\cong H_0(S/T,f^{\ast}M)$. \end{lemma}
{\em Proof:}
???
\eop

We have the following ``Hurewicz theorem''.

\begin{lemma}
Suppose $T\rightarrow S$ is a morphism of truncated weak flexible sheaves
of CW type, which is relatively $1$-connected. Suppose $H_i(T/S)=0$ for
$i<n$. Then $T\rightarrow S$ is relatively $n-1$-connected, and for any
morphism $p:S'\rightarrow T$ the Hurewicz map is an isomorphism $$
\pi _n (T\times _SS'/S', p)\cong H_n (T\times _SS'/S'), $$
\end{lemma}
{\em Proof:}
By an inductive argument, we may assume that $T\rightarrow S$ is relatively
$n-1$-connected, and it suffices to show that the Hurewicz map is an
isomorphism.
We may also suppose that $S=S'=\ast$, since a map of local systems which is
an isomorphism over each point is an isomorphism. Then $p$ becomes a point
in $\Gamma (\Xx , T)$. Let $q$ denote the image of $p$ in $\Gamma (\Xx ,
\tau ^{\rm pre}_{\leq n-1}(T)$, and let $T'$ be the fiber of the map
$T\rightarrow \tau ^{\rm pre}_{\leq n-1}(T)$ over $q$. The hypothesis that
$T$ is $n-1$-connected means that the morphism $T'\rightarrow T$ induces an
isomorphism on associated weak flexible sheaves. Therefore the horizontal
maps in the diagram $$
\begin{array}{ccc}
\pi _n ^{\rm pre}(T', p)&\rightarrow &\pi ^{\rm pre}_n(T,p)\\ \downarrow
&&\downarrow \\
H_n^{\rm pre}(T')&\rightarrow &H_n^{\rm pre}(T) \end{array}
$$
induce isomorphisms on associated sheaves. But for every $X\in \Xx$, $T'_X$
is $n-1$-connected, so the vertical map on the left is an isomorphism. This
implies that the vertical map on the right induces an isomorphism on
associated sheaves, as desired. \eop

These homology local systems behave functorially with respect to morphisms
$S'/T'\rightarrow S/T$: given such a morphism (which is really a diagram,
that is to say a flexible functor $\Xx \times I^2\rightarrow Top$) we get a
morphism $H_i (S'/T') \rightarrow f^{\ast}H_i (S/T)$, where $f$ is the
morphism from $T'$ to $T$.
For fixed $f$, this is invariant under homotopy of the rest of the diagram.
Under a composition of morphisms, this induced morphism for the composition
is equal to the composition of the induced morphisms for the components.
All of the constructions done so far are compatible with this induced
morphism on homology. There are similar statements for homology with
coefficients in a local system, and the Serre spectral sequence is
compatible with these pullbacks.

\subnumero{Unpointed Eilenberg-MacLane spaces}

Suppose $n\geq 2$ and
suppose $L$ is a local system of abelian groups on $S$. Suppose $T\rightarrow S$
is a morphism, with $S$ and $T$ truncated weak flexible sheaves of CW type.
Suppose $\psi :H_n(T/S)\rightarrow L$ is a morphism of local systems. We
say that the pair $(T/S, \psi )$
is a {\em relative unpointed Eilenberg-MacLane space $K^{\ast}(L/S,n)$} if
$T\rightarrow S$ is relatively $n$-truncated and relatively
$n-1$-connected, and if $\psi$ is an isomorphism.

Note that in the above case,
if $\sigma:S\rightarrow T$ is a section, then $(T/S, \sigma , \psi \circ
Hur )$ is a relative Eilenberg-MacLane space $K(L/S, n)$.

If $(T/S, \psi )$ is a relative unpointed Eilenberg-MacLane space
$K^{\ast}(L/S, n)$ and if $S'\rightarrow S$ is a morphism of truncated weak
flexible sheaves of CW type, then the induced map $\psi ':H_n (T\times
_SS'/S')\rightarrow L\times _SS'$ gives $(T\times _SS'/S', \psi ')$ a
structure of relative unpointed Eilenberg-MacLane space $K^{\ast}(L\times
_SS'/S', n)$.

We give a classification (analogue of the standard thing in topology) of
relative unpointed Eilenberg-MacLane spaces. Fix $n\geq 2$.
Suppose $S$ is a truncated weak flexible sheaf of CW type, and suppose $L$
is a local system of abelian groups on $S$. Then the morphism $$
S\stackrel{o}{\rightarrow} K(L/S , n+1)=:B $$
is a relative unpointed Eilenberg-MacLane space $K^{\ast}(L\times _SB/B, n)$.
The next lemma essentially says that relative unpointed $K^{\ast}(L/S,n)$
are classified by sections of $B/S$. This, of course, gives a similar
statement for spaces $K^{\ast}(L\times _SS'/S',n)$ using base changes.

\begin{lemma}
With this notation, suppose $T\rightarrow S$ is a relative unpointed
Eilenberg-MacLane space $K^{\ast}(L/S, n)$. Then the space of pairs $(f,g)$
where $f:S\rightarrow B= K(L/S, n+1)$ is a section and $g: S\times
_{o,B,f}S\sim T$ is an isomorphism over $S$ (the first being over $S$ via
the second projection), is weakly contractible. \end{lemma}
{\em Proof:}
???
\eop

Note that the space $T$ will have a section if and only if $f$ is homotopic
to the basepoint section $o$.

\subnumero{Postnikov towers}

By using the same definitions as in Section ???, we obtain notions of
towers of flexible functors, or of
(weak) flexible sheaves.

If $A$ is a tower of flexible functors from $\Xx$ to $Top$ then $\Gamma
(\Xx , A)$ is a tower of spaces.
\begin{lemma}
\mylabel{pi02} A tower of flexible functors may be considered as a flexible
functor on $\Xx \times \Nn$. We obtain $\Gamma (\Xx , A)$ on $\Nn$ and
$\Gamma (\Nn , A)$ on $\Xx$. These commute up to homotopy in the sense that
$$
\Gamma (\Xx , \Gamma (\Nn , A))\sim \Gamma (\Nn , \Gamma (\Xx , A))\sim
\Gamma (\Xx \times \Nn , A).
$$
\end{lemma}
{\em Proof:}
Left to the reader.
\eop

Recall that we have a truncation operation for a weak flexible sheaf of CW
type $S$,
yielding a weak flexible sheaf of CW type $\tau _{\leq n}(S)$ with a
morphism from $S$; and a relative truncation for a morphism $S\rightarrow
T$ yielding a weak flexible sheaf over $T$, denoted
$\tau _{\leq n}(S/T)\rightarrow T$,
with a morphism from $S$ relative to $T$.

Note that the absolute case is a special case of the relative case, where
$T=\ast$.

In what follows, we concentrate on these truncation operations for sheaves;
there are similar but easier statements for the operations $\tau _{\leq
n}^{\rm pre}$.

We have
morphisms $\tau _{\leq n}(S/T)
\rightarrow \tau _{\leq n-1}(S/T)$, so we obtain the {\em sheafified
Postnikov tower} $\tau (S/T)$ for $S$ over $T$. This is horizontally
relatively truncated and if $S$ and $T$ were
truncated to begin with then the sheafified Postnikov tower is completely
truncated.

\begin{lemma}
\mylabel{pi13}
Suppose $S$ is a truncated
weak flexible sheaf over another one $T$, and suppose $\sigma :T\rightarrow
S$ is a section. We obtain a section
also denoted by $\sigma $ of $\tau _{\leq n}(S/T)$. The natural morphism
induces isomorphisms of homotopy local systems $$
\pi _i
(S/T,\sigma )\cong \pi _i
(\tau _{\leq n}(S/T)/T,\sigma )
$$
for $i\leq n$.
Note that the homotopy local systems on the right are zero for $i>n$.
\end{lemma}
{\em Proof:}
The natural morphism induces isomorphisms of homotopy group pre-coverings
(the obvious modification of the definition of covering to allow things
which are not sheaves)
$\pi ^{pre}_i (S/T,\sigma )\cong \pi ^{pre}_i (\tau ^{pre}_{\leq
n}(S/T)/T,\sigma )$ for $i\leq n$, and hence of the associated
sheaves of groups.
As the homotopy local systems are preserved by the operation of taking
associated flexible sheaf (this is the local system version of the
statement that homotopy group sheaves are preserved by taking associated
flexible sheaves), we obtain the result. \eop

It follows from this lemma that
the morphism $\tau _{\leq n}(S/T)\rightarrow \tau _{\leq n-1}(S/T)$ is
relatively $n$-truncated and relatively $n-1$-connected.
Suppose $n\geq 2$.
Let
$$
Q_i := H_i (\tau _{\leq n}(S/T)/\tau _{\leq n-1}(S/T)) $$
be the relative homology sheaf.
The morphism from the $n$th truncation to the $n-1$st truncation is a
relative unpointed Eilenberg-MacLane
space of the form $K^{\ast}(Q_i/\tau _{\leq n-1}(S/T), n)$. Note that by
Lemma ???, the local system $Q_i$ is the pullback of a local system $P_i$
on $\tau _{\leq 1}(S/T)$. Also, by Lemma ???, the morphism $\tau _{\leq
n}(S/T)\rightarrow \tau _{\leq n-1}(S/T)$ is classified by a map
$$
\tau _{\leq n-1}(S/T) \rightarrow K(P_i / \tau _{\leq 1}(S/T), n+1) $$
over $\tau _{\leq 1}(S/T)$. In the terminology below, this is the
cohomology invariant of the $n$th stage in the Postnikov tower.

In general, if $T\rightarrow S$ is a relatively $1$-connected morphism then
we define $\pi _n(T/S)$ to be the local system on $S$ whose pullback to
$\tau _{\leq n-1}(T/S)$ is the homology
$$
H_n(\tau _{\leq n}(T/S)/\tau
_{\leq n-1}(T/S)).
$$
If $\sigma :S\rightarrow T$ is any section, then $\pi
_n(T/S, \sigma )$ is canonically isomorphic to $\pi _n(T/S)$ (but this
latter is defined even if no such section exists).

If $R\rightarrow T\rightarrow S$ is a pair of relatively $1$-connected
morphisms, then we write $\pi _n(R/T/S)$ for the local system on $S$ whose
pullback to $T$ is $\pi _n (R/T)$.

\begin{lemma}
Suppose $\{ T_n /S , T_n/S\rightarrow T_{n-1}/S\}$ is a vertically
truncated tower of truncated weak
flexible sheaves of CW type which are relatively $1$-connected over a base $S$.
Then the limit $T_{\infty}$ is truncated and relatively $1$-connected over
$S$, and
there is a spectral sequence
of local systems on $S$,
$$
E_2^{p,q}= \pi _{-p-q}( T_{-q}/T_{-q-1}/S)\Rightarrow \pi _{-p-q}(T_{\infty}/S).
$$
\end{lemma}
{\em Proof:}
???
\eop

Of course, this gives nothing new when applied to the Postnikov tower itself.

\subnumero{Cohomology}
Suppose $S\rightarrow T$ is a morphism of truncated weak flexible sheaves
of CW type. Suppose $L$ is a local system of abelian groups on $S$. For any
$n$ we obtain the relative Eilenberg-MacLane space $K(L/S , n) \rightarrow
S$, with base
section $o$. We have natural isomorphisms $\Omega ^o(K(L/S, n)/S)\cong
K(L/S, n-1)$.
Define the {\em cohomology local systems} $$
H^i(S/T, L):= \pi _0(Sect (S/T, K(L/S , i)/T), o) $$
over $T$. These have structures of local systems of abelian groups given by
the following lemma.

\begin{lemma}
For any $n$ there is an isomorphism
$$
H^i(S/T , L)\cong \pi _n (Sect (S/T , K(L/S , i+n)/T), o) $$
of local systems on $T$, induced by the isomorphisms $$
K(L/S, i)\cong (\Omega ^o)^n(K(L/S , i+n)/S); $$
and for various $n$ these induce the same structure of abelian group on the
cohomology local system. \end{lemma}
{\em Proof:}
???
\eop

Suppose
$$
\begin{array}{ccc}
S'&\stackrel{g}{\rightarrow} &S \\
\downarrow &&\downarrow \\
T' &\stackrel{f}{\rightarrow} & T
\end{array}
$$
is a diagram, and suppose
$L$ and $L'$ local systems on $S$ and $S'$ respectively. Suppose $\phi :
g^{\ast}L\rightarrow L'$ is a morphism of local systems. Then $\phi$
induces a morphism
$$
f^{\ast}H^i(S/T, L) \rightarrow H^i(S'/T', L'). $$
These morphisms satisfy a functoriality when there is a composable diagram
of diagrams; and they depend only on the homotopy class of the diagram
above a fixed $f$.

(construction---proof ???)

\begin{theorem}
{\rm (Universal coefficients theorem)}
(Suppose that $\Xx$ has enough points.)
Suppose $f:S\rightarrow T$ is a morphism, and suppose $L$ is a local system
on $T$. Then there is a spectral sequence of local systems on $T$, $$
E_2^{i,j}= Ext ^i (H_j (S/T), L)\Rightarrow H^{i+j}(S/T, f^{\ast}L). $$
This is compatible with the morphisms of functoriality established above.
\end{theorem}
{\em Proof:}
???
\eop

\subnumero{The Leray-Serre spectral sequence}

\begin{theorem}
Suppose $R\rightarrow S\rightarrow T$ is a composable pair of morphisms of
truncated weak flexible sheaves of CW type. Suppose $L$ is a local system
on $R$. Then there is a spectral sequence (called {\em Leray-Serre})
$$
E_2^{i,j}(R/S/T, L) = H^i (S/T, H^j(R/S , L))\Rightarrow H^{i+j}(R/T, L). $$
\end{theorem}
{\em Proof:}
???
\eop

We need to know that the Leray-Serre spectral sequence is compatible with
morphisms. Suppose
$$
\begin{array}{ccc}
R'&\stackrel{g}{\rightarrow }&R\\
\downarrow &&\downarrow \\
S'&\rightarrow &S \\
\downarrow &&\downarrow \\
T' &\stackrel{f}{\rightarrow} & T
\end{array}
$$
is a morphism of collections of data which go into the Leray-Serre spectral
sequence, and suppose $L$ and $L'$ are local systems of abelian groups on
$R$ and $R'$
respectively. Suppose $\phi : g^{\ast}L\rightarrow L'$ is a morphism of
local systems.

\begin{lemma}
In the above situation, there is a morphism of spectral sequences $$
f^{\ast}E^{i,j}_r(R/S/T, L)\rightarrow
E^{i,j}_r(R'/S'/T', L')
$$
which is the morphism induced by $\phi$ on the $E_2$ and the $E_{\infty}$ terms.
\end{lemma}
{\em Proof:}
???
\eop

\subnumero{Calculation of the space of global sections}

So far, everything we have done has stayed inside the realm of flexible
sheaves on $\Xx$. We now discuss what happens when we take global sections
to recover some information of topological spaces. Recall that the space of
morphisms $Mor (S,T)$ is equal to the global sections of the mapping space
$Map (S,T)$. Similarly, the space of sections of a morphism $T\rightarrow
S$ is the global sections of the section space $Sect (S, T)$. Thus it
suffices to treat the case of taking global sections of one particular
flexible sheaf. (Although there might be some other nicer statements in the
other cases...???.) We restrict to the absolute case, and speak of
Eilenberg-MacLane spaces instead of relative ones, for example.

\begin{corollary}
The fibers $F_i(\eta )$ in the sheafified Postnikov tower for $X$ are
locally Eilenberg MacLane sheaves locally of the form $K(\Gg _i , i)$ where
$\Gg_i$ are the homotopy group sheaves of $T$ based at a point $\eta \in
\Gamma (\Xx , \tau _{\leq i-1}T)$.
\end{corollary}
{\em Proof:}
Where basepoints exist, the long exact sequence for homotopy group sheaves,
and the previous lemma, imply that the homotopy group sheaves of the fibers
$F_i(\eta )$
vanish, except in degree $i$.
\eop

We obtain a tower of spaces
$$
P_i(\Xx , T):= \Gamma (\Xx , \tau _{\leq i}T), $$
with fibers $F_i (\Xx , T; \eta )=\Gamma (\Xx , F_i(\eta ))$. If $T$ is
truncated then this tower is completely truncated, and the limit is $\Gamma
(\Xx , T)$.
We obtain a spectral sequence
$$
E_2^{p,q}= \pi _{-p-q}(\Gamma (\Xx , F_{-q}(\eta )),\eta ) \Rightarrow \pi
_{-p-q}(\Gamma (\Xx , T),\eta ). $$
We use the
following lemmas to calculate the $E_2$ term.

{\em Remark:} In the Postnikov tower, the fibers are {\em locally} of the
form $K(\Gg , i)$. Thus there is no guarantee that the space $\Gamma (\Xx ,
F_i(\eta ))$ is nonempty. But if it is nonempty, then choosing a point will
give a basepoint $o$ for the Eilenberg-MacLane sheaf, so in this case
we can assume that the fiber is $K(\Gg , i)$ (note however that the sheaf
of groups $\Gg$ can change with choice of different connected component).

{\em Remark:} By taking the Postnikov tower of flexible functors (not
sheafifying), we obtain a spectral sequence which we call the {\em \v{C}ech
spectral sequence},
$$
E_2^{p,q}= \pi _{-p-q}(\Gamma (\Xx , F_{-q}^{\rm pre}(\eta )),\eta )
\Rightarrow \pi _{-p-q}(\Gamma (\Xx , T),\eta ). $$
Here $F_{-q}^{\rm pre}(\eta )$ denotes the fiber of the non-sheafified
Postnikov tower. It is a flexible functor whose value on $X\in \Xx$ is
$K(\pi _i(T_X, \eta _X ), i)$. We write this as $$
F_{i}^{\rm pre}(\eta )= K^{\rm pre}(\pi _i(T_X, \eta _X ), i). $$

\begin{lemma}
\mylabel{pi03}
Suppose $\Gg$ is a presheaf of groups on a site $\Xx$ and $F= K^{\rm
pre}(\Gg , i)$. Let $o$ be the basepoint section of $F$ over $\Xx$. Suppose
$\Bb$ is a sieve of $\Xx /X$ corresponding to a covering $\Uu$. Then $\pi
_j ( \Gamma (\Bb , F), o)$ is equal to the \v{C}ech cohomology group
$\check{H}^{i-j}(\Uu , \Gg ) $ on the category $\Xx$ (this can be
nonabelian cohomology if $i=1$, and is a pointed set for $i=0$).
\end{lemma}
{\em Proof:}
The description of $\Gamma (\Bb
, F)$ in terms of the covering family $\Uu$ (given toward the end of the
section on descent) immediately gives the \v{C}ech complex for the covering
$\Uu$ as calculating the homotopy groups.
\eop

{\em Remark:} The \v{C}ech cohomology group is the same as the right
derived functor of the global section functor on the category of
presheaves. We could have stated, in the lemma, that for any category $\Xx$
we have $$
\pi _j ( \Gamma (\Xx
, F), o)= R^{i-j}\Gamma (\Xx , \Gg ).
$$

\begin{lemma}
\mylabel{pi04}
Suppose $\Gg$ is a sheaf of groups which is \v{C}ech-acyclic (for sieves
and for global sections) in dimensions $\leq n$. Then the homotopy functor
$T$ defined by $T_X=K(\Gg (X), n)$ is a homotopy sheaf. In particular, $T$
is an Eilenberg-MacLane homotopy sheaf, and we get $$
K(\Gg , n)_X= K(\Gg (X), n).
$$
Furthermore, we have
$$
\Gamma (\Xx , K(\Gg , n))= K(\Gamma (\Xx , \Gg ), n). $$
\end{lemma}
{\em Proof:}
Suppose $X\in \Xx$ and $\Bb $ is a sieve over $X$ given by a covering
family $\Uu$. We apply the \v{C}ech spectral sequence to calculate the
homotopy groups of $\Gamma (\Bb , T|_{\Bb})$. All the terms except one
vanish, giving $$
\Gamma (\Bb , T|_{\Bb})= K(\Gamma (\Bb , \Gg ), n). $$
The sheaf condition for $\Gg$ insures that $\Gamma (\Bb , \Gg )= \Gg (X)$,
so the map
$$
T_X\rightarrow \Gamma (\Bb , T|_{\Bb})
$$
is a homotopy equivalence. The same argument works for global sections. \eop

\begin{lemma}
\mylabel{pi05}
If $\Ii$ is an injective sheaf on a site $\Xx$ then $\Ii$ is
\v{C}ech-acyclic (for sieves and for global sections).
Suppose $\{ \Ii ^{j}\} _{j=0,\ldots , n}$ is a complex of sheaves of groups
of length $n$ on a site $\Xx$. Suppose that it is exact except at the first
place, and suppose that $\Ii ^j$ are injective for $j<n$. Then $\Ii ^n$ is
\v{C}ech-acyclic (for sieves and for global sections) and acyclic for sheaf
cohomology, in degrees $\leq n$.
\end{lemma}
\eop

The following result is basically K. Brown's result of \cite{BrownThesis},
translated into our ``homotopy coherent'' situation.

\begin{lemma}
\mylabel{pi06}
{\rm (K. Brown \cite{BrownThesis})}
Suppose $\Gg$ is a sheaf of groups on a site $\Xx$ and $F= K(\Gg , n)$. Let
$o$ be the basepoint section of $F$ over $\Xx$. Then $\pi _j ( \Gamma (\Xx
, F), o)$ is zero for $j>n$ and for $j\leq n$ is equal to the sheaf
cohomology group $H_{\rm sh}^{n-j}(\Xx , \Gg ) $ on the site $\Xx$ (this
can be nonabelian cohomology if $n=1$, and is a pointed set for $n=0$).
\end{lemma}
{\em Proof:}
First we treat the case $n\geq 2$ where $\Gg$ is abelian. Then $\Gg$ admits
a possibly infinite injective resolution in the category of sheaves on
$\Xx$. By the previous lemma, we can truncate this to obtain a finite
resolution by sheaves which are \v{C}ech-acyclic and acyclic for sheaf
cohomology in degrees $\leq n$. Proceed by induction on the length of this
resolution (note that the inductive argument given next works also to start
the induction). By taking the first stage of the resolution and splitting
off an exact sequence, we obtain a short exact sequence $$
0\rightarrow \Gg \rightarrow \Ii \rightarrow \Ff \rightarrow 0 $$
where $\Ii$ is \v{C}ech-acyclic and acyclic for sheaf cohomology in degrees
$\leq n$, and where $\Ff$ admits a resolution of length shorter than that
of $\Gg$. By induction, we may suppose that the lemma is true for $\Ff$.
The lemma is true for $\Ii$ because of Lemma \ref{pi04} (giving the
vanishing of most of the homotopy groups) and since the same cohomology
groups vanish because $\Ii$ is acyclic. We obtain a fibration $$
K(\Gg , n)\rightarrow K(\Ii , n) \rightarrow K(\Rr , n), $$
yielding a long exact sequence of homotopy groups $$
\ldots \rightarrow \pi _j ( \Gamma (\Xx , K(\Gg , n)), o)\rightarrow $$
$$
\pi _j ( \Gamma (\Xx , K(\Ii , n)), o)
\rightarrow \pi _j ( \Gamma (\Xx , K(\Ff , n)), o)\rightarrow \ldots \; . $$
On the other hand, the exact sequence of sheaves gives the long exact
sequence of sheaf cohomology groups
$$
\ldots \rightarrow H_{\rm sh}^{n-j}(\Xx , \Gg ) \rightarrow $$
$$
H_{\rm sh}^{n-j}(\Xx , \Ii )
\rightarrow H_{\rm sh}^{n-j}(\Xx , \Ff ) \rightarrow \ldots \; . $$
By the inductive hypothesis and the acyclicity of $\Ii$ these exact
sequences coincide at all places involving $\Ff$ and $\Ii$ (the isomorphism
we are constructing inductively will be compatible with morphisms---this is
shown below). Comparing, we get the desired isomorphisms for $\Gg$. The
acyclicity of $\Ii$ means that the exact sequences give isomorphisms $$
H_{\rm sh}^{n-j}(\Xx , \Ff )\cong H_{\rm sh}^{n+1-j}(\Xx , \Gg ) $$ (and
similarly for the homotopy groups) when $n-j\geq 1$. Furthermore, we have a
map $$ H_{\rm sh}^{n-j}(\Xx , \Ii )
\rightarrow H_{\rm sh}^{n-j}(\Xx , \Ff ), $$
of which $H_{\rm sh}^{0}(\Xx , \Gg )$ is the kernel and $H_{\rm sh}^{1}(\Xx
, \Gg )$ the cokernel. Again, the same goes for the homotopy groups. Hence
the isomorphisms between the homotopy groups and the cohomology groups are
canonically determined by the isomorphisms for $\Ii$ and $\Ff$.

We have to check that our isomorphisms are compatible with morphisms of
sheaves of groups $\Gg$ (this was used in comparing the above exact
sequences). Given a morphism $\Gg \rightarrow \Gg '$
we may extend it to a morphism between injective
resolutions (and hence between their truncations used above). We get morphisms
$\Ii \rightarrow \Ii '$ and $\Ff \rightarrow \Ff '$ between the exact sequences
used above. By induction, we may suppose that $\Ff\rightarrow \Ff '$ induces the
same morphisms on homotopy and cohomology groups. In order to treat the first
case of the induction, we have to show that this is so for
$\Ii \rightarrow \Ii '$. But here the homotopy groups and cohomology groups are
both naturally equal to $\Gamma (\Xx , \Ii )$ or $\Gamma (\Xx , \Ii ')$, and the
maps are both the induced map here. This gives the statement for $\Ii
\rightarrow \Ii '$.

Finally we have to treat the cases $n=0$ and $n=1$. In the case $n=0$ we
may suppose that $F$ is equal to the sheaf of sets $\Gg$; then the
isomorphism is clear. Suppose $\Gg$ is a sheaf of nonabelian groups and
$n=1$. In terms of the operation $H$ considered in previous sections, we
have $$
K(\Gg , 1) = H \cdot K^{\rm pre}(\Gg ,1) $$
since the path spaces of $K^{\rm pre}(\Gg ,1)$ are of the form $K(\Gg , 0)$
which are already sheaves (see Theorem \ref{fibr07}). We may assume that
$\Xx$ embedds into a site with a final object $e$ and that $\Gg$ extends to
$\Gg '$ (and hence $F$ to $F'$); furthermore in this case $\Xx$ appears as
a sieve over $e$. Thus $$
\Gamma (\Xx , F)= F'_e= (H\cdot K^{\rm pre}(\Gg ',1))_e $$
$$
= \lim _{\rightarrow } \Gamma (\Bb , K^{\rm pre}(\Gg ',1)), $$
where the limit is taken over the sieves $\Bb$ over $e$. Thus $$
\pi _j ( \Gamma (\Xx , F), o)=
\lim _{\rightarrow } \pi _j (\Gamma (\Bb , K^{\rm pre}(\Gg ',1)), o). $$
But for sieves $\Bb^{\Uu}$ generated by coverings $\Uu$ (note that we may
suppose $\Bb \subset \Xx$) we have
$$
\pi _0 (\Gamma (\Bb ^{\Uu}, K^{\rm pre}(\Gg ',1)), o)= H^1(\Uu , \Gg ), $$
while
$$
\pi _1 (\Gamma (\Bb ^{\Uu}, K^{\rm pre}(\Gg ',1)), o)= H^0(\Uu , \Gg ). $$
The direct limits (over all $\Uu$) of these \v{C}ech cohomology groups are
the same as the sheaf cohomology groups (for $n=0$ and $n=1$ in the abelian
case, we define sheaf cohomology to be the same as \v{C}ech
cohomology---for want of a better definition). We get the desired
statements. \eop

\begin{corollary}
\mylabel{pi07}
Suppose $\zeta \in P_i(\Xx , T)$, and let $\eta \in \Gamma (\Xx , \tau
_{\leq i-1}T)$ be the image. Then $$
\pi _j (F_i(\Xx , T; \eta ), \zeta )= H^{i-j}_{\rm sh}(\Xx , \varpi
_i(T,\zeta )) .
$$
\end{corollary}
{\em Proof:}
This follows from the previous lemma because the fiber $F_i(\eta )$ with
basepoint section $\zeta $ is an Eilenberg-MacLane sheaf $K(\Gg , i)$ with
$\Gg = \varpi _i(T,\zeta )$.
\eop

\begin{corollary}
\mylabel{pi08}
The $E_2$ term in the spectral sequence for the tower $P_i(\Xx , T)$ at a
point $\zeta \in \Gamma (\Xx , T)$ is
$$
E_2^{p,q}=
H^{-p}_{\rm sh}(\Xx ,
\varpi _q(T,\zeta ))
\Rightarrow \pi _{p+q}(\Gamma (\Xx , T), \zeta ) $$
(in particular, it is concentrated in the second quadrant). \end{corollary}
\eop

\begin{corollary}
\mylabel{pi12}
Suppose $\Gg \rightarrow \Gg '$ is a morphism of groups (abelian for $n\geq
2$) with kernel $\Gg ''$. Then there exists a fibration triple of homotopy
sheaves $$
K(\Gg '' , n)\rightarrow K(\Gg , n)\rightarrow K(\Gg ' , n). $$
\end{corollary}
{\em Proof:}
There exists a morphism $K(\Gg , n)\rightarrow K(\Gg ' , n)$ by Theorem
\ref{pi09}; and by Lemma \ref{piA05} this can be completed to a fibration
triple. The long exact sequence of homotopy group sheaves implies that the
fiber is $K(\Gg '' , n)$.
(Note that if the total space in the fibration triple is pointed, then the
fiber is pointed.)
\eop

\subnumero{Homology and cohomology}

(CHUCK THIS ?????)

We can define the homology presheaves and sheaves of a flexible
functor in the same way as for homotopy. Put $H_i^{pre}(T)(X):= H_i(T_X,\zz )$.
This is a presheaf. Let $H_i(T)$ be the sheaf associated to the presheaf
$H^{pre}_i(T)$.

Suppose that $T$ is a functor. Then the singular chain complexes form a
functorial presheaf of complexes $C_{\cdot}(T)$ on $\Xx$. If $\Gg$ is a
sheaf of abelian groups on $\Xx$, choose an injective resolution
$\Gg\rightarrow \Ii ^{\cdot}$ in the category of sheaves over $\Xx$, and
put $$
H^i(T,\Gg )(X):= H^i(Hom _{\Xx /X}(C_{\cdot}(T),\Ii ^{\cdot }). $$
\begin{theorem}
\mylabel{pi14}
Suppose $\Gg$ is a sheaf of abelian groups on $\Xx$. Suppose $(K(\Gg ,
n),o)$ is a pointed flexible sheaf with homotopy groups $\Gg$ in degree $n$
and trivial otherwise. Then
$$
\pi _i (Mor (T, K(\Gg , n )), \underline{o})= H^{n-i}(T,\Gg ). $$
\end{theorem}
{\em Proof:}
???
\eop

\begin{corollary}
\mylabel{piA06}
There is a spectral sequence
$$
E^{i,j}_2=Ext ^i(H_j(S), \Gg )\Rightarrow \pi _{i+j+n}(Mor (S, K(\Gg ,
n)),\underline{o}).
$$
\end{corollary}
{\em Proof:}
By definition, cohomology is calculated by the double complex $Hom
(C_{\cdot}(S), \Ii ^{\cdot})$. Taking one of the spectral sequences of this
double complex, we get a spectral sequence $$
Ext ^i(H_j(S), \Gg )\Rightarrow H^{i+j}(S,\Gg ), $$
which gives the claimed one according to the theorem. \eop

We consider the following condition:
\newline
(Rat)\,\, For $i\geq 2$ and for any $X\in \Xx$, the homotopy groups $\pi _i
(S)$ are rational vector spaces.

We call a flexible presheaf satisfying this condition, {\em rational}.

Note that if $\Gg$ is a sheaf of rational vector spaces, then $K(\Gg , n)$
is rational.

We recall the calculation of the cohomology of a rational Eilenberg-MacLane
space.

\begin{lemma}
\mylabel{pi17}
Suppose $G$ is a rational vector space. Then if $n\geq 1$ we have $$
H_i(K(G,n))=\left\{
\begin{array}{ll}
Sym ^k_{\zz}(G) & i=kn, \;\; $n$\,\, \mbox{even} \\ \bigwedge ^k_{\zz}(G) &
i=kn, \;\; $n$\,\, \mbox{odd} \\ 0 & \mbox{otherwise} .
\end{array}
\right.
$$
\end{lemma}
{\em Proof:}
It is left to the reader to find the references. \eop

\begin{corollary}
\mylabel{pi18}
Suppose $\Gg$ is a sheaf of rational vector spaces. Then if $n\geq 1$ we have
$$
H_i(K(\Gg ,n))=\left\{
\begin{array}{ll}
Sym ^k_{\zz}(\Gg ) & i=kn, \;\; $n$\,\, \mbox{even} \\ \bigwedge
^k_{\zz}(\Gg ) & i=kn, \;\; $n$\,\, \mbox{odd} \\ 0 & \mbox{otherwise} .
\end{array}
\right.
$$
\end{corollary}
{\em Proof:}
The homology sheaves are the sheaves associated to the presheaves given by
the previous lemma; these are, by definition, the sheaves on the right.
\eop

Combining the above results, we obtain a spectral sequence converging to
the homotopy groups
$$
\pi _i(Mor(K(\Gg , n),K(\Gg ',m)),\underline{o}), $$
with the $E_2$ term consisting of groups of the form $$
Ext ^j(Sym ^k_{\zz}(\Gg ), \Gg ')
$$
or
$$
Ext ^j(\bigwedge ^k_{\zz}(\Gg ), \Gg '). $$
Below we will calculate these in the case where $\Xx$ is the category of
schemes over a field of characteristic zero, with the faithfully flat
finite type topology.

\newpage

\setcounter{section}{9}

\numero{Stacks}

In this section, we will make somewhat more precise the
statement that ``a stack of groupoids is the same thing as a $1$-truncated
flexible sheaf''.

Suppose that $\Xx$ is a site. Recall that a {\em fibered category over
$\Xx$} is a category $\Cc$ with a functor $a:\Cc \rightarrow \Xx$ such that
for each morphism $u:X_0\rightarrow X_1 $ in $\Xx$ and each object $Y\in
\Cc$ such that $a(Y)$ is the target $X_1$ of $u$, there exists a cartesian
lift of $u$ to a morphism
$v=v(u,Y)$ in $\Cc$---that is, a morphism $v$ with $a(v)=u$ and such that
any morphism $w$ with $a(w)=u$ factors uniquely as $$
w= v z
$$
for $z$ a morphism in $\Cc$ such that $a(z)$ is the
identity for the starting $X_0$ of $u$.

If we choose a lift $v(u,Y)$ for each pair $(u,Y)$ then we
obtain a functor $u^{\ast}: a^{-1}(X_1)\rightarrow a^{-1}(X_0)$, and a natural
isomorphism of $\xi _{u,v}: v^{\ast}u^{\ast} \cong (uv)^{\ast}$ satisfying a
pentagonal axiom. This gives a {\em pseudofunctor} from $\Xx$ to $Cat$ (cf
\cite{Jardine2}), and conversely given a pseudofunctor we obtain a fibered
category with choice of liftings.

A {\em functor of fibered categories} from
$a:\Cc \rightarrow \Xx$ to $a':\Cc '\rightarrow \Xx$ is a functor $b:\Cc
\rightarrow \Cc '$ such that $a'b=a$. A {\em natural transformation} between two
such functors is a natural transformation $\xi : b \rightarrow b'$ such that the
resulting natural endomorphism $a'\xi$ of $a$ is trivial. The set of fibered
categories over $\Xx$ is made into a strictly associative $2$-category with this
set of functors and natural transformations.

A {\em category fibered in groupoids over $\Xx$} is a
fibered category $a:\Cc \rightarrow \Xx$ such that the fibers are groupoids
(that is, their morphisms are all isomorphisms). A {\em stack of groupoids over
$\Xx$} is a category fibered in groupoids $a:\Cc \rightarrow \Xx$ which
satisfies a certain descent condition (see
\cite{GiraudThese} \cite{Artin} \cite{Deligne-Mumford}). We obtain strictly
associative $2$-categories of categories fibered in groupoids over $\Xx$, and of
stacks of groupoids over $\Xx$ (there are no extra conditions on the functors
and natural transformations).

\subnumero{Going from stacks to $1$-truncated weak sheaves}

Suppose that $a:\Cc \rightarrow \Xx$ is a category
fibered in groupoids, and suppose $v(u, Y)$ is a choice of cartesian lifting for
each $u, Y$. We will associate to this the {\em nerve} $N=N (\Cc /\Xx , v)$, a
flexible functor from $\Xx$ to $Top$, as follows. Let $N_X$ be the realization
of the nerve of the groupoid $a^{-1}(X)$. For $u:Y\rightarrow X$ let
$N(u):N_X\rightarrow N_Y$ be the morphism induced by the functor $u^{\ast}$. If
$u_1,\ldots , u_n$ is a composable sequence of morphisms, define $N(u_1, \ldots
, u_n, t_1, \ldots , t_{n-1})$ as follows ... ???

Note that the $N_X$ are $1$-truncated CW complexes---this
is generally true of the realization of the nerve of a groupoid. Hence the nerve
$N$ is a $1$-truncated flexible functor of CW type.

We leave the following verification
as an exercise for the reader.

\begin{lemma}
Suppose $\Cc /\Xx$ is a stack of groupoids, with choice of cartesian
liftings $v(u,Y)$. Then $\Cc / \Xx$ satisfies the descent
condition to be a stack of groupoids, if and only if the nerve $N(\Cc / \Xx ,
v)$ is a weak flexible sheaf. \end{lemma}
\eop

\subnumero{Going from $1$-truncated weak sheaves to stacks}

Suppose that $T$ is a $1$-truncated flexible functor of CW type over $\Xx$.
We will define a category fibered in groupoids $\Pp (T)\rightarrow \Xx$
with a distinguished choice of cartesian liftings $p(u,Y)$ ---this is a
generalization of the construction of the Poincare groupoid.

The category $\Pp (T)$ has one object for each pair $(X,s)$ where $X\in
\Xx$ and $s\in T_X$. The morphisms from $(X,s)$ to $(X',s')$ are the pairs
$(f, \gamma )$ where $f:X\rightarrow X'$ is a morphism in $\Xx$ and
$\gamma$ is a homotopy class of paths in $T_X$ from $s$ to $T(f)(s')$. The
composition of morphisms from $(X,s)$ to $(X',s')$ and then to $(X'' , s''
)$ is defined by
$$
(f',\gamma ,)(f, \gamma ) := (f'f, \gamma '' ) $$
where $\gamma ''$ is the homotopy class of paths obtained by starting with
$\gamma$ (from $s$ to $T(f)(s')$) then adding $T(f)(\gamma ')$ (from
$T(f)(s')$ to $T(f)T(f')(s'')$) and finally taking the path defined by
$T(f',f,t)(s'')$ as $t$ runs backward from $1$ to $0$ (which goes from
$T(f)T(f')(s'')$ to $T(f'f)(s'')$). The associativity of this composition
law is provided by the homotopies $T(f'', f', f, t_1,t_2)$.

\begin{lemma}
This $\Pp (T)$ is a category fibered in groupoids over $\Xx$ (via the
forgetful functor). The assignment $p(f,(X',s'))=(f,1)$ where $1$ denotes
the identity path from ??? to ??? is a cartesian section. The fibered
groupoid $\Pp (T)$ is a stack if and only if $T$ is a weak flexible sheaf.
\end{lemma}
{\em Proof:}
???
\eop

\subnumero{Comparison of these constructions}

These constructions $N(\Cc /\Xx , v)$ and $\Pp (T)$ are clearly not
inverses. However, they provide a sort of equivalence between the
$2$-category of stacks over $\Xx$
and the $\Delta$-category of $1$-truncated weak flexible sheaves of CW type
over $\Xx$.
We need to describe what this means.

A strictly associative $2$-category is a category with an additional
structure of set of natural transformations between the morphisms (we call
the morphisms functors) satisfying the same properties as satisfied in the
case of $Cat$, the set of categories.

\newpage

\setcounter{section}{10}

\numero{Classifying spaces}

Suppose $\Xx$ is a site.
Fix a universe $\Uu$ (basically, a cardinality bound). Then there is a
morphism of weak flexible sheaves of CW type $E\rightarrow B$ on $\Xx$ with
the following universal property: for any morphism $S\rightarrow T$ of
truncated weak flexible sheaves of CW type in the universe $\Uu$, the space
of pairs
$(f,\phi )$ where $f:T\rightarrow B$ is a morphism and $\phi : S/T\cong
E\times _B T/T$
is an isomorphism over $T$, is weakly contractible.

Note that $E\rightarrow B$ is not in $\Uu$ itself (this is ruled out by
Russell's paradox).

We call $B$ the {\em big classifying space}. It is essentially a union of
all possible classifying spaces.

The defining property of $E\rightarrow B$ implies that for any morphism
$f:T\rightarrow B$, the loop space $\Omega ^f(B\times T/T)$ is equivalent
to the relative space of equivalences $Equiv _{\infty} (E\times _BT/T,
E\times _BT/T)$
(this should even be true with the multiplicative structure---to state this
one would have to develop a loop space machinery for weak flexible
sheaves).

\subnumero{The construction of $E\rightarrow B$}

This section not finished...

\subnumero{Classifying spaces for local systems}

Suppose $\Ll$ is a type of local system (e.g. local system of abelian
groups, of rings, of groups, or just of sets). Then there is a classifying
space $B\Ll $ with a local system $E\Ll \rightarrow B\Ll$ of the required
type, such that for any $S$ and local system $L$ on $S$ of the required
type, the space of pairs $(f,\psi )$ where $f:S\rightarrow B\Ll$ and $\psi
: L\cong f^{\ast}(E\Ll )$ is an isomorphism of the required type of local
system over $S$, is weakly contractible.

Construction not finished...

In this case, $B\Ll$ is $1$-truncated (since the space of relative
equivalences of local systems of type$\Ll$ is $0$-truncated).

We describe $\pi _0(B\Ll )$. Let $I^{\rm pre}(\Ll )$ denote the presheaf on
$\Xx$ which associates to each object $Y$ the set of isomorphism classes of
sheaves of type $\Ll$ on $\Xx /Y$ (that is local systems over the flexible
sheaf represented by $Y$). Let $I(\Ll )$ denote the associated sheaf. We
claim that there is a natural isomorphism $I(\Ll )\cong \pi _0 (B\Ll )$.
Proof ???

Suppose $L\rightarrow S$ is a local system of type $\Ll$ over $S$, and
suppose that $\gamma \in \Gamma (\Xx , I(\Ll ))$. We say that $L$ is {\em
of similarity $\gamma$} if the image of the classifying map for $L/S$
projects to $\pi _0(B\Ll )$ by a map which factors through the point
$\gamma$.
Note that if $S$ is $0$-connected then there will always be such a $\gamma$.

Suppose that $G$ is a sheaf over $\Xx$ of type $\Ll$. Then we say that a
local system $L/S$ is {\em similar to $G$} if it is of similarity equal to the
point which is the image of the classifying map for $G/\ast$ (we will denote
this point by $\langle G\rangle $).

We have the following
characterization:
a local system $L/S$ is locally similar to $G$ if and only if,
for every object $X\in \Xx$ and every point $s\in S_X$, there exists a sieve
$\Bb$ over $X$ such that for any $g:Y\rightarrow X$ in $\Bb$, the fiber of
$L|_{\Xx /Y}$ over the point $g^{\ast}(s)$ is isomorphic to $G$.

(proof ???)

If $L$ and $M$ are two local systems over $S$, we say that $L$
is {\em similar to $M$} if the composition of the classifying
maps $S\rightarrow
\pi _0(B\Ll )$ for $L$ and $M$ are the same.

If $\gamma \in I(\Ll )$ let $B^{\gamma }\Ll$ denote the fiber of the morphism
$$
B\Ll \rightarrow \pi _0 (B\Ll )= I(\Ll ) $$
over the point $\gamma$.
If $G$ is a sheaf over $\Xx$ of type $\Ll$ then $Aut _{\Ll}(G)$
(or just $AUt (G)$ for short) is a sheaf of groups over $\Xx$ and we have $$
B^{\langle G \rangle} \Ll = K(Aut (G), 1) $$
with basepoint corresponding to the classifying map of the
local system $G\rightarrow \ast$.

In particular, if $L/S$ is a local system which is locally
similar to $G$ then the classifying map is a map
$$
f: S\rightarrow K(Aut (G), 1).
$$
There is a morphism
$$
\ast = \{ f, \psi \} \rightarrow \{ f\} = Map (S, K(Aut (G),1)) $$
where the fiber over a particular mapping $f$ is
the space of choices of $\psi$, that is the space of choices of isomorphism
$L/S\cong f^{\ast}E\Ll$. This is a principal homogeneous space for $Map (S, Aut
(G))$ which is the space of automorphisms of $f^{\ast}E\Ll $. This morphism is
just the path space fibration where the point is that corresponding to the
classifying map $f$ for $L$.

(Several unfinished sections from 1993 and 1995 supressed)

\end{document}